\documentclass{article}
\usepackage[verbose=true,letterpaper]{geometry}
\AtBeginDocument{
	\newgeometry{
		textheight=9in,
		textwidth=5.5in,
		top=1in,
		headheight=12pt,
		headsep=25pt,
		footskip=30pt
	}
}

\usepackage[utf8]{inputenc} 
\usepackage[T1]{fontenc}    
\usepackage{lmodern}

\usepackage{amsmath}
\usepackage{xfrac}
\makeatletter
\newcommand{\leqnomode}{\tagsleft@true\let\veqno\@@leqno}
\makeatother
\usepackage{amssymb}
\DeclareMathAlphabet\mathEulervmNormal{U}{zeur}{m}{n}
\DeclareMathAlphabet\mathEulervmBold{U}{zeur}{b}{n}
\DeclareMathAlphabet\mathEulervmScript{U}{zeus}{m}{n}
\DeclareMathAlphabet\mathEulervmScriptBold{U}{zeus}{b}{n}

\NewDocumentCommand{\randomVar}{m}{\mathEulervmScript{#1}}
\NewDocumentCommand{\topSpace}{m}{#1}
\NewDocumentCommand{\sigmaAlgebra}{m}{\mathcal{#1}}
\NewDocumentCommand{\measurableSpace}{m}{(\topSpace{#1}, \sigmaAlgebra{#1})}
\NewDocumentCommand{\sigmaAlgebraBorelRaw}{m}{\mathcal{B}(#1)}
\NewDocumentCommand{\sigmaAlgebraBorel}{m}{\mathcal{B}(\topSpace{#1})}
\NewDocumentCommand{\sigmaAlgebraBorelIdx}{mm}{\mathcal{B}(\topSpace{#1}_{#2})}
\NewDocumentCommand{\measurableSpaceBorel}{m}{(\topSpace{#1}, \sigmaAlgebraBorel{#1})}
\NewDocumentCommand{\measurableSpaceBorelIdx}{mm}{(\topSpace{#1}_{#2}, \sigmaAlgebraBorelIdx{#1}{#2})}
\NewDocumentCommand{\anyVar}{}{\randomVar{V}}
\NewDocumentCommand{\iidVarAny}{}{\randomVar{X}}

\NewDocumentCommand{\viewport}{}{\mathEulervmNormal{V}}
\NewDocumentCommand{\viewportEventual}{}{\mathEulervmBold{V}}
\NewDocumentCommand{\IGraph}{}{\mathEulervmNormal{G}_I}
\NewDocumentCommand{\AncDo}{}{\Anc_{\PearlDo}^L}
\NewDocumentCommand{\graphStructural}{}{G}
\NewDocumentCommand{\gUnderlying}{}{0}
\NewDocumentCommand{\gObs}{}{\txt{obs}}

\NewDocumentCommand{\gQuery}{}{\txt{query}}
\NewDocumentCommand{\graphUnderlying}{}{\graphStructural^\gUnderlying}
\NewDocumentCommand{\graphObs}{}{\graphStructural^\gObs}

\NewDocumentCommand{\graphQuery}{}{\graphStructural^\gQuery}
\NewDocumentCommand{\graph}{}{\mathcal{G}}
\NewDocumentCommand{\graphClassical}{}{\mathEulervmBold{G}}
\NewDocumentCommand{\graphClassicalDo}{}{\mathEulervmBold{G}_{\txt{do}}}
\NewDocumentCommand{\nodes}{}{\mathcal{N}}
\NewDocumentCommand{\edgesStructural}{}{\mathcal{E}}
\NewDocumentCommand{\nodesClassical}{}{\mathEulervmBold{N}}
\NewDocumentCommand{\edgesClassical}{}{\mathEulervmBold{E}}

\NewDocumentCommand{\nInner}{}{\txt{inner}}
\NewDocumentCommand{\nodesInner}{}{\nodes_{\txt{inner}}}
\NewDocumentCommand{\nodesInnerProper}{}{\nodes_{\txt{proper}}}
\NewDocumentCommand{\nodesOuterProper}{}{\nodes_{\txt{extern}}}
\NewDocumentCommand{\nodesOuterFixed}{}{\nodes_{\txt{pinned}}}
\NewDocumentCommand{\nodesOuter}{}{\nodes_{\txt{outer}}}

\NewDocumentCommand{\minLatentSets}{}{\mathfrak{L}}
\NewDocumentCommand{\ancestralStructure}{}{\mathcal{A}}
\NewDocumentCommand{\parentStructure}{}{\mathcal{P}}

\DeclareMathOperator{\absorbPa}{AbsorbPa}

\DeclareMathOperator{\absorbCh}{AbsorbCh}

\DeclareMathOperator{\disint}{disint}
\DeclareMathOperator{\marginalizeOp}{marg}

\NewDocumentCommand{\marginalize}{mm}{\marginalizeOp_{#2}(#1)}
\NewDocumentCommand{\marginalizeRight}{}{\marginalizeOp_{\leftarrow}}
\NewDocumentCommand{\disintLeft}{}{\disint_{\rightarrow}}

\DeclareMathOperator{\IdOp}{Id}
\NewDocumentCommand{\structureKernels}{}{\mathcal{F}}
\NewDocumentCommand{\modelKernels}{}{\mathcal{F}_{\txt{model}}}
\NewDocumentCommand{\knownKernels}{}{\mathcal{F}_{\txt{intervene}}}
\NewDocumentCommand{\pinnedKernels}{}{\mathcal{F}_{\txt{I}}}

\NewDocumentCommand{\forkAssign}{}{\overset{\forall}{:=}}

\NewDocumentCommand{\observedFunction}{}{\mathcal{F}_{\txt{obs}}}
\NewDocumentCommand{\knownFunction}{}{\mathcal{F}_{\txt{intervene}}}
\NewDocumentCommand{\shallowDistr}{}{P_\theta}
\NewDocumentCommand{\realizedDistr}{}{\shallowDistr^{\mathcal{D}}}

\NewDocumentCommand{\query}{}{$q=(\tilde{Y},\tilde{X},\theta)$}

\NewDocumentCommand{\knowledgeSet}{}{\mathcal{K}}
\NewDocumentCommand{\knowledgeSetBasic}{}{\knowledgeSet_{\txt{basic}}}

\NewDocumentCommand{\knowledgeSetBasicCompletion}{}{\bar\knowledgeSet_{\txt{basic}}}

\NewDocumentCommand{\kernelComposition}{mm}{#1\circ#2}

\NewDocumentCommand{\kernelCompound}{m}{[#1]}
\NewDocumentCommand{\bigKernelCompound}{m}{\big[#1\big]}
\NewDocumentCommand{\BigKernelCompound}{m}{\Big[#1\Big]}
\NewDocumentCommand{\kernelCompoundConfounded}{m}{\llbracket #1 \rrbracket}

\NewDocumentEnvironment{CopyClaim}{mmmm}
	{\theoremstyle{#1}%
	\newtheorem*{localcopy_#3}{#2 \ref{#3}}\begin{localcopy_#3}[#4]}%
	{\end{localcopy_#3}}
\NewDocumentEnvironment{CopyLemma}{mm}{\begin{CopyClaim}{plain}{Lemma}%
	{#1}{#2}}{\end{CopyClaim}}
\NewDocumentEnvironment{CopyDef}{mm}{\begin{CopyClaim}{definition}{Definition}%
{#1}{#2}}{\end{CopyClaim}}
\NewDocumentEnvironment{CopyThm}{mm}{\begin{CopyClaim}{plain}{Theorem}%
{#1}{#2}}{\end{CopyClaim}}

\usepackage{algorithm}
\usepackage{algpseudocode}

\usepackage{graphicx} 
\usepackage{caption}
\usepackage{placeins} 

\usepackage{booktabs}
\usepackage{makecell}

\usepackage{stmaryrd}	

\usepackage{tikz}
\usetikzlibrary{decorations.pathmorphing} 
\usetikzlibrary{calc}
\usetikzlibrary{backgrounds}

\usepackage{import}
\usepackage{graphicx} 
\usepackage{tabularx}
\usepackage{amsmath}
\usepackage{amsthm} 
\usepackage{amsfonts}
\usepackage{bbm}
\usepackage[utf8]{inputenc}
\usepackage{enumitem}
\usepackage{tikz}
\usetikzlibrary{arrows}
\usepackage{xfrac}		
\usepackage{enumitem}	
\usepackage{braket}
\usepackage[normalem]{ulem} 

\usepackage{todonotes}
\usepackage{dirtytalk} 

\usepackage{hyperref}
\hypersetup{allcolors=black, colorlinks}

\makeatletter
\newcommand*{\centernot}{%
	\mathpalette\@centernot
}
\def\@centernot#1#2{%
	\mathrel{%
		\rlap{%
			\settowidth\dimen@{$\m@th#1{#2}$}%
			\kern.5\dimen@
			\settowidth\dimen@{$\m@th#1=$}%
			\kern-.5\dimen@
			$\m@th#1\not$%
		}%
		{#2}%
	}%
}
\makeatother

\newcommand{\independent}{\perp\mkern-9.5mu\perp}
\newcommand{\notindependent}{\centernot{\independent}}

\newcommand{\eg}{e.\,g.\ }
\newcommand{\Eg}{E.\,g.\ }
\newcommand{\ie}{i.\,e.\ }
\newcommand{\Ie}{I.\,e.\ }

\newcommand{\asswlog}{w.\,l.\,o.\,g.\ }

\newcommand{\BreakingSmallSpace}{\hspace{.16667em}}
\newcommand{\Slash}{\,/\BreakingSmallSpace}

\DeclareMathOperator{\PearlDo}{do}

\DeclareMathOperator{\Pa}{Pa}
\DeclareMathOperator{\pa}{pa}
\DeclareMathOperator{\Ch}{Ch}
\DeclareMathOperator{\Anc}{Anc}
\DeclareMathOperator{\Dec}{Desc}

\DeclareMathOperator{\UniformDist}{Unif}

\DeclareMathOperator{\dx}{dx}
\newcommand{\Reals}{\mathbb{R}}
\NewDocumentCommand{\txt}{m}{\text{\normalfont{#1}}}

\DeclareMathOperator{\id}{id}

\DeclareMathOperator{\img}{img}
\DeclareMathOperator{\FinSubSets}{FiniteSubsetsOf}
\DeclareMathOperator{\ds}{ds}
\DeclareMathOperator{\dt}{dt}
\DeclareMathOperator{\du}{du}
\DeclareMathOperator{\dw}{dw}
\DeclareMathOperator{\dy}{dy}
\DeclareMathOperator{\dz}{dz}
\DeclareMathOperator{\dm}{dm}
\DeclareMathOperator{\intdp}{dp}

\NewDocumentCommand{\PaIdx}{m}{\Pa^{(#1)}}

\NewDocumentCommand{\approxsubset}{}{\substack{\subset\\[-0.15em]\sim}}

\NewDocumentCommand{\halfquad}{}{\mkern9mu}

\NewDocumentCommand{\val}{m}{\mathcal{#1}}

\theoremstyle{plain}
\newtheorem{lemma}{Lemma}[section]
\newtheorem{conjecture}{Conjecture} 
\newtheorem*{conjectureNamed}{Conjecture} 

\newtheorem{thm}{Theorem}

\newtheorem{cor}[lemma]{Corollary}

\theoremstyle{definition}

\newtheorem{definition}[lemma]{Definition}
\newtheorem{assumption}[lemma]{Assumption}

\newtheorem{lemmaDef}[lemma]{Lemma\Slash{}Definition}
\newtheorem{notation}[lemma]{Notation}

\newtheorem{example}[lemma]{Example}

\newtheorem{rmk}[lemma]{Remark} 

\usepackage{natbib}
\setcitestyle{numbers,square}

\usepackage{microtype}

\usepackage{hyperref}

\usepackage{graphicx}

\title{Symmetries and Causality:\\Causal Effect Identification Beyond IID Data}

\author{Martin Rabel\textsuperscript{a}, Jakob Runge\textsuperscript{a}}
\date{\today}

\NewDocumentCommand{\IDataset}{}{I_{\txt{dataset}}}
\NewDocumentCommand{\IVars}{}{I_{\txt{vars}}}
\NewDocumentCommand{\ISample}{}{I_{\txt{sample}}}

\DeclareMathOperator{\RelevantCh}{Ch\textsuperscript{relevant}}

\NewDocumentCommand{\IAux}{}{I^{\txt{aux}}}

\begin{document}
	\maketitle

	\noindent\begin{minipage}{\textwidth}
		\centering
		{\footnotesize\textsuperscript{a}%
		Department of Computer Science, University of Potsdam, Potsdam, Germany}
	\end{minipage}
	\begin{abstract}
		We develop theory for the identification of complex
		causal queries from non-IID data.
		We propose a formal description of stochastic systems
		based on symmetries in data leaving causal mechanisms invariant.
		These symmetries become the source
		of statistical repetition for identification
		and mandate knowledge-transfer from observations to queries.
		The result is an abstract, simple and very general language
		for causal reasoning.
		We give descriptions of models and queries,
		and the formal infrastructure for their
		mathematically rigorous identification.
		The proposed structure enables a number
		of simplifying principles:
		A shared description of mechanisms and interventions,
		a duality of observed and queried structure,
		and the sub-division of the
		identification task
		into three steps -- extraction, computation and query-expression.		
		Our approach reproduces and matches standard theoretical results on IID data
		and transport of experimental and non-experimental data.
		By design, it additionally applies beyond IID data and
		captures complex causal queries
		not expressible via do- or soft-interventions.
	\end{abstract}

	\section{Introduction} 

	The scientific method of exploring and understanding our reality seems to
	naturally evoke a notion of cause and effect in our reasoning.		
	We may not know \emph{why} causality arises in
	the structure of scientific study and knowledge, yet
	this evidence suggests a deep role of causal reasoning in
	\emph{efficient} world modeling.
	Efficient for example in the ability to gain agency.
	Causality in general, as motivating such claims,
	is a rather vaguely defined concept.
	Certain formalisms, like the do-operator and graphical-models
	\citep{PearlBook} or
	potential outcomes \citep{Rubin1974a, holland86},
	have had
	tremendous success in capturing causal aspects of stochastic systems.
	But many real-world problems and relevant questions cannot easily be captured
	by this notion of causality and form of model;
	an example is illustrated in Fig.\ \ref{fig:intro}.
	
	The main purpose of this paper is to make progress on these problems:
	Applicability of causal language to complex data-sets,
	and the extraction of answers to complex questions.
	Both problems are indications that we need a more
	expressive formal approach.
	Improving expressiveness will also bring
	the notion of causality actually captured by the formalism
	closer to the vague concept of causality that motivates its use in the first
	place.
	
	Many approaches to broaden the scope of causal inference have
	been proposed.
	Examples include
	the incorporation of experimental data and transportability
	between contexts \citep{Bareinboim2016TransportOverview, SelectionVars}
	for effect estimation or
	structure learning on multi-context data \citep{CD-NOD,JCI,JPCMCI},
	or more complex queries like soft-interventions \citep{correa2020general},
	interventions in time-series \citep{Runge2023},
	program-evaluation \citep[§4.4]{PearlBook}
	or mediation questions \citep{Pearl2012};
	many more have been studied.
	These approaches have brought many invaluable theoretical
	and practical insights. But at the same time,
	they seem to largely happen in isolation.
	
	We may be able to predict certain soft-interventions,
	and to predict do-interventions from experimental data,
	but already the prediction of anything from data obtained
	doing a experiment by soft-intervention
	would require additional technology.
	We may be able to combine IID data-sets
	\citep{Bareinboim2016TransportOverview, CD-NOD,JCI}, but for time-series data
	we suddenly need new ideas \citep{JPCMCI,gunther2024causal}.
	Of course any such combination could be solved individually,
	but this is not a theoretically or practically satisfactory answer.
	What is needed is
	a shared language, that can incorporate and combine all, or at least many,
	such ideas uniformly.
	
	Before we outline the contents of this paper,
	we first briefly analyze, on a conceptual level, why the problem at hand
	is difficult to handle in existing approaches,
	and what basic
	formal structure we seek to employ instead.

	\begin{figure}[!ht]
		\input{figures/intro_new}
		\caption{
			Real-world motivation (see
			example \ref{example:real_world_questions}):
			The step "$Z$" in an existing three-step production-process (1)
			shall be upgraded to "$Z_2$". Use qualification-testing
			results of $Z_2$ (2,3) to predict future performance (query).\\
			Nodes are data-points\Slash{}random variables,
			labeled with a capital letter indicating their causal mechanism.
			Edges point from causal parents to the node whose
			mechanism they affect.
			We draw observed nodes green (solid, rounded corners),
			latent nodes gray (dashed, rounded corners),
			nodes of interest are highlighted blue (circles)
			and interventions red (rectangles).
			We combine three datasets:
			(1) Univariate observations of a time-series.
			The model is three-periodic.
			(2) Repeated observations of two variables,			
			modeled as IID observations.			
			(3) Repeated observations
			obtained experimentally, \ie after an intervention,
			from the same system as (2). 
			A plausible question is: What happens if we replace
			the mechanism $Z$ by $Z_2$ in the time-series setup.
			E.\,g., given an initial state $Y_{\txt{init}}$,
			after four periods of evolution,
			what is the distribution of the final value (blue circle)?\\
			Intuitive solution:
			From repeated observation learn blue background boxes ("learn").
			Then piece this information together, gray background rectangle ("query").
			Example \ref{example:real_world_questions}\,(b) shows how to
			construct a sound estimator using the material of the paper,
			putting this intuitive idea into practice.
		}\label{fig:intro}
	\end{figure}
	
	\subsection{Abstraction Principles}
	\label{sec:intro:need_abstraction}
	
	One conceptual idea that practically all generalizing approaches
	seem to share is that of adding definitional structure to model and\Slash{}or query:
	Add selection-variables \citep{Bareinboim2016TransportOverview,SelectionVars}
	(often called context-variables) to grow the system into a "meta"-system
	\citep{JCI}, add new types of interventions, to describe queries that
	are not do-interventions, add proxies for hidden structure \citep{CD-NOD,JPCMCI}.
	Combining multiple formalisms will be expensive, if only because it
	requires a compatible redefinition of added structures.

	Our approach advances in the opposite direction, toward a more \emph{abstract} model.
	Instead of adding yet more definitional structure,
	the basic question we set out to ask is: What structure can be removed
	from model and query?
	A more abstract formalism
	applies more generally simply by virtue of having \emph{fewer constraints}
	to its application.
	
	The first structural aspect to be removed is IIDness of the data,
	which is implicit in the description of causal properties by structured
	causal models (SCMs) or potential outcomes (POs).
	For example mz-transport
	\citep{Bareinboim2016TransportOverview,SelectionVars}
	or JCI \citep{JCI} have to \emph{add} selection-\Slash{}context-variables
	to restore IIDness on a meta-SCM.
	If no IIDness must be restored, no additional structure is required.
	Complex questions in applied science involve combination of
	data from, and prediction for, different environments, sub-scale experiments etc.,
	thus inherently cannot be phrased in IID-setups.
	The ubiquity of IID setups in statistics is
	based primarily on the idea that a researcher curates suitable data-sets;
	but then the decision \emph{how} to curate data-sets
	is never phrased formally, thus also any subsequent combination
	of such data-sets must be performed explicitly and by hand.
	
	The second structural aspect to be removed is randomness
	-- the presence of random variables in the model --
	which is just as deeply rooted
	in existing formalisms.
	Removing randomness \emph{from the model and query}
	does not in any way impede the ability to describe stochastic systems.
	Indeed the perspective is similar to the one commonly taken on time-series.
	A random process can, under suitable conditions,
	be modeled by Markov-kernels -- certain probability-kernels,
	which are purely measure-theoretic: They do not reference
	a basic probability space $(\Omega,\sigmaAlgebra{A},P)$ or random elements.
	
	Writing a model in measure-theoretic rather than probability-theoretic terms
	is crucial not only for a clear and rigorous
	mathematical formalism but also to enable a number of
	far reaching structural principles.
	For example we can describe the model from probability-kernels
	as building-blocks. If a building block is a pre-existing mechanism
	or an intervention is insubstantial; thus describing soft- or stochasitic
	interventions and experimental data requires no additional work.
	Further, we can describe observations and queries uniformly in the same model;
	thus a query can be as complex as any observation might be,
	dramatically extending the scope of potential queries.
	This duality of observations and queries also allows
	the reuse of large parts of the formal structure.
	Finally, 
	we can also more easily sub-divide the problem of query-identification into
	multiple sub-tasks (see below).

	It is important that the simplified structure
	still \emph{retains the ability to describe the causal aspects} of the problem.
	One can always remove formal structure, if one
	is willing to give up on the ability to rigorously reason about,
	or automate, a problem at hand -- indeed it has sometimes
	been claimed that causal inference in general were unnecessary:
	it only translates causal problems to statistical ones,
	a task that could be performed by hand. This trivial
	type of removal of structure is \emph{not} what we aim for.
	Instead we show that the weak structure introduced below
	indeed suffices to tackle the problem of identification
	of causal queries from data.

	\subsection{Invariant Mechanisms}
	\label{sec:intro:invariant_mechaisms}
	
	We build our formalism on a concept that is already
	pervasive in causal study: Invariance of mechanisms.
	While formally somewhat elusive, invariance of mechanisms
	has played a pivotal role in causal inference in general
	\citep{PearlBook,Elements},
	for discovering causal structure
	\citep{CD-NOD,JCI, JPCMCI, peters2016causal}
	and for robust prediction
	\citep{RojasCarulla2018a, AnchorRegression}
	or transfer
	\citep{Bareinboim2016TransportOverview, SelectionVars}.
	Invariance of mechanisms is also baked into
	the SCM and PO frameworks, and while usually not made
	explicit has certainly played a core role in their motivation;
	for example	\citet[p.\,322]{PearlBook}
	states (slightly paraphrased) that
	"what we call 'causal knowledge'"
	are the "extra assumptions that identify
	what in the distribution remains invariant when the specified modification takes place".
	Invariance is the property of remaining unchanged under a given symmetry;
	a symmetry is a transformation of, in our case, the causal model
	(\eg, in Fig.\ \ref{fig:intro}, the model of (1) is invariant under
	time-translations by multiples of 3 time-steps, the model of (2) and (3)
	is "exchangeable", \ie invariant under permutations of samples).
	The symmetry Pearl seems to have in mind is the transformation
	switching mechanisms of observed and intervened model.
	Indeed this is the
	form of symmetry we will need to reproduce IID interventions as
	a special case (§\ref{apdx:iid}).
	\Eg in Fig.\ \ref{fig:intro} the mechanisms in the box "query"
	are the same as in some of the observations -- these mechanisms locally
	are invariant under a symmetry exchanging query and observation nodes.

	Invariant mechanisms in our formalism
	are probability kernels --
	inspired by the
	use of Markov-kernels in the
	process literature (see \eg \citep[§7]{kallenberg1997foundations}).
	Additionally  we explicitly attach a symmetry, leaving the mechanism
	invariant within a range of applicability.
	In our formalism, a single model can describe both
	the observed world and the query uniformly.
	This allows to make precise (in §\ref{sec:queries})
	the idea about invariance under switching
	mechanisms of observed and intervened model alluded to before.

	Some kernels are considered known a priori.
	These are conventionally called "interventions" and usually are \emph{not}
	symmetric in observation and query.
	Other kernels
	can be extracted, individually or only as part of larger confounded structures,
	from data by virtue of repetition from their symmetry properties.
	Of course statistics, at its very foundations, is about
	learning shared stochastic properties by combining many observations;
	so relating invariant structures built from mechanisms to
	suitable subsets of data is not a fundamentally new idea -- however its
	consistent, rigorous and explicitly formal execution to capture causality is novel.
	
	We elevate this weak notion of invariance of mechanisms
	to become the fundamental formal starting point of causal modeling.
	After establishing a canonical link to probabilistic observations
	(now including randomness), this captures
	a meaningful notion of causality. A notion of causality
	that contains the IID, interventional setup (SCMs, PO) as a special case,
	but readily extends to non-IID data and non-interventional queries,
	for example to Fig.\ \ref{fig:intro}.

	\subsection{Contributions and Content}
	
	An important reason why selection-\Slash{}context-variables
	have been so successful
	\citep{Bareinboim2016TransportOverview,SelectionVars,CD-NOD,JCI}
	is their ability to map a new problem to an established
	setup -- thus much of the required technology to work on the
	resulting mapped problem was already available.
	Our proposed formalism is a simpler (more low-level)
	language. Many ideas of modeling problems can be compiled
	into this simpler (and unified) language and thus remain intact;
	such ideas are among the core contributions of the existing frameworks
	\citep{pearl2022external,JCI}.
	However, getting information back out of a low-level modeled system is
	inherently more difficult as there is less definitional structure
	to work with. For example, identification-strategies like the Id-algorithm
	\citep{tian2002general,shpitser2006identification2} or mz-transport
	\citep{Bareinboim2013TransportAlgo, bareinboim2012TransportCompleteness}
	or even the do-calculus \citep{PearlBook}
	do not carry over.
	\textbf{The primary contribution of this paper is
	to provide the formal machinery to re-establish the link from
	modeled problem-statement to identification-strategy}.
	
	Constructions on simple structures are difficult.
	For example
	we will often encounter non-unique constructions (see
	\eg Rmk.\ \ref{rmk:non_uniqueness}).
	Despite this, we can keep the identification-problem tractable by dividing it into
	three logically independent parts
	-- largely by virtue of having deleted randomness from the model.
	These parts are the description and transformation of\Slash{}computation with
	knowledge about a model;
	the extraction of knowledge from observations; and the phrasing of
	causal queries as questions about what knowledge needs to be computed.
	After a formal definition of our modeling approach
	§\ref{sec:models}, these points also dictate the structure of this paper,
	with §\ref{sec:structured_kernels} introducing structure
	and graphical computations (Lemmas \ref{def:simplify_simp}, \ref{lemma:c_components:mt}, \ref{lemma:glue:mt}),
	§\ref{sec:extraction_from_data} extracting structured knowledge from data
	(Thm.\ \ref{thm:extract_from_backdoor_complete}),
	and §\ref{sec:queries} finally answering causal queries via
	the previously introduced material (Thm.\ \ref{thm:query_id}).
	Proofs are deferred to the appendix.
	
	We do \emph{not} discuss observational support,
	selection bias or correlated missingness.
	We also do \emph{not} focus on completeness of the proposed strategy
	-- it \emph{is} complete for IID-data and likely for mz-transport,
	see §\ref{sec:relate_to_iid} -- instead we formulate transparent
	conjectures \ref{conjecture:regular_computation},
	\ref{conjecture:regular_graphical}, \ref{conjecture:extraction},
	\ref{conjecture:decomp_first}
	that clarify what questions
	concerning completeness remain open.
	We provide algorithmic descriptions, to elucidate
	what systematic approach an applied researcher might take to
	write down a specific identification strategy, and to
	explain the formal structures.
	While our approach is amenable to automation, immediate
	implementation in code is not the primary goal of the presented
	algorithms.

	\subsection{Related Literature}	
	
	We focus on estimator-agnostic strategies for using causal information
	about a system to propose a concrete estimation-problem that answers a
	causal query, similar to \citep{tian2002general,shpitser2006identification2,
	bareinboim2012TransportCompleteness}. \Ie we focus on
	causal inference aspects \citep{PearlBook,Elements}, not on estimation theory.
	
	Some of the cornerstone results of causal inference in multi-context setups,
	like transportability of experimental and non-experiment information \citep{Bareinboim2013TransportAlgo,bareinboim2012TransportCompleteness,
		Bareinboim2016TransportOverview, SelectionVars},
	causal discovery by context-variable and similar approaches
	\citep{CD-NOD,JCI} and by invariant prediction \citep{peters2016causal}
	were already discussed above.
	Our results are similar in purpose to
	what would traditionally be called identification
	of interventional distributions or effect estimation.
	Among results in the literature transportability
	\citep{bareinboim2012TransportCompleteness}
	(for details see §\ref{apdx:mz_transport})
	seems closest in terms of goals
	(comparable results are available in potential outcome language,
	see \eg \citep{margueritte2026learning}).
	There are also many results for identification of interventional distributions
	in IID-data \citep{tian2002general,shpitser2006identification2,
	shpitser2006identification1} building on Pearls's
	do-calculus (see \eg \citep{PearlBook}). We encounter analogues for many of the
	structures like c-components \citep{tian2002general}, c-trees, c-forests
	\citep{shpitser2006identification2} and ultimately for do-interventional
	queries also hedges \citep{shpitser2006identification2}.
	A comprehensive comparison
	to the single-context IID-case
	\citep{shpitser2006identification2,shpitser2006identification1}
	and transportability \citep{bareinboim2012TransportCompleteness}
	is given in §\ref{apdx:iid}.
	Known results for the IID-case also include counter-factual
	queries \citep{IDAlgoCF,Shpitser2013}.
	We do not discuss counter-factuals specifically,
	but coincidentally our flexible query-formulation
	captures some mediation questions (natural direct effects)
	recovering known formulas \citep{Pearl2001}, see §\ref{apdx:iid:mediation}.
	
	More abstract formulations for causality
	have been proposed \eg by \citep{CausalSpaces, orujlu2026causal}
	to clarify the relation of causality to probability.
	Here, we primarily abstract to gain flexibility for
	non-IID applications and complex queries which leads to a mostly
	orthogonal to  \citep{CausalSpaces}
	(detached from IID and random elements) result;
	our primary focus is on re-establishing
	a pathway to identification of complex queries in non-IID systems
	from observations, the abstraction is guided by (and subordinate to) this goal.
	We also want to mention the perspective of \citep{janzing2023reinterpreting}
	of causality as the ability to predict unobserved joint distributions,
	which is philosophically close to our approach but focused on
	IID observations and statistical learning, and has
	very interesting consequences for analyzing finite-sample
	reliability of CD-methods \citep{faller2024self}.
	
	We do not investigate the discovery of the structures
	defining our models.
	Causal structure discovery (CD) results beyond
	the IID-case \citep{PCalgo,PCstable,spirtes2001causation}
	exist: time-series \citep{pcmci,pcmci_plus},
	multi-context data \citep{CD-NOD,JCI},
	proxies of hidden confounders \citep{CD-NOD, JPCMCI},
	context-specific structure
	\citep{gunther2024causal}, see also
	\citep{LDAG_definition,LDAG_logical,LDAG_learning},
	for learning both context-structure and graphical structure \citep{Saggioro2020,rahmani2023castor,BalsellsRodas2023,
		mameche2025spacetime,rabel2026contextspecificcausalgraphdiscovery},
	learning differences between data-sets only \citep{assaad2024causal}.
	Learning more flexible structure seems initially more difficult,
	but more flexible structure is also more likely to exist (which is nice
	for anything one wants to find). Further,
	the usefulness of non-IIDness for improving causal discovery
	\citep{CD-NOD,JCI,peters2016causal}, especially edge-orientations,
	suggests that a principled approach focused on invariance in modeling
	may be \emph{helpful} to extracting information.
	
	Symmetries and invariance also play an important role in modern probability
	theory (for a great overview, see \citep{kallenberg2005probabilistic}).
	The causal questions we ask sometimes touch on trivial cases of
	deep results like de Finetti's theorem (exchangeable is conditional IID;
	cf.\ Rmk.\ \ref{rmk:exchangable} or example \ref{example:multi_level_queries}),
	but the intricacies of such characterizations (\eg of exchangeability
	by conditional IIDness)
	and other topics commonly studied in probability-theory,
	while interesting problems of their own right, seem mostly disjoint from
	the causal questions we are interested in. 
	
	Also in machine-learning, especially for parameter-sharing,
	symmetries are of great relevance, for example in convolution-models
	(image-convolution, but also graph-convolutional networks etc.)
	or in attention-mechanisms (for example in transformers),
	see any textbook on the topic, \eg \citep{murphy2022probabilistic}.
	However, their use to guide the transfer to causal queries and the systematic
	structuring of world-models and -exploration seems to be
	largely unexplored.
	But toy-models for causal structure of world-models
	and exploration,
	specifically
	"causal bandits" \citep{lattimore2016causal},
	have received substantial attention lately.
	Causal models are believed to take a special role for invariant prediction,
	as justified for example by results like \citep{rodas2021causal,AnchorRegression}
	(and references therein). Finally, our results (Cor.\ \ref{cor:revealing_simp})
	have some unexpected
	connections to simple missing-data problems \citep{rubin1976inference}.
	A causal perspective on missing-data problems can be found for example
	in \citep{de2025causal} and references therein.

	\section{Models and Observations}\label{sec:models}
	
	We introduce the basic formalism,
	consisting of models
	and their relation to statistical observations,
	including observedness and asymptotic limits.
	
	\subsection{Preliminaries} 
	
	We are interested in causal properties, not
	in technical intricacies of measure-theory,
	hence we assume throughout this paper (see also §\ref{apdx:ptheo_basics}):
	
	\begin{assumption}\label{ass:standard_borel}
		Measurable spaces $\measurableSpaceBorel{X}$
		are standard Borel, \ie
		Polish spaces $\topSpace{X}$ with their associated Borel $\sigma$-algebra
		$\sigmaAlgebraBorel{X}$.
		Examples for Polish spaces are $\Reals^n$ (with standard topology) and
		open or closed subsets of $\Reals^n$ (with the induced topology),
		\eg finite sets with the discrete topology.
	\end{assumption}
	
	We describe causal models based
	on probability kernels. We give a preliminary definition and notation
	here and discuss technical details in §\ref{apdx:kernels}.
	
	\begin{definition}[Kernels, preliminary]\label{def:kernels_prelim}
		Given measurable spaces $\measurableSpaceBorel{S}$ and $\measurableSpaceBorel{T}$,
		a measurable
		mapping $f: \topSpace{S} \rightarrow \mathcal{P}(\sigmaAlgebraBorel{T})$
		into the probability measures on $\measurableSpaceBorel{T}$
		is called a probability kernel.
	\end{definition}
	
	Symmetry is an abstract concept, which we formally
	express via group-actions.
	Group-actions are ubiquitous in the natural sciences and in many cases
	intuitive (e.\,g., possibly $k$-periodic, time-translation symmetry,
	cf.\ Fig.\ \ref{fig:intro}).
	More abstract formal descriptions are possible §\ref{apdx:symmetries}.
	
	\begin{definition}[Symmetry]
		\label{def:symmetry}
		\emph{Group-Actions:}
		A (left) group action $G\curvearrowright I$ of a group $G$ on a set $I$ is 
		a mapping $\cdot:G \times I \rightarrow I, (g,i) \mapsto g\cdot i$
		such that $\forall i\in I: e \cdot i = i$ and
		$\forall g,h \in G: \forall i\in I: (g*h)\cdot i = g \cdot (h \cdot i)$.
		Given a sub-group $H \subset G$, there is an induced
		action $\cdot_H: H \times I \rightarrow I,
		(h,i) \mapsto h \cdot_H i := i_H(h) \cdot i$.
		
		\emph{Properties:}
		A group action $G \curvearrowright I$ is called: effective,
		if $(\forall i\in I: g \cdot i = i) \Rightarrow g=e$,
		free if $\forall i\in I: (g \cdot i = i \Rightarrow g=e)$,
		transitive if $\forall i,i' \in I:\exists g\in G: i' = g \cdot i$.
		Given $I_0\subset I$ we will call an action free\Slash{}transitive
		on $I_0$ if the corresponding condition holds $\forall i,i' \in I_0$.
		
		\emph{Orbits:}
		An orbit is a set $G\cdot i=\{i'\in I|\exists g\in G: g \cdot i = i'\}$,
		the set of orbits is $\sfrac{I}{G}$.
		
		\emph{Invariance\Slash{}Equivariance:}
		A mapping $f : I \rightarrow \topSpace{Z}$ is called
		$G \curvearrowright I$ invariant, if $\forall i\in I: \forall g\in G:
		f(g\cdot i) = f(i)$.
		A mapping $f' : I \rightarrow I'$ is called
		($G \curvearrowright I$, $G \curvearrowright I'$) equivariant,
		if $f'(g\cdot i) = g \cdot f'(i)$.
		
		\emph{Notation:}		
		We fix a group $G$ and an effective group action
		$G \curvearrowright I$.
		A symmetry is denoted as (induced action of) a subgroup $H\subset G$.
	\end{definition}

	\subsection{Model}	

	Our models are purely measure-theoretic constructions (no $(\Omega, \sigmaAlgebra{A},P)$
	or random elements appear).
	In contrast to the IID-setting, we cannot presuppose structure on
	the index-set (indexing the data points). 
	We fix only the following notation.
	
	\begin{notation}
		\label{notation:index_set_and_spaces}
		We fix a countable index-set $I$,
		and measurable spaces $\{ \measurableSpaceBorelIdx{X}{i} \}_{i\in I}$
		that are standard Borel (Ass.\ \ref{ass:standard_borel}).
		Given $J\subset I$, we denote
		$\topSpace{X}_J := \prod_{j\in J} \topSpace{X}_j$.
	\end{notation}
	
	Our models will specify mechanisms and their associated symmetries;
	associated in the sense of leaving this particular mechanism invariant.
	
	\begin{definition}[Mechanism]
		\label{def:mechanism}
		A mechanism $(f,J,H,\Pa)$ consists of probability kernel
		$f$ from
		$\topSpace{X}^{\Pa}=\prod_{k=1}^\kappa\topSpace{X}^{(k)}$,
		with $\kappa$ arguments, to
		$\topSpace{Y}$ (both standard Borel, Ass.\ \ref{ass:standard_borel})
		together with
		a subset $J\subset I$ which we call a region of applicability,
		a symmetry $H$, acting freely and transitively on $J$,
		and a $H$-equivariant mapping
		$\Pa : J \rightarrow I^\kappa$ of relative parents
		such that:
		
		Denoting the $k$th parent by $\PaIdx{k}$,
		these satisfy $\forall k: \PaIdx{k}(j)\neq j$ and
		$k\neq k'$ $\Rightarrow$ $\forall j\in J$:
		$\PaIdx{k}(j) \neq \PaIdx{k'}(j)$.
		And
		$\forall j\in J$,
		$\topSpace{Y} \subset \topSpace{X}_j$ and
		$\forall k$:
		$\topSpace{X}_{\PaIdx{k}(j)}\subset\topSpace{X}^{(k)}$.
	\end{definition}
	\begin{definition}[Model]
		\label{def:model}
		A model $\mathcal{M}=\{(f_J,J,H_J,\Pa_J)\}_{J\in\mathcal{J}}$
		is a collection
		of mechanisms, which we
		index by their region of applicability $J$,
		such that the index set $I$ (notation \ref{notation:index_set_and_spaces})
		is the disjoint union of regions of applicability,
		formally:
		$I = \cup_{J\in\mathcal{J}} J$
		and $J\neq J'\in \mathcal{J} \Rightarrow J\cap J'=\emptyset$.
		
		Mechanisms $f_J$ are in one of two (disjoint) categories
		specified as part of the model:		
		Either $f_J$ is considered initially unknown $f_J \in \observedFunction$
		or known $\tilde{f}_J \in \knownFunction$
		("interventions" transfer external, fixed knowledge
		into the model).
		Every $i\in I$ is by construction contained in exactly one $J\in\mathcal{J}$
		which we denote by $J(i)$ and we write $f_i := f_{J(i)}$.
	\end{definition}
	
	This definition contains IID and time-series SCMs as special cases
	(examples below). The single-context IID case is inherently simpler,
	and already at the end of §\ref{sec:structured_kernels}, in example
	\ref{example:iid_id_from_structured_only},
	we can identify do- and soft-interventions.
	Details on SCM notation (§\ref{apdx:iid:standard_results})
	and their "translation" into our formalism (§\ref{apdx:translate_IID_setup})
	are given in the appendix §\ref{apdx:iid}.
	
	\begin{example}[IID-Models from SCMs]
		\label{example:model_iid}
		Given is a SCM on variables indexed by
		some finite set $\IVars$,
		consisting of, for all $v\in\IVars$:
		a parent-set $\Pa_v$,
		a SCM-mechanisms $g_v$ and a noise-distributions $P(\eta_v)$.
		Then there is an associated model $\mathcal{M}$ on
		$I=\IVars \times \ISample$, where $\ISample=\mathbb{N}$ contains
		the sample-index of a data point. It is for the purpose of this example
		sufficient to consider a finite $\ISample=\{1,\ldots, N\}$,
		then $G=\mathfrak{S}_{\ISample}$ can be picked
		the group of permutations of $\ISample$, with $\pi\in\mathfrak{S}_{\ISample}$
		acting as $\pi\cdot(v,s)=(v,\pi(s))$;
		for $\ISample=\mathbb{N}$, this permutation-group can be replaced
		by the "symmetric group" of $\ISample$.
		Kernel-mechanisms are $f_v(\pa_v) := g_v(\pa_v,\cdot)_* P(\eta_v)$,
		the push-forward of $P(\eta_v)$ or, equivalently,
		the distribution\Slash{}law
		of the random variables $V_{\pa_v}(\omega):=g_v(\pa_v,\eta_v(\omega))$.
		Further,
		$H_v=G$, $J_v=\{v\}\times\ISample$ and for
		$j=(v,s)\in J_v$, set $\Pa_v(j)=\Pa_v\times\{s\}$.
		
		An intervention on a SCM replaces a SCM-mechanism $g_v$
		by a known function, for example a constant value for
		do-interventions. An interventional SCM
		$M^{\PearlDo}$ for $\PearlDo(X=x)$, where $X\subset\IVars$,
		induces a model $\mathcal{M}^{\PearlDo}$ by the same construction as above.
		Now, mechanisms associated to $v\in X$ are
		known: $\forall v\in X:f_v\in\knownFunction$.
		For do-interventions these $f_v\in\knownFunction$
		are singular measures $\delta[x_v]$ without parents,
		$\Pa^{\PearlDo}_v(j)=\emptyset$.
		The singular measure $\delta[x]$ is defined by: for all measurable $B$,
		$\delta[x](B) = 1$ if $x\in B$, $\delta[x](B) = 0$ otherwise.
		We will not need interventional models however,
		instead we can just define a unified model
		by modifying $\mathcal{M}$, for example on
		$I^{\txt{unified}}=\IVars \times (\ISample\sqcup\{*\})$,
		where $*$ is a disjoint point; mechanisms at $v\notin X$
		are symmetric under permutation of the full $\ISample\sqcup\{*\}$,
		while $v\in X$
		has $f_{(v,*)}=\delta[x_v]\in\knownFunction$.
		For $v\notin X$, this is the invariance
		under switching observed and intervened world
		informally introduced in §\ref{sec:intro:invariant_mechaisms}.
		We continue from here in
		example \ref{example:iid_do_intervention}.
		See also Fig.\ \ref{fig:do_query} or Fig.\ \ref{fig:mz_transport}.
		Further details are given in §\ref{apdx:iid}.
	\end{example}
	
	\begin{example}[Timeseries-Models]
		\label{example:timeseries}
		Given is a time-series SCM on variables indexed by
		a finite set $\IVars$. 
		For $v\in\IVars$, there are
		parent-sets $\Pa_v\subset \IVars \times \{-\tau,\ldots, 0\}$,
		where $\tau$ is some maximal lag,
		SCM-mechanisms $g_v$ and noise-distributions $P(\eta_v)$.
		Given an initial state $\vec\nu$ (for example a stationary distribution
		if one exists)
		for $\tau + 1$ time-steps,
		there is a model $\mathcal{M}$ with
		$I=\IVars \times T$, where $T=\mathbb{N}_0$, 
		with the first $\tau$ steps fixed to $\vec\nu\in\knownFunction$.
		At later time-steps, we have time-translation invariant mechanisms, \ie
		$H_v=G=\mathbb{Z}$, acting by $a\cdot (v,t):= (v,t+a)$ for $a\in G$.
		Such mechanisms are
		often called time-homogeneous \citep{kallenberg1997foundations},
		or sometimes "causally stationary" \citep{Runge2023}.
		As before, define kernel-mechanisms
		$f_v(\pa_v) := g_v(\pa_v,\cdot)_* P(\eta_v)$.
		The $H_v$-equivariance of $\Pa_v : J \rightarrow I^{\kappa_v}$
		amounts to fixing pairs $(v,\Delta t)$ of a variable and a lag,
		automatically matching the intuition and hand-crafted form of
		a time-series SCM.
		The choice of an initial state
		is made explicit in our formalism.
		While sometimes obfuscated a little (usually by defaulting
		to the use of a stationary distribution), this cannot be avoided:
		In cases where no or multiple stationary distributions exist,
		this choice cannot be made implicitly. Even in cases where a unique
		stationary distribution does exist, there seems to be little reason
		to needlessly restrict the formalism to this particular choice,
		as other initial states may make practical sense.
	\end{example}

	\subsection{Shallow Distribution and Realization} 

	We fix a model $\mathcal{M}$.
	Via the parent-mappings $\Pa_J$ of its mechanisms,
	$\mathcal{M}$ defines a single "large" graph on nodes $I$.
	With $I$ indexing data-points, this graph has one node for every data-point;
	in the IID case with $N$ samples, it consists of $N$ disjoint copies of
	the standard causal graph (see illustration in Fig.\ \ref{fig:intro}).
	\begin{definition}[$I$-Graph]
		\label{def:Igraph}
		The $I$-graph $\IGraph$ has nodes $I$ and
		an edge $i\rightarrow i'$ if $i\in \Pa_I(i')$, where
		$\Pa_I(i') := \{i\in I|\exists k:i = \PaIdx{k}_{J(i')}(i')\}$.
	\end{definition}
	
	We assume the $I$-graph is acyclic, as we want independent mechanisms.
	In principle, indices contained in a finite cycle can be replaced by a single
	index (consistent with Ass.\ \ref{ass:standard_borel}, which is closed under finite
	products).
	
	\begin{assumption}[Acyclic $I$-Graph]\label{ass:acyclic}
		We assume that the $I$-graph is acyclic,
		and fix a total order $\pi_I$.
	\end{assumption}
	
	For formal simplicity (to avoid issues with non-unique solvability of
	equations and filtrations of $\sigma$-algebras),
	we assume the \emph{modeled} past is finite.
	
	\begin{assumption}[Finite Past]\label{ass:finite_past}
		For all $i\in I$, $|\Anc_I(i)|<\infty$ is finite.
	\end{assumption}
	\begin{rmk}[Boundary Construction]
		It is possible, rather generally, to construct
		a "boundary" (\eg an initial state for a time-series\Slash{}Markov-process)
		that translates infinite past models to finite past models,
		see example \ref{example:timeseries}.		
		Assumption \ref{ass:finite_past} is primarily a statement about
		decomposition of the problem, where we discuss the causal aspects
		and detach the logically distinct question
		of choice of (non-unique) solutions into the choice
		of a (non-unique) initial state;
		see §\ref{sec:queries_beyond_basic}.
	\end{rmk}
	
	The model $\mathcal{M}$ defines a distribution, jointly over $\topSpace{X}_i$
	associated to $i\in I$
	by composing mechanisms.
	The result is
	given by mechanisms only, formally $\shallowDistr = \otimes^{\pi_I}_{i} f_{J(i)}$
	(using notation from §\ref{sec:structured_kernels}).
	We can avoid constructions on an infinite number of mechanisms by the following
	formulation.
	
	\begin{lemmaDef}[Shallow Distribution]
		\label{def:shallow_distr}
		Given acyclicity (Ass.\ \ref{ass:acyclic})
		and finite past (Ass.\ \ref{ass:finite_past}),
		then there exists a unique
		probability measure $P_\theta$ on
		$\prod_{i\in I} \topSpace{X}_i$
		with the following property:
		For all finite $I' \subset I$, the marginalization
		to $\topSpace{X}_{I'}$ satisfies
		for all tuples of measurable $B_i\in\sigmaAlgebraBorelIdx{X}{i}$ for $i\in I'$
		the equation
		\begin{equation*}
			\shallowDistr
			\Big(
			\prod_{i'\in I'} B_{i'}
			\Big)
			=
			\bigKernelCompound{
				\otimes^{\pi_I}_{i'\in \Anc_I(I')} f_{J(i')} 
			}
			\Big(
			\big(\prod_{i'\in I'} B_{i'}\big)
			\times				
			\big(\prod_{i\in \Anc_I(I')\setminus I'} \topSpace{X}_i \big)
			\Big)
			\txt.
		\end{equation*}
		We call $\shallowDistr$ the shallow distribution (Rmk.\ \ref{rmk:shallowness}).
		$\shallowDistr$ depends on the fixed $\mathcal{M}$,
		we make the parameters $\theta=\knownFunction$ explicit in the notation.
	\end{lemmaDef}
	
	While the model is not probabilistic, the observable
	world it describes is:
	
	\begin{definition}[Observable World]\label{def:obs_world}
		An observable world realizing the model $\mathcal{M}$ is
		a family of random variables
		$\{\anyVar_i: \Omega \rightarrow \topSpace{X}_i\}_{i\in I}$,
		such that there are measurable maps
		$f_i:\mathcal{X}_{\Pa_I(i)}\times[0,1] \rightarrow \topSpace{X}_i$
		and jointly independent random variables $\eta_i\sim U([0,1])$ called noises,
		with
		\begin{equation*}
			\anyVar_i = f_i(\anyVar_{\Pa_I(i)},\eta_i)
			\quad\text{and}\quad
			P(\{\anyVar_i\}_{i\in I}) = \shallowDistr
			\txt.
		\end{equation*}
		Mandated by the second equality, we will often simply write
		$\shallowDistr(\{\anyVar_i\}_{i\in I'})$ for $I'\subset I$.
	\end{definition}
	\begin{rmk}[Shallowness]
		\label{rmk:shallowness}
		The reason for calling $\shallowDistr$
		"shallow" is conceptually important:
		The way we model our data, for each $i\in I$, there is at most
		one observation.
		Thus each random variable $\anyVar_i$ is observed
		at most once.
		If we had more observations, we would instead extend
		$I$ and enhance the model's symmetry to capture this repetition.
		So $\shallowDistr$, while formally well-defined, does not allow
		for a statistical analysis without additional use of
		the model's symmetries.
		It becomes, in this language, an important aspect of identifiability,
		to ensure that any statistical analysis of the data
		must be able to draw from an
		invariant subset of data that becomes large as $N\rightarrow \infty$,
		thus is "deep" as opposed to "shallow".
	\end{rmk}
	
	\begin{lemma}[Observable World Existence]
		\label{lemma:obs_world_existence}
		Given acyclicity (Ass.\ \ref{ass:acyclic})
		and finite past (Ass.\ \ref{ass:finite_past}),
		an observable world exists.
		Furthermore,
		it is unique up to equality in distribution
		and it is locally Markov, \ie
		$\anyVar_i \independent \anyVar_K | \Pa_I(i)$ for any
		$K\subset I$ with $K\cap\Dec_I(i)=\emptyset$.
	\end{lemma}

	\begin{example}[IID-Models from SCMs]
		\label{example:model_iid_distr}
		Continuing example \ref{example:model_iid},
		which given an SCM $M$ constructed an associated model $\mathcal{M}$.
		For $\mathcal{M}$ and $\forall s\in\ISample$,
		a simple computation yields
		$\shallowDistr(\{\anyVar_{v,s}\}_{v\in\IVars}) = P^{\txt{obs}}$,
		where $P^{\txt{obs}}$ is the observed distribution of $M$.
		In particular the random model described by $M$
		is an observable world of $\mathcal{M}$.
		Similarly, for $\mathcal{M}^{\PearlDo}$,
		one obtains
		$\shallowDistr^{\PearlDo}(\{\anyVar_{v,s}\}_{v\in Y})=P(Y|\PearlDo(X=x))$.
		In the unified model, the first equality holds
		if $s\in\ISample$, the second equation holds if $s=*$.
		For details, see Lemma \ref{def:translate_model}.
	\end{example}

	\subsection{Observedness and Asymptotics} 

	So far, we have no notion of observed vs.\ hidden variables.
	We use a simple and flexible definition that only fixes which indices
	are observed and how limits are taken.
	
	\begin{definition}[Viewport]
		\label{def:viewport}
		Given a model, a viewport is given by an infinite
		subset $\mathcal{N}\subset\mathbb{N}$
		together with a mapping
		\begin{align*}
			\viewport: \mathcal{N} \rightarrow \FinSubSets(I),
			\halfquad
			N\mapsto \viewport(N)\halfquad
			&\txt{such that }
			|\viewport(N)|=N
			\txt{ and}\\
			&N \leq N' \in\mathcal{N}
			\Rightarrow
			\viewport(N)\subset\viewport(N')
			\txt.
		\end{align*}
		We write $\viewportEventual := \cup_{N\in\mathcal{N}} \viewport(N)\subset I$.
	\end{definition}
	\begin{assumption}[Independent Missingness]\label{ass:independent_missingness}
		The viewport, as defined in \ref{def:viewport},
		is not a random object. Thus it is independent of values taken by
		variables in the observable world.
		So our setup implicitly assumes (describes only) \emph{independent}
		missingness \citep{rubin1976inference}.
	\end{assumption}
	\begin{example}
		\label{example:viewport}
		Consider an IID model with three variables
		$\IVars=\{x,y,l\}$, where $x$ and $y$ are always observed,
		while $l$ is never observed.
		In this case we can choose $\mathcal{N} := 2\mathbb{N}$
		the even numbers (we could try to say what is observed
		first, $x$ or $y$, but we do not have to, this is why
		$\mathcal{N}\neq\mathbb{N}$ is allowed).
		Then for $N=2n\in \mathcal{N}$,
		define $\viewport(N) = \{x,y\} \times \{1,\ldots, n\}
		\subset \IVars \times \mathbb{N}$ (we chose
		$\ISample$ as the natural numbers here)
		to capture the described missingness pattern.
		In the asymptotic limit $n\rightarrow\infty$
		we observe $\viewportEventual=\{x,y\}\times\mathbb{N} \subset I$.
	\end{example}
	
	It only remains to combine observable worlds and viewports,
	and to formally encode the assumption that
	each variable of an observable world is observed at most once
	(Rmk.\ \ref{rmk:shallowness}).
	
	\begin{definition}[Realized World]\label{def:realized_world}
		Given an observable world and a viewport,
		we call a sample
		$\mathcal{D}^N(\omega) := \{\anyVar_i(\omega)\}_{i\in \viewport(N)}$
		a realized data-set\Slash{}observation
		and for $\mathcal{D}:=\cup_{N}\mathcal{D}^N(\omega)$ define
		$\realizedDistr:=
		\shallowDistr((\anyVar_i)_{i\in I}|(\anyVar_i)_{i\in \viewportEventual}
			=\mathcal{D})$
		a realized world (distribution).
	\end{definition}
	\begin{rmk}
		\label{rmk:realized_world_interpretation}
		While $\theta=\knownFunction$ are treated as
		fixed parameters, we may think of
		mechanisms $\observedFunction$ as latents.
		If we study properties of the realized world
		$\realizedDistr$, we study, by definition,
		$\shallowDistr$
		with each $\anyVar_i$ for $i\in\viewportEventual$ fixed to a value.
		Thus we input exactly one value (observation) for each
		element of the viewport.
		We use this observed information, plus the structure
		provided by symmetries, to draw conclusions about
		$\observedFunction$ (Thm.\ \ref{thm:extract_from_backdoor_complete}),
		then transform knowledge gained this way (§\ref{sec:graphical_operations})
		to answer queries (Thm.\ \ref{thm:query_id}) about $\realizedDistr$.
		The realized world is a formal means of asking
		questions under the restriction of only one observation
		per variable. It formally enforces the idea of shallowness
		(Rmk.\ \ref{rmk:shallowness}) on a mathematical object.
	\end{rmk}

	We cannot generally extrapolate to
	arguments we have never observed.
	The choice of the formal assumption made below
	is explained in detail in Rmk.\ \ref{rmk:kernels_uniqueness};
	this assumption is usually not made explicit, but formally necessary
	and implicitly made in other approaches.
	\begin{assumption}[Valid Support for Transfer]
		\label{ass:support}
		We assume that queries and other computations
		on kernels identified from observations
		remain within the observational support.
		Formally, any disintegration $\nu_x$
		of a observed joint distribution $\mu\otimes\nu_x$
		is only transferred to
		products $\xi\otimes\nu_x$ where $\xi\ll\mu$
		is dominated by $\mu$ (\ie $\forall$ measurable $B$:
		$\mu(B)=0 \Rightarrow \xi(B)=0$).
		For practical purposes it is often enough if
		the computed result can be interpreted as if this assumption
		were satisfied;
		\eg by a regularity
		assumption like smoothness of densities.
	\end{assumption}

	\section{Structured Kernels} 
	\label{sec:structured_kernels}
	
	We need formal objects to capture pieces of knowledge
	that will be obtained from observations in §\ref{sec:extraction_from_data}.
	Further, we need tools to transform these formal objects in
	a mathematically rigorous yet practically applicable calculus.
	
	\subsection{Preliminaries} 
	\label{sec:structured_kernels_preliminaries}

	Probability kernels are in wide spread use especially
	in probability theory and time-series statistics.
	We provide a very brief property-oriented overview
	of standard constructs here
	and a more detailed discussion connecting this language
	to SCMs in §\ref{apdx:kernels}.
	
	\begin{figure}[ht]
		\begin{minipage}{0.6\textwidth}
			\colorlet{blockcolor}{black!10!white}
\colorlet{flipcolor}{black!10!blue}

\begin{tikzpicture}
	\draw (0,0) node (mu)[rectangle, fill=blockcolor, rounded corners, inner sep=0.4em]
		{$\mu$};	
	\draw (1,0) node (nu)[rectangle, fill=blockcolor, rounded corners, inner sep=0.4em]
		{$\nu$};	
	\draw (3.5,0) node (prod)[rectangle, fill=blockcolor, rounded corners, inner sep=0.4em]
		{$\mu\otimes\nu$};
	
	\draw (2.25,0.75) node (opProd) {$\otimes$};
	\draw [-, thick] (mu) edge[bend left] (opProd);
	\draw [-, thick] (nu) edge[bend left] (opProd);
	\draw [->, thick] (opProd) edge[bend left] (prod);
	
	\draw [->, thick] (prod) edge[bend left]
		node[pos=0.4,above]{$\disintLeft$} (nu);
	\draw (prod.south) node(edgeStart){};
	\draw [->, thick] (edgeStart) edge[bend left]
		node[pos=0.5,below]{$\marginalizeRight$} (mu);
%
	\draw (5,0.3) node (anticausal)
		[anchor=west, rectangle, fill=blockcolor, rounded corners, inner sep=0.4em]
		{$(\nu\circ\mu) \otimes (\mu|\nu)$};	
	\draw [<->, thick,flipcolor] (prod) -- (anticausal);
	\draw (4.5,-0.1) node (mapping)[anchor=north west,align=left,flipcolor]{
			$X \times Y \rightarrow Y \times X$,\\
			$(x,y) \mapsto (y,x)$.
		};
\end{tikzpicture}
		\end{minipage}
		\hfill
		\begin{minipage}{0.35\textwidth}
			\caption{Basic kernel operations.
			Helpful analogues are:}
			\label{fig:basic_kernels}
			\vspace*{-1.85em}
			\begin{align*}
				\mu &\leftrightarrow P(X|Z=z)\\
				\nu &\leftrightarrow P(Y|X=x,Z=z)\\
				\mu\otimes\nu &\leftrightarrow P(X,Y|Z=z)
			\end{align*}
		\end{minipage}
	\end{figure}
	
	Some basic operations on kernels (illustrated in Fig.\ \ref{fig:basic_kernels})
	are: $\otimes$-products, given kernels $\mu$ from $\topSpace{Z}$ to
	$\topSpace{X}$ and $\nu$ from $\topSpace{Z}\times\topSpace{X}$
	to $\topSpace{Y}$,
	construct a kernel $\mu\otimes\nu$
	from $\topSpace{Z}$	to $\topSpace{X}\times\topSpace{Y}$;
	note the causal ordering (asymmetry in $\mu\leftrightarrow\nu$).
	Marginalizations (on the right), $\marginalizeRight$,
	drop a term at the end (in causal order),
	disintegrations (on the left), $\disintLeft$,
	turn a term at the beginning (in causal order)
	into an argument.
	If a suitable joint distribution including $\topSpace{Z}$ exists,
	these expressions can be thought of as conditional
	distributions.
	Reordering terms, for example switching $(x,y)\mapsto(y,x)$
	is a Borel-isomorphism -- given $P(X,Y)$ we know what $P(Y,X)$ is.
	We can marginalize after switching to get a composition $\nu\circ\mu$ of kernels.
	Hiding a cause $X$ of $Y$ does not change the distribution of $Y$,
	it still depends on the mechanisms at both $X$ and $Y$.
	Hence our expression $\nu\circ\mu$ in terms of mechanisms
	must be a \emph{binary} operator, \ie depend on both $\mu$ and $\nu$.
	Similar by disintegration after switching we get an "anti-causal"
	disintegration $(\mu|\nu)$.
	
	
	With $\xi := \mu\otimes\nu$ being a kernel itself,
	operations like
	$\xi \otimes \vartheta=\kernelCompound{\mu\otimes\nu}\otimes \vartheta$,
	and thus larger products are immediately defined (and turn out to be associative).
	It is however not immediately obvious, how reordering terms,
	marginalizations on the right and disintegrations on the left
	commute.
	In analogy to "$P(X,Y) = P(Y|X) P(X)$", which holds for tuples\Slash{}multi-variate
	$X$, $Y$ regardless of causal ordering, arbitrary marginalizations
	and disintegrations can be defined uniquely with these characterizing
	properties:
	\begin{lemmaDef}[Order Independent Operations]
		\label{lemma:order_independent_ops:mt}
		Given $\mu = \mu_1 \otimes \ldots \otimes \mu_n$ together with
		$L\subset \{1,\ldots,n\}$ and $C=\{1,\ldots,n\}\setminus L$,
		there are
		unique (Rmk.\ \ref{rmk:kernels_uniqueness})
		kernels $\marginalize{\mu}{L}$
		and $\disint_C(\mu)$
		such that $\forall$ measurable $B^C_k\in\sigmaAlgebraBorelIdx{X}{k}$
		where $B^C_k=\topSpace{X}_k$ if $k\in L$:
		\begin{align*}
			\marginalize{\mu}{L}(\textstyle\prod_{k\in C} B_i)
			&= \mu(\textstyle\prod_{k=1}^n B^C_k)\\
			\marginalize{\mu}{L} \otimes \disint_C(\mu)
			&= \mu\halfquad\txt{(up to reordering).}
		\end{align*}
	\end{lemmaDef}

	These operations all produce suitably unique
	results. Combining
	finitely many such operations again produces
	a unique result. We call such finite computations
	with these operations \emph{regular} functionals or simply regular.
	For generic kernels, without known internal structure like linearity,
	there do not seem to be any more well-defined (finite; excluding limits
	and fixed-points) computations. Thus we conjecture:
	\begin{conjecture}[Completeness of Regular Computation]
		\label{conjecture:regular_computation}
		If a functional uniquely computes a result from finitely many generic
		kernels then it is regular (Def.\ \ref{def:regular_functionals},
		Rmk.\ \ref{rmk:regular_limits}).
	\end{conjecture}

	\subsection{Causality and Sparsity} 
	\label{sec:causality_and_sparsity}

	So far, $\otimes$-products were defined with shared arguments in $\topSpace{Z}$
	plus the value-space of the left-hand-side as additional argument to the
	right-hand-side.
	Next, we define a product of many terms in "causal order",
	meaning any arguments that depend on another kernel's value are such
	that this parent is to the left.
	This product is defined inductively, term by term from left to right.
	To append the next kernel $\mu_x$,
	extend its argument space such that it
	includes all ancestors $a=(x,x')$;	
	simply by ignoring non-parents $\mu_{(x,x')} := \mu_x$.
	Here, "ancestors"
	are the kernels in the already defined product
	plus all shared arguments.
	
	Causal reasoning oftentimes relies precisely on
	the sparsity of parents among the ancestors.	
	As is common practice, both for tensor networks and causal
	modeling \citep{PearlBook}, we encode this sparsity in a graph
	(Def.\ \ref{def:structural_graph_simpl})
	attached to the product constructed above.
	Shared arguments for kernels are fundamentally
	different from their values, the definition
	thus uses two types of nodes (outer for shared arguments,
	inner for values).
	See also §\ref{apdx:structured_kernels}.
	
	\begin{definition}[Structural Graph]\label{def:structural_graph_simpl}
		A structural graph $\graphStructural$ is given by a finite set of nodes
		$\nodes := \nodesInner\,\dot\cup\,\nodesOuter$
		with $\nodesInner \cap \nodesOuter=\emptyset$,
		together with a set of directed edges
		$\edgesStructural\subset \nodes\times\nodesInner$,
		such that edges never point to an element of $\nodesOuter$.
	
		\emph{Alignment:}
		A model aligned structural graph is an acyclic (no directed cycles)
		structural graph $\graphStructural$ together
		with
		a probability kernel $\mu^n$ for each inner node $n\in\nodesInner$
		and
		a bijective assignment (cf.\ \ref{def:contraction})
		of the $\kappa$ arguments of $\mu^n$
		to parents of $n$:
		\begin{equation*}
			\pa^n: \{1, 2, \ldots, \kappa\}
			\xrightarrow{1:1} \Pa_{\graphStructural}(n)
			\txt.
		\end{equation*}
		Value-spaces of parents are subspaces of corresponding factors
		of the domain of children.
	\end{definition}
	
	The mappings $\pa^n$ tell us how to "wire up" each kernel $\mu$,
	\ie which shared arguments or values of kernels to the left of $\mu$
	are plugged into which of the $\kappa$ arguments of $\mu_x$.
	Previously, when defining a product in causal order,
	we implicitly assumed to know the order of parents in the argument $x$.
	Now we can formally clarify this point.
	To append the next term $\mu_x$ in causal order,
	reorder the (inductively obtained) product over ancestors
	such that $\pa^n(1), \ldots, \pa^n(\kappa)$ are the $\kappa$ first
	terms. Then add trivial dependence
	on all other non-parent ancestors $x'$, now all to the right of $x$,
	$\mu^n_{(x,x')}:=\mu^n_x$
	and write down the standard product.
	Thereby we obtain a structured, by a model-aligned graph, kernel.
	
	Later on, we want to describe
	knowledge about the model obtained from observations. This will
	require reasoning about latents.
	But, at the moment, we are talking about products of kernels,
	not causal models, observable worlds or viewports; there is
	no implied latent structure. Thus a subset $L\subset\nodesInner$,
	which will later be related to latents, is fixed as part of the
	definition.

	\begin{definition}[Structured Kernel]\label{def:structured_kernel_simp}
		Given a tuple $(\graphStructural,L)$, called a structural model below,
		of a
		model aligned
		structural graph $\graphStructural$ together with a subset
		$L\subset\nodesInner$,
		we define
		\begin{equation*}
			\mu(\graphStructural,L)_{x}
			\halfquad:=\halfquad
			\marginalize{
				\otimes_{n\in\nodesInner}^{\txt{wired}}
				\mu^n
			}{ L }_{x}
			\txt,
		\end{equation*}
		where $x=(x_i)_{i\in\nodesOuter}$
		is the tuple of shared arguments in $\nodesOuter$.
		The super-script "wired" on the product indicates that the product
		is constructed as described above, using the mappings $\pa^n$
		to "wire up" kernels.
		Given a set $\structureKernels$ of probability kernels,
		a $\structureKernels$-structured kernel is a kernel $\mu_x$ that can be
		written in the form $\mu_x = \mu(\graphStructural,L)_x$
		with $\forall n\in\nodesInner: \mu^n\in\structureKernels$.
	\end{definition}
	
	This definition implicitly depends on a choice $\pi$ of causal order.
	Because we usually need not be concerned about reordered results,
	in the sense explained in §\ref{sec:structured_kernels_preliminaries},
	we suppress $\pi$ in the notation (cf.\ \ref{notation:uniqueness_transposition}).
	Finally, the next subsection will require a formal notion of sub-graphs:
	
	\begin{definition}[Structural Subgraphs]
		\label{def:subgraphs_simp}
		A structural subgraph $\graphStructural^A\leq\graphStructural^B$
		of a structural graph $\graphStructural^B$
		is a structural graph $\graphStructural^A$ such that:
		\begin{enumerate}[label=(\roman*)]
			\item\label{def:subgraphs:nodesets}
			\emph{Node Sets:}
			$\nodes^A\subset\nodes^B$
			and	$\nodesInner^A\subset\nodesInner^B$.
			\item\label{def:subgraphs:edges}
			\emph{Edge Sets:}
			edges are exactly those in $\graphStructural^B$
			from nodes in $\nodes^A$
			to inner nodes of $\nodesInner^A$.
			\item\label{def:subgraphs:inner_parents}
			\emph{Inner Parents:}
			$\forall n\in\nodesInner^A$:
			if $p\in\Pa_{\graphStructural^B}(n)$,
			then $p\in\nodes^A$.
			\item\label{def:subgraphs:alignment}
			\emph{Alignment:} for model-aligned graphs
			$\forall n\in\nodesInner^A$:
			$(\mu^A)^n = (\mu^B)^n$.
		\end{enumerate}
	\end{definition}
	
	\subsection{Graphical Operations}
	\label{sec:graphical_operations}

	\begin{figure}[ht]
		\begin{minipage}{\textwidth}
			\centering
			\colorlet{dblue}{blue!50!black}
\colorlet{dgreen}{green!50!black}
\colorlet{dyellow}{yellow!50!black}
\colorlet{dpurple}{purple!50!black}
\colorlet{c1}{dblue!20!white}
\colorlet{c2}{dgreen!20!white}
\colorlet{c3}{dyellow!20!white}
\colorlet{c4}{dpurple!20!white}
\colorlet{dgray}{gray!50!black}
\colorlet{c0}{dgray!20!white}
\colorlet{b0}{c0!40!white}

\usetikzlibrary{decorations.pathmorphing}

\begin{tikzpicture}
	\begin{scope}[xshift=-0.5cm,yshift=-0.45cm, xscale=1.1, yscale=0.75]
	\draw[dgray, line width=0.2em, fill=white] (0,0) -- (4.5,0)
			-- (5.5,-0.5) -- (12,-0.5)
			-- (12, 3.5) -- (5.5, 3.5)
			-- (4.5, 3) -- (0,3) -- cycle;
	\end{scope}
	\draw[dgray, line width=0.2em, fill=white, opacity=0.85]
		(5.5, -0.5)  -- (5.5,1.8) -- (10.2,1.8) -- (10.2,-0.5)
		-- (12.5,-1.25)
		-- (12.5,-3.4) -- (6,-3.4) -- (6,-1.25) -- cycle;
		
	\draw[dgray, line width=0.2em, fill=white, opacity=0.85]
		(0, -3.25) -- (5, -3.25) -- (5.725,-3.15) -- (9.25,-3.15)
		-- (9.25, -1.25) -- (5.725,-1.25) -- (5,-1.1) -- (0,-1.1) -- cycle;
		
	\draw[->, very thick, line join=round,
	decorate, decoration={
		zigzag,
		segment length=8,
		amplitude=.9,post=lineto,
		post length=2pt
	}] (4.3,0.25) -- (5.25,0.25) node[above,pos=0.5, align=center]
	{\colorbox{white}{c-components}\\(L\,\ref{lemma:c_components:mt})};
	\draw[->, very thick, line join=round,
	decorate, decoration={
		zigzag,
		segment length=8,
		amplitude=.9,post=lineto,
		post length=2pt
	}] (8.75,-0.4) -- (9,-1.1) node[left,pos=0.5, align=center,xshift=-0.3em,yshift=-0.1em]
	{simplify\\(L\,\ref{def:simplify_simp})};
	\draw[->, very thick, line join=round,
	decorate, decoration={
		zigzag,
		segment length=8,
		amplitude=.9,post=lineto,
		post length=2pt
	}] (5.85,-2.5) -- (4.85,-2.5) node[above,pos=0.5, xshift=1.4em, align=center]
	{\colorbox{white}{c-components}\\(L\,\ref{lemma:c_components:mt})};
	\draw[->, very thick, line join=round,
	decorate, decoration={
		zigzag,
		segment length=8,
		amplitude=.9,post=lineto,
		post length=2pt
	}] (0.5,-4.75) -- (1.5,-4.75) node[above,pos=0.5,yshift=0.2em, align=center]
	{$2\times$glue\\(L.\,\ref{lemma:glue:mt})};
	
	\fill[rounded corners, fill=c0] (2.5,-5.75) rectangle (3.5,-3.65);
	\draw[dgray] (2.5,-3.65) node[anchor=north east] {\textbf{A}};	
	\fill[rounded corners, fill=c0] (5.625,-5.75) rectangle (6.625,-3.65);
	\draw[dgray] (5.625,-3.65) node[anchor=north east] {\textbf{B}};
	
	\begin{scope}[xscale=1.25]
		\draw (0,0) node (W)[rectangle, fill=c1, rounded corners, inner sep=0.4em,
		draw=dblue,	thick]
		{$W$};
		\draw (1,0) node (Z)[rectangle, fill=c2, rounded corners, inner sep=0.4em,
		draw=dgreen, thick]
		{$Z$};
		\draw (2,0) node (X)[rectangle, fill=c1, rounded corners, inner sep=0.4em,
		draw=dblue,	thick]
		{$X$};
		\draw (3,0) node (Y)[rectangle, fill=c1, rounded corners, inner sep=0.4em,
		draw=dblue,	thick]
		{$Y$};
		
		\draw [->, very thick] (W) -- (Z);	
		\draw [->, very thick] (Z) -- (X);	
		\draw [->, very thick] (X) -- (Y);
		
		\draw [<->, very thick,dblue] (W) edge[bend left=2.5cm]
		node[pos=0.6,above][rectangle, fill=c1, rounded corners, inner sep=0.4em,
		draw=dblue, thick, dashed]{$L_2$} (Y);
		\draw [<->, very thick, dblue] (W) edge[bend left=1.75cm]
		node[pos=0.5,above][rectangle, fill=c1, rounded corners, inner sep=0.4em,
		draw=dblue, thick, dashed]{$L_1$} (X);
	\end{scope}
	\begin{scope}[xshift=6.0cm]		
		\begin{scope}[xscale=1.25]
			\draw (0,0) node (W)[rectangle, fill=c1, rounded corners, inner sep=0.4em,
			draw=dblue,	thick]
			{$W$};
			\draw (1,0) node (Z)[circle,draw, very thick, inner sep=0.2em]
			{};
			\draw (2,0) node (X)[rectangle, fill=c1, rounded corners, inner sep=0.4em,
			draw=dblue,	thick]
			{$X$};
			\draw (3,0) node (Y)[rectangle, fill=c1, rounded corners, inner sep=0.4em,
			draw=dblue,	thick]
			{$Y$};
			
			\draw [->, very thick] (Z) -- (X);	
			\draw [->, very thick] (X) -- (Y);
			
			\draw [<->, very thick,dblue] (W) edge[bend left=2.5cm]
			node[pos=0.6,above][rectangle, fill=c1, rounded corners, inner sep=0.4em,
			draw=dblue, thick, dashed]{$L_2$} (Y);
			\draw [<->, very thick, dblue] (W) edge[bend left=1.75cm]
			node[pos=0.5,above][rectangle, fill=c1, rounded corners, inner sep=0.4em,
			draw=dblue, thick, dashed]{$L_1$} (X);
		\end{scope}		
		\draw (4.5,0.5) node{\textbf{\&}};
		\begin{scope}[xshift=5cm, yshift=0.5cm, xscale=1.25]			
			\draw (0,0) node (W)[circle,draw, very thick, inner sep=0.2em]
			{};
			\draw (1,0) node (Z)[rectangle, fill=c2, rounded corners, inner sep=0.4em,
			draw=dgreen, thick]
			{$Z$};
			
			\draw [->, very thick] (W) -- (Z);	
		\end{scope}
	\end{scope}
	\begin{scope}[xshift=7.25cm, yshift=-2.75cm]
		\begin{scope}[xshift=0.25cm, xscale=1.25]		
			\draw (-1,0) node (Z)[circle,draw, very thick, inner sep=0.2em]
			{};
			\draw (0,0) node (X)[rectangle, fill=c1, rounded corners, inner sep=0.4em,
			draw=dblue,	thick] {$X$};
			\draw (0,1) node (L1)[rectangle, fill=c1, rounded corners, inner sep=0.4em,
			draw=dblue, thick, dashed] {$L_1$};
			\draw (1,0) node (Y)[rectangle, fill=c3, rounded corners, inner sep=0.4em,
			draw=dyellow, thick] {$Y$};
			\draw (1,1) node (L2)[rectangle, fill=c3, rounded corners, inner sep=0.4em,
			draw=dyellow, thick, dashed] {$L_2$};
			
			\draw [->, very thick] (Z) -- (X);
			\draw [->, very thick] (X) -- (Y);	
			\draw [->, very thick] (L1) -- (X);	
			\draw [->, very thick] (L2) -- (Y);	
		\end{scope}
		\draw (2.5,0.5) node{\textbf{\&}};
		\begin{scope}[xshift=4cm, xscale=1.25]		
			\draw (0,0) node (W)[rectangle, fill=c1, rounded corners, inner sep=0.4em,
			draw=dblue,	thick]
			{$W$};
			\draw (0.5,1) node (L1)[rectangle, fill=c1, rounded corners, inner sep=0.4em,
			draw=dblue, thick, dashed] {$L_1$};
			\draw (-0.5,1) node (L2)[rectangle, fill=c1, rounded corners, inner sep=0.4em,
			draw=dblue, thick, dashed] {$L_2$};
			
			\draw [->, very thick] (L1) -- (W);	
			\draw [->, very thick] (L2) -- (W);	
		\end{scope}		
	\end{scope}	
	\begin{scope}[xshift=1.25cm, yshift=-2.75cm]
		\begin{scope}[xshift=0.25cm, xscale=1.25]		
				\draw (-1,0) node (Z)[circle,draw, very thick, inner sep=0.2em]
			{};
			\draw (0,0) node (X)[rectangle, fill=c1, rounded corners, inner sep=0.4em,
			draw=dblue,	thick] {$X$};
			\draw (0,1) node (L1)[rectangle, fill=c1, rounded corners, inner sep=0.4em,
			draw=dblue, thick, dashed] {$L_1$};			
			
			\draw[dgray] (L1.north west) node[anchor=north east] {\textbf{A}};
			
			\draw [->, very thick] (Z) -- (X);
			\draw [->, very thick] (L1) -- (X);	
		\end{scope}
		\draw (1.25,0.5) node{\textbf{\&}};
		\begin{scope}[xshift=3cm, xscale=1.25]		
			\draw (-1,0) node (X)[circle,draw, very thick, inner sep=0.2em]
			{};
			\draw (0,0) node (Y)[rectangle, fill=c3, rounded corners, inner sep=0.4em,
			draw=dyellow,	thick] {$Y$};
			\draw (0,1) node (L2)[rectangle, fill=c3, rounded corners, inner sep=0.4em,
			draw=dyellow, thick, dashed] {$L_2$};
						
			\draw[dgray] (L2.north west) node[anchor=north east] {\textbf{B}};
			
			\draw [->, very thick] (X) -- (Y);
			\draw [->, very thick] (L2) -- (Y);	
		\end{scope}		
	\end{scope}
	\begin{scope}[xshift=3cm, yshift=-5.2cm, xscale=1.25]	
		\draw (-1,0) node (Z)[circle,draw, very thick, inner sep=0.2em]
		{};
		\draw (0,0) node (X)[rectangle, fill=c1, rounded corners, inner sep=0.4em,
		draw=dblue,	thick] {$X$};
		\draw (0,1) node (L1)[rectangle, fill=c1, rounded corners, inner sep=0.4em,
		draw=dblue, thick, dashed] {$L_1$};
		
		\draw [->, very thick] (Z) -- (X);
		\draw [->, very thick] (L1) -- (X);	
		
		\draw (2.5,0) node (Y)[rectangle, fill=c3, rounded corners, inner sep=0.4em,
		draw=dyellow,	thick] {$Y$};
		\draw (2.5,1) node (L2)[rectangle, fill=c3, rounded corners, inner sep=0.4em,
		draw=dyellow, thick, dashed] {$L_2$};
		
		\draw (1.25,0) node (Xsoft)[rectangle, fill=c4, inner sep=0.4em,
		draw=dpurple,	thick] {$\tilde{X}_{\txt{soft}}$};
		
		\draw[dgray] (Xsoft.north west) node[anchor=south west]{\textbf{C}};
		
		\draw [->, very thick] (X) -- (Xsoft);
		\draw [->, very thick] (Xsoft) -- (Y);
		\draw [->, very thick] (L2) -- (Y);	
	\end{scope}
\end{tikzpicture}
		\end{minipage}\\[-0.5em]
		\caption{Examples for graphical operations.
		All shown structured kernels can be computed from
		Pearl's "napkin"-graph (top-left) and the (know) soft-intervention
		$\tilde{X}_{\txt{soft}}$ (a single inner node plus a
		single outer node for its argument).
		The computation is step-by-step, for example
		$\{W\}\neq\Dec(\{W\})$ in the initial graph,
		but after taking c-components $\{W\}=\Dec(\{W\})$ 
		in the graph at the right of the first $\rightsquigarrow$,
		so Lemma \ref{def:simplify_simp} applies in the next step.
		Inner nodes are colored by
		structural c-component, labeled (with a capital roman letter)
		by mechanism,
		hidden nodes $L_i$ have dashed border, outer nodes
		are draw as circles.}
		\label{fig:graphical_ops}
	\end{figure}
	
	We have seen certain regular computations ($\otimes$-products,
	$\marginalizeRight$, $\disintLeft$, reordering)
	on kernels,
	and we have structured kernels by graphs.
	Next we show that many (maybe all) relevant regular computations on
	kernels can be represented by simple and intuitive operations
	performed on \emph{model-aligned} graphs, see Fig.\ \ref{fig:graphical_ops}.
	The first plausible candidate are sub-graphs.
	Given $\mu(\graphStructural,L)$ we can
	compute $\mu(\graphStructural', L')$ for trivial cases
	where we simply discard hidden structure irrelevant to observed nodes:
	\begin{lemmaDef}[Simplify]
		\label{def:simplify_simp}
		Given a structured kernel $\mu(\graphStructural, L)$
		and a sub-set $B \subset \nodesInner$
		with $\Dec_\graphStructural(B) \subset B$,
		then
		$\mu(\graphStructural',L \setminus (B\cap L))$
		can be computed as a regular functional of
		$\mu(\graphStructural, L)$,
		where we call the sub-graph
		$\graphStructural'\leq \graphStructural$
		with inner nodes $\nodesInner'=\nodesInner\setminus B$
		and outer nodes $\nodesOuter'=\Pa_{\graphStructural}
		(\nodesInner')\setminus\nodesInner'$
		a simplification.
		Restricting to $B\subset L$,
		there is a unique maximal $B\subset L$
		and we call the corresponding
		(minimal) $\mu(\graphStructural',L\setminus B)$ simplified.
	\end{lemmaDef}
	
	Also a class of sub-graphs well known from the IID-case \citep{tian2002general}
	appear: c-components.
	\begin{definition}[C-Components]
		\label{def:c_components_local:mt}
		Given a structural model $(\graphStructural, L)$,
		we call $l\in L$ a hidden confounder of
		$y,w\in\nodesInner\setminus L$
		if there are directed paths $\gamma_y$ from $l$ to $y$ and $\gamma_w$
		from $l$ to $w$ such that all non-endpoint nodes of $\gamma_y, \gamma_w$
		are in $L$.
		On $\nodesInner\setminus L$ define a relation
		$y \sim w :\Leftrightarrow \exists$ hidden confounder $l$ of $y, w$.
		Define $\sim_L$ as the equivalence relation generated by $\sim$.
		We call equivalence-classes of $\sim_L$ on $\nodesInner\setminus L$
		c-components.
		
		We can consistently extend this notion to $L$
		by defining $l\approx y$ for $y\in \nodesInner\setminus L$ if there
		exists a directed path $\gamma_y$ from $l$ to $y$
		such that all non-endpoint nodes of $\gamma_y$ are in $L$.
		We call a sub-graph $\graphStructural^c\leq \graphStructural$
		with $\nodesInner^c$
		an equivalence-class of $\approx_L$ and
		$\nodesOuter^c = \Pa_{\graphStructural}(\nodesInner^c)\setminus\nodesInner^c$
		together with $L^c=L\cap\nodesInner^c=L\cap\nodes^c$
		a structural c-component.
	\end{definition}
	\begin{lemma}[Computation of C-Components]
		\label{lemma:c_components:mt}
		Given a structural c-component $\graphStructural^c\leq \graphStructural$, then
		there is a regular functional computing
		$\mu(\graphStructural^c,L^c)$ from $\mu(\graphStructural,L)$.
	\end{lemma}
	
	On the other hand, we can regularly compute (larger) products from their
	(smaller) constituents.
	A graphical analogue to obtain a larger graph from smaller ones
	is a gluing operation:
	
	\begin{lemma}[Gluing]
		\label{lemma:glue:mt}
		Given a structural kernel $\mu(\graphStructural, L)$,
		and subgraphs $\graphStructural^A, \graphStructural^B \leq \graphStructural$,
		such that
		each structural c-component $\graphStructural^c \leq \graphStructural$
		is a sub-graph $\graphStructural^c\leq \graphStructural^A$
		or $\graphStructural^c\leq \graphStructural^B$,		
		then there is a regular functional computing
		$\mu(\graphStructural,L)$
		from $\mu(\graphStructural^A,L\cap\nodes^A)$ and
		$\mu(\graphStructural^B,L\cap\nodes^B)$.
	\end{lemma}
	
	Our notion of observed and hidden variables is flexible,
	and it is possible to end up in situations where
	structural models on a graphical overlap differ by their latent sets.
	Almost by accident, we end up with a (non-parametric) causal
	missing-value imputation technique:
	
	\begin{cor}[Revealing]
		\label{cor:revealing_simp}
		Given a structural kernel $\mu(\graphStructural, L)$,
		a structural c-component
		$\graphStructural^c \leq \graphStructural$
		and $L'\subset L^c$,
		then there is a regular functional computing
		$\mu(\graphStructural,(L\setminus L^c) \cup L')$
		from $\mu(\graphStructural,L)$
		and $\mu(\graphStructural^c,L')$.
	\end{cor}
	
	Provided knowledge of a set of structured kernels $\knowledgeSet$,
	we can systematically find regular computations to learn
	a (unknown) structured kernel of interest.
	It will be practically useful to order gluing operations
	to the end, \ie to first decompose
	elements of $\knowledgeSet$
	into smaller pieces (which can be done algorithmically, see
	Algo.\ \ref{algo:decomp}), then glue the kernel of interest from
	these pieces (which can also be done algorithmically, see Algo.\ \ref{algo:svs}).
	Finally, for details on how the results of
	these operations can be computed in practice,
	see Rmk.\ \ref{rmk:compute_atoms_practice}, \ref{rmk:compute_gluing}.
	It seems plausible that all regular computations of \emph{structured} results
	(see §\ref{sec:queries_beyond_basic} and nested queries however)
	can be achieved by these graphical operations.
	\begin{conjecture}[Graphical Regular Computation]
		\label{conjecture:regular_graphical}
		If a structured kernel $\mu(\graphStructural,L)$ can be
		regularly computed from a set $\knowledgeSet$ of
		structured kernels,
		then there is a finite sequence of graphical operations
		computing $\mu(\graphStructural,L)$ from $\knowledgeSet$.
	\end{conjecture}
	Our main
	concern regarding this conjecture
	are potential generalizations of the revealing operation
	(Cor.\ \ref{cor:revealing_simp}). Revealing currently is a simple corollary
	of other operations, but might be possible more generally
	(see §\ref{apdx:reveal})
	
	Finally we conclude this section by the observation that at
	this point we can already identify interventions in single-context IID setups:
	\begin{example}[IID Interventions]
		\label{example:iid_id_from_structured_only}
		In example \ref{example:model_iid_distr},
		the relations $\shallowDistr = P^{\txt{obs}}$
		and $\shallowDistr^{\PearlDo}=P(\dots|\PearlDo(X=x))$ were found.
		By definition of $\shallowDistr=\otimes_{v\in\IVars} f_v$,
		and assuming we observe $\IVars \setminus L$ (the latents $L$ are fixed
		uniformly, for all samples identically), then
		$P^{\txt{obs}}=\mu(\graphClassical,L)$, where $\graphClassical$
		is the structural graph with $\nodesInner=\nodes=\IVars$
		aligned to $\mu^v = f_v$. Here $f_v := g_v(\pa,\cdot)_*P_{\eta^v}$
		as in example \ref{example:model_iid}.
		Thus, if we know $P^{\txt{obs}}$, we can decompose
		$\mu(\graphClassical,L)$ as illustrated in Fig.\ \ref{fig:graphical_ops},
		\ie we can compute, from $P^{\txt{obs}}$, the smaller pieces
		(like $Y_{x,[l_2]}\circ L_2$).
		
		With $\shallowDistr^{\PearlDo}=P(\dots|\PearlDo(X=x))$,
		to compute $P(Y|\PearlDo(X=x))$, it is enough
		to compute $\shallowDistr^{\PearlDo}=\mu(\graphClassicalDo, L')$,
		where $\graphClassicalDo$ is $\graphClassical$ with the nodes in $X$
		aligned to $\delta[x_i]\in\knownFunction$.
		Importantly $L' := \IVars\setminus Y$ can be large, we are \emph{only}
		interested in the distribution of elements of $Y$.
		Thus, it may be possible to simplify (Lemma \ref{lemma:simplify})
		$\mu(\graphClassicalDo, L')$ and glue a \emph{smaller}
		graph from the pieces computed above (marginal identification is
		easier than joint identification).
		As illustrated in Fig.\ \ref{fig:graphical_ops}, the intervention
		can be a soft intervention.
		While extremely simple -- we only have to apply graphical operations
		to decompose $(\graphClassical,L)$, then glue
		a simplification of $(\graphClassicalDo,L')$ --
		this procedure is complete in the case of IID do-interventions
		Cor.\ \ref{cor:iid:complete_single_context_uncond}, see
		§\ref{apdx:iid:single_context_empty_cond} (for conditioning,
		§\ref{apdx:iid_conditional}).
	\end{example}

	\section{Extraction from Data}
	\label{sec:extraction_from_data}

	\RenewDocumentCommand{\nodesInnerProper}{}{\nodesInner}
	\RenewDocumentCommand{\nodesOuterProper}{}{\nodesOuter}
	\RenewDocumentCommand{\nodesOuterFixed}{}{\nodes_{\txt{pinned}}}

	Having introduced a language of structured kernels,
	we still need to connect this notion to our model
	and to analyze which structured kernels can be learned from data.
	The reader interested in the reasons behind
	formal choices made in the definitions of this section
	may want to analyze
	the proof of Thm.\ \ref{thm:extract_from_backdoor_complete},
	especially step 1,
	in §\ref{apdx:id_from_embeddings}.
	We focus on the single-level case (cf.\ §\ref{sec:multi_level}),
	not to be confused with the single-\emph{context} \citep{CD-NOD,JCI,JPCMCI} case
	(see also §\ref{apdx:translate_mz_transport}).
	A genuinely multi-\emph{level} example is given in
	§\ref{apdx:multi_level_examples}.
	All results and proofs are given for the multi-level case in the appendix.
	Multi-level statistics (see \eg \citep{Gelman2006}) have not played
	a major role in causal literature so far and the specialization
	to a single level allows for a simplified notation.
	We start with a simple example to help elucidate how
	the present section generalizes beyond example
	\ref{example:iid_id_from_structured_only}.
	We will return to this question in example
	\ref{example:iid_extraction_direct}.
	\begin{example}[Necessity of Extraction]
		\label{example:extraction_frontdoor}
		Consider the following setup inspired by \citep{JPCMCI}.
		We have multiple time-series on $t\in T$ recorded
		at multiple sites\Slash{}locations $s\in S$
		in space. There are three observed variables
		$X \rightarrow Z^s \rightarrow Y^t$. Their mechanisms
		are in part shared, $X$ is the same at all times on all sites,
		$Z^s$ is the same at all times, but differs by site,
		$Y^t$ is the same at all sites, but changes over time.
		Transport \citep{bareinboim2012TransportCompleteness},
		JCI \citep{JCI} and also \citep{JPCMCI} would introduce
		context-variables $C_S$ (a parent of $Z$) and $C_T$
		(a parent of $Y$).
		If we want to estimate $P(Y|\PearlDo(X=x))$,
		the extended graph (including contexts) suggests
		$P(Y|\PearlDo(X=x))=P(Y|X=x)$, or maybe
		$P(Y|\PearlDo(X=x),C_T=t)=P(Y|X=x,C_T=t)$
		if we want to make explicit that we ask for a specific time $t$.
		Maybe even 
		$P(Y|\PearlDo(X=x),C_T=t, C_S=s)=P(Y|X=x,C_T=t,C_S=s)$?
		It is not formally clear which of these should be used.
		Even worse, the last (and from a practical perspective most useful)
		one $P(Y|\PearlDo(X=x),C_T=t, C_S=s)$ is seemingly identified from
		$P(Y|X=x,C_T=t,C_S=s)$; but we have precisely one data-point for this
		distribution, no meaningful estimation is possible.
		
		By separating extraction, including a formal account to
		availability of data, from computation we find:
		We can, for any $s\in S$,
		learn $\circ \rightarrow Z^s$ from data at $s$, with $|T|$ data-points.
		We can, for any $t\in T$,
		learn $\circ \rightarrow Y^t$ from data at $t$, with $|S|$ data-points.	
		If both $|T|$ and $|S|$ are large, reasonable estimation of both is possible.	
		Then we can compute (by gluing, Lemma \ref{lemma:glue:mt})
		$X_{\PearlDo} \rightarrow Z^s \rightarrow Y^t$.
		In particular $P(Y|\PearlDo(X=x),C_T=t, C_S=s)$ is identifiable
		(in a meaningful sense, with data), just not via
		$P(Y|X=x,C_T=t,C_S=s)$, but by a "frontdoor-like" argument instead.
	\end{example}

	\subsection{Families of Embeddings} 

	We relate model-aligned graphs to individual local structures in
	the model (and $I$-graph) by a suitable notion of embeddings.
	For statistical reasoning, we need certain (for example sufficiently
	non-degenerate) repeated observations, thus we study
	families of embeddings additionally to individual embeddings.
	In the single-level case we can directly embed
	structured graphs, we nevertheless call them "local" graphs
	(and denote them $\graph$ instead of $\graphStructural$)
	in this section
	to emphasize that conceptually local graphs are not the same as structured graphs.

	\begin{figure}[ht]
		\colorlet{dblue}{blue!50!black}
\colorlet{dgreen}{green!50!black}
\colorlet{dyellow}{yellow!50!black}
\colorlet{dpurple}{purple!50!black}
\colorlet{dorange}{orange!50!black}
\colorlet{c4}{dblue!20!white}
\colorlet{c3}{dgreen!20!white}
\colorlet{c1}{dyellow!20!white}
\colorlet{c2}{dpurple!20!white}
\colorlet{c5}{dorange!20!white}

\pgfdeclarelayer{bg}
\pgfsetlayers{bg,main}

\begin{tikzpicture}
	\draw (0,1) node (X0) {};
	\draw (0,0) node (Y0) {};
	\foreach \t  in {1,...,12}{
		\draw (\t,1) node (X\t) {$X$};
		\draw (\t,0) node (Y\t) {$Y$};
		
		\pgfmathsetmacro{\prev}{\t-1}
		\draw[->, very thick] (X\prev) -- (X\t);
		\draw[->, very thick] (Y\prev) -- (Y\t);
		\draw[->, very thick] (X\prev) + (0,-0.4em) -- (Y\t);
	}
	\draw (-0.5,0.5) node {$\cdots$};
	\draw (12.7,0.5) node {$\cdots$};

	\draw[dyellow] (2,1.5) node {OK};
	\draw[dpurple] (5,1.5) node {not injective};
	\draw[dgreen] (7.5,1.5) node {not (i)};
	\draw[dblue] (9.5,1.5) node {not (ii)};
	\draw[dorange] (11.5,1.5) node {not (iii)};
	
	\begin{scope}[yshift=-2cm]
		\draw (1,1) node[circle,draw, very thick, inner sep=0.2em](x1) {};
		\draw (2,1) node[circle,draw, very thick, inner sep=0.2em](x2) {};
		\draw (1,0) node[circle,draw, very thick, inner sep=0.2em](y1) {};
		
		\draw (3,1) node(x3) {$X$};
		\draw (2,0) node(y2) {$Y$};
		\draw (3,0) node(y3) {$Y$};
		
		\draw[->, very thick] (x2) -- (x3);
		\draw[->, very thick] (x1) -- (y2);
		\draw[->, very thick] (x2) -- (y3);
		\draw[->, very thick] (y1) -- (y2);
		\draw[->, very thick] (y2) -- (y3);
		
		\begin{pgfonlayer}{bg}			
			\draw[c1, line width=0.1cm] (x1) edge[bend right] (X1);
			\draw[c1, line width=0.1cm] (x2) edge[bend right] (X2);
			\draw[c1, line width=0.1cm] (x3) edge[bend right] (X3);
			\draw[c1, line width=0.1cm] (y1) edge[bend right] (Y1);
			\draw[c1, line width=0.1cm] (y2) edge[bend right] (Y2);
			\draw[c1, line width=0.1cm] (y3) edge[bend right] (Y3);			
		\end{pgfonlayer}
	\end{scope}
	\begin{scope}[xshift=3cm, yshift=-2cm]
		\draw (1,1) node[circle,draw, very thick, inner sep=0.2em](x1) {};
		\draw (2,1) node[circle,draw, very thick, inner sep=0.2em](x2) {};
		\draw (1,0) node[circle,draw, very thick, inner sep=0.2em](y1) {};
		
		\draw (3,1) node(x3) {$X$};
		\draw (2,0) node(y2) {$Y$};
		\draw (3,0) node(y3) {$Y$};
		
		\draw[->, very thick] (x2) -- (x3);
		\draw[->, very thick] (x1) -- (y2);
		\draw[->, very thick] (x2) -- (y3);
		\draw[->, very thick] (y1) -- (y2);
		\draw[->, very thick] (y2) -- (y3);
		
		\begin{pgfonlayer}{bg}
			\draw[c2, line width=0.1cm] (x1) edge[bend right] (X4);
			\draw[c2, line width=0.1cm] (x2) edge[bend right] (X5);
			\draw[c2, line width=0.1cm] (x3) edge[bend right] (X5);
			\draw[c2, line width=0.1cm] (y1) edge[bend right] (Y4);
			\draw[c2, line width=0.1cm] (y2) edge[bend right] (Y5);
			\draw[c2, line width=0.1cm] (y3) edge[bend right] (Y6);		
		\end{pgfonlayer}
	\end{scope}
	\begin{scope}[xshift=6cm, yshift=-2cm]
		\draw (1,1) node[circle,draw, very thick, inner sep=0.2em](x1) {};
		
		\draw (2,0) node(y2) {$Y$};
		
		\draw[->, very thick] (x1) -- (y2);
		
		\begin{pgfonlayer}{bg}
			\draw[c3, line width=0.1cm] (x1) edge[bend right] (X7);
			\draw[c3, line width=0.1cm] (y2) edge[bend right] (Y8);	
		\end{pgfonlayer}
	\end{scope}
	\begin{scope}[xshift=8cm, yshift=-2cm]
		\draw (1,1) node[circle,draw, very thick, inner sep=0.2em](x1) {};
		\draw (1,0) node[circle,draw, very thick, inner sep=0.2em](y1) {};
		
		\draw (2,1) node(x2) {$X$};
		\draw (2,0) node(y2) {$Y$};
		
		\draw[->, very thick] (x1) -- (x2);
		\draw[->, very thick] (y1) -- (y2);
		
		\begin{pgfonlayer}{bg}
			\draw[c4, line width=0.1cm] (x1) edge[bend right] (X9);
			\draw[c4, line width=0.1cm] (y2) edge[bend right] (Y10);	
			\draw[c4, line width=0.1cm] (x2) edge[bend right] (X10);
			\draw[c4, line width=0.1cm] (y1) edge[bend right] (Y9);	
		\end{pgfonlayer}
	\end{scope}
	
	\begin{scope}[xshift=10cm, yshift=-2cm]
		\draw (1,0) node[circle,draw, very thick, inner sep=0.2em](y1) {};
		
		\draw (2,0) node(y2) {$Y$};
		
		\draw[->, very thick] (y1) -- (y2);
		
		\begin{pgfonlayer}{bg}
			\draw[c5, line width=0.1cm] (y1) edge[bend right] (X11);
			\draw[c5, line width=0.1cm] (y2) edge[bend right] (X12);
		\end{pgfonlayer}
	\end{scope}
\end{tikzpicture}\\[-1.5em]
		\caption{Examples illustrating conditions on local graph embeddings.
		Alignment is indicated by capital letters,
		outer nodes are drawn as circles. On the top
		is an $I$-graph (in this example a time series), below are
		multiple local graphs. The nodes of local graphs are mapped to the $I$-graph,
		colored curved lines, but only the first of the shown ones is embedded
		(Def.\ \ref{def:local_graph_embedding}).}
		\label{fig:embeddings}
	\end{figure}
	
	\begin{definition}[Local Graph Embedding]\label{def:local_graph_embedding}
		An embedding of
		a model aligned local graph $\graph$ in the index set $I$,
		denoted by
		$\psi:\graph \hookrightarrow I$, is an injective mapping
		$\psi:\nodes \hookrightarrow I$
		such that the following conditions are satisfied
		(see Fig.\ \ref{fig:embeddings}):
		\begin{enumerate}[label=(\roman*)]
			\item\label{def:local_graph_embedding:inner_parents_incl}
			Inner-Node $I$-Parents are Included:
			If $n\in \nodesInnerProper$, then
			$\forall i' \in \Pa_I(\psi(n))$, $\exists n' \in \nodes$
			such that $i' = \psi(n')$.
			\item\label{def:local_graph_embedding:proper_edges_a}
			Edges equal $I$-Graph Edges:
			For $n' \in \nodesInnerProper$,
			there is a edge $n \rightarrow n'$
			in $\graph$ if and only if
			$\psi(n) \in \Pa_I(\psi(n'))$.
			\item\label{def:local_graph_embedding:applicable}
			Model Alignment:
			$\forall n\in\nodesInnerProper$:
			$\mu^n = f_{J(\psi(n))}$ is a model mechanism (Def.\ \ref{def:model}).
		\end{enumerate}
	\end{definition}

	\begin{definition}[Families of Embeddings]\label{def:local_graph_embedding_family}
		A family of local graph embeddings
		\begin{equation*}
			(J_0, \{\psi_j\}_{j\in J_0}, \graph, H, y_0)
		\end{equation*}
		is a collection of local graph embeddings
		$\psi_j: \graph\hookrightarrow I$
		of a single model-aligned local graph $\graph$.
		Further there is a fixed element $y_0 \in \nodes$,
		which will be called the anchor\allowbreak{}\mbox{(-}node),
		a symmetry $H\subset G$ and a range of applicability $J_0\subset I$.
		For $n\in\nodes$ denote:
		\begin{equation*}
			\psi_*(n): J_0 \rightarrow I,
			\halfquad
			j\mapsto \psi_j(n)
			\txt.
		\end{equation*}
		Finally, we will 
		require the following conditions to be satisfied:
		\begin{enumerate}[label=(\Roman*)]
			\item\label{def:local_graph_embedding_family:trivial_anchor}
			Trivial on Anchor:
			$\forall j\in J_0$: $\psi_j(y_0) = j$.
			\item\label{def:local_graph_embedding_family:rigidity}
			Rigidity:
			$\forall n \in \nodes$:
			$\psi_*(n)$
			is $H$-equivariant.
			\item\label{def:local_graph_embedding_family:freeness}
			Freeness:
			$\forall n\in\nodesInnerProper$:
			$\psi_*(n)$
			is injective.
		\end{enumerate}
	\end{definition}	
	\begin{rmk}
		\label{rmk:families_and_symmetries}
		The intuition behind using \emph{families} of embeddings is
		to collect repeated occurrences of the same structure, as a
		prerequisite of statistical learning.
		The properties can be understood as follows:
		\ref{def:local_graph_embedding_family:trivial_anchor} simply
		says that indexing $J_0$ within the family is by occurrences
		$J_0=\{\psi_j(y_0)|j\in J_0\}$
		of the anchor node $y_0$ plus correct neighborhood in the $I$-graph.
		\ref{def:local_graph_embedding_family:freeness} will
		avoid degeneracy of data-sets (for an injective map, the
		image contains as many elements as the domain $J_0$;
		\emph{different nodes} may overlap, but for finitely many nodes
		this will not affect limits).
		
		\ref{def:local_graph_embedding_family:rigidity}
		is specific to our formal description of symmetries by group-actions.
		Orbit-based approaches might be more general (§\ref{apdx:symmetries}),
		but for formal analysis group-actions have the nice and
		intuitive property that absorbing a $H'$-symmetric
		mechanism into a $H$-symmetric structure yields
		a $H\cap H'$-symmetric structure, \ie symmetry of larger
		structures (as intuitively expected) is reduced to the intersection
		of the symmetries of their parts.
		There are other advantages to the choice of group-actions, especially
		with regard to potential future work on structure discovery.
		See §\ref{apdx:symmetries}.
		
		\emph{Example:} One also finds for $H_1$-symmetric $X$ and
			$H_2$-symmetric $Y_x$ the product will appear in the $I$-graph
			$H_1\cap H_2$-symmetric. The causal disintegration 
			is simply $Y_x$, thus $H_2$-symmetric,
			while the anti-causal disintegration
			$(X|Y)_y$ is only $H_1\cap H_2$-symmetric,
			which provides a more concrete connection
			to invariance-based edge-orientation
			in structure-dis\-co\-ve\-ry,
			compare \eg to \citep{peters2016causal,CD-NOD,JCI}.
	\end{rmk}

	\subsection{Decorated Families} 

	These objects, so far, do not account for hidden variables.
	Further we need an analogue of the well-known "backdoor-paths" \citep{PearlBook}
	from the IID setting.
		
	\begin{definition}[Observedness]\label{def:observed}
		We call $I_0\subset I$ observed
		if $|I_0\cap \viewport(N)| \rightarrow \infty$
		for $N\rightarrow\infty$.		
		Given a
		family  of local graph embeddings $\{\psi_j\}_{j\in J_0}$,
		we call a subset $O \subset \nodes$ observed
		if
		\begin{align*}
			&J^{\text{obs}}_0(N,O)
			\halfquad:=\halfquad
			\{j \in J_0| \psi_j(O)\subset \viewport(N)\}
			\txt,\\
			&\txt{satisfies }
			|J^{\text{obs}}_0(N,O)|
			\rightarrow \infty
			\quad\txt{for } N\rightarrow\infty
			\txt.
		\end{align*}
	\end{definition}
	
	\begin{definition}[Minimal Latent Subsets]\label{def:minimal_hidden}
		Given a
		family of local graph embeddings $\{\psi_j\}_{j\in J_0}$,
		we call	an observed subset $O \subset \nodes$
		maximal observed
		if $O' \supsetneq O \Rightarrow$ $O'$ is not observed.
		In this case, we call the complement
		$L=\nodes-O$ a minimal latent subset.
		
		\emph{Remark:} A structured kernel $\mu(\graph,L)$
		encodes information on the distributions of $\nodesInner \setminus L$.
		Smaller $L$ means more information and we want to pick $L$ small.
		However, we also need data for its extraction. Hence $L$ must be
		large enough that its complement is observed.
	\end{definition}
	
	Similarly, some of the ancestral structure of the $I$-graph will
	be relevant for identifiability.
	
	\begin{definition}[Ancestral Structure]
		\label{def:ancestral_structure}
		Given a local graph $\graph$,
		an ancestral structure $\ancestralStructure$ on $\graph$
		is a set of directed edges, denoted $\rightsquigarrow$.
		We call edges in $\edgesStructural$ proper,
		and edges in $\ancestralStructure$ ancestral.
		Each ancestral edge starts at an inner node and ends at an outer node.
		We require that $\graph$ with proper \emph{and} ancestral edges is acyclic.
		
		\emph{Validity:}
		An ancestral structure $\ancestralStructure$ on $\graph$ is
		valid for an embedding $\psi:\graph\hookrightarrow I$
		if
		$\forall y\in\nodesInnerProper, \forall x\in\nodesOuterProper$
		the following condition is satisfied:
		If there is directed path $\gamma$ from $\psi(y)$ to $\psi(x)$ in
		$\IGraph \setminus \img(\psi)$, the $I$-graph with
		$\psi$-image of proper edges of $\graph$ removed,
		then
		there is an ancestral edge $y\rightsquigarrow x$ in $\ancestralStructure$.
		See Fig.\ \ref{fig:extraction} for examples of valid
		ancestral structures.
		
		We call $\ancestralStructure$ valid for a family of embeddings
		$\{\psi_j\}_{j\in J_0}$ if $\forall j\in J_0$ it is valid for $\psi_j$.
	\end{definition}

	\begin{definition}[Minimal Ancestral Structures]		
		Given are a local graph $\graph$,
		and an ancestral structure $\ancestralStructure$
		on $\graph$ valid for a family of embeddings $\{\psi_j\}_{j\in J_0}$.
		Further $J^{\txt{obs}}(N)$
		is a given, $\subset$-increasing, sequence of subsets of $I$.
		We say $\ancestralStructure$ is $J^{\txt{obs}}(N)$-minimal,
		if for all strict subsets $\ancestralStructure'\subsetneq\ancestralStructure$,
		and for all sequences of subsets $J_0'(N)\subset J^{\txt{obs}}(N)$
		such that $\ancestralStructure'$ is
		valid for $\{\psi_j\}_{j\in J'_0}$,
		the number of elements in this sub-sequence is bounded,
		\ie $\forall N:|J_0'(N)|\leq C<\infty$.
		
		\emph{Remark:} If $|J_0'(N)|$ is bounded, then $\{\psi_j\}_{j\in J'_0}$
		cannot be observed (Def.\ \ref{def:observed}) for any $L$.
	\end{definition}
	\begin{rmk}[Non-Uniqueness]
		\label{rmk:non_uniqueness}
		Neither one of these minimal structures is unique.
		For example consider a local graph containing
		the structure $X_1 \rightarrow Y \leftarrow X_2$
		such that for each $j\in J_0$ of a associated family of
		embeddings exactly one of $X_1$ or $X_2$ is observed
		(and either one asymptotically infinitely often).
		Then there are \emph{two} minimal
		latent subsets: $L_1=\{X_1\}$ and $L_2=\{X_2\}$.
	\end{rmk}
	
	These attached structures "decorate" a family of embeddings
	(non-uniquely):
	\begin{definition}[Decorated Families]
		A decorated family of embeddings
		is a family of embeddings $\{\psi_j\}_{j\in J_0}$
		together with a minimal latent subset $L\subset\nodes$
		and a $J^{\txt{obs}}_0(N,O)$-minimal
		ancestral structure $\ancestralStructure$.
	\end{definition}
	
	For identification, the interaction of backdoor-paths and latents
	will be important.
	
	\begin{definition}[Backdoor-Freeness]\label{def:backdoor_free}
		Given a decorated family of embeddings $(\{\psi_j\}_{j\in J_0}$,
		$L$, $\ancestralStructure)$, we call an ancestral edge
		$l \rightsquigarrow x \in \ancestralStructure$
		a backdoor if it starts at $l\in L$.
		We call the family backdoor-free, if
		$L\cap\nodesOuterProper=\emptyset$
		and there are no backdoors.
	\end{definition}

	\begin{figure}[ht]
		\hspace*{2em}
		\begin{minipage}{0.9\linewidth}
			\colorlet{dblue}{blue!50!black}
\colorlet{dgreen}{green!50!black}
\colorlet{dyellow}{yellow!50!black}
\colorlet{dpurple}{purple!50!black}
\colorlet{dorange}{orange!50!black}
\colorlet{c4}{dblue!20!white}
\colorlet{c3}{dgreen!20!white}
\colorlet{c1}{dyellow!20!white}
\colorlet{c2}{dpurple!20!white}
\colorlet{c5}{dorange!20!white}
\colorlet{dgray}{gray!50!black}
\colorlet{c0}{dgray!20!white}

\colorlet{b0}{c0!40!white}
\colorlet{b5}{c5!40!white}
\colorlet{b3}{c3!40!white}
\colorlet{b4}{c4!40!white}

\usetikzlibrary{decorations.pathmorphing}
\pgfdeclarelayer{bg}
\pgfsetlayers{bg,main}

\begin{tikzpicture}
	\draw (0,0) node (Igraph)
		{
			\begin{tikzpicture}[xscale=1.25]
				\draw (0,0) node (X)[rectangle, fill=c0, rounded corners, inner sep=0.4em,
				draw=dgray,	thick]
				{$X$};
				\draw (1,0) node (Z)[rectangle, fill=c0, rounded corners, inner sep=0.4em,
				draw=dgray,	thick]
				{$Z$};
				\draw (2,0) node (Y)[rectangle, fill=c0, rounded corners, inner sep=0.4em,
				draw=dgray,	thick]
				{$Y$};
				\draw (1,1) node (L)[rectangle, fill=c0, rounded corners, inner sep=0.4em,
				draw=dgray,	thick, dashed]
				{$L$};
				
				\draw [->, very thick] (X) -- (Z);	
				\draw [->, very thick] (Z) -- (Y);	
				\draw [->, very thick] (L) -- (X);	
				\draw [->, very thick] (L) -- (Y);	
			\end{tikzpicture}
		};
		\draw (Igraph.north west) + (0,0.15) node (lI) [anchor=south west] {in $I$-graph:};
		\draw (Igraph.north east) + (0.5,0.15) node (lX) [anchor=south west] {start at $X$:};
		\draw (Igraph.north east) + (0.5,0) node (atX1) [anchor=north west] {			
			\begin{tikzpicture}[xscale=1.25]
				\draw (1,1) node(L)
				[circle, draw, very thick, inner sep=0.15em, outer sep=0.15em, anchor=center]{};
				
				\draw (0,0) node (X)[rectangle, rounded corners, inner sep=0.4em,
				draw=dgray, fill=c0, thick, anchor=center]
				{$X$};	
				
				\draw [->, very thick] (L) -- (X);
				
				\draw[rounded corners, dpurple, dotted, very thick]
				(0.7,0.7) rectangle (1.3,1.3);	
			\end{tikzpicture}
		};
		\draw [->,line width=0.4em, c0] (atX1.east) + (0.2,0) --
			node[above, yshift=0.5em, dgray,align=center]{Absorb\\Parent}
			++(1.5,0);
		\draw (atX1.north east) + (1.9,0) node (atX2) [anchor=north west] {						
			\begin{tikzpicture}[xscale=1.25]
				\draw (1,1) node (L)[rectangle, rounded corners, inner sep=0.4em,
				draw=dgray, dashed, fill=c0, thick]
				{$L$};
				
				\draw (0,0) node (X)[rectangle, rounded corners, inner sep=0.4em,
				draw=dgray, fill=c0, thick]
				{$X$};	
				
				\draw [->, very thick] (L) -- (X);
			\end{tikzpicture}
		};		
		\draw (atX2.north east) + (0.5,0.15) node (lZ) [anchor=south west] {start at $Z$:};
		\draw (atX2.north east) + (0.5,-0.15) node (atZ) [anchor=north west] {	
			\begin{tikzpicture}[xscale=1.25]
				\draw (0,0) node(X)
				[circle, draw, very thick, inner sep=0.15em, outer sep=0.15em, anchor=center]
				{};
				
				\draw (1,0) node (Z)[rectangle, rounded corners, inner sep=0.4em,
				draw=dgray, fill=c0, thick, anchor=center]
				{$Z$};	
				
				\draw (1,1) node[anchor=center, inner sep=0.4] {\vphantom{$L$}}; 
				
				\draw [->, very thick] (X) -- (Z);
			\end{tikzpicture}
		};
		\draw (Igraph.south west) + (0,-0.4) node (lY) [anchor=north west] {start at $Y$:};
		\draw (Igraph.south west) + (0,-0.75) node (atY1) [anchor=north west] {
			\begin{tikzpicture}[xscale=1.25]
				\draw (1,0) node(Z)
				[circle, draw, very thick, inner sep=0.15em, outer sep=0.15em, anchor=center]{};
				\draw (1,1) node(L)
				[circle, draw, very thick, inner sep=0.15em, outer sep=0.15em, anchor=center]{};
				
				\draw (2,0) node (Y)[rectangle, fill=c0, rounded corners, inner sep=0.4em,
				draw=dgray,	thick, anchor=center]
				{$Y$};
				
				\draw [->, very thick] (Z) -- (Y);	
				\draw [->, very thick] (L) -- (Y);	
				
				\draw[rounded corners, dpurple, dotted, very thick]
				(0.7,0.7) rectangle (1.3,1.3);	
			\end{tikzpicture}
		};		
		\draw [->,line width=0.4em, c0] (atY1.east) + (0.2,0) --
			node[above, yshift=0.5em, dgray,align=center]{Absorb\\Parent}
		++(1.5,0);
		\draw (atY1.north east) + (1.9,0.3) node (atY2) [anchor=north west] {
			\begin{tikzpicture}[xscale=1.25]
				\draw (1,0) node(Z)
				[circle, draw, very thick, inner sep=0.15em, outer sep=0.15em,
				anchor=center]
				{};
				
				\draw (1,1) node (L)[rectangle, rounded corners, inner sep=0.4em,
				draw=dgray, dashed, fill=c0, thick,
				anchor=center]
				{$L$};
				\draw (2,0) node (Y)[rectangle, fill=c0, rounded corners, inner sep=0.4em,
				draw=dgray,	thick,
				anchor=center]
				{$Y$};
				
				\draw [->, very thick] (Z) -- (Y);
				\draw [->, very thick] (L) -- (Y);
				\draw [->, very thick, line join=round,
				decorate, decoration={
					zigzag,
					segment length=4,
					amplitude=.9,post=lineto,
					post length=2pt
				}] (L) -- (Z);

				\draw[rounded corners, dgreen, dashed, very thick]
				(0.6,-0.4) rectangle (1.4,1.5);	
			\end{tikzpicture}
		};
		\draw [->,line width=0.4em, c0] ([yshift=-0.2em] atY2.east) + (0.2,0) --
			node[above, yshift=0.5em, dgray,align=center]{Absorb\\Child}
			++(1.5,0);
		\draw (atY2.north east) + (1.9,0) node (atY3) [anchor=north west] {			
			\begin{tikzpicture}[xscale=1.25]
				\draw (1,0) node(Z)
				[circle, draw, very thick, inner sep=0.15em, outer sep=0.15em,
				anchor=center]
				{};
				
				\draw (0,0) node (X)[rectangle, rounded corners, inner sep=0.4em,
				draw=dgray, fill=c0, thick,
				anchor=center]
				{$X$};	
				\draw (1,1) node (L)[rectangle, rounded corners, inner sep=0.4em,
				draw=dgray, dashed, fill=c0, thick,
				anchor=center]
				{$L$};
				\draw (2,0) node (Y)[rectangle, fill=c0, rounded corners, inner sep=0.4em,
				draw=dgray,	thick,
				anchor=center]
				{$Y$};
				
				\draw [->, very thick] (Z) -- (Y);
				\draw [->, very thick] (L) -- (Y);
				\draw [->, very thick] (L) -- (X);
				\draw [->, very thick, line join=round,
				decorate, decoration={
					zigzag,
					segment length=4,
					amplitude=.9,post=lineto,
					post length=2pt
				}] (X) -- (Z);
			\end{tikzpicture}
		};
%
%
	\begin{pgfonlayer}{bg}			
		\draw[fill=b0, rounded corners, draw=c0] (lI.north west) rectangle
			(Igraph.south east);
		\draw[fill=b5, rounded corners, draw=c5] (lX.north west) rectangle
			(atX2.south east);
		\draw[fill=b3, rounded corners, draw=c3] (lZ.north west) rectangle
			(atZ.south east);
		\draw (atY3.south east) + (0,-0.4) node (pY4) {};
		\draw[fill=b4, rounded corners, draw=c4] (lY.north west) rectangle
			(pY4);
	\end{pgfonlayer}
\end{tikzpicture}
		\end{minipage}
		\caption{Systematic construction of minimal backdoor-free families.
		While there are: hidden outer nodes (red-dotted boxes) absorb parents;
		backdoor-paths (green-dashed boxes) absorb children.
		Note that $\mu(\graphStructural,L)$ extracted by
		Thm.\ \ref{thm:extract_from_backdoor_complete}
		will in the end drop the remaining ancestral edge (it depends
		only on $\graph$ and $L$, not on $\ancestralStructure$),
		and can be further decomposed (cf.\ Fig.\ \ref{fig:graphical_ops}).
		A single gluing-operation then yields the well-known
		frontdoor formula \citep{PearlBook}.
		}\label{fig:extraction}
	\end{figure}
	
	There is a simple algorithmic approach (cf.\ Algo.\ \ref{algo:fcs},
	§\ref{apdx:algo_bd_free})
	to find (enough, Lemma \ref{lemma:c_conn_extr_enough}, "minimal")
	backdoor-free families by "growing"
	embeddings successively.
	As a starting point we use \emph{direct embeddings} of mechanisms:
	A mechanism $f_J$ (Def.\ \ref{def:mechanism}) with range of applicability
	$J$ defines an embedded family indexed by $J$ with
	a single inner node $n$ aligned to $f_J$, one outer node per parent and
	$\psi(j,n)=j$, $\psi(j,k)=\PaIdx{k}(j)$
	for $k\in\nodesOuterProper$.
	In the next step, as illustrated in Fig.\ \ref{fig:extraction},
	we "absorb" latent outer nodes as inner
	nodes and close backdoors by "absorbing" their
	starting node (a child of a latent node) as an inner node.	
	This step must be repeated until the result is backdoor-free.
	Here \emph{absorbing} a node means checking overlaps with direct
	embeddings (see above) to find a family of embeddings for
	an enlarged local graph.
	Note however, that the same mechanism
	can appear in the $I$-graph with different neighborhoods, so absorbing
	parents or children of $L$ is not unique.
	
	There are at least three practical challenges to this approach:
	(i) Many constructions (like minimal latent-sets or candidates for absorption)
	are non-unique,	so we usually have to deal with sets of results.
	(ii) We can start this construction from any directly embedded mechanism,
	but even univariate queries
	will require intermediate results on joint descriptions (products of kernels),
	so besides starting-points of backdoor-paths, additional
	children of latents have to be absorbed (depending on the query).
	(iii) Minimal backdoor-free families can, in principle,
	have arbitrarily many nodes which means algorithms
	either terminate after a finite number of steps or they are
	complete, but never both; in the IID case this number is
	bounded by the number of system variables. In practice
	large backdoor-free families require the estimation of
	kernels with many parameters and use (often) small subsets of
	the full data-set, so if we are interested in
	plausibly estimable results only, we can bail on large graphs
	and algorithms \emph{do} effectively terminate.

	\subsection{Notions of Identifiability} 
	\label{sec:identifiability}

	In the IID-case, identifiability is conventionally
	defined as being uniquely determined by the observational
	distribution $P_{\txt{obs}}$ \citep{tian2002general,shpitser2006identification2,
		shpitser2006identification1,bareinboim2012TransportCompleteness}.
	Our closest analogue of $P_{\txt{obs}}$
	is the shallow distribution $\shallowDistr$ and the observable world.
	But each variable in the observable world is observed at most once,
	so we cannot "learn" $\shallowDistr$ (at least not
	immediately as for $P_{\txt{obs}}$), so uniquely expressing
	anything relative to $\shallowDistr$ does not identify it
	in any reasonable sense.
	In fact, a suitable "joint" distribution to learn from
	will often not exist (a similar, but less obvious
	problem occurs in example \ref{example:extraction_frontdoor}):	
	\begin{example}[Random Walk]\label{example:random_walk}
		Let $\eta_t \sim \mathcal{N}(b,1)$ be normal distributed
		and $X_n := \sum_{t=1}^n \eta_t$ a random walk.
		The kernel that describes $X_{n+1}$ relative
		to the preceding time-step $X_{n+1}(x_n) = x_n + \mathcal{N}(b,1)$
		is perfectly well-defined. Yet, already the \emph{marginal} distributions
		even of the unbiased, $b=0$, one-dimensional case 
		of the "joint distribution" of $(X_n, X_{n+1})$ are
		(by CLT) of the form $X_n \sim \mathcal{N}(0,\sigma_n^2=n)$.
		They not only depends on $n$, even worse, the density pointwise
		approaches zero for large $n$.
	\end{example}

	The fundamental idea of statistics is to combine many observations
	of the same thing.
	We employ a notion of identifiability that is
	a rather direct formalization of this idea.
	This approach helps to separate causal ideas from problems of 
	statistical estimation theory (which is not a main topic of the present paper).
	Further discussion can be found in §\ref{apdx:identifiablity}.

	\begin{definition}[Data-Set]\label{def:data_set}
		A data-set is a tuple
		$\mathcal{D} = (\randomVar{X}_j, \randomVar{Y}_j)_{j\in J_0}$,
		where $\randomVar{X}_j$ and $\randomVar{Y}_j$
		are tuples
		of $\anyVar_i$, \ie there exist $i_{X,j}^{(1)}, \ldots,
		i_{X,j}^{(n)} \in I$,
		such that $\randomVar{X}_j = (\anyVar_{i_{X,j}^{(1)}},
		\ldots, \anyVar_{i_{X,j}^{(n)}})$ and analogously for $\randomVar{Y}_j$.
		
		We call the data-set valid if $\randomVar{X}$ and $\randomVar{Y}$
		are observed that is
		\begin{equation*}
			J_0^{\txt{obs}}(N)
			\halfquad:=\halfquad
			\{\halfquad
			j\in J_0
			\halfquad|\halfquad
			\forall m:
			i_{X,j}^{(m)} \in \viewport(N),
			\forall m':
			i_{Y,j}^{(m')} \in \viewport(N)
			\halfquad\}
		\end{equation*}
		satisfies $|J_0^{\txt{obs}}(N)|\rightarrow\infty$ as $N\rightarrow\infty$
		and is non-degenerate
		$j\neq j' \Rightarrow \forall m: i_{Y,j}^{(m)}\neq i_{Y,j'}^{(m)}$
		(this last condition is imposed only on $Y$, not on $X$;
		see Rmk.\ \ref{rmk:datasets_and_support_multi_level}).
	\end{definition}
	
	\begin{example}[Data-Set of Decorated Families]
		\label{example:data_set_of_embedding}
		Given a decorated family of embeddings
		$(\{\psi_j\}_{j\in J_0},L, \ancestralStructure)$ with
		$L\cap\nodesOuterProper=\emptyset$,
		there is a valid data-set
		$(\randomVar{X}_j, \randomVar{Y}_j)_{j\in J_0}$
		(Def.\ \ref{def:data_set})
		defined for $j\in J_0$ as the tuples
		\begin{align*}
			\randomVar{X}_j
			\halfquad&=\halfquad
			(\anyVar_{\psi_j(n)})_{n\in\nodesOuterProper}\\
			\randomVar{Y}_j
			\halfquad&=\halfquad
			(\anyVar_{\psi_j(n)})_{n\in\nodesInnerProper\setminus L}
			\txt.
		\end{align*}
	\end{example}
	
	\begin{definition}[Direct Identifiability]\label{def:identification_direct_mt}
		Given a valid data-set (Def.\ \ref{def:data_set})
		$\mathcal{D} = (\randomVar{X}_j, \randomVar{Y}_j)_{j\in J_0}$
		and kernels $\{X^{(j)}\}_{j\in J_0}$, $Y_{x}$
		(we write the argument as index $x$ to indicate that
		the corresponding parent is $X$, notation \ref{notation:causal_wiring}),
		such that the shallow distribution $\shallowDistr$
		(Def.\ \ref{def:obs_world}) satisfies
		\begin{equation*}
			\forall j\in J_0:\halfquad
			\shallowDistr(\randomVar{X}_j, \randomVar{Y}_j)
			= X^{(j)} \otimes Y_{x}
			\txt,
		\end{equation*}
		then we call $Y_{x}$ directly identifiable. 
		Formally we also consider the elements of
		$\knownFunction$, which is
		known a priori (Def.\ \ref{def:model}), to be
		directly identifiable.
	\end{definition}
	\begin{example}
		\label{example:direct_id_iid}
		If the data are IID, then $X^{(j)} \equiv X$ do not depend on $j$
		and given (as\-ymp\-tot\-i\-cal\-ly infinite)	data for the joint 
		$\shallowDistr(\randomVar{X}_j, \randomVar{Y}_j) = X \otimes Y_x$,
		we consider $Y_x = P(\randomVar{Y}|\randomVar{X}=x)$ directly identifiable.
		Indeed in this case our notion of direct identifiability is the same
		as considering a SCM's observed distribution $P_{\txt{obs}}$
		and its conditionals known
		(Lemma \ref{lemma:iid:direct_id}).
	\end{example}
	
	\begin{definition}[Identifiability]\label{def:identification}		
		We call a kernel $\mu$ identifiable if
		$\mu$ can be uniquely computed from
		finitely many directly identifiable kernels.
		We call a collection of directly identifiable kernels
		$\mu_1, \ldots, \mu_n$ (plus their datasets)
		together with a regular functional $F$ computing
		$\mu=F[\mu_1,\ldots,\mu_n]$ an identification strategy for $\mu$.
	\end{definition}

	\subsection{Identifiability from Embeddings} 

	Having fixed a formal notion of identifiability
	from data-sets, we continue by giving
	formal statements on the identifiability of structured kernels
	from decorated families of embeddings.

	The main result of this section is that backdoor-free families
	have identifiable $\mu(\graph,L)$.
	There is a simple intuition for this:
	For each $j$, the image of $\psi_j(\graph)$
	depends on its $I$-graph ancestors only through
	$\nodesOuterProper$ (by Def.\ \ref{def:local_graph_embedding}~%
	\ref{def:local_graph_embedding:inner_parents_incl}),
	if there are no backdoor-paths $L \rightsquigarrow \nodesOuterProper$,
	then inner nodes $\nodesInner$ form a union of c-components
	of the $I$-graph, thus $\psi_j(\graph)$ is essentially
	a union of structural c-components of the $I$-graph and
	a variant of Lemma \ref{lemma:c_components:mt} applies.

	\begin{thm}\label{thm:extract_from_backdoor_complete}
		Given a backdoor free family of embeddings 
		$(\{\psi_j\}_{j\in J_0},L,\ancestralStructure)$
		(Def.\ \ref{def:backdoor_free}),
		then
		\begin{equation*}
			\mu(\graph,L)
			\txt{ (Def.\ \ref{def:structured_kernel_simp}) is identifiable
				(Def.\ \ref{def:identification}).}
		\end{equation*}				
	\end{thm}

	For the computation of $\mu(\graph,L)$ in practice
	see Rmk.\ \ref{rmk:practical_computation_of_extraction}.
	Backdoor-free families of embeddings
	are also rather generic means of curating data-sets with invariant properties,
	thus the following seems plausible:
	
	\begin{conjecture}[Completeness of Extraction]
		\label{conjecture:extraction}
		In the absence of selection bias,
		if a structured kernel $\mu$ is semi-directly identifiable
		(a $\otimes$-product of directly identifiable kernels),
		then there is a backdoor-free family of embeddings
		on $(\graph,L)$ with $\mu = \mu(\graph, L)$,
		at least for an orbit-based notion of symmetry,
		see §\ref{apdx:symmetries} and after including
		multi-level structure (see appendix, §\ref{apdx:extraction_from_data}).
	\end{conjecture}
	
	At this point, we capture for example IID multi-context
	cases (often called observational and experimental transport
	\citep{bareinboim2012TransportCompleteness,
		Bareinboim2016TransportOverview,pearl2022external}):
	\begin{example}[IID Extraction]
		\label{example:iid_extraction_direct}
		For any individual IID context $c$,
		\ie an IID data-set $\{c\}\times\IVars\times\ISample$, but not necessarily
		\emph{all} available data $\IDataset\times\IVars\times\ISample$,
		with graph $\graphClassical^c$
		(possibly involving nodes aligned to do- or soft-interventions)
		there is a backdoor-free family $\{\psi^c_j\}_{j\in\ISample}$,
		embedding $\graphClassical^c$ with $\nodesInner^c=\nodes^c=\IVars$
		as $\psi^c_j(v) = (c,v,j)\in \IDataset\times\IVars\times\ISample$.
		This embedding, by Thm.\ \ref{thm:extract_from_backdoor_complete},
		extracts $\mu(\graphClassical^c,L^c)$ (where $L^c$ are the latents in 
		context $c$.
		Interventional distributions can be recovered
		as in example \ref{example:iid_id_from_structured_only}.
		Gluing can combine pieces obtained as decomposition-results
		from multiple contexts.
		This approach seems to be complete for mz-transport
		\citep{bareinboim2012TransportCompleteness} questions
		(§\ref{apdx:iid:mz_transport}), but can additionally
		combine data-sets with different latent-sets, different graphs,
		or even completely different underlying structure (like time-series).
	\end{example}

	\subsection{Summary: Extracted Knowledge}
	
	We briefly summarize what kind of knowledge will be available
	for reasoning about causal queries.
	
	\begin{definition}[Extracted Knowledge Set]
		\label{def:knowledge_set_extracted}
		Let $\mathcal{B}$ be the set of all backdoor-free
		families of embeddings.	
		For $(\psi,L,\ancestralStructure)\in\mathcal{B}$
		define a structured model
		$k(\psi,L,\ancestralStructure) := (\graph, L)$.
		By Thm.~\ref{thm:extract_from_backdoor_complete},
		$\mu(\graph,L)$ is identifiable
		(we simply say $k$ is identifiable).
		
		For $\tilde{f} \in \knownFunction$, define
		a structured model $k(\tilde{f})=(\graphStructural(\tilde{f}),L=\emptyset)$,
		where $\graphStructural(\tilde{f})$ is a graph with a single inner node $y$
		aligned to $\mu^y=\tilde{f}$, an outer node
		for each of the $\kappa$ parents
		$\nodesOuter = \{x_1, \ldots, x_\kappa\}$
		and a proper edge from each $x_k \rightarrow y$.
		By definition (Def.\ \ref{def:model} and \ref{def:identification_direct_mt}),
		elements of $\knownFunction$ are considered known,
		thus $k(\tilde{f})$ is identifiable for all $\tilde{f} \in \knownFunction$.

		Define the, thus (element-wise)
		identifiable, single-level basic knowledge set as
		\begin{equation*}
			\knowledgeSet_{\txt{basic}}
			\halfquad=\halfquad
			\bigcup_{(\psi,L,\ancestralStructure)\in\mathcal{B}}
				k(\psi,L,\ancestralStructure)
			\halfquad\cup\halfquad
			\bigcup_{\tilde{f}\in\knownFunction}
				k(\tilde{f})
			\txt.
		\end{equation*}
	\end{definition}

	\subsection{Extension: Multi-Level Statistics}
	\label{sec:multi_level}
	
	Multi-level statistics \citep{Gelman2006} models collections
	of datasets hierarchically (see example below).
	Our models do not explicitly fix a hierarchy on
	elements of $I$ (nodes of the $I$-graph).
	But, while parent-sets in the $I$-graph are finite by definition,
	arbitrary numbers of children are allowed.
	Thus parts of the $I$-graph can be described meaningfully as
	hierarchical systems; in our formalism hierarchy is
	an emerging phenomenon -- it is explicit only in the machinery
	used for identification of queries, not in the
	queries or models themselves.
	
	\begin{example}\label{example:multi_level_illustrate}
		Consider a case, similar to the motivating example of \citep{JPCMCI},
		where data was collected on water-levels (and similar properties)
		at multiple river-sites over time.
		Some physical mechanisms may be shared, but could depend on
		site-specific, constant in time,
		"contextual" properties like slope or form of the riverbed.
		It is possible to model these mechanisms per site
		(equivalently, a unique site-id is treated as "context-variable"
		\citep{JCI,JPCMCI}).
		There are relevant questions, \eg about edge-orientations or interventions
		at a particular site that rely on such information.
		
		However, we might want to introduce a second level of modeling
		that inter-relates contextual properties:
		The form (\eg depth vs.\ width) of a riverbed may
		be driven by slope when comparing different sites.
		In that case, there should
		be hidden variables (\eg slope) per site, with infinitely
		many children (time-points at that site).
		Such structure is allowed by our models implicitly.
		The argument above, and choice of hierarchy,
		can be made dynamically:
		The same variable may, with regard to some families of embeddings,
		take the role of a context,
		while appearing as an ordinary variable in others.
	\end{example}

	\section{Queries and Prediction} 
	\label{sec:queries}

	At this point, we introduced a formalism to describe
	a model, machinery to extract structured kernels from observations,
	and intuitive graphical transformations to work with such structured kernels.
	What is left to do, is to provide a way to phrase causal questions, including
	interventional examples \ref{example:iid_id_from_structured_only} and
	\ref{example:iid_extraction_direct}, but also more complex
	queries as the one in Fig.\ \ref{fig:intro}.
	We do so in two steps: First we define simple formal objects
	(called basic queries below),
	which we make then accessible through a more intuitive graphical language.
	Finally, the identifiability of queries will be a simple consequence 
	of the results obtained in the previous sections.
	We will additionally phrase the question of query-identifiability
	in a slightly more abstract way. This facilitates and organizes
	the incorporation of internal structure,
	for example via instrumental variable arguments,
	for future work.
	Finally, we compare to standard formalisms in the IID-case.

	We start with an example, see also Fig.\ \ref{fig:do_query},
	to clarify how IID interventions can be mapped to these queries as a special case.
	
	\begin{figure}[tb]
		\colorlet{dblue}{blue!50!black}
\colorlet{dgreen}{green!50!black}
\colorlet{dyellow}{yellow!50!black}
\colorlet{dpurple}{purple!50!black}
\colorlet{dorange}{orange!50!black}
\colorlet{c4}{dblue!20!white}
\colorlet{c3}{dgreen!20!white}
\colorlet{c1}{dyellow!20!white}
\colorlet{c2}{dpurple!20!white}
\colorlet{c5}{dorange!20!white}
\colorlet{dgray}{gray!50!black}
\colorlet{c0}{dgray!20!white}

\pgfdeclarelayer{bg}
\pgfsetlayers{bg,main}

\begin{tikzpicture}[xscale=1.25]
	\draw (0,0) node (W)[rectangle, fill=c3, rounded corners, inner sep=0.4em,
	draw=dgreen,	thick]
	{$W$};
	\draw (1,0) node (Z)[rectangle, fill=c3, rounded corners, inner sep=0.4em,
	draw=dgreen, thick]
	{$Z$};
	\draw (2,0) node (X)[rectangle, fill=c3, rounded corners, inner sep=0.4em,
	draw=dgreen,	thick]
	{$X$};
	\draw (3,0) node (Y)[rectangle, fill=c3, rounded corners, inner sep=0.4em,
	draw=dgreen,	thick]
	{$Y$};
	
	\draw [->, very thick] (W) -- (Z);	
	\draw [->, very thick] (Z) -- (X);	
	\draw [->, very thick] (X) -- (Y);
	
	\draw [<->, very thick,dgray] (W) edge[bend left=2.5cm]
	node[pos=0.6,above][rectangle, fill=c0, rounded corners, inner sep=0.4em,
	draw=dgray, thick, dashed]{$L_2$} (Y);
	\draw [<->, very thick, dgray] (W) edge[bend left=1.75cm]
	node[pos=0.5,above][rectangle, fill=c0, rounded corners, inner sep=0.4em,
	draw=dgray, thick, dashed]{$L_1$} (X);
	
	\foreach \s in {0,...,3}{		
		\draw (4.25+1.78*\s, 0.75) node {\scalebox{0.4}[0.4]{
		\begin{tikzpicture}[xscale=1.25]
			\draw (0,0) node (W\s)[rectangle, fill=c3, rounded corners, inner sep=0.4em,
			draw=dgreen,	thick]
			{$W$};
			\draw (1,0) node (Z)[rectangle, fill=c3, rounded corners, inner sep=0.4em,
			draw=dgreen, thick]
			{$Z$};
			\draw (2,0) node (X)[rectangle, fill=c3, rounded corners, inner sep=0.4em,
			draw=dgreen,	thick]
			{$X$};
			\draw (3,0) node (Y)[rectangle, fill=c3, rounded corners, inner sep=0.4em,
			draw=dgreen,	thick]
			{$Y$};
			
			\draw [->, very thick] (W\s) -- (Z);	
			\draw [->, very thick] (Z) -- (X);	
			\draw [->, very thick] (X) -- (Y);
			
			\draw [<->, very thick,dgray] (W\s) edge[bend left=2.5cm]
			node[pos=0.6,above][rectangle, fill=c0, rounded corners, inner sep=0.4em,
			draw=dgray, thick, dashed]{$L_2$} (Y);
			\draw [<->, very thick, dgray] (W\s) edge[bend left=1.75cm]
			node[pos=0.5,above][rectangle, fill=c0, rounded corners, inner sep=0.4em,
			draw=dgray, thick, dashed]{$L_1$} (X);
		\end{tikzpicture}
		}};	
		\pgfmathsetmacro{\idx}{int(\s+2)}
		\draw (4.25+1.78*\s, 0.1) node {$s_{\idx}$};
	}
	\draw (10.65,0.75) node {$\cdots$};
	\draw (1.5,-0.65) node {$s_1$};
	
	\begin{scope}[xshift=4cm, yshift=-1.65cm]
		\draw (3,0) node (Y)[rectangle, fill=c0, rounded corners, inner sep=0.4em,
		draw=dgray,	thick, dashed]
		{$Y$};
		\draw (3,1) node (L2)[rectangle, fill=c0, rounded corners, inner sep=0.4em,
		draw=dgray,	thick, dashed]
		{$L_2$};
		\draw (2,0) node (X)[rectangle, fill=c2, inner sep=0.4em,
		draw=dpurple,	thick, dashed]
		{$X^{\txt{do}}$};
		
		\draw [->, very thick] (X) -- (Y);
		\draw [->, very thick] (L2) -- (Y);
	\end{scope}
	\draw (9,-1.15) node  [align=left] {Query (at $s=*$) is\\also part of $I$-graph.};
	
	\draw (2,-1.5) node(V)[align=center]
	{"observed" $i\in\viewportEventual$\\(solid\,/\,green)};
	\draw[line width=0.2cm, ->, c3] (V) edge[bend left] ([yshift=-0.1cm]W.south);
	\draw[line width=0.2cm, ->, c3] (V) edge[bend right] (5,0.25);
\end{tikzpicture}\\[-1.5em]
		\caption{Illustration of $I$-graph for do-intervention and IID-data
			(unified model). Nodes correspond to variables of the observable world,
			labels (capital letters) indicate mechanisms (Def.\ \ref{def:mechanism}).}
		\label{fig:do_query}
	\end{figure}
	
	\begin{example}[IID Do-Intervention]		
		\label{example:iid_do_intervention}
		Recall example \ref{example:model_iid}:
		Given a SCM $M$ on $\IVars$ and $X, Y\subset \IVars$,
		we want to predict $Y$ under intervention on $X$.
		There is an associated unified model $\mathcal{M}$.
		This model $\mathcal{M}$ has index-set
		$I=\IVars \times (\ISample\cup \{*\})$,
		with model-mechanisms induced by SCM-mechanisms,
		for $(v,*)$ with $v\in X$ the mechanism is
		an intervention $\delta[x_v]\in\knownFunction$.		
		We will denote by $\theta=x$ the interventional value
		$\PearlDo(X=x)$, which fully parametrizes $\knownFunction$.
		
		Note that the $I$-graph has no edges between nodes at
		different samples in $s\neq s' \in \ISample \cup \{*\}$,
		in particular we are in a single-level setup.
		We fix $\tilde{Y}:=Y\times \{*\}$. In this setup
		$\realizedDistr(\tilde{Y})=\shallowDistr(\tilde{Y})$,
		\ie
		shallow, Def.\ \ref{def:shallow_distr}, and realized,
		Def.\ \ref{def:realized_world}, distribution
		agree for $\anyVar_{\tilde{Y}}$; in the spirit of
		Rmk.\ \ref{rmk:realized_world_interpretation}
		one may want to clarify explicitly
		$\realizedDistr(\tilde{Y}|\observedFunction=\mathEulervmScript{F}_{\txt{obs}})
		=\shallowDistr(\tilde{Y}|\observedFunction=\mathEulervmScript{F}_{\txt{obs}})$.
		This follows from the local Markov condition
		(Lemma \ref{lemma:obs_world_existence}).
		
		From example \ref{example:model_iid_distr} we recall
		$\shallowDistr(\tilde{Y})=P(Y|\PearlDo(X=x))$.
		The dataset $\mathcal{D}(\omega)$ corresponds
		directly to the observations from the (unintervened)
		SCM $M$.
		Therefore, identifying $\realizedDistr(\tilde{Y})$
		from $\mathcal{D}(\omega)$ is the same
		as identifying $P(Y|\PearlDo(X=x))$
		from $P^{\txt{obs}}$ (the observational distribution).
		
		We will formalize questions about $\realizedDistr$
		(like "find $\realizedDistr(\tilde{Y})$") as "basic queries".
		Identification of interventional distributions
		(and thus causal effects etc.)
		in the IID\Slash{}SCM sense is a special case of such a
		query by the construction described above.
		Under weak assumptions always satisfied
		for do-interventions on IID-models,
		basic queries can be transformed into
		structured queries which can be answered soundly
		using §\ref{sec:structured_kernels} and §\ref{sec:extraction_from_data}
		(see Thm.\ \ref{thm:query_id}).
		This approach is complete in the IID case
		§\ref{apdx:iid:standard_results}.
	\end{example}
	
	In §\ref{sec:intro:invariant_mechaisms} of the introduction,
	we reviewed the informal idea of causality as study of invariances between
	observations and query.
	Example \ref{example:model_iid} made this intuition formally precise
	by introducing a unified model for observed and intervened setup
	(see also Fig.\ \ref{fig:do_query}).
	This uniform modeling is not just an interesting viewpoint,
	it further makes the description of complex queries much simpler:
	Given a complex practical question, we could, in analogy to
	§\ref{sec:models} and §\ref{sec:extraction_from_data} build formal
	infrastructure to phrase this question as a formal statement.
	By uniform modeling, we can simply reuse
	§\ref{sec:models} and §\ref{sec:extraction_from_data}:
	Observations and query can be described dually\Slash{}uniformly.
	
	A query asks what parts of a world of a given form \emph{would} look like.
	The identification of a query
	then amounts to "running a simulation" of the query-part of the world.
	For interventional effects this simulation replaces mechanisms
	in the observed part, but retains the graphical form.
	The extension to queries as presented here can be thought of
	as an extension to more general simulations,
	for example the query in Fig.\ \ref{fig:intro}
	(see also example \ref{example:simulation}).

	\subsection{Basic Queries} 

	A basic query is formulated directly from a realized world.
	Query-identifiability is thus a tight statement; it directly relates
	simple and abstract structures, without referencing any of the
	technology like embeddings or structured kernels in its hypothesis or claim.

	\begin{definition}[Query]\label{def:query}
		A basic query \query{} consists of finite subsets
		$\tilde{Y},\tilde{X} \subset I$
		with $\tilde{X}\cap\tilde{Y} = \emptyset$
		plus a set of parameters $\theta=\knownFunction$ (of the shallow distribution),
		with target
		\begin{equation*}
			\realizedDistr(\tilde{Y}|\tilde{X}=x)
			\txt.
		\end{equation*}
		On (standard) Borel spaces (Ass.\ \ref{ass:standard_borel}),
		a basic query has a regular version, that is, there is a
		unique probability kernel $\tilde{Y}_x$ such that
		$\realizedDistr(\tilde{Y}|\tilde{X}=x)=\tilde{Y}_x$.
		We call a basic query identifiable,
		if $\tilde{Y}_x$ is identifiable (Def.\ \ref{def:identification})
		from the dataset $\mathcal{D}(\omega)$.
		
		\emph{Remark:} For interventional queries $\tilde{X}=\emptyset$,
		the intervention is parameterized by $\theta$.
		Other $\tilde{X}$ allow queries for intermediate
		results (cf.\ Fig.\ \ref{fig:graphical_ops}) and are necessary
		\eg in example \ref{example:random_walk},
		see also example \ref{example:query_argument}.
	\end{definition}

	\subsection{Structured Queries} 

	The illustrations in Fig.\ \ref{fig:intro} and \ref{fig:do_query}
	of queries really show (embedded) local graphs with a few simple
	properties, not abstract basic queries.
	This can be made formally clearer.
	
	\begin{definition}[Structured Query]
		\label{def:structured_query}
		A structured query is an embedded local graph
		$\psi:\graph\hookrightarrow I$ together with a set of interest
		$Y \subset \nodesInnerProper$,
		such that $\ancestralStructure=\emptyset$ is valid
		(Def.\ \ref{def:ancestral_structure}).
		
		\emph{Underlying Query:}
		Given a structured query $(\psi, Y)$,
		there is an underlying basic query $q(\psi,Y):=(\tilde{Y}=\psi(Y),
		\tilde{X}=\psi(\nodesOuterProper), \theta)$.
		We call $(\psi, Y)$ identifiable, if the underlying
		basic query is identifiable.
	\end{definition}
	
	For many basic queries
	there is a canonical structured query
	constructed by, starting from $\tilde{Y}$, adding $I$-graph parents
	until a parent is either in $\tilde{X}$ (and added as outer node),
	or no more parents are available.
	By Ass.\ \ref{ass:finite_past} (finite past),
	this construction terminates after a finite number of steps with
	a finite graph. Embed this graph by the identity-mapping
	$\psi=\id$ and identify $Y:=\tilde{Y}$ with its $\psi$-image.
	\begin{lemma}[Associated Structured Query]
		\label{lemma:assoc_structured_query}
		Given is a basic query \query{}
		about a single-level statement,
		formally: $\tilde{Y}$ is disconnected (not reachable 
		by directed or undirected paths in the $I$-graph)
		from $\viewportEventual$.
		Further, $\tilde{X}$ cannot be bypassed in the sense of:
		For $y\in\tilde{Y}$
		and $w \in \Anc_I(\tilde{X})\setminus\tilde{X}$, there
		is no directed path $\gamma: w \rightsquigarrow y$ in the $I$-graph
		with $\gamma\cap\tilde{X} = \emptyset$.
		
		Then 
		there is a structured query $(\psi(q),\tilde{Y})$
		with underlying query $q$.
	\end{lemma}
	\begin{rmk}[Assumptions for Structured Representation]
		\label{rmk:obs_before_query} 
		"Standard" interventions, 
		as in example \ref{example:iid_do_intervention},
		and even transport-questions (§\ref{apdx:translate_mz_transport})
		can be described as single-level statements in the sense
		of the hypothesis of Lemma \ref{lemma:assoc_structured_query}.
		Further, for interventional formulations (like the do-calculus),
		$\tilde{X}=\emptyset$ (see §\ref{apdx:iid}),
		so the "no bypassing" assumption is always trivially
		satisfied. 
		In other cases (parts of) $\tilde{X}$
		can be included into $\tilde{Y}$ and disintegrated
		at the end, see §\ref{sec:queries_beyond_basic}.
		However nodes might be allowed in $\tilde{X}$ but not
		in $\tilde{Y}$, see example \ref{example:random_walk}.
		In the multi-level case the hypothesis of
		Lemma \ref{lemma:assoc_structured_query}
		can be weakened (Lemma \ref{lemma:assoc_structured_query:apdx}),
		and it becomes an if and only if statement.
	\end{rmk}
	
	Structured queries can be identified using the previously
	introduced machinery.
	
	\begin{thm}
		\label{thm:query_id}
		Given a structured query
		$(\psi,Y)$,
		using
		$L:=\nodesInnerProper\setminus Y$,
		if there is a regular functional $F$ computing
		$\mu(\graph,L)=F[\mu_1,\ldots,\mu_n]$
		with $\mu_1,\ldots,\mu_n$ all identifiable,
		then the underlying query
		$q(\psi,Y)$ is identifiable.
		
		\emph{Remark:}
		$F$ can be constructed as any finite sequence of graphical
		operations from §\ref{sec:graphical_operations},
		$\mu_1,\ldots,\mu_n$ can be obtained from any back-door free families
		by
		Thm.\ \ref{thm:extract_from_backdoor_complete}.
	\end{thm}

	\begin{algorithm}[tb]
		\renewcommand{\thealgorithm}{EDAIdentify}
		\caption{Extract--Decompose--Assemble Identification}
		\label{algo:id}
		\textbf{Input:} A query $\mu(\graphStructural,L)$.\\
		\textbf{Output:} A (potentially empty) set of identification
		strategies.
		\begin{enumerate}[label=\arabic*)]
			\item\label{algo:id:E}
			\emph{Extract:}\\
			$\mathcal{B}$ := \ref{algo:fcs}.\\
			$\knowledgeSet$ := Apply Thm.\ \ref{thm:extract_from_backdoor_complete}
			to elements of $\mathcal{B}$.
			\item\label{algo:id:D}
			\emph{Decompose:}\\
			$\knowledgeSet'$ := Apply \ref{algo:decomp} to elements of $\knowledgeSet$.
			\item\label{algo:id:A}
			\emph{Assemble:}\\
			\textbf{Return}
			\ref{algo:svs}($\knowledgeSet' \cup_{\tilde{f}\in\knownFunction} k(\tilde{f}) \cup_{i\in\viewportEventual} k(i)$,
			$(\graphStructural,L)$).
			\Comment{$k(i)$ are included here for completeness, they are relevant
				only to the multi-level case (cf.\ 
				Def.\ \ref{def:knowledge_set_extracted:apdx}).}
		\end{enumerate}
	\end{algorithm}
	
	A systematic approach to find identification strategies
	for a given query is sketched in Algo.\ \ref{algo:id}.
	First, minimal backdoor-free families are constructed.
	Then $\mu(\graph, L)$ (obtained by Thm.\ \ref{thm:extract_from_backdoor_complete})
	are decomposed into smaller parts (repeated application of
	Lemmas \ref{def:simplify_simp}, \ref{lemma:c_components:mt})
	and finally assembled by gluing (Lemma \ref{lemma:glue:mt})
	into the query-graph of interest.
	This procedure collects all information required
	for the application of Thm.\ \ref{thm:query_id}.
	This approach is evidently only complete if graphical
	operations are complete (Conj.\ \ref{conjecture:regular_graphical})
	and additionally:
	\begin{conjecture}[Decomposition Before Assembly]
		\label{conjecture:decomp_first}
		If a regular computation is possible by graphical operations,
		then it is possible to first perform all required decompositions,
		via Lemmas \ref{def:simplify_simp}, \ref{lemma:c_components:mt},
		followed by all required gluing operations, Lemma \ref{lemma:glue:mt},
		afterwards.
	\end{conjecture}
	
	We leave a detailed analysis of algorithms and their properties to future work,
	we do however give multiple partial results
	(like Lemmas
	\ref{lemma:monotonicity_of_L}, \ref{lemma:monotonicity_of_Anc},
	\ref{lemma:c_conn_extr_enough})
	that clarify why
	some simplifications in the algorithms presented in the appendix
	are possible without loss of generality.
	Summarizing, the main open questions concerning completeness
	are in the already highlighted conjectures:
	
	\begin{conjectureNamed}[Completeness]
		If Conjectures \ref{conjecture:regular_computation},
		\ref{conjecture:regular_graphical}, \ref{conjecture:extraction},
		\ref{conjecture:decomp_first}
		hold true, then
		in absence of selection-bias, a structured query is identifiable
		if and only if Algo.\ \ref{algo:id} returns a non-empty set.
	\end{conjectureNamed}

	\subsection{Knowledge-Closures and Identification} 
	\label{sec:knowledge_closures}

	At the end of §\ref{sec:extraction_from_data}, we summarized
	the extracted knowledge as a set $\knowledgeSetBasic$
	(Def. \ref{def:knowledge_set_extracted}).
	Next, we ask, from an abstract perspective, when is a query identifiable
	relative to $\knowledgeSetBasic$?
	Algorithmic completeness, relative to extraction-results
	(\ie modular completeness of computation only),
	can be phrased as
	the recovery of a closure under a given regularity condition.
	
	\begin{definition}[Closures of Knowledge Sets]
		Given a knowledge-set $\knowledgeSet$ and
		a class $\mathcal{R}$ ("regularity class")
		of functionals closed under finite compositions,
		the $\mathcal{R}$-closure of $\knowledgeSet$,
		denoted $\bar{\knowledgeSet}^{\mathcal{R}}$,
		is the set of elements of knowledge
		$\mathcal{R}$-identified from $\knowledgeSet$:
		\begin{equation*}
			\bar{\knowledgeSet}^{\mathcal{R}}
			\halfquad:=\halfquad
			\{\halfquad
			k'
			\halfquad|\halfquad	
			\exists F\in\mathcal{R},
			\exists k_1,\ldots,k_n \in \knowledgeSet:
			k'=F[k_1,\ldots,k_n]
			\halfquad\}\txt.
		\end{equation*}
		If $\mathcal{R}$ is the set of (standard-)regular functionals in the sense
		of Def.\ \ref{def:regular_functionals},
		then we simply write $\bar{\knowledgeSet}$ (dropping the superscript),
		and $\mu$ is identifiable (Def. \ref{def:identification})
		if $\mu\in\knowledgeSetBasicCompletion$.
	\end{definition}
	\begin{example}[Instrumental Variables]
		\label{example:instrumental_variables}
		Assume there are subsets
		$\mathcal{F}_{\txt{mono}}, \mathcal{F}_{\txt{epi}}
		\subset \knownFunction \cup \observedFunction$
		of left- and right-cancelative kernels respectively.
		Left-cancelative means
		$\mu \in \mathcal{F}_{\txt{mono}}$ iff
		$\nu \circ \mu = \nu' \circ \mu$ $\Rightarrow$ $\nu = \nu'$
		for all $\nu \in \knownFunction \cup \observedFunction$;
		similar for right-cancelative.
		Define $\mathcal{R}^{\txt{IV}}$ as the class of
		regular functionals (Def.\ \ref{def:regular_functionals}),
		functionals of the form $\nu \circ \mu \mapsto \nu$ if
		$\mu\in\mathcal{F}_{\txt{mono}}$
		and $\nu \circ \mu \mapsto \mu$ if 
		$\nu\in\mathcal{F}_{\txt{epi}}$
		and compositions of such functionals.
		These operations can be iterated together
		with Lemmas \ref{def:simplify_simp}, \ref{lemma:c_components:mt}
		in the decomposition-step of Algo.\ \ref{algo:id}.
		
		For example if all elements of $\knownFunction \cup \observedFunction$
		are linear, then $\mathcal{F}_{\txt{mono}} = \mathcal{F}_{\txt{epi}}
		=\knownFunction \cup \observedFunction$ and this deconvolution is
		called an "instrumental variable" argument \citep[p.\,247]{PearlBook}.
		Formally, the distinction should be made between
		a priori knowledge of $\mathcal{F}_{\txt{mono}}, \mathcal{F}_{\txt{epi}}$
		and their existence (since not every $f_J$ is identifiable,
		$f_J\in\mathcal{F}_{\txt{mono}}, \mathcal{F}_{\txt{epi}}$
		may not be testable).
		
		Since (by definition) regular functionals are in $\mathcal{R}^{\txt{IV}}$,
		automatically
		$\bar{\knowledgeSet} \subset \bar{\knowledgeSet}^{\mathcal{R}^{\txt{IV}}}$
		for all $\knowledgeSet$.
	\end{example}

	\subsection{Relation to the IID-Case} 
	\label{sec:relate_to_iid}

	The do-calculus applies,
	if data are IID, interventions are do-interventions,
	and missingness is uniform (an IID-variable $v\in\IVars$ is either observed
	in all samples $s\in\ISample$ or never).
	Under these conditions,
	the id-algorithm \citep{tian2002general} is known to be complete
	\citep{shpitser2006identification2}.
	In §\ref{apdx:iid:standard_results}, we show in detail
	that Algo.\ \ref{algo:id} returns an empty set only if there is a "hedge".
	By a criterion of \citep{shpitser2006identification2}
	this implies that our algorithm is complete,
	at least for unconditional do-interventions, thus matching the id-algorithm.
	For conditional queries, some care should be taken.
	Rule 2 of the do-calculus exploits the internal structure
	of do-interventions, which we intentionally exclude.
	Conditional queries are treated in detail by \citep{shpitser2006identification1}
	in the IID case. We can reproduce the results of
	\citep{shpitser2006identification1},
	if we match this internal structure assumption.
	This
	can readily be achieved (see §\ref{sec:queries_beyond_basic})
	without substantial overhead, see §\ref{apdx:iid_conditional}.
	The easiest technical approach is by "conditioning" on the intervened
	variable by adding it to the otherwise empty $\tilde{X}$ of the basic query,
	which is allowed and makes sense in our formalism.
	
	We also compare our approach in detail to the mz-transportability
	setup \citep{pearl2022external, Bareinboim2016TransportOverview, 
		bareinboim2012TransportCompleteness} in §\ref{apdx:iid:mz_transport}.
	This setup aims to transport experimental and non-experimental
	information between IID-contexts (with a number of restrictions,
	for example all contexts share a graph, and hiddenness is cross-context
	uniform). It uses "selection-variables" (often called context-variables in
	other sources), which (in similar form) are at the foundation of many influential
	multi-context approaches to causal questions \citep{CD-NOD, JCI, JPCMCI}.
	We find that also here Algo.\ \ref{algo:id} seems to be complete
	(by a criterion of \citep{bareinboim2012TransportCompleteness}),
	matching the algorithm of
	\citep{Bareinboim2013TransportAlgo,bareinboim2012TransportCompleteness}.
	
	Finally, somewhat surprisingly, we can also readily reproduce
	the mediation formula for natural direct effects (NDE)
	(or rather: the formula for underlying
	conditional probabilities).
	We do not actually
	consider counter-factuals.
	However, our query-formulation is general enough to capture
	the question motivating NDE directly, without the detour
	through counterfactuals, see §\ref{apdx:iid:mediation}.
	
	We can readily reproduce many relevant results from the literature,
	importantly we can do so in a unified approach (and algorithm).
	While literature on the topic of causal effect-estimation
	is dominated by IID-setups and do-interventions
	(and hence so is this comparison-section),
	the main purpose and strength of our formalism is its 
	applicability far from IID and its flexible query formulation.
	While we do not know if our algorithm is complete,
	it does not lead to a regression below the state of the art
	in standard setups.

	\section{Conclusion} 
	
	We introduced a simple and widely applicable abstract language to
	describe causal reasoning §\ref{sec:models}.
	This weak imposed structure proves sufficient to learn
	and reason about a model §\ref{sec:structured_kernels},
	§\ref{sec:extraction_from_data} in order to ultimately
	identify a notion of causal queries §\ref{sec:queries},
	that includes many conventionally
	considered questions as special cases.
	
	While some questions remain open, concerning selection-bias,
	correlated missingness,
	counterfactuals and also about the completeness of our identification strategy
	in general, we do re-discovery many core
	structures like c-components, c-trees, c-forests and hedges
	known to be pivotal to IID do-interventional queries
	\citep{tian2002general,shpitser2006identification2}
	from a new perspective. This not only sheds new light on these
	objects, but it also validates our results
	against known special cases.
	Also, the ability to (correctly) reproduce the mediation-formula
	for natural direct effects, and mz-transport results
	is quite encouraging.
	
	From a broader perspective, it would certainly be interesting
	to understand how and to what degree the models employed by
	our formalism can be learned from data. Initial
	ideas and feasibility are explored in §\ref{apdx:structure_discovery}.
	However, specifying the right structure to learn
	must necessarily precede finding ways of learning it.
	
	Returning to the example from the introduction (Fig.\ \ref{fig:intro}),
	note that the panels highlighted as "learn" can indeed be
	identified by Thm.\ \ref{thm:extract_from_backdoor_complete},
	and the query can be glued (Lemma \ref{lemma:glue:mt}) from these
	partial results and the known $Y_{\txt{init}}$. In a last step,
	variables not of interest can be marginalized.
	A detailed construction of a sound statistical estimator
	is given in example \ref{example:real_world_questions}\,b.
	
	Aspirationally, these results may also contribute to
	better and better explainable AI. 
	For example in reinforcement learning,
	an agent that requires curated datasets is not autonomous.
	So IIDness makes little conceptual sense for learning
	complex world-models. At the same time exploration requires
	efficient offline evaluation of possible future trajectories,
	brought into the scope of causal reasoning by our flexible queries.
	Finally, understanding and explaining AI-behavior strongly hinges on
	human-accessible structure given to learned solutions.

	\RenewDocumentCommand{\nodesInnerProper}{}{\nodes_{\txt{proper}}}
	\RenewDocumentCommand{\nodesOuterProper}{}{\nodes_{\txt{extern}}}
	\RenewDocumentCommand{\nodesOuterFixed}{}{\nodes_{\txt{pinned}}}

	\section*{Acknowledgments}
	
	The authors would like to thank Wiebke Günther
	for very helpful and detailed feedback on the paper.
	J.R. has received funding from the European Research Council (ERC) Starting Grant CausalEarth under the European Union’s Horizon 2020 research and innovation program (Grant Agreement No. 948112).
	J.R. and M.R. have received funding from the European Union’s Horizon 2020 research and innovation program under grant agreement No 101003469 (XAIDA).
	
	\bibliography{./BibTex}

\begin{thebibliography}{63}
\providecommand{\natexlab}[1]{#1}
\providecommand{\url}[1]{\texttt{#1}}
\expandafter\ifx\csname urlstyle\endcsname\relax
  \providecommand{\doi}[1]{doi: #1}\else
  \providecommand{\doi}{doi: \begingroup \urlstyle{rm}\Url}\fi

\bibitem[Assaad(2024)]{assaad2024causal}
C.~K. Assaad.
\newblock Causal reasoning in difference graphs.
\newblock \emph{arXiv preprint arXiv:2411.01292}, 2024.

\bibitem[Avin et~al.(2005)Avin, Shpitser, and Pearl]{avin2005identifiability}
C.~Avin, I.~Shpitser, and J.~Pearl.
\newblock Identifiability of path-specific effects.
\newblock In \emph{Proceedings of the 19th International Joint Conference on
  Artificial Intelligence}, IJCAI'05, page 357–363, San Francisco, CA, USA,
  2005. Morgan Kaufmann Publishers Inc.

\bibitem[Balsells-Rodas et~al.(2023)Balsells-Rodas, Wang, and
  Li]{BalsellsRodas2023}
C.~Balsells-Rodas, Y.~Wang, and Y.~Li.
\newblock On the identifiability of markov switching models, 2023.

\bibitem[Bareinboim and Pearl(2012)]{bareinboim2012TransportCompleteness}
E.~Bareinboim and J.~Pearl.
\newblock Transportability of causal effects: Completeness results.
\newblock In \emph{Proceedings of the AAAI Conference on Artificial
  Intelligence}, volume~26, pages 698--704, 2012.

\bibitem[Bareinboim and Pearl(2013)]{Bareinboim2013TransportAlgo}
E.~Bareinboim and J.~Pearl.
\newblock A general algorithm for deciding transportability of experimental
  results.
\newblock \emph{Journal of Causal Inference}, 1\penalty0 (1):\penalty0
  107--134, 2013.

\bibitem[Bareinboim and Pearl(2016)]{Bareinboim2016TransportOverview}
E.~Bareinboim and J.~Pearl.
\newblock Causal inference and the data-fusion problem.
\newblock \emph{Proceedings of the National Academy of Sciences}, 113:\penalty0
  7345 -- 7352, 2016.

\bibitem[Bongers et~al.(2021)Bongers, Forr{\'e}, Peters, and
  Mooij]{BongersCyclic}
S.~Bongers, P.~Forr{\'e}, J.~Peters, and J.~M. Mooij.
\newblock Foundations of structural causal models with cycles and latent
  variables.
\newblock \emph{The Annals of Statistics}, 49\penalty0 (5):\penalty0
  2885--2915, 2021.

\bibitem[Chang and Pollard(1997)]{chang1997conditioning}
J.~T. Chang and D.~Pollard.
\newblock Conditioning as disintegration.
\newblock \emph{Statistica Neerlandica}, 51\penalty0 (3):\penalty0 287--317,
  1997.

\bibitem[Colombo et~al.(2014)Colombo, Maathuis, et~al.]{PCstable}
D.~Colombo, M.~H. Maathuis, et~al.
\newblock Order-independent constraint-based causal structure learning.
\newblock \emph{J. Mach. Learn. Res.}, 15\penalty0 (1):\penalty0 3741--3782,
  2014.

\bibitem[Corander et~al.(2019)Corander, Hyttinen, Kontinen, Pensar, and
  V{\"a}{\"a}n{\"a}nen]{LDAG_logical}
J.~Corander, A.~Hyttinen, J.~Kontinen, J.~Pensar, and J.~V{\"a}{\"a}n{\"a}nen.
\newblock A logical approach to context-specific independence.
\newblock \emph{Annals of Pure and Applied Logic}, 170\penalty0 (9):\penalty0
  975--992, 2019.

\bibitem[Correa and Bareinboim(2020)]{correa2020general}
J.~Correa and E.~Bareinboim.
\newblock General transportability of soft interventions: Completeness results.
\newblock \emph{Advances in Neural Information Processing Systems},
  33:\penalty0 10902--10912, 2020.

\bibitem[de~Aguas et~al.(2025)de~Aguas, Henckel, Pensar, and
  Biele]{de2025causal}
J.~de~Aguas, L.~Henckel, J.~Pensar, and G.~Biele.
\newblock Causal inference amid missingness-specific independencies and
  mechanism shifts.
\newblock \emph{arXiv preprint arXiv:2506.15441}, 2025.

\bibitem[Faller et~al.(2024)Faller, Vankadara, Mastakouri, Locatello, and
  Janzing]{faller2024self}
P.~M. Faller, L.~C. Vankadara, A.~A. Mastakouri, F.~Locatello, and D.~Janzing.
\newblock Self-compatibility: Evaluating causal discovery without ground truth.
\newblock In \emph{International Conference on Artificial Intelligence and
  Statistics}, pages 4132--4140. PMLR, 2024.

\bibitem[Gelman and Hill(2006)]{Gelman2006}
A.~Gelman and J.~Hill.
\newblock \emph{Data analysis using regression and multilevel/hierarchical
  models}.
\newblock Cambridge university press, 2006.

\bibitem[G{\"u}nther et~al.(2023)G{\"u}nther, Ninad, and Runge]{JPCMCI}
W.~G{\"u}nther, U.~Ninad, and J.~Runge.
\newblock Causal discovery for time series from multiple datasets with latent
  contexts.
\newblock In \emph{Uncertainty in Artificial Intelligence}, pages 766--776.
  PMLR, 2023.

\bibitem[G{\"u}nther et~al.(2024)G{\"u}nther, Popescu, Rabel, Ninad, Gerhardus,
  and Runge]{gunther2024causal}
W.~G{\"u}nther, O.-I. Popescu, M.~Rabel, U.~Ninad, A.~Gerhardus, and J.~Runge.
\newblock Causal discovery with endogenous context variables.
\newblock \emph{Advances in Neural Information Processing Systems},
  37:\penalty0 36243--36284, 2024.

\bibitem[Guo and Perkovi{\'c}(2022)]{guo2022efficient}
F.~R. Guo and E.~Perkovi{\'c}.
\newblock Efficient least squares for estimating total effects under linearity
  and causal sufficiency.
\newblock \emph{Journal of Machine Learning Research}, 23\penalty0
  (104):\penalty0 1--41, 2022.

\bibitem[Holland(1986)]{holland86}
P.~W. Holland.
\newblock Statistics and causal inference.
\newblock \emph{Journal of the American Statistical Association}, 81\penalty0
  (396):\penalty0 945--960, 1986.
\newblock ISSN 01621459, 1537274X.

\bibitem[Huang et~al.(2020)Huang, Zhang, Zhang, Ramsey, Sanchez-Romero,
  Glymour, and Sch{\"o}lkopf]{CD-NOD}
B.~Huang, K.~Zhang, J.~Zhang, J.~Ramsey, R.~Sanchez-Romero, C.~Glymour, and
  B.~Sch{\"o}lkopf.
\newblock Causal discovery from heterogeneous/nonstationary data.
\newblock \emph{The Journal of Machine Learning Research}, 21\penalty0
  (1):\penalty0 3482--3534, 2020.

\bibitem[Hyttinen et~al.(2018)Hyttinen, Pensar, Kontinen, and
  Corander]{LDAG_learning}
A.~Hyttinen, J.~Pensar, J.~Kontinen, and J.~Corander.
\newblock Structure learning for bayesian networks over labeled dags.
\newblock In \emph{International conference on probabilistic graphical models},
  pages 133--144. PMLR, 2018.

\bibitem[Janzing et~al.(2023)Janzing, Faller, and
  Vankadara]{janzing2023reinterpreting}
D.~Janzing, P.~M. Faller, and L.~C. Vankadara.
\newblock Reinterpreting causal discovery as the task of predicting unobserved
  joint statistics.
\newblock \emph{arXiv preprint arXiv:2305.06894}, 2023.

\bibitem[Kallenberg(1997)]{kallenberg1997foundations}
O.~Kallenberg.
\newblock \emph{Foundations of modern probability}.
\newblock Springer, 1997.

\bibitem[Kallenberg(2005)]{kallenberg2005probabilistic}
O.~Kallenberg.
\newblock \emph{Probabilistic symmetries and invariance principles}.
\newblock Springer, 2005.

\bibitem[Lattimore et~al.(2016)Lattimore, Lattimore, and
  Reid]{lattimore2016causal}
F.~Lattimore, T.~Lattimore, and M.~D. Reid.
\newblock Causal bandits: Learning good interventions via causal inference.
\newblock \emph{Advances in neural information processing systems}, 29, 2016.

\bibitem[Mameche et~al.(2025)Mameche, Cornanguer, Ninad, and
  Vreeken]{mameche2025spacetime}
S.~Mameche, L.~Cornanguer, U.~Ninad, and J.~Vreeken.
\newblock Spacetime: Causal discovery from non-stationary time series.
\newblock \emph{arXiv preprint arXiv:2501.10235}, 2025.

\bibitem[Margueritte et~al.(2026)Margueritte, Balc{\i}o{\u{g}}lu, Krijthe,
  Zachariah, and Johansson]{margueritte2026learning}
A.-U. Margueritte, A.~Z. Balc{\i}o{\u{g}}lu, J.~Krijthe, D.~Zachariah, and
  F.~D. Johansson.
\newblock Learning plug-in surrogate endpoints for randomized experiments.
\newblock \emph{arXiv preprint arXiv:2605.12051}, 2026.

\bibitem[Mooij et~al.(2020)Mooij, Magliacane, and Claassen]{JCI}
J.~M. Mooij, S.~Magliacane, and T.~Claassen.
\newblock Joint causal inference from multiple contexts.
\newblock \emph{The Journal of Machine Learning Research}, 21\penalty0
  (1):\penalty0 3919--4026, 2020.

\bibitem[Murphy(2022)]{murphy2022probabilistic}
K.~P. Murphy.
\newblock \emph{Probabilistic machine learning: an introduction}.
\newblock MIT press, 2022.

\bibitem[Orujlu et~al.(2026)Orujlu, Willig, and
  Ze{\v{c}}evi{\'c}]{orujlu2026causal}
T.~Orujlu, M.~Willig, and M.~Ze{\v{c}}evi{\'c}.
\newblock Causal substrates: A set-theoretic foundation of causal inference.
\newblock In \emph{UAI 2026 Workshop on Causality for Decision Making}, 2026.

\bibitem[Park et~al.(2023)Park, Buchholz, Sch\"{o}lkopf, and
  Muandet]{CausalSpaces}
J.~Park, S.~Buchholz, B.~Sch\"{o}lkopf, and K.~Muandet.
\newblock A measure-theoretic axiomatisation of causality.
\newblock In A.~Oh, T.~Naumann, A.~Globerson, K.~Saenko, M.~Hardt, and
  S.~Levine, editors, \emph{Advances in Neural Information Processing Systems},
  volume~36, pages 28510--28540. Curran Associates, Inc., 2023.
\newblock \doi{10.52202/075280-1239}.

\bibitem[Pearl(2000)]{PearlBook}
J.~Pearl.
\newblock \emph{Causality: Models, reasoning and inference}.
\newblock Cambride University Press, 2000.

\bibitem[Pearl(2001)]{Pearl2001}
J.~Pearl.
\newblock Direct and indirect effects.
\newblock \emph{Proceedings of the Seventeenth Conference on Uncertainty in
  Artificial intelligence}, 2001.

\bibitem[Pearl(2012)]{Pearl2012}
J.~Pearl.
\newblock The causal mediation formula—a guide to the assessment of pathways
  and mechanisms.
\newblock \emph{Prevention science}, 13:\penalty0 426--436, 2012.

\bibitem[Pearl and Bareinboim(2014)]{SelectionVars}
J.~Pearl and E.~Bareinboim.
\newblock External validity: From do-calculus to transportability across
  populations.
\newblock \emph{Statistical Science}, 29\penalty0 (4):\penalty0 579--595, 2014.

\bibitem[Pearl and Bareinboim(2022)]{pearl2022external}
J.~Pearl and E.~Bareinboim.
\newblock External validity: From do-calculus to transportability across
  populations.
\newblock In \emph{Probabilistic and causal inference: The works of Judea
  Pearl}, pages 451--482. 2022.

\bibitem[Pensar et~al.(2015)Pensar, Nyman, Koski, and
  Corander]{LDAG_definition}
J.~Pensar, H.~Nyman, T.~Koski, and J.~Corander.
\newblock Labeled directed acyclic graphs: a generalization of context-specific
  independence in directed graphical models.
\newblock \emph{Data mining and knowledge discovery}, 29:\penalty0 503--533,
  2015.

\bibitem[Perkovi{\'c} et~al.(2015)Perkovi{\'c}, Textor, Kalisch, and
  Maathuis]{Perkovic2015}
E.~Perkovi{\'c}, J.~Textor, M.~Kalisch, and M.~H. Maathuis.
\newblock A complete generalized adjustment criterion.
\newblock \emph{arXiv preprint arXiv:1507.01524}, 2015.

\bibitem[Peters et~al.(2016)Peters, B{\"u}hlmann, and
  Meinshausen]{peters2016causal}
J.~Peters, P.~B{\"u}hlmann, and N.~Meinshausen.
\newblock Causal inference by using invariant prediction: identification and
  confidence intervals.
\newblock \emph{Journal of the Royal Statistical Society Series B: Statistical
  Methodology}, 78\penalty0 (5):\penalty0 947--1012, 2016.

\bibitem[Peters et~al.(2017)Peters, Janzing, and Sch{\"o}lkopf]{Elements}
J.~Peters, D.~Janzing, and B.~Sch{\"o}lkopf.
\newblock \emph{Elements of causal inference: foundations and learning
  algorithms}.
\newblock The MIT Press, 2017.

\bibitem[Rabel and Runge(2026)]{rabel2026contextspecificcausalgraphdiscovery}
M.~Rabel and J.~Runge.
\newblock Context-specific causal graph discovery with unobserved contexts:
  Non-stationarity, regimes and spatio-temporal patterns, 2026.

\bibitem[Rabel et~al.(2025)Rabel, Günther, Runge, and Gerhardus]{EndoTheory}
M.~Rabel, W.~Günther, J.~Runge, and A.~Gerhardus.
\newblock Causal modeling in multi-context systems: Distinguishing multiple
  context-specific causal graphs which account for observational support.
\newblock \emph{arXiv preprint arXiv:2410.20405}, 2025.

\bibitem[Rahmani and Frossard(2023)]{rahmani2023castor}
A.~Rahmani and P.~Frossard.
\newblock Castor: Causal temporal regime structure learning.
\newblock 2023.

\bibitem[Robins and Richardson(2010)]{robins2010alternative}
J.~M. Robins and T.~S. Richardson.
\newblock Alternative graphical causal models and the identification of direct
  effects.
\newblock \emph{Causality and psychopathology: Finding the determinants of
  disorders and their cures}, 84:\penalty0 103--158, 2010.

\bibitem[Rodas et~al.(2021)Rodas, Tu, and Kjellstrom]{rodas2021causal}
C.~B. Rodas, R.~Tu, and H.~Kjellstrom.
\newblock Causal discovery from conditionally stationary time-series.
\newblock \emph{arXiv preprint arXiv:2110.06257}, 2021.

\bibitem[Rojas-Carulla et~al.(2018)Rojas-Carulla, Sch{\"o}lkopf, Turner, and
  Peters]{RojasCarulla2018a}
M.~Rojas-Carulla, B.~Sch{\"o}lkopf, R.~Turner, and J.~Peters.
\newblock Invariant models for causal transfer learning.
\newblock \emph{The Journal of Machine Learning Research}, 19\penalty0
  (1):\penalty0 1309--1342, 2018.

\bibitem[Rothenh{\"a}usler et~al.(2021)Rothenh{\"a}usler, Meinshausen,
  B{\"u}hlmann, and Peters]{AnchorRegression}
D.~Rothenh{\"a}usler, N.~Meinshausen, P.~B{\"u}hlmann, and J.~Peters.
\newblock Anchor regression: Heterogeneous data meet causality.
\newblock \emph{Journal of the Royal Statistical Society Series B: Statistical
  Methodology}, 83\penalty0 (2):\penalty0 215--246, 2021.

\bibitem[Rubin(1974)]{Rubin1974a}
D.~B. Rubin.
\newblock Estimating causal effects of treatments in randomized and
  nonrandomized studies.
\newblock \emph{J. Educ. Psychol.}, 66\penalty0 (5):\penalty0 688--701, 1974.
\newblock ISSN 00220663.
\newblock \doi{10.1037/h0037350}.

\bibitem[Rubin(1976)]{rubin1976inference}
D.~B. Rubin.
\newblock Inference and missing data.
\newblock \emph{Biometrika}, 63\penalty0 (3):\penalty0 581--592, 1976.

\bibitem[Runge(2020)]{pcmci_plus}
J.~Runge.
\newblock Discovering contemporaneous and lagged causal relations in
  autocorrelated nonlinear time series datasets.
\newblock In \emph{Conference on Uncertainty in Artificial Intelligence}, pages
  1388--1397. PMLR, 2020.

\bibitem[Runge(2021)]{OSets}
J.~Runge.
\newblock Necessary and sufficient graphical conditions for optimal adjustment
  sets in causal graphical models with hidden variables.
\newblock \emph{Advances in Neural Information Processing Systems},
  34:\penalty0 15762--15773, 2021.

\bibitem[Runge et~al.(2019)Runge, Nowack, Kretschmer, Flaxman, and
  Sejdinovic]{pcmci}
J.~Runge, P.~Nowack, M.~Kretschmer, S.~Flaxman, and D.~Sejdinovic.
\newblock Detecting and quantifying causal associations in large nonlinear time
  series datasets.
\newblock \emph{Science advances}, 5\penalty0 (11):\penalty0 eaau4996, 2019.

\bibitem[Runge et~al.(2023)Runge, Gerhardus, Varando, Eyring, and
  Camps-Valls]{Runge2023}
J.~Runge, A.~Gerhardus, G.~Varando, V.~Eyring, and G.~Camps-Valls.
\newblock Causal inference for time series.
\newblock \emph{Nature Reviews Earth \& Environment}, pages 1--19, 2023.

\bibitem[Saggioro et~al.(2020)Saggioro, de~Wiljes, Kretschmer, and
  Runge]{Saggioro2020}
E.~Saggioro, J.~de~Wiljes, M.~Kretschmer, and J.~Runge.
\newblock Reconstructing regime-dependent causal relationships from
  observational time series.
\newblock \emph{Chaos: An Interdisciplinary Journal of Nonlinear Science},
  30\penalty0 (11), 2020.

\bibitem[Shah and Peters(2020)]{shah2020hardness}
R.~D. Shah and J.~Peters.
\newblock The hardness of conditional independence testing and the generalised
  covariance measure.
\newblock \emph{The Annals of Statistics}, 48\penalty0 (3):\penalty0 1514,
  2020.

\bibitem[Shpitser(2013)]{Shpitser2013}
I.~Shpitser.
\newblock Counterfactual graphical models for longitudinal mediation analysis
  with unobserved confounding.
\newblock \emph{Cognitive science}, 37\penalty0 (6):\penalty0 1011--1035, 2013.

\bibitem[Shpitser(2023)]{shpitser2023does}
I.~Shpitser.
\newblock When does the id algorithm fail?
\newblock \emph{arXiv preprint arXiv:2307.03750}, 2023.

\bibitem[Shpitser and Pearl(2006{\natexlab{a}})]{shpitser2006identification1}
I.~Shpitser and J.~Pearl.
\newblock Identification of conditional interventional distributions.
\newblock In \emph{Uncertainty in Artificial Intelligence}, volume~22, pages
  437--444, 2006{\natexlab{a}}.
\newblock arXiv:1206.6876.

\bibitem[Shpitser and Pearl(2006{\natexlab{b}})]{shpitser2006identification2}
I.~Shpitser and J.~Pearl.
\newblock Identification of joint interventional distributions in recursive
  semi-markovian causal models.
\newblock In \emph{AAAI}, pages 1219--1226, 2006{\natexlab{b}}.

\bibitem[Shpitser and Pearl(2012)]{IDAlgoCF}
I.~Shpitser and J.~Pearl.
\newblock What counterfactuals can be tested.
\newblock \emph{arXiv preprint arXiv:1206.5294}, 2012.

\bibitem[Spirtes and Glymour(1991)]{PCalgo}
P.~Spirtes and C.~Glymour.
\newblock An algorithm for fast recovery of sparse causal graphs.
\newblock \emph{Social Science Computer Review}, 9:\penalty0 62--72, 1991.

\bibitem[Spirtes et~al.(2001)Spirtes, Glymour, and
  Scheines]{spirtes2001causation}
P.~Spirtes, C.~Glymour, and R.~Scheines.
\newblock \emph{Causation, prediction, and search}.
\newblock MIT press, 2001.

\bibitem[Tian and Pearl(2002)]{tian2002general}
J.~Tian and J.~Pearl.
\newblock A general identification condition for causal effects.
\newblock In \emph{Aaai/iaai}, pages 567--573, 2002.

\bibitem[Vapnik(2006)]{VapnikEstimation}
V.~Vapnik.
\newblock Estimation of dependences based on empirical data.
\newblock In M.~Jordan, J.~Kleinberg, and B.~Sch\"{o}lkopf, editors,
  \emph{Vapnik: Estimation of Dependences Based on Empirical Data}, Information
  Science and Statistics, pages 1--399. Springer Science + Business Media, New
  York, second edition, 2006.
\newblock Reprint of 1982 edition.

\end{thebibliography}

	\appendix

	\section*{Reading Guide}
	\addcontentsline{toc}{section}{Reading Guide to the Appendix}
	
	This appendix is rather long and in places very technical, so we
	start with a brief overview to help the reader orient themselves.
	A table of contents is included at the very end of the paper.
	
	The two technically most interesting and relevant results are
	the proofs of graphical operations (§\ref{sec:graphical_operations})
	and of Thm.\ \ref{thm:extract_from_backdoor_complete} (extraction
	from backdoor-free families).
	It should be possible to start reading the respective proofs,
	then reference back to other material of the appendices
	where required (or of interest).	
	Graphical operations are proved in §\ref{apdx:tranform_graphs}.
	The main proof technique is a reduction to statements about
	certain univariate kernels called "atoms" in §\ref{apdx:tranform_graphs}.
	The proof of Thm.\ \ref{thm:extract_from_backdoor_complete}
	in §\ref{apdx:id_from_embeddings}
	relies on results about these "atoms" as well, hence it is
	recommendable to get an overview of results shown in
	§\ref{apdx:tranform_graphs} before starting to read
	the proof of Thm.\ \ref{thm:extract_from_backdoor_complete}
	(this is especially true for proof-steps 2 and 3).
	Some technical details on the relation of data-sets to
	models relevant to the proof of Thm.\ \ref{thm:extract_from_backdoor_complete}
	have been detached into §\ref{apdx:models_obs_worlds}
	to make the proof more accessible.
	
	We removed multi-level results from the main-text for improving
	readability. Where differences arise (in §\ref{apdx:extraction_from_data}
	and §\ref{apdx:queries}) an initial sub-section introduces
	the multi-level setup as a simple modification of the main text.
	It should also be possible to read the corresponding
	appendices (especially §\ref{apdx:extraction_from_data})
	for the single-level case only, without much confusion
	(see introduction to §\ref{apdx:extraction_from_data};
	for the proof of Thm.\ \ref{thm:extract_from_backdoor_complete}
	all multi-level aspects are contained within a separate
	proof-step 3).
	For queries (§\ref{apdx:queries}),
	multi-level aspects are relevant for phrasing more
	general questions: In the single-level case,
	after marginalizing to the query
	$\realizedDistr(\tilde{X},\tilde{Y})=\shallowDistr(\tilde{X},\tilde{Y})$
	and many technical problems vanish. For an intuition
	why additional non-trivial statements can be made
	with multi-level structures see
	the example \ref{example:learn_pinned}.
	Further examples of multi-level queries and their
	usefulness are given in §\ref{apdx:multi_level_examples}.
	Especially to the reader familiar with
	results for the IID-case
	\citep{PearlBook,tian2002general,shpitser2006identification1}
	it may also be helpful to have a look at §\ref{apdx:iid},
	to better understand the query-formulation and IID-completeness.

	By separating the random elements of an observable world
	from the purely measure-theoretical model our approach
	can be made mathematically rigorous.
	To clarify notation and technical details,
	§\ref{apdx:ptheo_basics} builds the detailed technical
	connection from citeable text-book results to the
	notation and computational means used throughout
	the remainder of this paper.

	\section{Examples and Applications}
	
	This section provides detailed examples
	to help the reader in making the connection
	from abstract theory to concrete computation and estimation from data.
	
	\subsection{Example from the Introduction}
	
	We start with details on Fig.\ \ref{fig:intro} from the introduction.
	First (example \ref{example:real_world_questions}a)
	we provide a concrete real-world scenario where such a problem may arise.
	In applications, statistical questions are proxies to help the expert
	understand the true dynamics of a stochastic system better.
	A more flexible query formulation allows proxies closer to, and thus
	more informative on, the actual point of interest.
	Second (example \ref{example:real_world_questions}b)
	we exercise in detail how the framework presented in this paper
	is employed to write down a statistical estimator for the indicated
	query.
	Finally (example \ref{example:real_world_questions}c),
	we use this example to give further intuition
	on observational support problems
	(cf.\ Ass.\ \ref{ass:support}, Rmk.\ \ref{rmk:kernels_uniqueness}).

	\begin{example}[Real-World Questions]
		\label{example:real_world_questions}	
		\textbf{(a) Motivation\Slash{}Illustration:}
		
		We give a toy-example, further illustrated in Fig.\ \ref{fig:intro},
		motivated by production processes in factories.
		A robotic arm performs three separate actions in sequence on each
		workpiece, say welding $X$, active cooling $Z$ (of the workpiece),
		then moving $Y$
		the workpiece to the next station.
		We measure for example operating temperature (of the robotics)
		as a uni-variate time-series; for simplicity,
		on each performed step one data-point is recorded.
		There may be latent confounders, for example both welding
		and moving may be affected by variations in workpiece weight $L_1$.
		See Fig.\ \ref{fig:intro}, (1).
		
		We want to replace the active cooling module $Z$ by a newer
		version $Z_2$. We do not want to do so in the full factory setup,
		rather we want to qualify the part in a unit-test.
		We put a dummy heating module $X_2$ simulating the welder
		together with $Z_2$ on a robotic arm and repeatedly
		test the activation of $X_2$ followed by $Z_2$.
		However $X_2$ and $Z_2$ are confounded by a latent variable,
		for example the power-supply $L_2$ feeding $X_2$ and $Z_2$
		or the ambient temperature of the lab $L'_2$ where the
		test-setup runs is different from the factory.
		See Fig.\ \ref{fig:intro}, (2).
		We could easily modify $X_2$ from heating given a fixed
		wattage to heating to a specific fixed known temperature $T$,
		ensuring $P(X_{\PearlDo})=\delta[T]$ (a do-intervention
		in the sense of Pearl).
		See Fig.\ \ref{fig:intro}, (3).
		
		Does this analysis provide enough data to understand
		at what temperatures the continuously operating system
		will run if $Z$ is replaced by $Z_2$?
		We could analyze what stationary temperature is reached
		(§\ref{sec:queries_beyond_basic}), but for simplicity
		let's say we ask for the temperature after some
		finite number of cycles, say 4, as a proxy.
		See "predict" in Fig.\ \ref{fig:intro}.
		
		The intuitive answer seems simple: From each data-set
		we can learn building blocks (see "learn" in Fig.\ \ref{fig:intro})
		and piece them together to answer the question.
		Clearly, we want a systematic and mathematical approach however,
		that tells us for example what "building blocks" are as mathematical objects,
		what learning them requires, and when\Slash{}how "piecing them together"
		is allowed.
		Also we need tools to understand
		when and how such questions can be answered soundly.
		Indeed there are some subtleties: For example
		the transportability of the latent $L_2$ (or $L_2'$) from the unit-test
		to the integrated setup is relevant.
		A dependence on ambient temperature $L_2'$ would clearly require measurements
		to adjust for different temperature-distributions in lab and factory.
		A dependence on a power-supply $L_2$ may be transportable,
		\eg if the same power-supply is used.
		Hence, the latent $L_2$ in the query, despite not being
		a confounder, must be included to obtain a correct result.
		Reasoning about its relation to the test-setup is necessary.
		
		Standard frameworks do not seem to provide much guidance on how
		to even start on this problem.
		Some of the points made above (like transportability of $L_2$)
		are related to multi-context setups
		\citep{pearl2022external,bareinboim2012TransportCompleteness,JCI,CD-NOD,JPCMCI}.
		But combining such extensions to \emph{specific} setups,
		for complex problems, is a very non-trivial task (see main text),
		for deep conceptual reasons 
		§\ref{sec:intro:need_abstraction}.
		An example necessarily is very specific, our formalism is not
		in any way designed for, or limited to, this particular problem.
		This example only serves to illustrate how the conventional
		SCM\Slash{}PO perspectives struggle to
		capture even simple and plausible examples.
		
		\textbf{(b) Estimation\Slash{}Computation:}
		
		Still based on Fig.\ \ref{fig:intro}.
		We start from indices
		$I=(T\sqcup 3T) \sqcup (\{x,z,l\}\times I) \sqcup (\{x',z',l'\}\times I')$
		as a disjoint union of terms for the three datasets.
		The viewport is $\viewportEventual
		=T \sqcup (\{x,z\}\times I) \sqcup (\{x',z'\}\times I')$.
		We first focus on the first, dataset (1), 
		with indices in $t\in T=\mathbb{N}_0 \cup \{-1\}$ (the univariate time-series)
		and $t_L \in 3T$ (the latent every three time-steps);
		we will use the notations $t$ and $t_L$ to distinguish
		elements of both subsets.
		We formally prepend some initial state $Y_{\txt{init}}^{\txt{obs}}$
		at $t=-1$ (we assume we do not know what the mechanism
		$Y_{\txt{init}}^{\txt{obs}}$ is).
		We observe values
		$\randomVar{A}_0, \randomVar{A}_1, \ldots$ (starting at $t_0=0$)
		associated to $t\in T$.
		These observations have different associated mechanisms,
		$\randomVar{A}_t$ for $t\equiv 0 \mod 3$ the mechanism
		is $X$, for $t\equiv 1 \mod 3$ it is $Z$, for
		$t\equiv 2 \mod 3$ it is $Y$.
		
		There is a backdoor-free family of embeddings,
		embedding the local graph $\graph$ shown in the box
		"learn" (Fig.\ \ref{fig:intro}, dataset 1):
		Call the three inner nodes $n_x, n_l, n_y$ (aligned to $X, L_1, Y$)
		and the two outer nodes $n_1$ (the parent of $x$)
		and $n_2$ (the parent of $y$).
		Define $J_0 = 3T + 2 := \{n\in\mathbb{N}|n\equiv 2 \mod 3\}$.
		The anchor node will be $n_y$.
		For $j = 3t + 2 \in J_0$, define
		\begin{align*}
			\psi_j(n_x) = 3t
			\txt,
			\qquad&\psi_j(n_1) = 3t - 1\\
			\psi_j(n_y) = 3t + 2 = j
			\txt,
			\qquad&\psi_j(n_2) = 3t + 1\\
			\psi_j(n_l) = (3t)_L
			\txt.
			\qquad&
		\end{align*}
		The set $L=\{n_l\}$ is a minimal latent set and the
		ancestral structure $\ancestralStructure=\{n_x \rightsquigarrow n_2\}$
		is valid.
		Note that the outer nodes are not in $L$ and the one
		ancestral edge $n_x \rightsquigarrow n_2$ does \emph{not}
		start at a latent node, so this is indeed backdoor-free.
		Thus by Thm.\ \ref{thm:extract_from_backdoor_complete},
		$\mu(\graph,L)$ is identified.
		By Rmk.\ \ref{rmk:practical_computation_of_extraction}
		(or alternatively: by inspecting the proof of Thm.\ \ref{thm:extract_from_backdoor_complete}),
		$\mu(\graph,L)$ can be computed\Slash{}estimated in practice as follows.
		
		First, curate a dataset $\mathcal{D}^1=\{ \randomVar{X}^1_j,
		\randomVar{Y}^1_j\}_{j\in J_0}$,
		given as $\randomVar{X}^1_j=\psi_j(n_1)$ and
		$\randomVar{Y}^1_j=\psi_j(n_x)$;
		\ie data-points are of the form $(x,y)=(\randomVar{A}_{3t-1},
		\randomVar{A}_{3t})$.
		Estimate a probability kernel $\hat\mu^1_x$
		(a distribution parametrized by $x$)
		from $(x,y)\in\mathcal{D}^1$ with
		$y\sim\hat\mu^1_x$
		(\eg choose bins $X_1, \ldots, X_b$, then for each bin $X_k$ estimate a density
		$\hat\mu^1_k$
		from the values $y$ of pairs $(x,y)$ with $x\in X_k$).
		We write $\hat{X}^1_{y^{\txt{past}}} = \hat\mu^1_{y^{\txt{past}}}$.		
		Note that this is an estimator of the atom Def.\ \ref{def:atoms}
		at $n_x$, not of the mechanism $X$ (indeed it estimates
		the composition $X_{y^{\txt{past}},[l_1]}\circ L_1$,
		using notation \ref{notation:causal_wiring},
		for example it has a single argument, while $X$ has two).
		
		Second, curate a dataset $\mathcal{D}^2=\{ \randomVar{X}^2_j,
		\randomVar{Y}^2_j\}_{j\in J_0}$,
		given as the multi-variate $\randomVar{X}^2_j=(\psi_j(n_x),
		\psi_j(n_2))$ -- $\psi_j(n_x)$ must be included, because it is an
		observed ancestor within the c-component -- and
		$\randomVar{Y}^2_j=\psi_j(n_y)$;
		\ie data-points are of the form $(\vec{x},y)=((\randomVar{A}_{3t},\randomVar{A}_{3t+1}),
		\randomVar{A}_{3t+2})$.		
		Estimate a probability kernel $\hat\mu^2_{\vec{x}}$
		(a distribution parametrized by $\vec{x}$)
		from $(\vec{x},y)\in\mathcal{D}^2$ with
		$y\sim\hat\mu^2_{\vec{x}}$.
		We write $\hat{Y}^2_{x,z} = \hat\mu^2_{x,z}$.
		Note that this is an estimator of the atom Def.\ \ref{def:atoms}
		at $n_y$, not of the mechanism $Y$.
		
		Then (see Rmk.\ \ref{rmk:practical_computation_of_extraction})
		$\mu(\graph,L)_{y^{\txt{past}},z} = 
		\hat{X}^1_{y^{\txt{past}}}
		\otimes
		\hat{Y}^2_{x,z}$, \ie
		for measurable $B$
		\begin{align*}
			\mu(\graph,L)_{y^{\txt{past}},z}(B)
			\halfquad&=\halfquad
			\kernelCompound{
				\hat{X}^1_{y^{\txt{past}}}
				\otimes
				\hat{Y}^2_{x,z}
			}(B)
			\\
			\halfquad&=\halfquad
			\int 
			\hat{X}^1_{y^{\txt{past}}}(\dx)
			\int
			\hat{Y}^2_{x,z}(\dy)
			1_B(x,y)\txt.
		\end{align*}
		However, for the present example we do not actually have to compute this,
		knowing the "atoms"
		$\hat{X}^1_{y^{\txt{past}}}$
		and
		$\hat{Y}^2_{x,z}$ will be enough
		(see example \ref{example:frontdoor} however).

		Next, we focus on the dataset (3), 
		with indices in $i=(v,i)\in\{x',z'\}\times I'$,
		we assume $I' = \mathbb{N}_0$.
		We observe values
		$\randomVar{B}^X_0, \randomVar{B}^X_1, \ldots$
		and
		$\randomVar{B}^Z_0, \randomVar{B}^Z_1, \ldots$.
		
		Similarly as before,
		there is a backdoor-free family of embeddings,
		embedding the local graph $\graph'$ shown in the box
		"learn" (Fig.\ \ref{fig:intro}, dataset 3):
		Call the two inner nodes $n_l, n_z$ (aligned to $L_2, Z_2$)
		and the single outer node $n_1$ (the parent of $n_z$).
		Define $J'_0 = \{z'\}\times I'$.
		The anchor node will be $n_z$.
		For $j = (z',i) \in J'_0$, define
		\begin{align*}
			\psi_j(n_l) = (l',i)
			\txt,
			\qquad&\psi_j(n_1) = (x',i)\\
			\psi_j(n_z) = (z',i) = j
			\txt.
			\qquad&
		\end{align*}
		The set $L'=\{n_l\}$ is a minimal latent set and the
		ancestral structure $\ancestralStructure=\emptyset$
		is valid.
		Note that the outer node is not in $L$ and
		there are no ancestral edges, thus no backdoors,
		so this family is indeed backdoor-free.
		Thus by Thm.\ \ref{thm:extract_from_backdoor_complete},
		$\mu(\graph',L')$ is identified.
		By Rmk.\ \ref{rmk:practical_computation_of_extraction}
		(or alternatively: by inspecting the proof of Thm.\ \ref{thm:extract_from_backdoor_complete}),
		$\mu(\graph',L')$ can be computed\Slash{}estimated in practice as follows.
		
		Curate a dataset $\mathcal{D}^3=\{ \randomVar{X}^3_j,
		\randomVar{Y}^3_j\}_{j\in J'_0}$,
		given as $\randomVar{X}^3_j=\psi_j(n_1)$ and
		$\randomVar{Y}^3_j=\psi_j(n_z)$;
		\ie data-points are of the form $(x,z)=(\randomVar{B}^X_{i},
		\randomVar{B}^Z_{i})$.
		Estimate a probability kernel $\hat\mu^3_x$
		(a distribution parametrized by $x$)
		from $(x,z)\in\mathcal{D}^3$ with
		$y\sim\hat\mu^3_x$.
		We write $\hat{Z}^3_{x} = \hat\mu^3_{x}$.
		Then $\mu(\graph',L')_x=\hat{Z}^3_{x}$.
		Again, this is an estimator of the atom Def.\ \ref{def:atoms}
		at $n_z$, not of the mechanism $Z$.
		
		Finally, the query (see Fig.\ \ref{fig:intro})
		is introduced into the $I$-graph,
		it can be structured by Lemma \ref{lemma:assoc_structured_query},
		which uses the identity on this subset of $I$ to embed the
		graph shown in the panel "query" of Fig.\ \ref{fig:intro}.
		By Thm.\ \ref{thm:query_id}, the query can be identified
		via the structured kernel associated to this query.
		
		We glue this query (as structured kernel, initially with
		$L$ as indicated in Fig.\ \ref{fig:intro})
		by using Lemma \ref{lemma:glue:mt}.
		We start by gluing a single period of the time-series.
		The gluing target $\graph''$ is:\\[0.2em]
		\begin{minipage}{\textwidth}
			\colorlet{dblue}{blue!50!black}
			\colorlet{dgreen}{green!50!black}
			\colorlet{dyellow}{yellow!50!black}
			\colorlet{dpurple}{purple!50!black}
			\colorlet{dorange}{orange!50!black}
			\colorlet{c4}{dblue!20!white}
			\colorlet{c3}{dgreen!20!white}
			\colorlet{c1}{dyellow!20!white}
			\colorlet{c2}{dpurple!20!white}
			\colorlet{c5}{dorange!20!white}
			\colorlet{dgray}{gray!50!black}
			\colorlet{c0}{dgray!20!white}
			\colorlet{b0}{c0!40!white}
			\colorlet{b5}{c5!40!white}
			\colorlet{b3}{c3!40!white}
			\colorlet{b4}{c4!40!white}
			\centering
			\begin{tikzpicture}
				\draw (-1,0) node(past) [circle,draw, very thick, inner sep=0.2em] {};
				
				\draw (0,0) node (X0)[rectangle, fill=c3, rounded corners, inner sep=0.4em,
				draw=dgreen,	thick]
				{$X$};
				\draw (1,0) node (Z0)[rectangle, fill=c3, rounded corners, inner sep=0.4em,
				draw=dgreen, thick]
				{$Z_2$};
				\draw (2,0) node (Y0)[rectangle, fill=c3, rounded corners, inner sep=0.4em,
				draw=dgreen,	thick]
				{$Y$};
				\draw (1,-1) node (L0)[rectangle, fill=c0, rounded corners, inner sep=0.4em,
				draw=dgray, dashed, thick]
				{$L_2$};
				\draw (0,1) node (L10)[rectangle, fill=c0, rounded corners, inner sep=0.4em,
				draw=dgray, dashed, thick]
				{$L_1$};

				\draw [->, very thick] (past) -- (X0);
				\draw [->, very thick] (X0) -- (Z0);	
				\draw [->, very thick] (Z0) -- (Y0);	
				\draw [->, very thick] (L0) -- (Z0);	
				\draw [->, very thick] (L10) -- (X0);	
				\draw [->, very thick] (L10) edge[bend left] (Y0);
			\end{tikzpicture}
		\end{minipage}\\[0.2em]
		This graph has two c-components each of which is contained (actually
		identical to) one of the two learned graphs above --
		and these are also subgraphs in the sense of Def.\ \ref{def:subgraphs_simp}.
		By Rmk.\ \ref{rmk:compute_gluing},
		we can compute $\mu(\graph'',L'')$ (where $L''=\{l_1,l_2\}$)
		by "interlacing" the atoms of these c-components:
		\begin{equation*}
			\mu(\graph'',L'')
			\halfquad=\halfquad
			\hat{X}^1_{y^{\txt{past}}}
			\otimes
			\hat{Z}^3_{x}
			\otimes
			\hat{Y}^2_{x,z}
			\txt.
		\end{equation*}
		Similarly, the known initial state
		of the query $Y_{\txt{init}}$ and further copies of
		$\mu(\graph'',L'')$ can be glued to obtain the query
		$\mu^{\txt{query}}$
		as shown in Fig.\ \ref{fig:intro}.
		\begin{equation*}
			\hat\mu^{\txt{query}}
			\halfquad=\halfquad
			Y_{\txt{init}}
			\otimes
			\mu(\graph'',L'')
			\otimes
			\mu(\graph'',L'')
			\otimes
			\mu(\graph'',L'')
			\otimes
			\mu(\graph'',L'')
			\txt.
		\end{equation*}
		Technically, the query asked for the distribution
		of the final $Y$ only, so we may marginalize
		$\marginalizeOp_{L^q}\hat\mu^{\txt{query}}$
		where $L^q$ contains all nodes of the query other than the
		finial $y$.		
		
		In the end we combine three estimates on three datasets,
		further
		every $\otimes$-product and the final marginalization
		come with an integral (multiple integrals for the final marginalization),
		so $\hat\mu^{\txt{query}}$ is actually a rather complex estimator.
		However,
		through the separation of extraction
		and computation
		we easily draft a very intuitive plan
		(boxes in Fig.\ \ref{fig:intro}).
		Then, by the graphical operations for
		the computations involved, 
		we immediately know how to write down this sound estimator.
		
		\textbf{(c) Observational Support:}		
		
		Recall that we assumed away (Ass.\ \ref{ass:support},
		cf.\ Rmk.\ \ref{rmk:kernels_uniqueness})
		observational support issues.
		It may of course happen, particularly for modifications of
		time-series as indicated above,
		that a particular choice of initial state
		or modified Markov-kernel lead us outside of our observational support.
		In practice, the detection of such problems may require
		explicit annotation on estimators (if they extrapolate by default,
		\eg linear regression) or may occur from explicitly failing prediction
		(if extrapolation is undefined or bogus, \eg nearest neighbor estimators).
		Helpful workarounds might involve comparing
		different accessible (in the sense of staying in the observational
		support) partial solutions:
		If we want to understand long term behavior, we may estimate
		the state after some fixed duration
		(\eg four periods in the example (a,b) above),
		if at four periods (12 time-steps) we have left observational
		support, we may instead try to understand what
		happens after one, two or three periods
		to try and extrapolate the general behavior from there.
		It may also be possible to use such partial results
		to guide the choice of parameters for experiments
		to do \emph{few} (rather than no) experiments.
		
		There is of course no general strategy
		to deal with lack of observational support.
		However flexible queries can help find the most
		helpful available results.
	\end{example}

	\subsection{Computation: Frontdoor Formula}
	
	The previous example \ref{example:real_world_questions}
	was focused on illustrating the beyond-interventional
	query capability on a concrete example and with
	a sound identification strategy.
	Here, we further illustrate the use and computational
	realization of graphical operations §\ref{sec:graphical_operations}
	on the standard example of the frontdoor-formula (see \eg
	\citep[§3.3.2 (p.\,81ff)]{PearlBook});
	a similar computation for the napkin-graph
	is given in Rmk.\ \ref{rmk:compute_atoms_practice}.
	The reader should notice how, on a level of graphical
	operations, the approach is very intuitive.
	The complexity of the final estimator only arises once
	the computational details of the graphical operations
	are inlined into a single formula.
	Slight modifications can be made for the multi-level case,
	as motivated by example \ref{example:learn_pinned}.
	For computations with kernels in general, see §\ref{apdx:kernels}--\ref{apdx:sparsity}.

	\begin{figure}[tb]
		\begin{minipage}{\textwidth}
			\centering
			\colorlet{dblue}{blue!50!black}
\colorlet{dgreen}{green!50!black}
\colorlet{dyellow}{yellow!50!black}
\colorlet{dpurple}{purple!50!black}
\colorlet{dorange}{orange!50!black}
\colorlet{c4}{dblue!20!white}
\colorlet{c3}{dgreen!20!white}
\colorlet{c1}{dyellow!20!white}
\colorlet{c2}{dpurple!20!white}
\colorlet{c5}{dorange!20!white}
\colorlet{dgray}{gray!50!black}
\colorlet{c0}{dgray!20!white}

\pgfdeclarelayer{bg}
\pgfsetlayers{bg,main}

\begin{tikzpicture}[xscale=1.25]
	\begin{scope}
		\draw (-0.15,1) node {(i)};
		
		\draw (0,0) node (X)[rectangle, fill=c3, rounded corners, inner sep=0.4em,
		draw=dgreen,	thick]
		{$X$};
		\draw (1,0) node (Z)[rectangle, fill=c3, rounded corners, inner sep=0.4em,
		draw=dgreen, thick]
		{$Z$};
		\draw (2,0) node (Y)[rectangle, fill=c3, rounded corners, inner sep=0.4em,
		draw=dgreen,	thick]
		{$Y$};
		
		\draw [->, very thick] (X) -- (Z);	
		\draw [->, very thick] (Z) -- (Y);	
		
		\draw [<->, very thick,dgray] (X) edge[bend left=2.5cm]
		node[pos=0.5,above,yshift=-0.4em]
		[rectangle, fill=c0, rounded corners, inner sep=0.4em,
		draw=dgray, thick, dashed]{$L$} (Y);
	\end{scope}
	\begin{scope}[xshift=3cm]
		\draw (-0.15,1) node {(ii)};
		
		\draw (0,0) node (X)[rectangle, fill=c3, rounded corners, inner sep=0.4em,
		draw=dgreen,	thick]
		{$X$};
		\draw (1,0) node (Z)
		[circle, draw, very thick, inner sep=0.2em, outer sep=0.15em]{};
		\draw (2,0) node (Y)[rectangle, fill=c3, rounded corners, inner sep=0.4em,
		draw=dgreen,	thick]
		{$Y$};
		
		\draw [->, very thick] (Z) -- (Y);	
		
		\draw [<->, very thick,dgray] (X) edge[bend left=2.5cm]
		node[pos=0.5,above,yshift=-0.4em]
		[rectangle, fill=c0, rounded corners, inner sep=0.4em,
		draw=dgray, thick, dashed]{$L$} (Y);
	\end{scope}
	\begin{scope}[xshift=6cm]
		\draw (-0.15,1) node {(iii)};
		
		\draw (0,0) node (X)[rectangle, fill=c0, rounded corners, inner sep=0.4em,
		draw=dgray,	thick, dashed]
		{$X$};
		\draw (1,0) node (Z)
		[circle, draw, very thick, inner sep=0.2em, outer sep=0.15em]{};
		\draw (2,0) node (Y)[rectangle, fill=c3, rounded corners, inner sep=0.4em,
		draw=dgreen,	thick]
		{$Y$};
		
		\draw [->, very thick] (Z) -- (Y);	
		
		\draw [<->, very thick,dgray] (X) edge[bend left=2.5cm]
		node[pos=0.5,above,yshift=-0.4em]
		[rectangle, fill=c0, rounded corners, inner sep=0.4em,
		draw=dgray, thick, dashed]{$L$} (Y);
	\end{scope}
	\begin{scope}[xshift=9cm]
		\draw (-0.15,1) node {(iv)};
		
		\draw (0,0) node (Z)
		[circle, draw, very thick, inner sep=0.2em, outer sep=0.15em]{};
		\draw (1,0) node (Y)[rectangle, fill=c3, rounded corners, inner sep=0.4em,
		draw=dgreen,	thick]
		{$Y$};
		
		\draw (1,1) node (L)[rectangle, fill=c0, rounded corners, inner sep=0.4em,
		draw=dgray,	thick, dashed]
		{$L$};
		
		\draw [->, very thick] (Z) -- (Y);	
		\draw [->, very thick] (L) -- (Y);	
	\end{scope}
	\begin{scope}[yshift=-2cm]
		\draw (-0.15,1) node {(v)};
		
		\draw (0,0) node (X)
		[circle, draw, very thick, inner sep=0.2em, outer sep=0.15em]{};
		\draw (1,0) node (Z)[rectangle, fill=c3, rounded corners, inner sep=0.4em,
		draw=dgreen,	thick]
		{$Z$};
		
		\draw [->, very thick] (X) -- (Z);	
	\end{scope}
	\begin{scope}[yshift=-2cm, xshift=2.5cm]
		\draw (-0.15,1) node {(vi)};
		
		\draw (0,0) node (X)[rectangle, fill=c2, inner sep=0.4em,
		draw=dpurple,	thick]
		{$X_{\txt{do}}$};
	\end{scope}
	\begin{scope}[yshift=-2cm, xshift=4cm]
		\draw (-0.15,1) node {(vii)};
		
		\draw (0,0) node (X)[rectangle, fill=c2, inner sep=0.4em,
		draw=dpurple,	thick]
		{$X_{\txt{do}}$};
		\draw (1,0) node (Z)[rectangle, fill=c3, rounded corners, inner sep=0.4em,
		draw=dgreen,	thick]
		{$Z$};
		\draw (2,0) node (Y)[rectangle, fill=c3, rounded corners, inner sep=0.4em,
		draw=dgreen,	thick]
		{$Y$};
		
		\draw (2,1) node (L)[rectangle, fill=c0, rounded corners, inner sep=0.4em,
		draw=dgray,	thick, dashed]
		{$L$};
		
		\draw [->, very thick] (X) -- (Z);	
		\draw [->, very thick] (Z) -- (Y);	
		\draw [->, very thick] (L) -- (Y);	
	\end{scope}
	\begin{scope}[yshift=-2cm, xshift=7.5cm]
		\draw (-0.15,1) node {(viii)};
		
		\draw (0,0) node (X)[rectangle, fill=c2, inner sep=0.4em,
		draw=dpurple,	thick]
		{$X_{\txt{do}}$};		
		\draw (1,0) node (Z)[rectangle, fill=c0, rounded corners, inner sep=0.4em,
		draw=dgray,	thick, dashed]
		{$Z$};
		\draw (2,0) node (Y)[rectangle, fill=c3, rounded corners, inner sep=0.4em,
		draw=dgreen,	thick]
		{$Y$};
		
		\draw (2,1) node (L)[rectangle, fill=c0, rounded corners, inner sep=0.4em,
		draw=dgray,	thick, dashed]
		{$L$};
		
		\draw [->, very thick] (X) -- (Z);	
		\draw [->, very thick] (Z) -- (Y);	
		\draw [->, very thick] (L) -- (Y);	
	\end{scope}
\end{tikzpicture}
		\end{minipage}
		\caption{Steps taken in example \ref{example:frontdoor}\,a.}
		\label{fig:frontdoor}
	\end{figure}
	
	\begin{example}[Frontdoor-Formula]
		\label{example:frontdoor}
		As an illustration of the use of graphical computations
		we (a) derive the frontdoor formula, then (b) put it in
		relation to the "pseudo-frontdoor" argument
		of example \ref{example:extraction_frontdoor}.
		
		\textbf{(a) The frontdoor formula:}
		
		In the spirit of example \ref{example:iid_id_from_structured_only},
		we can skip the extraction procedure for IID data (see (b) below however).
		We start from the structured kernel
		Fig.\ \ref{fig:frontdoor}\,(i), with
		$\mu(\graphStructural^{\txt{(i)}},L)=P^{\txt{obs}}$;
		we write $L$ both for the latent $L$ and the set containing $L$.
		Note that $\graphStructural^{\txt{(i)}}$ has
		two structural c-components, shown in Fig.\ \ref{fig:frontdoor}\,(ii,\,v).
		
		Thus by Lemma \ref{lemma:c_components:mt} we can compute
		$\mu(\graphStructural^{\txt{(ii)}},L)$
		and $\mu(\graphStructural^{\txt{(v)}},\emptyset)$
		from $P^{\txt{obs}}$.
		To do so in practice, we consult
		Rmk.\ \ref{rmk:compute_atoms_practice}
		(see also Lemmas \ref{lemma:atom_properties},
		\ref{lemma:properties_c_components}).
		First we can write
		$\mu(\graphStructural^{\txt{(i)}},L) = A^x \otimes A^z_x \otimes A^y_{x,z}$,
		where the argument $x$ of $A^z_x$ is "wired" to $A^x$
		etc.\ (see §\ref{sec:causality_and_sparsity}).
		Generally the atoms $A^v$ can be computed by marginalization and disintegration
		(disintegrations from the left are a basic measure-theoretic
		operation on kernels,
		either our practical encoding of kernels directly supports their computation
		or we can estimate them \eg as a conditional distribution relative
		to the result for $\mu$, here $P^{\txt{obs}}$).
		In the present case,
		by uniqueness of disintegrations
		(see Lemma \ref{lemma:order_independent_ops:mt}),
		we can read off
		\begin{align*}
			A^x&=P^{\txt{obs}}(X)\\
			A^z_x&=P^{\txt{obs}}(Z|X=x)\\
			A^y_{x,z}&=P^{\txt{obs}}(Y|X=x,Z=z)
			\txt.
		\end{align*}
		By Lemmas \ref{lemma:atom_properties},
		\ref{lemma:properties_c_components}, thus
		\begin{align*}
			\mu(\graphStructural^{\txt{(ii)}},L)_z
			&=
			A^x\otimes A^y_{x,z}\\			
			\mu(\graphStructural^{\txt{(v)}},L)_x
			&=
			A^z_{x}
			\txt.
		\end{align*}
		Next, we simplify $\mu(\graphStructural^{\txt{(ii)}},L)$
		(by Lemma \ref{def:simplify_simp}),
		by marginalizing $X$ Fig.\ \ref{fig:frontdoor}\,(iii),
		then dropping $X$ Fig.\ \ref{fig:frontdoor}\,(iv).		
		\begin{align*}
			\mu(\graphStructural^{\txt{(iii)}},L)_z
			&=
			\marginalizeOp_{\{X\}}
			(\mu(\graphStructural^{\txt{(ii)}},L)_z)\\
			\mu(\graphStructural^{\txt{(iv)}},L)_z
			&=
			\mu(\graphStructural^{\txt{(iii)}},L)_z
			\txt.
		\end{align*}
		Formally, marginalization of $X$ is
		the same as plugging the full space $\topSpace{X}$
		into the measure (see Lemma \ref{lemma:order_independent_ops:mt}),
		\ie for measurable $B\in \sigmaAlgebraBorel{\topSpace{Y}}$
		(using the definition of the $\otimes$-product, Def.\ \ref{def:kernel_product}):
		\begin{align*}
			\mu(\graphStructural^{\txt{(iii)}},L)_z(B)
			&=
			\int A^x(\dx) \int A^y_{x,z} (\dy) 1_{\topSpace{X}\times B}(x,y)
			\\
			&=
			\int P^{\txt{obs}}(X\in\dx)
			\int P^{\txt{obs}}(Y\in\dy|X=x,Z=z) 1_{B}(y)
			\txt.
		\end{align*}
		Or using density notation (assuming densities exist),
		writing $\rho(\ldots)$ for the density of $\mu(\ldots)$,
		$p^{\txt{obs}}$ for the density of $P^{\txt{obs}}$:
		\begin{equation*}			
			\rho(\graphStructural^{\txt{(iv)}},L)_z(y)
			=	
			\rho(\graphStructural^{\txt{(iii)}},L)_z(y)
			=
			\int p^{\txt{obs}}(x)
			p^{\txt{obs}}(y|x,z) \dx
			\txt,
		\end{equation*}
		where the second integral has vanished, because it only
		integrated out the indicator function of $B$
		(the left-hand-side is now a density evaluated at $y$,
		rather than a measure evaluated on $B$).
		Similarly we can use parametrized densities and the results above to write
		\begin{equation*}			
			\rho(\graphStructural^{\txt{(v)}},\emptyset)_x(z)
			=	
			p^{\txt{obs}}(z|x)
			\txt.
		\end{equation*}
		Next the intervention Fig.\ \ref{fig:frontdoor}\,(vi)
		is known a priori by definition,
		and is given by $\mu(\graphStructural^{\txt{(vi)}},\emptyset) = \delta[x_0]$
		for a intervention target $x_0$.
		By applying gluing Lemma \ref{lemma:glue:mt} twice
		we can obtain $\mu(\graphStructural^{\txt{(vii)}},L)$.
		In practice, we follow Rmk.\ \ref{rmk:compute_gluing},
		which asks us to write the atoms of the glued parts
		in causal order as a $\otimes$-product.
		Note that all parts have a single observed node,
		thus a single atom given by the structured kernel itself
		(Lemma \ref{lemma:atom_properties}).
		In general we may have to compute them as marginalization
		and disintegrations of our intermediate results.
		Thus we find
		\begin{equation*}
			\mu(\graphStructural^{\txt{(vii)}},L)
			=
			\mu(\graphStructural^{\txt{(vi)}},\emptyset)^{x'}
			\otimes
			\mu(\graphStructural^{\txt{(v)}},\emptyset)^z_{x'}
			\otimes
			\mu(\graphStructural^{\txt{(iv)}},\emptyset)_z
			\txt,
		\end{equation*}
		where upper indices clarify what is computed for "wiring"
		purposes, \ie $\mu(\graphStructural^{\txt{(vi)}},\emptyset)^{x'}$
		outputs "$x'$" which is used as input for
		$\mu(\graphStructural^{\txt{(v)}},\emptyset)^z_{x'}$ and so on.
		Finally, we are interested only in the distribution of
		$Y$, thus we marginalize $Z$ to get
		$\mu(\graphStructural^{\txt{(viii)}},L)$.
		As a density:
		\begin{align*}
			\rho(\graphStructural^{\txt{(viii)}},L)(y)
			&=
			\int
			\delta[x_0](\dx')
			\int
			p^{\txt{obs}}(z|x')
			\left(
			\int p^{\txt{obs}}(x)
			p^{\txt{obs}}(y|x,z) \dx
			\right) \dz\\
			&=
			\int
			p^{\txt{obs}}(z|x_0)
			\left(
			\int p^{\txt{obs}}(x)
			p^{\txt{obs}}(y|x,z) \dx
			\right) \dz
			\txt.
		\end{align*}
		
		\textbf{(b) The extraction frontdoor:}
		A similar two-step estimation was
		used in example \ref{example:extraction_frontdoor}
		to combine spatially and temporally repeating
		mechanisms.
		Both arguments can be combined. However
		the two pieces Fig.\ \ref{fig:frontdoor}\,(ii,\,v)
		are actually the smallest parts that possess backdoor-free
		embeddings (Def.\ \ref{def:backdoor_free}).
		For (v) this is no restriction, it has symmetry $H_Z$
		of its only inner node.
		For (ii) however, the symmetry is $H_X \cap H_L \cap H_Y$,
		\ie the \emph{shared symmetry} of all inner nodes.
		In the setup of example \ref{example:extraction_frontdoor},
		thus all three of $X$, $L$ and $Y$ have to simultaneously be
		either all spatially invariant or all temporally invariant.
		For example we could have $Z$ spatially invariant (but changing over time)
		and $X$, $L$ and $Y$ all temporally invariant (but changing over space),
		but not $X$ changing over time and $L$ changing over space.
		Once $\mu(\graphStructural^{\txt{(v)}},\emptyset)$
		and $\mu(\graphStructural^{\txt{(ii)}},L)$ are extracted
		via Thm.\ \ref{thm:extract_from_backdoor_complete},
		we may proceed as in part (a) above.
	\end{example}

	\subsection{Complex Queries and Simulation}
	
	In §\ref{sec:queries} we briefly introduce the viewpoint
	of queries as simulations.
	The general idea is as follows.
	The query shares its building blocks (in our case structured kernels
	modeling \eg physical processes)
	with the observations or they are know a priori (interventions).
	Thus, given perfect knowledge of the model,
	we could imagine the query as the result of running a
	simulation of a imaginary setup pieced together
	from building-blocks.
	In the interventional setup the arrangement (graph)
	of these building blocks must be identical to their arrangement
	in the observed setup, only the replacement of individual mechanisms
	with interventions (known a priori building blocks) is allowed.
	In our case, we only insist on knowing the building-blocks
	but allow any arrangement for the query.
	This is a more general notion of "simulation"
	as further illustrated by the following example.
	This perspective on query-construction extends
	canonically to multi-level setups,
	see example \ref{example:multi_level_queries}.

	\begin{example}[Simulation]
		\label{example:simulation}
		We equip a sheet made of material A with temperature-sensors
		in a regular 2 dimensional grid and take time-series data.
		We model the system as time- and space-translation invariant,
		where each temperature-measurement $T_{x,y;t}$ is
		given by its own and its neighbors' previous, at $t-1$, values and
		a local noise-term.
		The boundary of the sheet is modeled by separate nodes (we need not
		be concerned about them for the purpose of this example;
		\eg one edge of the sheet may be actively heated and the opposite
		one connected to a heat-sink).
		The $I$-graph is a "tower" of two-dimensional grids,
		where all non-boundary nodes are aligned to the same mechanism $\mu^A$,
		with parents in the time-slice below.
		Since everything is observed, we can easily learn $\mu^A$.
		
		Let's say we do the same experiment (in the intuitive sense, not
		specifically in the sense this word is often used in causal literature)
		with a sheet made of material B to learn $\mu^B$.
		
		We want to know: If we produce a tube of finite length made
		of rings of materials A and B and heat one end (\eg
		intervene to a fixed temperature), how is heat
		conducted to the other end.
		Obviously simulating heat-diffusion is a comparatively simple
		task, when thought about as simulation.
		However, we can now also model the query by exactly the same
		ideas we used to model the observations; for modeling purposes,
		it is insubstantial
		that we are currently formulating a question, not a description
		of observations.
	\end{example}
	
	For this example, the causal reasoning does of course not contribute much.
	It is an analogue of an SCM model for $X \rightarrow Y$ telling us
	that for $E[Y|\PearlDo(X=x)]$ we have to plug $x$ into a regressor.
	The point is that (a) "plug in x" (SCM) has become "run this simulation" (ours),
	and (b) given more variables at each grid point with more complex causal
	relations and possibly latents, the identification of the simulation
	(on \emph{what} regressors to run on and how) is no longer trivial
	(in analogy to how, on graphs larger than $X\rightarrow Y$,
	causal inference informs
	us which regressors to learn and how to combine them).

	\section{Details on Structured Kernels}
	\label{apdx:structured_kernels}
	
	Basic notions, especially concerning kernels,
	are recalled and explained with examples in SCM language
	in §\ref{apdx:ptheo_basics}.
	Here we continue with non-standard constructions on kernels,
	designed for our use-case of encoding (and computing with)
	the sparsity-structure of a (causal) structural graph while
	respecting latent structure.

	\subsection{Structural Graphs}
	\label{apdx:structural_graphs}
	
	Large $\otimes$-products of kernels with sparse inter-dependencies,
	are essentially a type of tensor network. Tensor networks are often
	represented graphically. Also the representation of causal
	relationships by graphical models is very common.
	It therefore is not particularly surprising, that
	a suitable form of graphical representation will also help
	substantially in organizing the knowledge we encounter in our
	reasoning here. We start by defining a precise notion of graphical
	representation that can encapsulate the causal sparsity information
	of §\ref{apdx:sparsity} in a more accessible form.
	
	First, we give a slightly more detailed version of
	Def.\ \ref{def:structural_graph_simpl} split into
	a graph-definition \ref{def:structural_graph}
	and model-alignment \ref{def:model_alignment_and_kernel}.
	
	\begin{definition}[Structural Graph And Causal Orders]\label{def:structural_graph}
		A structural graph $\graphStructural$ is a finite set of nodes
		$\nodes := \nodesInner\,\dot\cup\,\nodesOuter$
		with $\nodesInner \cap \nodesOuter=\emptyset$,
		together with a set of directed edges
		$\edgesStructural\subset\nodes\times\nodesInner\setminus\Delta$,
		where $\Delta=\{(n,n)|n\in\nodesInner\}$ is the diagonal,
		(proper edges)
		between nodes in
		$\nodesInner$ and from nodes in
		$\nodesOuter$ to nodes in $\nodesInner$.
		
		We call a structural graph $\graphStructural$ acyclic, if there are no directed cycles.
		If $\graphStructural$ is acyclic, then there is a partial order on (all) nodes $\nodes$
		defined by $n_1 \leq n_2 :\Leftrightarrow $ there is a directed path
		from $n_1$ to $n_2$.
		A $\graphStructural$-causal order on $\nodes$ is a
		total order extending this partial order.
		Such a total order always exists, but it is in general not unique.
	\end{definition}
	
	\begin{rmk}
		Edges are relevant (only) for parent relation-ships
		(and properties deduced thereof, like ancestral relation-ships),
		it thus is irrelevant if more than one edge is allowed between nodes
		for each orientation and that there are no edges from a node to itself.
		We did adopt the convention of allowing at most one edge
		(per orientation) between any two nodes (the definition of edges
		$\edgesStructural\subset\nodes\times\nodesInner\setminus\Delta$
		enforces this),
		which is in line with standard causal analysis \citep{PearlBook,Elements}.
		
		With the same mechanism $f_Y$ potentially featuring in different
		contexts,
		it may occur in practice, with the same variable acting as
		multiple parents, \eg $f_Y(x_1, x_2)$ could be evaluated as
		$x\mapsto f_Y(x,x)$. This can be modeled easily, by defining
		a diagonal map $\Delta_{\topSpace{X}}: \topSpace{X}\rightarrow\topSpace{X}\times\topSpace{X}, x\mapsto (x,x)$
		with $\Delta_{\topSpace{X}}\in\knownFunction$ and
		considering $f_y\circ \Delta_{\topSpace{X}}$ as a soft-intervention
		(which is allowed both in observations an queries, which we do not distinguish
		formally);
		technically the reader may have noticed, that $\Delta_{\topSpace{X}}$
		is not actually a probability-kernel, however we can replace
		it by one which has singular measure $1$ at the value taken by
		$\Delta_{\topSpace{X}}$.
		
		Note that this interpretation does not seem to affect identifiability:
		If $f_Y$ is identifiable, then the composition with
		an intervention $\Delta_{\topSpace{X}}\in\knownFunction$
		will be considered identifiable (if we know $f_Y(x_1,x_2)$ we can
		in particular plug in the same argument twice), but knowing
		(from observations) only $f_Y \circ \Delta_{\topSpace{X}}$
		will only retain this composition (\ie, observing only
		$f_Y(x,x)$ we cannot make predictions for queries
		that involve the full two-dimensional $f_Y(x_1,x_2)$, at least not
		without substantial knowledge about internal structure, cf.\ %
		§\ref{sec:knowledge_closures}).
	\end{rmk}
	
	\begin{definition}[Model Alignment]\label{def:model_alignment_and_kernel}
		A model aligned structural graph is an acyclic structural graph $\graphStructural$ together
		with
		\begin{enumerate}[label=(\roman*)]
			\item 
			a probability kernel $\mu^n$ for each inner node $n\in\nodesInner$.
			and
			a bijective assignment (cf.\ \ref{def:contraction}) of the $\kappa$ arguments of $\mu^n$
			to proper parents of $n$:
			\begin{equation*}
				\pa^n: \{1, 2, \ldots, \kappa\}
				\xrightarrow{1:1} \Pa_{\graphStructural}(n)
				\txt.
			\end{equation*}
			\item
			$\mu^n$ is a kernel on $\kappa^n$ arguments,
			from $\topSpace{X}^n_1, \ldots, \topSpace{X}^n_{\kappa^n}$
			to $\topSpace{Y}^n$, we require $\forall n$:
			$\topSpace{Y}^{\pa^n(1)}\subset\topSpace{X}^n_1$,
			\dots, 
			$\topSpace{Y}^{\pa^n(\kappa^n)}\subset\topSpace{X}^n_{\kappa^n}$.
		\end{enumerate}
	\end{definition}
	\begin{notation}[Structural Graph Contraction]	
		\label{notation:graphical_sparse_contraction}	
		Given a model aligned acyclic
		structural graph $\graphStructural$ together with a $\graphStructural$-causal order $\pi$
		and a subset
		$\nodes' \subset \nodesInner$
		(which is totally ordered by $\leq_\pi$),
		we write
		\begin{align*}
			&\otimes^\pi_{n\in\nodes'}\mu^n
			:=			
			\otimes_{n\in\nodes'}^{\pa^*} \mu^n
			\txt{ (cf.\ Def.\ \ref{def:contraction})}\\
			&\txt{which is a kernel from }
			\prod_{x\in\Pa_{\graphStructural}(\nodes')\setminus\nodes'}
			\Big(
				\cap_{\{(n,k)|\pa^{n}(k)=x\}}
				\topSpace{X}^n_k
			\Big)
			\txt{ to }
			\prod_{n\in\nodes'}
			\topSpace{Y}^n\txt,
		\end{align*}
		for the product of the $\mu^n$ in $\pi$-order, contracted
		by the mapping $\pa^n$ of Def.\ \ref{def:model_alignment_and_kernel},
		\ie using notation \ref{notation:causal_wiring}
		if nodes in the graph are named
		$\randomVar{X}, \randomVar{Y}, \randomVar{Z}, \ldots$
		with associated kernels $X$, $Y$, $Z$, then for example 
		$\randomVar{Y}$ having parents $\randomVar{X},\randomVar{Z}$
		will have arguments\Slash{}indices $x,z$, $Y_{x,z}$
		(the order is relevant, by convention we write indices in $\pi$-order,
		that is the kernels formally depends on $\pi$, albeit only by
		reordering its arguments).
		More formally
		this means we set $w^n:=\pa^n$ in Def.\ \ref{def:contraction}.	
	\end{notation}
	
	Next we use these structural graphs to attach sparsity-information
	to kernels, detailing Def.\ \ref{def:structured_kernel_simp} of the main text;
	\emph{the main text additionally uses notation
		\ref{notation:uniqueness_transposition}}:
	
	\begin{definition}[Structured Kernel]\label{def:structured_kernel}
		A structural model $(\graphStructural,L)$ is a tuple of a
		model aligned
		structural graph $\graphStructural$ together with a subset
		$L\subset\nodesInner$.
		Given a structural model $(\graphStructural,L)$
		and a $\graphStructural$-causal order $\pi$,
		we define (using notation \ref{notation:graphical_sparse_contraction};
		"in causal order compute the value at node $n$,
		by plugging into the $k$-th argument of $\mu^n$,
		the value of the node $\pa^n(k)$, where outer nodes $x$ are left open
		as arguments of the resulting $\mu(\graphStructural,L,\pi)_{x}$")
		\begin{equation*}
			\mu(\graphStructural,L,\pi)_{x}
			\halfquad:=\halfquad
			\marginalize{
				\otimes_{n\in\nodesInner}^{\pi}
				\mu^n
			}{ L }_{x}
			\txt,
		\end{equation*}
		where $x=(x_i)_{i\in\nodesOuter}$
		summarizes arguments in $\nodesOuter$.
		
		Given a set $\structureKernels$ of probability kernels
		(later typically $\structureKernels=\modelKernels\cup\knownKernels\cup\{\delta_i|i\in I\}$),
		a $\structureKernels$-structured kernel is a kernel $\mu_x$ that can be
		written in the form $\mu_x = \mu(\graphStructural,L,\pi)_x$,
		where all nodes $n\in\nodesInner$
		are aligned to an element $\mu^n\in\mathcal{F}$.
	\end{definition}
	
	This definition depends on the $\graphStructural$-causal order $\pi$,
	but only up to transposition (Def.\ \ref{def:transposition}),
	see \ref{notation:uniqueness_transposition}.
	It is helpful to briefly recall (slightly more formal)
	the precise definition of subgraphs introduced in the main text,
	as it is relevant to multiple graphical operations:
	\begin{CopyDef}{def:subgraphs_simp}{Structural Subgraphs}
		A structural subgraph $\graphStructural^A\leq\graphStructural^B$
		of a structural graph $\graphStructural^B$
		is a structural graph $\graphStructural^A$ such that:
		\begin{enumerate}[label=(\roman*)]
			\item\label{def:subgraphs:nodesets:apdx}
				\emph{Node Sets:}
				$\nodes^A\subset\nodes^B$
				and	$\nodesInner^A\subset\nodesInner^B$.
			\item\label{def:subgraphs:edges:apdx}
				\emph{Edge Sets:}
				$\edgesStructural^A=\edgesStructural^B \cap (\nodes^A\times\nodesInner^A)$,
				\ie edges are exactly those in $\graphStructural^B$
				from nodes in $\graphStructural^A$
				to inner nodes of $\graphStructural^A$.
			\item\label{def:subgraphs:inner_parents:apdx}
				\emph{Inner Parents:}
				$\forall n\in\nodesInner^A$:
				if $p\in\Pa_{\graphStructural^B}(n)$,
				then $p\in\nodes^A$.
			\item\label{def:subgraphs:alignment:apdx}
				\emph{Alignment:} For model-aligned graphs
				$\forall n\in\nodesInner^A$:
				$(\mu^A)^n = (\mu^B)^n$.
		\end{enumerate}
		We write		
		$\graphStructural^A < \graphStructural^B$,
		if $\graphStructural^A \leq \graphStructural^B$
		and
		$\nodesInner^A\subsetneq\nodesInner^B$.
	\end{CopyDef}

	\subsection{Regularity}
	\label{apdx:regularity}
	
	This subsection is conceptually important to understanding the limits
	of computation and identifiability within the formalism proposed
	in this paper and beyond it.
	However, it is somewhat technical and, on a first reading,
	may be skipped, as it may disrupt the logical flow of the presentation.	

	Statistical and causal inference fundamentally ask:
	What questions about a system are "identifiable", that is well-defined,
	\ie have answers that are uniquely determined by observations
	(at least in principle and asymptotically),
	and how can they be "estimated", that is how can answers be
	approximated (from finite data).
	In this section we clarify what kind of computations are
	possible with kernels within this philosophy.
	
	However, first a few remarks about uniqueness of kernels
	and the relevance of distributional support to transfer
	problems are required.
	\begin{rmk}[Support]
		\label{rmk:kernels_uniqueness}
		Disintegrations are typically only
		"$\vartheta\otimes\mu$-almost everywhere"
		unique (see Lemma \ref{lemma:disintegration}).
		In statistics, "almost everywhere" is often read essentially
		as "good enough".
		From the perspective of learning \emph{transferable}
		properties of a model, this hides a very fundamental
		aspect of empirical science:
		We can only ever learn things that we actually see
		(encounter with non-zero-probability).
		If we observe $P(\randomVar{X},\randomVar{Y})$ 
		of a causal (say IID) model $\randomVar{X} \rightarrow \randomVar{Y}$
		with marginal $\mu = P(\randomVar{X})$, then whatever we try
		to learn about (a regular version of) the conditional distribution
		$\nu_x := P(\randomVar{Y}|\randomVar{X}=x)$, we can
		learn at best in a $\mu$-almost everywhere sense.
		This is not a problem of our formalism, it is a deeply rooted
		issue of modeling
		transfer of knowledge (for example to an intervened system).
		
		Before we go into any detail, it seems important
		to point out, that this is (usually) not actually a problem in practice.
		When \citet{PearlBook} (and others) talk about
		do-interventions to a value $x_0$ then
		either $\randomVar{X}$ is categorical (and $P(\randomVar{X}=x_0)\neq 0$,
		in which case $\delta[x_0]\otimes \nu_x$ is well-defined, even if
		$\nu_x$ is defined only $\mu$-almost everywhere)
		or the actual real-world problem
		formally expressed as "$\PearlDo(\randomVar{X}=x)$"
		is understood as
		"once $\randomVar{X}$ is forced to a value very close to $x$";
		this last notion is actually perfectly fine,
		for example if "very close to $x$" is a measure dominated by
		(absolutely continuous w.\,r.\,t.) $\mu$.
		
		Essentially, the usual argument is "the formal problem is just
		a proxy for the real-world problem".
		Despite this, there is the (inherently formal) aspect
		of mathematical rigor. For (most) practical computations,
		the interpretation of do-interventions given above makes sense.
		But what about proofs? A proof-technique that can "prove" both
		true and false statements does not seem particularly convincing.
		Overly simplistic notions of conditional probability
		(like quotients of joints)
		have long since been known for their pitfalls
		(see \eg \citep{chang1997conditioning}).
		The problem here is as follows:
		Consider two IID models $M_1$ and $M_2$ with
		$\randomVar{X} \rightarrow \randomVar{Y}$
		and $X\sim\mathcal{N}(0,1)$,
		in $M_1$ the mechanism at $Y$ is
		$f_Y(x) = x + \mathcal{N}(0,1)$ if $x\neq 0$
		but $f_Y(x) = 42 + \mathcal{N}(0,1)$ if $x=0$,
		while in $M_2$ the mechanism at $Y$ is
		simply $g_Y(x) = x + \mathcal{N}(0,1)$.
		Both models produce the same joint observational distribution
		$P(\randomVar{X},\randomVar{Y})$, yet
		$P_1(\randomVar{Y}|\PearlDo(\randomVar{X}=0))
		\neq P_2(\randomVar{Y}|\PearlDo(\randomVar{X}=0))$.
		Thus the do-intervention is (formally) not identifiable 
		from observations or even from perfect knowledge of
		the observational joint distribution.
		The problem is, that by $P(\randomVar{X}=0)=0$, so
		we will almost never (with probability $0$)
		figure out the difference between both models.
		Note that $p(\randomVar{X})>0$, so a positivity-assumption
		cannot fix this issue.
		This may often not be a problem in computations (see above),
		but obviously, if it is not actually formally true
		that this singular intervention can be identified (the example
		above is a counter-example),
		then any rigorous mathematical approach must \emph{see}
		this formal problem.
		
		Using a rigorous language for conditional probabilities,
		on standard Borel spaces (Ass.\ \ref{ass:standard_borel})
		disintegrations are available, the problem actually becomes apparent:
		The general idea of identifying $P(\randomVar{Y}|\PearlDo(\randomVar{X}=0))$
		is to write down a Markov-factorization (actually a disintegration)
		$P(\randomVar{X},\randomVar{Y})=\int P(X=x) P(Y|X=x) \dx$,
		then to replace $P(X)$ by the intervention $\delta[0]$
		singular at $0$, and compute the new joint distribution.
		Writing $\mu=P(X)$, $\nu_x = P(Y|X=x)$ (a regular version),
		then $P(\randomVar{X},\randomVar{Y})=\mu\otimes\nu$.
		However, it also becomes clear (cf.\ Cor.\ \ref{cor:disint_measure})
		that $\nu$ (computed form the joint) is unique only $\mu$-almost everywhere.
		So the second step, replacing $\mu$ by $\xi$
		(more in line with our symmetry-guided perspective:
		transfer $\nu_x$ from the observational system to the query;
		for example $\xi=\delta[0]$) produces a well-defined (unique)
		result $\xi\otimes\nu$ iff $\xi \ll \mu$ is dominated by $\mu$
		(\ie iff $\forall$ measurable $B$: $\mu(B)=0\Rightarrow\xi(B)=0$).
		For the example above $\{0\}$ is a measurable set
		and $\mu(\{0\})=P(\randomVar{X}\in\{0\})=0$,
		but $\xi(\{0\})=\delta[0](\{0\})=1$, in particular
		$\delta[0] \not\ll \mu$ and our formalism correctly alerts us
		to the fact
		that the transfer is not well-defined.
		
		Practical interpretation is often possible by replacing
		$\xi\ll\mu$ by a regularity assumption, for example by
		existence of smooth densities.
		Given a smooth density, a limit over a Dirac-sequence
		can be exchanged with integration, thus $\delta[0]$
		can be safely expressed as
		\begin{align*}
			\delta[0] \otimes \nu_x
			&=
			\big(\lim_{n\rightarrow \infty} \delta_n[0]\big)\otimes \nu_x & \\
			&=
			\lim_{n\rightarrow \infty} \kernelCompound{\delta_n[0]\otimes \nu_x}
			&\txt{by smoothness.}
		\end{align*}
		So in this case as long as a Dirac-sequence with
		$\forall n: \delta_n[0] \ll \mu$ can be chosen,
		the computation is "safe".
		
		The uniqueness "$\vartheta\otimes\mu$"-almost everywhere
		in our formalism can mostly
		(for cases like example \ref{example:random_walk} this
		problem is more complicated, \eg it depends on $b\neq 0$
		or dimension $\geq 3$)
		be traced back to
		some data-set on which we learned our initial information.
		So it is a formal means of making this reliance on
		support in data explicit.
		On the other hand, reliance on observational support
		is itself a very complicated topic, and not the main topic
		of this paper.
		So we will usually "forget about" observational support
		and simply write $\mu = \nu$ if this is true
		almost surely on our data.
		We make this transparent however in Ass.\ \ref{ass:support},
		and want to emphasize that our approach does
		provide a formal means of making explicit and analyzing this problem,
		thus sets the stage for future work on this topic.
		
		While often ignored, observational support has been understood
		to be relevant for learning and interpreting causal structure on 
		multi-context systems \citep{gunther2024causal, EndoTheory}.
	\end{rmk}
	
	Another, less deep, aspect of our notation so far is that
	it tracks order of factors even though
	this order can be trivially changed post hoc, thus is
	largely irrelevant (cf.\ §\ref{apdx:kernels_order_independent_ops} for details).
	We will usually suppress this order in notation:
	
	\begin{notation}[Equality and Transposition]
		\label{notation:uniqueness_transposition}
		We call two structured kernels $\mu, \nu$ essentially equivalent
		$\mu \approx \nu$,
		if $\mu^T = \nu$ (Def.\ \ref{def:transposition}).
		For example, given two $\graphStructural$-causal orders
		$\pi, \pi'$ then
		$\mu(\graphStructural,L,\pi) \approx \mu(\graphStructural,L,\pi')$
		by Lemma \ref{lemma:causal_transformations}.
		We will write $\mu(\graphStructural,L)$ for the $\approx$-equivalence
		class of any (thus all) causal orders.
		
		\emph{Remark:} The reader may think of this as
		a writing $\mu(\graphStructural,L)$ for $\mu(\graphStructural,L,\pi)$
		in cases where $\pi$ is not of interest, with slight abuse of
		notation (as equality "$=$" on equivalence classes technically
		means $\approx$ on actual kernels).
	\end{notation}

	Knowing something about a system, in form of one or multiple
	structured kernels, we seek to understand what else can be concluded
	or "computed" from this information.
	We will revisit this question in §\ref{apdx:identifiablity}
	when discussing identifiability more generally.
	In the IID case, most treatments
	\citep{tian2002general,shpitser2006identification1, shpitser2006identification2},
	allow for "general probability-theoretical transformations" (and the do-calculus,
	which is itself proved by probability-theoretical transformations
	employing causally mandated independencies \citep{PearlBook}).
	This is in practice taken to mean "quotients" (or formally more
	rigorous and more generally applicable: disintegrations
	\citep{chang1997conditioning}) to compute conditional distributions,
	marginalizations, and certain integrals
	(corresponding to $\otimes$-products in our case).
	This does, in the IID-case, not seem to lead to any confusion,
	as in practice there are very few available statements (§\ref{apdx:iid}).

	In the present case, with a purely measure-theoretic model
	(and correspondingly no troubles from trying to define suitable
	random elements) we are in a position where we can actually
	state, in a mathematically rigorous way, what exactly
	we mean by computation.
	Indeed, there is not a unique meaningful notion:
	It is well-understood that for example in the linear
	case instrumental-variable arguments can improve identification-results.
	This is one of many examples that can be understood as an extended notion
	of computation and is discussed in §\ref{sec:knowledge_closures}.
	We will, for the purposes of the present section, start with
	what in the extended notion of §\ref{sec:knowledge_closures} is called
	standard-regular (computations), \ie we make, for now, no assumptions
	about the internal structure (like linearity) of mechanisms.
	The following notion of regular computations
	will in the IID case reproduce the results
	of \citep{tian2002general,shpitser2006identification1, shpitser2006identification2,
		bareinboim2012TransportCompleteness} (see §\ref{apdx:iid}),
	but applies much more generally:
	
	\begin{definition}[Regularity]\label{def:regular_functionals}
		We call a functional $F$ of probability kernels
		(standard\mbox{-)}\allowbreak{}regular 
		if $F$ is the identity or projection to one of its arguments, or if
		$F$ consists of a finite combination of $\otimes$-products,
		left-disintegrations, right-marginalizations and transpositions.
		
		More formally
		$F[\mu_1, \ldots, \mu_n]$ is regular if $\exists F^1, \ldots, F^m$
		such that
		defining $\mu^0_1 := \mu_1$, $\ldots$, $\mu^0_n=\mu_n$,
		then inductively for $i+1=1$ to $m$,
		let $k\leq n+i$: $\mu^{i+1}_k = \mu^i_k$ and
		$\mu^{i+1}_{n+i+1} = F^{i+1}[\mu^{i}_1, \ldots, \mu^i_{n+i}]$
		where each
		$F^i$ is either of the form $\mu^i_l \otimes \mu^i_j$
		or a right-marginalization of an argument $\mu^i_j$ or a left-disintegration
		of an argument $\mu^i_j$
		or a transposition of an argument $\mu^i_j$
		and finally $F[\mu_1, \ldots, \mu_n] = \mu^m_{n+m}$.
	\end{definition}
	\begin{example}
		\label{example:regular_functionals_simple_rules}
		Any finite
		composition of regular functionals is thus also regular.
		Anti-causal disintegrations can be written as a combination
		of a transposition and a left-dis\-in\-te\-gra\-tion, thus are regular.
		The general marginalization and disintegration operators
		can be written as a combination of left-disintegrations,
		anti-causal disintegrations, transpositions and right-marginalizations,
		thus are regular. Compositions are marginalizations of products, thus are
		regular.
	\end{example}
	\begin{rmk}[Regularity]
		We call these functionals regular, without actually fixing a regularity condition.
		The idea is as follows: In practice, given some assumptions (like smoothness etc)
		about mechanisms and an estimator (for example of a kernel\Slash{}"conditional density") suitable (converging) under those assumptions, it is usually possible
		to give an actual regularity condition (like smoothness in some sense)
		such that convergence of the estimator(s) in the arguments entails
		convergence of the functional.
		For the chosen estimator to be applicable with the formalism outlined,
		this estimator-specific regularity must at least include
		Def.\ \ref{def:regular_functionals}; in practice, additionally data-support
		is relevant -- we cannot intervene to a range of values we have never seen
		(without extremely strong parametric assumptions like linearity),
		cf.\ Rmk.\ \ref{rmk:kernels_uniqueness}. 
	\end{rmk}
	
	On generic kernels (\ie without internal structure)
	no further computations should be possible, motivating
	Conjecture \ref{conjecture:regular_computation}, a relevant details to
	this problem is however:
	\begin{rmk}[Fixed-Points and Limits]
		\label{rmk:regular_limits}
		Given a kernel $\mu_x$ such that
		$\kernelCompound{\mu_x \otimes \mu_{x'}}_x$
		is well-defined (\eg the Markov-kernel of a time-series),
		it is reasonable to ask for a limit
		$\lim_{n\rightarrow\infty} \kernelCompound{\otimes_{k=1}^n \mu}_x$
		(\eg a stationary distribution of a time-series).
		There is no reason why this computation
		would be regular (in the sense above).
		However, there is also no reason why such limits
		would be well-defined (exist uniquely) for \emph{generic}
		kernels, so Conjecture \ref{conjecture:regular_computation}
		still seems plausible, albeit for non-trivial reasons.
	\end{rmk}	
	While this clarifies what computations we consider on general kernels,
	it does not yet tell us how these computations on kernels, when applied
	to structured kernels, correspond to operations on graphs of structured kernels.

	\subsection{Graphical Operations}
	\label{apdx:tranform_graphs}
	
	In this subsection, we describe how structured kernels,
	can be transformed and combined.
	To this end, we study which graphical operations correspond
	to regular computations (see §\ref{apdx:regularity})
	on the described structured kernels.
	Note that the main text uses notation \ref{notation:uniqueness_transposition}
	to suppress causal orders $\pi$ in these results.
	The main technical result is Lemma \ref{lemma:properties_c_components}~%
	\ref{lemma:properties_c_components:atoms}.
	This is actually a result about "atoms", factors associated
	to individual observable nodes of a structured kernel.
	The full power of Lemma \ref{lemma:properties_c_components}~%
	\ref{lemma:properties_c_components:atoms} is a little obscured by its
	very simple and very technical claim, but it immediately implies
	Lemma \ref{lemma:c_components:mt} (c-components, cf.\ Lemma \ref{lemma:properties_c_components}~%
	\ref{lemma:properties_c_components:regular}),
	makes Lemma \ref{lemma:glue:mt} (gluing, cf.\ Lemma \ref{lemma:glue:apdx})
	astonishingly easy to state and proof
	and is a major ingredient in the proof of Thm.\ \ref{thm:extract_from_backdoor_complete}
	(extraction from backdoor-free families, cf.\ §\ref{apdx:id_from_embeddings}).
	Given its central importance, the proof of Lemma \ref{lemma:properties_c_components}~%
	\ref{lemma:properties_c_components:atoms} is stated with great detail;
	the reader may want to first go through Rmk.\ \ref{rmk:compute_atoms},
	which explains the general idea behind the structures studied 
	carefully in the proof of Lemma \ref{lemma:properties_c_components}~%
	\ref{lemma:properties_c_components:atoms}.
	
	We will start with sub-graphs: Given $\mu(\graphStructural,L,\pi)$, which
	sub-graphs $\graphStructural'\leq \graphStructural$
	encode a kernel $\mu(\graphStructural',L',\pi')$ that
	can be computed from $\mu(\graphStructural,L,\pi)$?
	Besides trivial simplifications that only delete hidden parts
	of the graph (possibly after first hiding them by marginalization),
	there is also an analogue of c-components
	which play an important role in the IID case \citep{tian2002general},
	and their simplifications, which form an analogue of c-trees and c-forests
	\citep{shpitser2006identification1,shpitser2006identification2}.
	Sub-graphs (at least as they are formulated here)
	are local, in the sense that information
	about the kernel of a sub-graph $\graphStructural'\leq\graphStructural$
	is obtained from the kernel of $\graphStructural$, without any reference to other
	possibly known information.
	The next construction glues two graphs with controlled overlap,
	that is two "smaller" graphs' kernels provide information about
	the kernel associated to a larger graph containing both as sub-graphs.
	If both smaller graphs agree on which variables are hidden,
	this operation is semi-local, in the sense that it requires
	two, but only two, graphs as inputs.
	Finally, we conclude with a construction that reveals\Slash{}imputates
	certain	hidden nodes in a larger graph, by virtue of using a second
	smaller graph with smaller latent set.
	This construction seemingly is also semi-local (as is gluing).
	A more general revealing operation might be possible
	(see §\ref{apdx:reveal} for further discussion),
	in which case gluing on revealed parts might become
	truly non-local.
	
	We start by a more detailed version of Lemma\Slash{}Def.\ \ref{def:simplify_simp}.
	
	\begin{lemmaDef}[Simplified Graphs]
		\label{lemma:simplify}
		Given a structural model $(\graphStructural, L)$,
		then for any $L' \subset\nodesInner$,
		$\mu(\graphStructural, L \cup L') =
		\marginalize{\mu(\graphStructural, L)}{L'\setminus L}$
		can be regularly computed from $\mu(\graphStructural, L)$.
		
		Further, for a sub-set $B \subset L$
		with $\Dec_\graphStructural(B) \subset B$,
		there is a $\graphStructural$-causal order $\pi$
		putting $B$ at the end, \ie
		such that $\forall n\in\nodes \forall b\in B$:
		$(n\geq b \Rightarrow n\in B)$.
		Using this $\pi$,
		the sub-graph $\graphStructural'\leq\graphStructural$ with
		nodes $\nodesInner'	=\nodesInner\setminus B$
		and $\nodesOuter'=\Pa_{\graphStructural}(\nodesInner')\setminus\nodesInner'
		\subset \nodesOuter$ (with $B\subset L \subset \nodesInner$
		removing $B$ only removes inner nodes, however, some of the original
		arguments, \ie outer nodes, may only be relevant to nodes in $B$,
		in which case we can discard them)
		satisfies
		\begin{equation*}
			\mu(\graphStructural',L\setminus B, \pi)
			=
			\mu(\graphStructural, L, \pi)
			\txt.
		\end{equation*}
		We call $(\graphStructural',L\setminus B)$ a direct simplification of $(\graphStructural, L)$.
		There is a unique maximal such $B$ (in the sense that any $B'$ satisfying the hypothesis is a subset $B'\subset B$),
		and we call the corresponding
		(minimal) $\graphStructural'$ simplified.
		The maximal $B$ can be constructed in practice as described by the proof.
		
		Combining both, given a sub-set $B \subset \nodesInner$
		with $\Dec_\graphStructural(B) \subset B$,
		define $L':= B\setminus L$, and again
		the sub-graph $\graphStructural'\leq\graphStructural$ with
		nodes $\nodesInner'	=\nodesInner\setminus B$
		and $\nodesOuter'=\Pa_{\graphStructural}(\nodesInner')\setminus\nodesInner'
		\subset \nodesOuter$,
		then $\mu(\graph', L \setminus B)$ can be regularly computed
		from $\mu(\graph, L)$, we call this construction a simplification.
	\end{lemmaDef}
	\begin{proof}
		The first claim follows directly by definition, the last claim
		follows from the first two claims,
		thus we proof the remaining second claim.
		
		For the existence of $\pi$, first choose a
		$\pi^A$ of $\nodes\setminus B$ (compatible with the ancestral
		partial order) and $\pi^B$ of $B$ (compatible with the ancestral
		partial order).
		Define the total order $\pi$ (\ie $\leq_\pi$) as follows:
		Given arbitrary $n,n'\in\nodes$, if both $n,n'\in B$,
		define $n \leq_\pi n' :\Leftrightarrow n \leq_{\pi^B} n'$,
		if both $n,n'\notin B$,
		define $n \leq_\pi n' :\Leftrightarrow n \leq_{\pi^A} n'$,
		otherwise (one of $n,n'$ is in $B$, one is not),
		define $n \leq_\pi n' :\Leftrightarrow n'\in B$
		(\ie order nodes not in $B$ before those in $B$).
		This total order is compatible with the ancestral
		partial order on $\nodes\setminus B$ and $B$
		because $\pi^A$ and $\pi^B$ are, on mixed terms because
		$\Dec_\graphStructural(B) \subset B$ (\ie if $n,n'$ with exactly one of them
		in $B$ are ancestral-comparable, then they must compare
		as the one not in $B$ ordered after the one in $B$).
		By construction $\pi$ satisfies
		$\forall n\in\nodes \forall b\in B$:
		$(n\geq_\pi b \Rightarrow n\in B)$. We will drop the subscript $\pi$
		on comparison operators for the remainder of this proof.
		
		Using $N := \nodesInner\cup(\nodesOuter\cap L)$ we compute:
		\begin{align*}
			\mu(\graphStructural,L, \pi)
			&=
			\marginalize{\otimes^\pi_{n\in N}\mu^n}{L}
			&\txt{(by definition)}
			\\
			&=
			\marginalize{\kernelCompound{\otimes^\pi_{n\in N\setminus B}\mu^n}
			\otimes
			\kernelCompound{\otimes^\pi_{n\in B}\mu^n}
			}{L}
			&\txt{(by property of $\pi$)}
			\\
			&=
			\marginalize{\kernelCompound{\otimes^\pi_{n\in N\setminus B}\mu^n}}
			{L\setminus B}
			&\txt{(trivial marg.\ on rhs)}
			\\
			&=
			\mu(\graphStructural',L\setminus B, \pi)
			&\txt{(by definition).}
		\end{align*}
		
		To construct a maximal $B$, filter $L$ by longest descendant-chains,
		\ie $B^{(0)}:=\{n\in L|\Ch(n)=\emptyset\}$ and inductively
		$B^{(k+1)}:=\{n\in L|\Ch(n)\subset B^{(k)}\}$
		(note that $B^{(k)}\subset B^{(k+1)}$).
		Define $B:= B^{(\infty)}:=\cup_k B^{(k)}$ (since $L$ is finite, obviously there is
		a finite number $m$ where this filtration becomes stationary and
		$B^{(\infty)}=B^{(m)}$ is not a "real" limit, just a convenient notation).
		
		Let $B'$ satisfying the hypothesis be arbitrary.
		We show $B'\subset B^{(\infty)}$.
		Since $B'$ satisfies the hypothesis $B'\subset L$ and
		$\Dec(B') \subset B' \subset L$.
		Define $(B')^{(0)}:=\{n\in B'|\Ch(n)=\emptyset\}$,
		then $(B')^{(0)}\subset B^{(0)}$ (because $B'\subset L$).
		Define inductively
		$(B')^{(k+1)}:=\{n\in B'|\Ch(n)\subset (B')^{(k)}\}$,
		then (inductively, \ie by $(B')^{(k)} \subset B^{(k)}$)
		we have for arbitrary $n\in (B')^{(k+1)}$ that $n\in B' \subset L$
		and that $\Ch(n)\subset (B')^{(k)} \subset B^{(k)}$,
		thus $n\in B^{(k+1)}$, showing $(B')^{(k+1)}\subset B^{(k+1)}$;
		after a finite number of steps $(B')^{(\infty)} = (B')^{(m)} \subset B^{(m)}
		\subset B^{(\infty)}$.
		Finally, by $\Dec(B') \subset B'$ (and by finiteness of $\nodes$)
		$B'=(B')^{(\infty)}$.
	\end{proof}
	
	The main technical ingredient for proving the more sophisticated
	statements below is a
	decomposition of the product describing $\mu(\graphStructural,L,\pi)$
	in per-observable-node factors or "atoms":
	For each observed inner node $n$,
	we define $A^n$ for a causal order $\pi$ on nodes
	by disintegrating from the left until we reach $n$
	and marginalizing from the right until we reach $n$
	(by Lemma \ref{lemma:product_operator_properties}
	\ref{lemma:product_operator_properties:disint_marg} the order in
	which this is executed does not matter).
	
	\begin{definition}[Atoms]\label{def:atoms}
		Given a structural model $(\graphStructural, L)$, a
		\graphStructural-causal order $\pi$, then
		for $n\in\nodesInner\setminus L$,
		we define
		the atom at $n$ as
		\begin{equation*}
			A^n(\graphStructural, L,\pi)
			:= \disint_{\{n'\in\nodesInner\setminus L|n'<_\pi n\}}
			\marginalize{\mu(\graphStructural,L,\pi)}{L\cup\{n'\in\nodesInner\setminus L|n'>_\pi n\}}
			\txt.
		\end{equation*}
	\end{definition}
	
	These atoms by construction are regularly computable for
	a given structured kernel and turn out to contain all relevant information
	to reconstruct the structured kernel:
	
	\begin{lemma}[Properties of Atoms]
		\label{lemma:atom_properties}
		Given a structural model $(\graphStructural, L)$
		and a $\graphStructural$-causal order $\pi$
		then:
		\begin{enumerate}[label=(\alph*)]
			\item\label{lemma:atom_properties:atoms_regular}
			$A^n(\graphStructural, L,\pi)$ is a regular functional of $\mu(\graphStructural,L,\pi)$.
			
			\item\label{lemma:atom_properties:product}
			$
				\mu(\graphStructural,L,\pi) =
				\otimes^\pi_{n\in\nodesInner\setminus L} A^n(\graphStructural,L,\pi)
				\txt.
			$
			
			\item\label{lemma:atom_properties:regular_from_atoms}
			$\mu(\graphStructural,L,\pi)$ is a regular functional of its atoms.
			
			\item\label{lemma:atom_properties:simplification}
			Given $B \subset \nodesInner$
			with $\Dec_\graphStructural(B) \subset B$,
			and $\pi$ the $\graphStructural$-causal order and $\graphStructural'\leq
			\graphStructural$ are the results produced by
			Lemma \ref{lemma:simplify}, then
			$\forall n\in \nodesInner'\setminus L'$:
			$A^n(\graphStructural,L,\pi) = A^n(\graphStructural',L',\pi')$.
		\end{enumerate}
	\end{lemma}
	\begin{proof}		
		\textbf{Part \ref{lemma:atom_properties:atoms_regular}:}
		Disintegrations and marginalizations are regular, so are their compositions.

		\textbf{Part \ref{lemma:atom_properties:product}:}
		Inductively over $n\in\nodesInner \setminus L$ along $\pi$
		we show:
		Using $L^n := L \cup \{ n' \in \nodesInner | n' >_\pi n \}$
		and $\nodes^{\leq n} := \{ n' \in \nodesInner\setminus L | n' \leq_\pi n \}$,
		\begin{equation*}\mu(\graphStructural,L^n,\pi) =
			\otimes^\pi_{n'\in \nodes^{\leq n}} A^{n'}(\graphStructural,L,\pi)
			\txt.
		\end{equation*}
		
		For the start of induction (the first, in $\pi$-order,
		$n\in\nodesInner \setminus L$): By definition
		$A^n(\graphStructural,L,\pi)\allowbreak{}=\marginalize{
			\kernelCompound{\otimes_{l\in\Anc_\graphStructural(n)}\mu^l}
			\otimes \mu^n }{L}$,
		where ancestors of $n$ are in $L$ (by $n$ being the first observed
		inner node) by trivial marginalization on the right this is
		$\mu(\graphStructural,L^n,\pi)$.
		
		Inductive step ($n$ to "$n+1$", where $n+1$ is the next, in $\pi$-order,
		node in $\nodesInner \setminus L$ after $n$):
		By the characterizing property of the disintegration (Lemma \ref{def:disint_product}),
		$\mu(\graphStructural,L^n,\pi) \otimes A^{n+1}(\graphStructural,L,\pi)
		= \mu(\graphStructural,L^{n+1},\pi)$.
		By inductive hypothesis $\mu(\graphStructural,L^n,\pi) 
		= \otimes^\pi_{n'\nodes^{\leq n}} A^{n'}(\graphStructural,L, \pi)$,
		thus proving the claim.
		
		\textbf{Part \ref{lemma:atom_properties:regular_from_atoms}:}
		Immediate by \ref{lemma:atom_properties:product} (products are regular).

		\textbf{Part \ref{lemma:atom_properties:simplification}:}
		Intuitively, since $\pi$ puts $B$ at the end of the product
		in \ref{lemma:atom_properties:product}, and
		$\mu(\graph',L',\pi')=\marginalize{\mu(\graph,L,\pi)}{B\setminus L}$
		by Lemma \ref{lemma:simplify} is given by trivial marginalizations on the right
		this follows by comparing terms in the representation by \ref{lemma:atom_properties:product}.
		Formally, for example use
		$\mu(\graph',L',\pi')=\marginalize{\mu(\graph,L,\pi)}{B\setminus L}$
		by Lemma \ref{lemma:simplify}, then
		Lemma \ref{lemma:kernel_marginalization_safe} and
		Lemma \ref{lemma:product_operator_properties}~%
		\ref{lemma:product_operator_properties:marg_prod}
		and the definition of atoms.		
	\end{proof}
	\begin{rmk}[Computation of Atoms I]\label{rmk:compute_atoms}
		We briefly illustrate how atoms are related to the
		kernels $\mu^n$ associated to nodes $n$ of the graph
		(a more formal version can be found in the proof of
		Lemma \ref{lemma:properties_c_components}).
		First, we associate the $\mu^n$ into groups $L_i$ (if they are in $L$)
		and $Y_i$ (otherwise) and write (up to transposition, denoted "$\approx$", notation
		\ref{notation:uniqueness_transposition})
		\begin{equation*}
			\otimes_{n} \mu^n
			\halfquad\approx\halfquad
			\kernelCompound{ L_1 \otimes Y_1 }
			\otimes
			\ldots
			\otimes			
			\kernelCompound{ L_m \otimes Y_m }
		\end{equation*}
		We can, without loss of generality, assume that each
		$Y_i$-term is exactly a single node (by injecting
		empty $L_i$-terms with the understanding that $\mu\otimes\emptyset=\mu$
		by convention).
		We reorder this successively by moving $L$-terms further to the right.
		First, $L_1 \otimes Y_1 = (Y_1 \circ L_1) \otimes (L_1|Y_1)$.
		We write $L_1^\leq := L_1$, then this reads
		$L_1 \otimes Y_1 = (Y_1 \circ L_1^\leq) \otimes (L_1^\leq|Y_1)$
		The next term in the expression for $\otimes_{n} \mu^n$ is $L_2$,
		we collect this term into $L_2^\leq := L_1^\leq \otimes L_2$
		and swap it with $Y_2$ to get
		\begin{align*}
			\kernelCompound{ L_1 \otimes Y_1 }
			\otimes
			\kernelCompound{ L_2 \otimes Y_2 }
			&\approx
			(Y_1 \circ L_1^\leq)
			\otimes
			L_2^\leq
			\otimes
			Y_2\\
			&\approx
			(Y_1 \circ L_1^\leq)
			\otimes
			(Y_2 \circ L_2^\leq)
			\otimes
			(L_2^\leq|Y_2)
			\txt.
		\end{align*}
		Then inductively define $L_{k+1}^\leq := L_k^\leq \otimes L_{k+1}$
		and swap to get
		\begin{equation*}
			\otimes_{n} \mu^n
			\halfquad\approx\halfquad
			(Y_1 \circ L_1^\leq)
			\otimes\ldots\otimes
			(Y_m \circ L_m^\leq)
			\otimes L_{m+1}^\leq
			\txt.
		\end{equation*}
		This expression can be trivially marginalized at the right
		(all elements of $L$ were transposed\Slash{}swapped to the right
		and were collected into the product $L^\leq_{m+1}$) to compute
		\begin{equation*}
			\mu(\graphStructural,L,\pi)
			=
			\marginalize{\otimes_{n} \mu^n}{L}
			\approx
			(Y_1\circ L^\leq_1)
			\otimes\ldots\otimes
			(Y_m\circ L^\leq_m)
			\txt.
		\end{equation*}
		Disintegrating observable nodes from left to right is now
		a trivial disintegration on the left, thus we can read of
		\begin{equation*}
			A^n(\graphStructural,L,\pi)
			=
			Y_n\circ L^\leq_n\txt.
		\end{equation*}
		We have not yet considered which arguments these
		terms actually carry, we will investigate this
		in the proof of
		Lemma \ref{lemma:properties_c_components},
		and continue in Rmk.\ \ref{rmk:compute_atoms_practice}.
	\end{rmk}
	
	Similar to the IID-case \citep{tian2002general},
	c-components of graphs (in the sense of the next definition)
	play an important role for computing results,
	detailing Def.\ \ref{def:c_components_local:mt}:
	
	\begin{definition}[Local C-Component]\label{def:c_component_local}
		Given a structural model $(\graphStructural, L)$,
		we call $l\in L$ a hidden confounder of
		$y,w\in\nodesInner\setminus L$
		if there are directed paths $\gamma_y$ from $l$ to $y$ and $\gamma_w$
		from $l$ to $w$ such that all non-endpoint nodes of $\gamma_y, \gamma_w$
		are in $L$.
		On $\nodesInner\setminus L$ define a relation
		$y \sim w :\Leftrightarrow \exists$ hidden confounder $l$ of $y, w$.
		Define $\sim_L$ as the equivalence relation generated by $\sim$
		(the "smallest" equivalence relation with $y \sim w \Rightarrow y\sim_L w$,
		since $\sim$ is already symmetric,
		this amounts to
		trivially making $\sim$ reflective and
		transitively closing it, \ie define $y\sim_L w$ if
		$y=w$ or if $y\sim z_1 \sim \ldots \sim z_n \sim w$).
		We call equivalence-classes $\sfrac{\nodesInner\setminus L}{\sim_L}$
		c-components of $(\graphStructural,L)$
		and say $(\graphStructural,L)$ is c-connected if there is only a single c-component,
		$|\sfrac{\nodesInner\setminus L}{\sim_L}|=1$.
		
		We can extend this notion to $\nodesInner$
		by defining $l\approx y$ for $y\in \nodesInner\setminus L$ if there
		exists a directed path $\gamma_y$ from $l$ to $y$
		such that all non-endpoint nodes of $\gamma_y$ are in $L$.
		Then each $l\in\Anc_\graphStructural(\nodesInner\setminus L)\cap L$ is
		equivalent to at least one $y$, and if it is equivalent
		via $\gamma_w$ to another $w\in \nodesInner\setminus L$,
		then $y\sim w$ via $(\gamma_y,\gamma_w)$.
		Thus $\approx$ and $\sim$ together (\ie two nodes $n$, $n'$ are
		related if $n\approx n'$ or $n\sim n'$),
		then generates an equivalence relation $\approx_L$
		on $\Anc_\graphStructural(\nodesInner\setminus L)$
		with the same number of equivalence-classes as $\sim_L$
		(indeed the intersections of $\approx_L$ equivalence-classes
		of $\Anc_\graphStructural(\nodesInner\setminus L)$
		with $\nodesInner\setminus L$ are exactly the $\sim_L$
		equivalence-classes; every c-component thus has a unique associated
		$\approx_L$ equivalence class, containing additional hidden ancestors).
		Finally, we formally make all elements of
		$L\setminus\Anc_\graphStructural(\nodesInner\setminus L)$ equivalent,
		this adds zero or one additional classes
		(depending on whether $L\setminus\Anc_\graphStructural(\nodesInner\setminus L)=\emptyset$)
		witch we call a trivial c-component.
		
		A structural c-component is a subgraph
		$\graphStructural^c\leq\graphStructural$
		containing as inner nodes $\nodesInner^c$ the
		elements of a non-trivial $\approx_L$-equivalence class
		and as outer nodes
		$\nodesOuter^c = \Pa_{\graphStructural}(\nodesInner^c)\setminus\nodesInner^c$.
		By slight abuse of notation, we also refer to the
		structural model $(\graphStructural^c,L^c)$ with
		$L^c:= L\cap\nodes^c$ as structural c-component;
		note that by construction of $\approx_L$,
		$\nodesInner^c$ is closed under latent parents,
		thus $L^c\subset\nodesInner^c$.
	\end{definition}
	
	These c-components are a direct analogue of c-components in the IID-case
	which are well-understood to be of great importance to the
	question of identifiability of interventional distributions \citep{tian2002general}.
	Similarly we find also here (note that part (b) is
	Lemma \ref{lemma:c_components:mt}):
	
	\begin{lemma}[Properties of C-Components]		
		\label{lemma:properties_c_components}
		Given a structural model $(\graphStructural, L)$,
		a $\graphStructural$-causal order $\pi$
		and a c-component
		$(\graphStructural^c, L^c) \leq (\graphStructural, L)$ with
		$\graphStructural^c$-causal order $\pi^c=\pi|$
		then:
		\begin{enumerate}[label=(\alph*)]
			\item \label{lemma:properties_c_components:atoms}
				If $n\in \nodesInner^c\setminus L^c$,
				then
				$A^n(\graphStructural^c,L^c,\pi^c) = A^n(\graphStructural,L,\pi)$.
			\item \label{lemma:properties_c_components:regular}
				$\mu(\graphStructural^c,L^c,\pi^c)$ is a regular functional
				of $\mu(\graphStructural,L,\pi)$.
			\item 
				$\mu(\graphStructural,L,\pi)$ is a regular functional of
				jointly all non-trivial c-components'
				$\mu(\graphStructural^c,L^c,\pi^c)$.
		\end{enumerate}
	\end{lemma}
	\begin{proof}
		We follow the idea outlined in Rmk.\ \ref{rmk:compute_atoms},
		but carefully track kernel-arguments.
		For any fixed structural c-component,
		we will be able split the $L^\leq$ terms of the remark into
		terms $L^\leq$ relevant to that particular c-component
		and terms $W^\leq$ irrelevant to that particular c-component.
		We write $C = \nodesInner^c \setminus L^c \subset\nodesInner\setminus L$
		for the c-component with associated structural c-component $\graphStructural^c$.
		We may assume \asswlog that there are only non-trivial
		structural c-components (otherwise simplify, Lemma \ref{lemma:simplify}).
		
		We may assume, without loss of generality
		that the $\pi$-ordered $\nodesInner=\{1,\ldots,m\}$.
		Thus $\mu(\graphStructural,L,\pi)=\marginalize{\mu^1 \otimes\ldots\otimes\mu^m}{L}$.
		To simplify notation, we define "empty" kernels $L_0^\leq=\emptyset$ and $W_0^\leq=\emptyset$
		(with the understanding that $\forall\mu$:
		$\emptyset\otimes\mu=\mu\otimes\emptyset=\mu$ and
		$\mu \circ \emptyset=\mu$; this is merely for notational convenience),
		then we inductively for $k=1,\ldots,m$
		construct $L_k^\leq, W_k^\leq$ with the properties (to be shown below):
		\begin{enumerate}[label=(\alph*)]
			\item
				$L_k^\leq$ and $W_k^\leq$ have arguments only in observed nodes
				(nodes not in $L$).
			\item
				$L_k^\leq$ and $W_k^\leq$ contain exactly those
				factors associated to nodes $n\leq k$ with
				$n\in L^c$ and $n\in L\setminus L^c$ respectively
				(\ie the transposition, Def.\ \ref{def:transposition},
				$T$ in (c) is such that it puts these terms to the right,
				see also proof below).
			\item
				$(\otimes_{i\leq k}\mu^i)^T
				= \kernelCompound{\otimes_{i \notin L, i\leq k}
					(\mu^i \circ Z_{i-1}^\leq)}
					\otimes
					L_k^\leq \otimes W_k^\leq$,
					where $Z_{i-1}^\leq=L_{i-1}^\leq$ if $i\in C$,
					and $Z_{i-1}^\leq=W_{i-1}^\leq$ otherwise.
		\end{enumerate}
		For the inductive start at $k=0$, define $L^\leq_0 =\emptyset$,
		$W^\leq_0 =\emptyset$, then (a,b) are trivial, defining the empty
		$\otimes$-product as $\emptyset$ also (c) is trivial
		(and this will be consistent, as for $k=1$, the product
		containing only $\mu^1$ by definition is $\mu^1=\emptyset\otimes\mu^1$).
		
		For the inductive step ($k-1 \mapsto k\geq 1$), we note that
		for each $k=1,\ldots, m$ we are in exactly one of the following cases:
		\begin{enumerate}[label=\arabic*)]
			\item
				Case $k\in L^c$: Define $W_k^\leq := W_{k-1}^\leq$ and
				$L_k^\leq := L_{k-1}^\leq \otimes \mu^k$.
				Arguments of $L_k^\leq$ are those of $L_{k-1}^\leq$ plus
				those of $\mu^k$, \ie $\Pa(k)$ (parents of inner nodes
				in $\graphStructural^c$ agree with those in $\graphStructural$,
				cf.\ Def.\ \ref{def:subgraphs_simp}).
				Note that $\Pa(k)\cap L \subset L^c$:
				Let $p\in \Pa(k)\cap L$ be arbitrary.
				Since $k\in L^c\subset \nodesInner^c$,
				$k\approx_L y$ for some $y\in C$, \ie
				there is a path $\gamma:k\rightsquigarrow y$ with non-endpoints in $L$.
				Prefixing this path with the edge $p\rightarrow k$
				we get $p\approx_L y$, thus $p\in \nodesInner^c$.
				Since also $p\in L$, and $L^c=L\cap\nodesInner^c$ by
				definition, $p\in L^c$.
				
				We show (a):
				Arguments of $W_k^\leq$ and $L_{k-1}^\leq$ already satisfy (a) by inductive
				hypothesis (a).
				Given $p\in \Pa(k)\cap L$,
				it remains to show $p$ is not an argument of $L_k^\leq$,
				\ie it is contracted with $L_{k-1}^\leq$.			
				Clearly parents are ancestors, thus $\pi$-before
				$k$, \ie $p<_\pi k$,
				thus this parent $p\leq k$ with $p\in L^c$ is contracted into $L_{k-1}^\leq$ by inductive hypothesis (b).
				
				We show (b): By inductive hypothesis (b) this is true
				for $n<_\pi k$. We absorbed $\mu^k$ into
				$L_k^\leq := L_{k-1}^\leq \otimes \mu^k$, thus it remains
				true for $n\leq_\pi k$.
				
				We show (c): Plugging in the inductive hypothesis (c),
				$(\otimes_{i\leq k}\mu^i)^T
				= \kernelCompound{\otimes_{i \notin L, i\leq k-1}
					(\mu^i \circ Z_{i-1}^\leq)}
				\otimes
				L_{k-1}^\leq \otimes W_{k-1}^\leq \otimes \mu^k$.
				The new term $\mu^k$ has only observed arguments
				and arguments in $L^c$ (a; shown above)
				while $W_k^\leq$ contains only hidden node factors in 
				$L \setminus L^c$ (b; shown above),
				so $W_{k-1}^\leq \otimes \mu^k = \mu^k \otimes W_{k-1}^\leq$
				by sparsity, Lemma \ref{lemma:causal_transformations}.
				Thus (c) follows with the definition of
				$L_k^\leq = L_{k-1}^\leq \otimes \mu^k$.
			\item 
				Case $k\in C$: Define $W_k^\leq := W_{k-1}^\leq$ and
				$L_k^\leq := (L_{k-1}^\leq|\mu^k)$.
				
				We show (a):
				By inductive hypothesis arguments of $W_k^\leq$
				and $L_{k-1}^\leq$ are in observed nodes.
				$L_k^\leq := (L_{k-1}^\leq|\mu^k)$ has additional
				arguments $k$ and $\Pa(k)\setminus L^c$ (from the disintegration,
				Lemma \ref{lemma:causal_transformations}~%
				\ref{lemma:causal_transformations:disint_args}
				and using the inductive hypothesis (b)). $k$ itself
				is observed by case-hypothesis.
				$\Pa(k)$ are either observed or in $L^c$ by definition of $\graphStructural^c$, so $\Pa(k)\setminus L^c$ are observed.
				
				We show (b):
				Since $k$ is observed by case-hypothesis
				the claim of (b) is already given by the inductive
				hypothesis (b).
				
				We show (c):
				Plugging in the inductive hypothesis (c),
				$(\otimes_{i\leq k}\mu^i)^T
				= \kernelCompound{\otimes_{i \notin L, i\leq k-1}
					(\mu^i \circ Z_{i-1}^\leq)}
				\otimes
				L_{k-1}^\leq \otimes W_{k-1}^\leq \otimes \mu^k$.
				All hidden arguments of
				the new term $\mu^k$ are in $L^c$ (by definition
				of $\approx_L$ and $\graph^c$)
				while $W_k^\leq$ contains
				only hidden node factors in $L\setminus L^c$ (b, shown above),
				so $W_{k-1}^\leq \otimes \mu^k = \mu^k \otimes W_{k}^\leq$
				by sparsity, Lemma \ref{lemma:causal_transformations}.
				By anti-causal disintegration
				$\kernelCompound{\otimes_{i \notin L, i\leq k-1}
					(\mu^i \circ Z_{i-1}^\leq)}
				\otimes
				L_{k-1}^\leq \otimes \mu^k \otimes W_{k-1}^\leq
				=
				\kernelCompound{\otimes_{i \notin L, i\leq k-1}
					(\mu^i \circ Z_{i-1}^\leq)}
				\otimes
				(\mu^k \circ L_{k-1}^\leq) \otimes
				(L_{k-1}^\leq|\mu^k) \otimes W_{k-1}^\leq$.
				Thus (c) follows with the definition of $L_k^\leq$.
				
			\item 
				Case $k\in L\setminus L^c$,
				define $W_k^\leq := W_{k-1}^\leq \otimes \mu^k$
				and $L_k^\leq = L_{k-1}^\leq$.
				Arguments of $W_k^\leq$ are those of $W_{k-1}^\leq$ plus
				those of $\mu^k$, \ie $\Pa(k)$.
				Note that $\Pa(k)\cap L^c = \emptyset$:
				By contradiction. If there were a $p\in \Pa(k)\cap L^c$,
				then there is a latent path $\gamma_y: p\rightsquigarrow y$ for some
				observed $y\in C$.
				Since $k\notin \nodesInner^c$ is in a non-trivial (\asswlog
				via simplify, see above) structural c-component
				associated to $k\in C' \neq C$,
				there is a latent path $\gamma_w: k \rightsquigarrow w$ for
				some observed $w\in C'\setminus L$.
				Prefix this path with the edge $p\rightarrow k$
				to obtain the latent path $\gamma_w': p\rightsquigarrow w$.
				Using $(\gamma_y, \gamma_w')$, thus $y \sim_L w$ and
				$C' = C$ contradicting $C' \neq C$.
				
				We show (a):
				Arguments of $W_{k-1}^\leq$ and $L_k^\leq$ are
				observed by inductive hypothesis (a).
				Clearly parents are ancestors, thus $\pi$-before
				$k$, \ie $p<_\pi k$ and if $p\in L$, then
				$p\notin L^c$ (see above),
				thus such arguments are contracted into $W_{k-1}^\leq$.
				
				We show (b):
				For factors $n < k$ this is true by inductive hypothesis (b).
				For $n=k$ it remains true by definition of $W_k^\leq$.
				
				We show (c):
				Property (c) holds immediately by inductive hypothesis
				and definition of $W_k^\leq$.
			\item 
				Case $k\notin C\cup L$: Define $W_k^\leq := (W_{k-1}^\leq|\mu^k)$
				and $L_k^\leq:=L_{k-1}^\leq$.
				
				We show (a):
				By inductive hypothesis arguments of $W_k^\leq$
				and $L_{k-1}^\leq$ are observed nodes.
				$W_k^\leq := (W_{k-1}^\leq|\mu^k)$ has a additional
				arguments $k$ and $\Pa(k)\setminus (L\setminus L^c)$ (from the disintegration,
				Lemma \ref{lemma:causal_transformations}~%
				\ref{lemma:causal_transformations:disint_args}
				and using the inductive hypothesis (b)).
				$k$ itself is observed by case-hypothesis.
				$\Pa(k)$ are either observed or in $L\setminus L^c$ by definition of $\graphStructural^c$, so $\Pa(k)\setminus (L\setminus L^c)$ are observed.
				
				We show (b):
				Since $k$ is observed by case-hypothesis
				the claim of (b) is already given by the inductive
				hypothesis (b).
				
				We show (c):
				Plugging in the inductive hypothesis (c),
				$(\otimes_{i\leq k}\mu^i)^T
				= \kernelCompound{\otimes_{i \notin L, i\leq k-1}
					(\mu^i \circ Z_{i-1}^\leq)}
				\otimes
				L_{k-1}^\leq \otimes W_{k-1}^\leq \otimes \mu^k$.
				By anti-causal disintegration
				$\kernelCompound{\otimes_{i \notin L, i\leq k-1}
					(\mu^i \circ Z_{i-1}^\leq)}
				\otimes
				L_{k-1}^\leq \otimes W_{k-1}^\leq \otimes \mu^k 
				= 
				\kernelCompound{\otimes_{i \notin L, i\leq k-1}
					(\mu^i \circ Z_{i-1}^\leq)}
				\otimes
				L_{k-1}^\leq \otimes (\mu^k\circ W_{k-1}^\leq)
				\otimes (W_{k-1}^\leq|\mu^k)
				= 
				\kernelCompound{\otimes_{i \notin L, i\leq k-1}
					(\mu^i \circ Z_{i-1}^\leq)}
				\otimes
				L_{k}^\leq \otimes (\mu^k\circ W_{k-1}^\leq)
				\otimes W_{k}^\leq$
				by definitions of $L_k^\leq$ and $W_k^\leq$.
				
				$W_{k-1}^\leq$ has no hidden arguments (a; shown above)
				and all hidden arguments of
				the new term $\mu^k$ are in $L\setminus L^c$ (by definition
				of $\approx_L$ and $\graph^c$)
				thus all hidden arguments of 
				$\mu^k\circ W_{k-1}^\leq$ are in $L\setminus L^c$.
				On the other hand $L_k^\leq$ contains
				only hidden node factors in $L^c$ (b, shown above),
				so $L_{k}^\leq \otimes (\mu^k\circ W_{k-1}^\leq)
				= (\mu^k\circ W_{k-1}^\leq) \otimes L_{k}^\leq$
				by sparsity, Lemma \ref{lemma:causal_transformations}.
				Thus (c) follows.
		\end{enumerate}
		In the end, (after $m$-steps), we obtain
		by property (c):		
		$\otimes_{i}\mu^i
		= \kernelCompound{\otimes_{i \notin L}
			(\mu^i \circ Z_{i-1}^\leq)}
		\otimes
		L_m^\leq \otimes W_m^\leq$,
		where $Z_{i-1}^\leq=L_{i-1}^\leq$ if $i\in C$,
		and $Z_{i-1}^\leq=W_{i-1}^\leq$ otherwise.
		By property (b), $L_m^\leq \otimes W_m^\leq$
		contains (as factors) exactly the term in $L$,
		thus
		\begin{equation*}
			\mu(\graphStructural, L, \pi)
			=
			\marginalize{\otimes_{i}\mu^i}{L}
			= \otimes_{i \notin L}
			(\mu^i \circ Z_{i-1}^\leq)
			\txt.
		\end{equation*}
		By trivially marginalizing on the right, then
		trivially disintegrating from the left we get:
		\begin{equation*}
			\tag{$*$}
			A^i(\graphStructural,L,\pi)
			=
			(\mu^i \circ Z_{i-1}^\leq)
			\txt.
		\end{equation*}
		Next we apply the same process to $\graphStructural^c$
		instead of $\graphStructural$.
		By definition, $\graphStructural^c$ contains exactly the $k\in \graph^c$,
		thus cases 1) and 2) are invoked exactly as before,
		we obtain the same $L_k^\leq$, while
		steps 3) and 4) cannot occur, thus all $W_k^\leq$ are trivial.
		We thereby obtain:
		\begin{equation*}
			\mu(\graphStructural^c, L^c, \pi^c)
			=
			\marginalize{\otimes_{i\in C\cup L^c}\mu^i}{L^c}
			= \otimes_{i \in C}
			(\mu^i \circ Z_{i-1}^\leq)
			\txt.
		\end{equation*}
		By trivially marginalizing on the right, then
		trivially disintegrating from the left we get for $i\in C$:
		\begin{equation*}
			A^i(\graphStructural^c,L^c,\pi^c)
			=
			(\mu^i \circ L_{i-1}^\leq)
			\overset{(*), i\in C}{=}
			A^i(\graphStructural,L,\pi)
			\txt.
		\end{equation*}
		This proves part (a) of the lemma.
		
		Parts (b) and (c) follow immediately with
		Lemma \ref{lemma:atom_properties}.
	\end{proof}

	\begin{rmk}[Computation of Atoms II]
		\label{rmk:compute_atoms_practice}
		Fundamentally, atoms are computable essentially
		by definition via disintegrations from the left
		and marginalizations from the right.
		Marginalizations from the right are trivial in practice:
		They just drop the last factor.
		Disintegrations on the left can be complicated
		(they are like estimating conditional probabilities).
		
		If we want to compute an atom $A^n(\graph,L,\pi)$ in practice,
		then Lemma \ref{lemma:properties_c_components}
		\ref{lemma:properties_c_components:atoms}
		and its proof
		tells us something about the \emph{actual} arguments.
		That is, we drop terms strictly $\pi$-after $n$,
		then disintegrate those strictly $\pi$-before $n$,
		but since we know that the result is actually
		a function (or kernel) in arguments 
		$C \cup \Pa(C)$
		(where $C$ is the c-component containing $n$),
		we can simply ignore all other disintegrated nodes and shared
		arguments.
		Thus the kernel $A^n(\graph,L,\pi)$ in practice gets harder to
		estimate if the c-component is large, but at least it
		depends only on the c-component (and its parents).
		Estimating each atom in $C$ relative to its ancestors in $C$
		and parents of $C$
		is of course very closely related to estimating
		the joint distribution of nodes in $C$
		by the (not simplifiable in this case, cf.\ \citep{tian2002general})
		Markov-factorization
		of $P(C|\Pa(C)\setminus C)$.
		
		\emph{Example:}
		For the napkin graph (cf.\ Fig.\ \ref{fig:do_query}) in the IID case,
		with causal order $\pi=(L_1,L_2,W,Z,X,Y)$
		the atom $A^w$ is simply
		$P(W)=W_{[l_1,l_2]}\circ\kernelCompound{L_1\otimes L_2}$,
		the atom $A^z$ is simply $P(Z|W=w)=Z_w$,
		the atom $A^x$ is $P(X|W=w,Z=z)$ (because $Z$ is a parent and
		$W$ is an ancestor in the same c-component),
		but note that
		\begin{equation*}
			P(X|W=w,Z=z)=\big(X_{z,[l_1]}\circ (L_1|W)_w\big)_{w,z}
		\end{equation*}
		actually contains an anti-causal term!
		The atom at $A^y$ is $P(Y|W=w,Z=z,X=x)$ (because $W,X$ are ancestors
		in the same c-component, $Z$ is additionally a parent of the c-component).
		The dependence is non-trivial in $Z$,
		because $Y \notindependent Z | W,X$, due to the colliders at
		$W,X$ there is an open path
		$Z \rightarrow X \leftrightarrow W \leftrightarrow Y$.
		To estimate $A^y$ from data, we only have to know its arguments (see above),
		yet the question remains: What is $A^y$ as expressed by model-kernels?
		Tracing Rmk.\ \ref{rmk:compute_atoms} and the proof of the lemma we find
		(cf. $\pi$ given above)
		$L_1^\leq = L_1$,
		$L_2^\leq = L_1 \otimes L_2$,
		from $W\in C$ at index $3$ we get
		$L_3^\leq = (L_2^\leq|W)_w = (L_1 \otimes L_2|W)_w$
		from $Z\notin C$ at index $4$ simply $L_4^\leq=L_3^\leq$,
		from $X\in C$ at index $5$ we get
		and $L_5^\leq = ((L_4^\leq)_w|X_{z,[l_1]})_{w,z,x}
		=((L_1 \otimes L_2|W)_w|X_{z,[l_1]})_{w,z,x}$
		so finally:
		\begin{equation*}
			P(Y|W=w,Z=z,X=x)=
				Y_{x,[l_2]} \circ L_5^\leq
			=
				Y_{x,[l_2]}\circ 
				((L_1 \otimes L_2|W)_w|X_{z,[l_1]})_{w,z,x}
			\txt.
		\end{equation*}
		Computing $\mu(\graph^c,L^c,\pi)$ for the c-component
		containing $W, X, Y$ from these atoms in the form
		$A^w \otimes A^x \otimes A^y$, through the definition
		of the $\otimes$-product (Def.\ \ref{def:kernel_product})
		contains an integral over $W$ (if we then marginalize $W$,
		this simply "integrates out" $W$)
		which corresponds to the well-known integrating out of the
		adjustment set $\{W\}$ for the effect $P(\randomVar{X},\randomVar{Y}|\PearlDo
		(\randomVar{Z}=z))$ (we study an IID example after all).
		
		In case the reader was wondering, why we bothered introducing
		graphical operations, they might take note that
		in Fig.\ \ref{fig:graphical_ops} it is almost immediately
		evident that the c-subgraph marginalizing $W$ from $\graph^c$
		has two c-components and can be split further,
		\ie we actually \emph{can} compute the much simpler
		$Y_{x,[l_2]} \circ L_2$.
		Indeed as pointed out above, we know the arguments of atoms,
		thus we may estimate $A^x = P(X|W=w,Z=z)$ and $A^y=P(Y|W=w,Z=z,X=x)$,
		then marginalizing $W$ means computing the integral
		$P'_z(X,Y)=\int P(W=w) P(X|W=w,Z=z)P(Y|W=w,Z=z,X=x)$,
		finally we know the argument of the (new) atoms,
		and $(A^y)'$ has argument $X$ (cf.\  Fig.\ \ref{fig:graphical_ops}),
		thus we can estimate $(A^y)' = P'_z(Y|X=x)$ for any $z$.
		In practice to use data optimally we will probably want to integrate $z$
		over the observational distribution $P(Z)$.
		Note that from Fig.\ \ref{fig:graphical_ops} it is
		also clear that we can compute $P(Y|\PearlDo(X=x))$
		by gluing a $\delta[x]\in\knownFunction$ (which,
		cf.\ proof of Lemma \ref{lemma:glue:apdx} and Rmk.\ \ref{rmk:compute_gluing},
		is simply a product
		with $A^y$), the product with $\delta[x]$ in turn is the same
		as plugging $x$ in for $X$, thus
		(if we are willing to write the disintegration as a quotient,
		which is common but may not be a good idea in general
		\citep{chang1997conditioning}):
		\begin{equation*}
			P(Y|\PearlDo(X=x))
			=
			\frac{
				\int P(W=w) P(X|W=w,Z=z)P(Y|W=w,Z=z,X=x) \dw
			}{
				\int P(W=w) P(X|W=w,Z=z) \dw
			}
		\end{equation*}
		This is of course the same result commonly obtained using
		do-calculus.
	\end{rmk}
	
	As illustrated in Fig.\ \ref{fig:graphical_ops},
	simplifications of c-components (and more generally
	repeated identification of c-components and their simplifications)
	is often relevant:
	
	\begin{definition}[c-Subgraphs]
		Given a structural model $(\graphStructural, L)$,
		and a c-component $(\graphStructural^c, L^c) \leq (\graphStructural, L)$,
		we call a simplification (Def.\ \ref{lemma:simplify})
		$(\graphStructural',L^c\setminus B)$ of $(\graphStructural^c, L^c)$
		a c-subgraph.
	\end{definition}
	\begin{cor}
		\label{cor:c_subgraphs}
		Given a structural model $(\graphStructural, L)$,
		and a c-subgraph $(\graphStructural', L') \leq (\graphStructural, L)$,
		then
		$\mu(\graphStructural',L',\pi')$ is a regular functional
		of $\mu(\graphStructural,L,\pi)$.
	\end{cor}
	\begin{proof}
		Apply Lemma \ref{lemma:properties_c_components}
		\ref{lemma:properties_c_components:regular},
		then Lemma \ref{lemma:simplify}.
	\end{proof}
	
	Given suitable overlap, two graphs (with know kernels)
	can be glued to a larger graph with regularly computable kernel.
	The next result is Lemma \ref{lemma:glue:mt}:
	
	\begin{lemma}[Gluing]
		\label{lemma:glue:apdx}
		Given a structural model $(\graphStructural, L)$, a $\graphStructural$-causal order $\pi$
		(with restrictions $\pi^A$, $\pi^B$ respectively)
		and subgraphs $\graphStructural^A, \graphStructural^B \leq \graphStructural$,
		such that
		each structural c-component $\graphStructural^c$
		is contained in at least one subgraph as
		$\graphStructural^c\leq \graphStructural^A$
		or $\graphStructural^c\leq \graphStructural^B$
		then there is a regular functional
		$F_\cup[\mu(\graphStructural^A,L\cap\nodes^A,\pi^A),
		\mu(\graphStructural^B,L\cap\nodes^B,\pi^B)]$ which we denote as
		$\graphStructural^A \cup_{(\graphStructural,L)} \graphStructural^B$ computing $\mu(\graphStructural,L,\pi)$,
		\begin{equation*}
			\mu(\graphStructural,L,\pi)
			\halfquad=\halfquad
			\graphStructural^A \cup_{(\graphStructural,L)} \graphStructural^B
			\txt.
		\end{equation*}
	\end{lemma}
	\begin{proof}		
		Let $n\in \nodesInner\setminus L$ be arbitrary.
		Let  $\graphStructural^c\leq\graphStructural$  be the
		structural c-component associated to the c-component of $n$
		(\ie the unique structural c-component with $n\in\nodesInner^c$).
		By Lemma \ref{lemma:properties_c_components}~%
		\ref{lemma:properties_c_components:atoms},
		$A^n(\graphStructural, L, \pi)=A^n(\graphStructural^c, L^c, \pi^c)$.
		By hypothesis $\graphStructural^c$ is contained in at least one of $\graphStructural^A$ or $\graphStructural^B$.
		If $\graphStructural^c\leq \graphStructural^A$,
		then Lemma \ref{lemma:properties_c_components}~%
		\ref{lemma:properties_c_components:atoms} applies to
		$\graphStructural^A$ and yields
		$A^n(\graphStructural^A, L^A, \pi^A)=A^n(\graphStructural^c, L^c, \pi^c)$.
		By Lemma \ref{lemma:atom_properties}~%
		\ref{lemma:atom_properties:atoms_regular},
		$A^n(\graphStructural, L, \pi)$ is thus a regular functional of
		$\mu(\graphStructural^A, L^A, \pi^A)$.
		If $\graphStructural^c \leq \graphStructural^B$,
		the same argument yields
		$A^n(\graphStructural, L, \pi)$ as a regular functional of
		$\mu(\graphStructural^B, L^B, \pi^B)$.
		
		Since $n\in \nodesInner\setminus L$ was arbitrary,
		we thus have regular functionals computing all the
		(finitely many) $A^n(\graphStructural, L, \pi)$ from
		$\mu(\graphStructural^A, L^A, \pi^A)$ or $\mu(\graphStructural^B, L^B, \pi^B)$.
		By Lemma \ref{lemma:atom_properties}~%
		\ref{lemma:atom_properties:regular_from_atoms},
		$\mu(\graphStructural,L,\pi)$ is a regular functional of these
		$A^n(\graphStructural, L, \pi)$.		
		Since composition of regular functionals is regular,
		thus $\mu(\graphStructural,L,\pi)$ is a regular functional of
		$\mu(\graphStructural^A,L\cap\nodes^A,\pi^A)$
		and $\mu(\graphStructural^B,L\cap\nodes^B,\pi^B)$.
	\end{proof}
	
	\begin{rmk}[Computation of Gluings]
		\label{rmk:compute_gluing}
		As pointed out for the example in Rmk.\ \ref{rmk:compute_atoms_practice},
		the gluing lemma is not only easy to proof once we
		have atoms of c-components Lemma \ref{lemma:properties_c_components}~%
		\ref{lemma:properties_c_components:atoms},
		but are also easy to compute given atoms:
		For each (non-trivial) structural c-component $\graph^c$,
		by hypothesis either $\graph^c\leq \graph^A$
		or $\graph^c\leq\graph^B$,
		compute atoms in $C= \nodesInner^c\setminus L^c$
		in this respective graph as described in Rmk.\ \ref{rmk:compute_atoms_practice}.
		Once we have all atoms (from all c-components),
		use Lemma \ref{lemma:atom_properties}~%
		\ref{lemma:atom_properties:product},
		\ie compute $\mu(\graphStructural, L, \pi)$ as a product over atoms.
		The product is essentially a Markov factorization of the joint distribution.
	\end{rmk}
	
	The revealing corollary (Cor.\ \ref{cor:revealing_simp}) follows immediately:
	\begin{cor}[Revealing]
		\label{cor:revealing:apdx}
		Given a structural kernel $\mu(\graphStructural, L)$,
		a structural c-component
		$\graphStructural^c \leq \graphStructural$
		and $L'\subset L^c$,
		then there is a regular functional computing
		$\mu(\graphStructural,(L\setminus L^c) \cup L')$
		from $\mu(\graphStructural,L)$
		and $\mu(\graphStructural^c,L')$.
	\end{cor}
	\begin{proof}
		Use Lemma \ref{lemma:properties_c_components}~%
		\ref{lemma:properties_c_components:regular}
		to extract all (other) c-components of
		$(\graphStructural, L)$,
		then glue (by repeated application of Lemma \ref{lemma:glue:apdx})
		$(\graphStructural, (L\setminus L^c) \cup L')$
		from these c-components and $(\graphStructural^c,L')$.
	\end{proof}
	The question as to whether such a revealing operation
	can be computed more generally is quite subtle and discussed in
	\ref{apdx:reveal}.

	\subsection{Algorithmic Construction}
	\label{apdx:structural_algos}
	
	We give a "simplest viable search"
	algorithm (Algo.\ \ref{algo:svs}) that constructs,
	by a regular computation,
	a structural kernel of interest
	by gluing elements from a set of known
	structural kernels.
	We do not know if this algorithm is complete
	(cf.\ also \ref{apdx:reveal}).
	But \emph{first} decomposing knowledge, then 
	gluing (only growing the result)
	makes the approach rather simple, while
	remaining complete in the IID case.
	This is a "search" in the sense that other than
	for the local operations of the last sub-section,
	we have to search a knowledge-set $\knowledgeSet$
	for pieces to glue.
	See also §\ref{apdx:algorithm_details}.	
		
	\begin{algorithm}
		\renewcommand{\thealgorithm}{SVSearch}
		\caption{\texttt{simplest\_viable\_search}} 
		\label{algo:svs}
		\textbf{Input:} A structural model $(\graphStructural, L)$
			to construct, 
			a set structured kernels $\knowledgeSet = \{ \mu(\graph^k, L^k, \pi^k) \}$,
			\\\hphantom{\textbf{Input:} }%
			optionally $(\graphStructural',L')\leq(\graphStructural,L)$
			and a regular functional $F'$ computing $\mu(\graphStructural', L', \pi')$
			\\\hphantom{\textbf{Input:} }%
			from finitely many elements of $\knowledgeSet$.
		\\
		\textbf{Output:}
		A (possibly empty) set of regular functionals
		$F$ computing $\mu(\graphStructural, L, \pi)$
		from finitely
		\\\hphantom{\textbf{Output:} }%
		many elements of $\knowledgeSet$.
		\begin{algorithmic}
			\State $R := \emptyset$.
			\For{$k\in \knowledgeSet$, $(\graphStructural^\cup, L^\cup)$
				\textbf{with} 
					$\big($\,\texttt{matching}$_{(\graphStructural^\cup, L^\cup)}$%
					($k$,$(G, L)$)
					\textbf{and}
					\texttt{gluable}$_{(\graphStructural^\cup, L^\cup)}$%
					($k$, $(\graph',L')$)\,$\big)$ }
				\State $R := R$ $\cup$
					\texttt{SVS}(\,$(G, L)$, $\knowledgeSet$,
						$(\graphStructural^\cup, L^\cup)$,
						$F^\cup[F', k]$\,).
					\Comment{Use Lemma \ref{lemma:glue:mt}.}
			\EndFor
			\If{$(G',L')$==$(G,L)$}
				\State\Return $\{F'\}$.
			\Else
				\State\Return $R$.
			\EndIf
		\end{algorithmic}
		\vspace{0.3em}\hrule\vspace{0.3em}
		\emph{Remark:} $\graphStructural^\cup$ with
			$\graphStructural, \graphStructural^k \leq \graphStructural^\cup$
			is used to describe where and how these graphs overlap.
			If no $(\graph',L'), F'$ are given, then
			\texttt{gluable} is true and $F^\cup[F,k] := \IdOp[k]$.
	\end{algorithm}

	The "opposite" of searching for glued structures,
	namely an algorithmic decomposition	into c-subgraphs, is in many ways much simpler
	(it is "local" in the sense of not having to search over a set $\knowledgeSet$
	of potential building-blocks).
	It is described in Algo.\ \ref{algo:decomp} and returns the set of
	smaller sub-graphs plus regular functionals computing their structured kernels.
	
	\begin{algorithm}
		\renewcommand{\thealgorithm}{CGDecomp}
		\caption{\texttt{csubgraph\_decomposition}}
		\label{algo:decomp}
		\textbf{Input:} A structured kernel $\mu=\mu(\graphStructural, L, \pi)$.\\
		\textbf{Output:} A set of pairs of a structural model
			$(\graph',L')$ and a regular functional
			$F'$ computing $\mu(\graph',L',\pi')=F'[\mu]$,
			such that the projection to the first component is
			a bijection to the set of c-subgraphs of $(\graphStructural, L)$.
		\begin{enumerate}[label=\arabic*)]
			\item
			\emph{C-Components:}\\
			Construct structural c-components (as equivalence-classes of $\approx_L$).\\
			\textbf{For each} structural c-component $(\graph^c,L^c)$,
			denote by $F^c$ the regular functional
			from Lemma \ref{lemma:properties_c_components}~%
			\ref{lemma:properties_c_components:regular}
			computing $\mu(\graph^c,L^c,\pi^c)=F^c[\mu]$ and
			perform the steps below.
			
			\item
			\emph{Latent-Sets:}\\
			Enumerate $L'\subsetneq\nodesInner$ with $L^c\subset L'$.\\
			\textbf{For each} $L'$, denote by $F^m$ the (regular) marginalization
				of $L'\setminus L^c$ computing $\mu(\graph^c,L',\pi^c)
				=F^m[\mu(\graph^c,L^c,\pi^c)]=F^m\circ F^c[\mu]$,
				and perform the steps below.
			\item 
			\emph{Simplify:}\\
			Construct the maximal $B$ for $L'$ as described in the proof
			of Lemma \ref{lemma:simplify}.\\
			Denote the $B$-simplification of $(\graph^c,L')$ by
			$(\graph',L'\setminus B)$, denote by $F^s$ the
			regular functional (actually $F^s=\IdOp$) from Lemma \ref{lemma:simplify},
			computing $\mu(\graph',L'\setminus B,\pi')=
			F^s[\mu(\graph^c,L',\pi^c)]=F^s\circ F^m \circ F^c[\mu]$.
			\item 
			\emph{Finalize or Iterate:}
			\\
			\textbf{If} $(\graph',L'\setminus B)$ is c-connected:
			\textbf{Yield} $((\graph',L'\setminus B), F^s\circ F^m \circ F^c)$.\\
			\textbf{Else}: \textbf{Repeat} from 1).
		\end{enumerate}
		\vspace{0.3em}\hrule\vspace{0.3em}
		\emph{Remark:} Yield adds the given expression to the output-set.
		For an example, see Fig.\ \ref{fig:graphical_ops}
		(ignoring the last gluing operation, which is not part of decomposition).
		Importantly, also the trivial, $L'=L^c$, simplification is in the output,
		so for example the first "large" c-component in Fig.\ \ref{fig:graphical_ops}
		(which cannot be glued from its simplifications!) is in the output.
	\end{algorithm}

	\subsection{Revealing Operations}
	\label{apdx:reveal}
	
	The gluing operation (Lemma \ref{lemma:glue:apdx}) was formulated
	enforcing the same missingness: nodes are latent in one
	subgraph, if and only if they are latent in the other one.
	To account for different missingness from different sub-graphs,
	we formulated a revealing operation (Cor.\ \ref{cor:revealing:apdx})
	-- this was however a simple combinations of other graphical operations.
	Intuitively it seems plausible that a more elementary (and more general)
	revealing operation, \eg on c-subgraphs (instead of c-components)
	or by combining gluing and revealing
	might be possible; a careful analysis of the situation suggests however
	that this might be an illusion
	(at least if one first decomposes into c-connected pieces).
	The difficulties are
	quite subtle, so we illustrate them on a simple example.
	
	\begin{example}[C-Subgraph Revealing]
		\label{example:revealing}
		Consider we are given information about structured kernels
		with the following structural graphs (mechanisms in square
		brackets are hidden), where we want to "reveal"\Slash{}imputate some
		of the hidden structure on the left graph from one on the right:\\
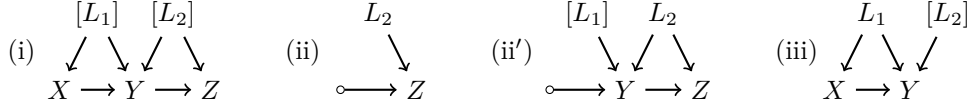
		\begin{minipage}{\textwidth}
			\centering
			\hfill~
			\begin{tikzpicture}
				\draw (-0.5,0.5) node {(i)};
				
				\draw (0,0) node (X){$X$};
				\draw (1,0) node (Y){$Y$};
				\draw (2,0) node (Z){$Z$};
				\draw (0.5,1) node (L1){$[L_1]$};
				\draw (1.5,1) node (L2){$[L_2]$};
				
				\draw [->,thick] (X) -- (Y);
				\draw [->,thick] (Y) -- (Z);
				
				\draw [->,thick] (L1) -- (X);
				\draw [->,thick] (L1) -- (Y);
							
				\draw [->,thick] (L2) -- (Y);
				\draw [->,thick] (L2) -- (Z);
			\end{tikzpicture}
			\hfill~
			\begin{tikzpicture}
				\draw (0.5,0.5) node {(ii)};
				
				\draw (1,0) node (Y)[circle, inner sep=0.1em, draw]{};
				\draw (2,0) node (Z){$Z$};
				\draw (1.5,1) node (L2){$L_2$};
				
				\draw [->,thick] (Y) -- (Z);		
				
				\draw [->,thick] (L2) -- (Z);
			\end{tikzpicture}
			\hfill~
			\begin{tikzpicture}
				\draw (-0.5,0.5) node {(ii$'$)};
				
				\draw (0,0) node (X)[circle, inner sep=0.1em, draw]{};
				\draw (0.5,1) node (L1){$[L_1]$};
				\draw (1,0) node (Y){$Y$};
				\draw (2,0) node (Z){$Z$};
				\draw (1.5,1) node (L2){$L_2$};
				
				\draw [->,thick] (X) -- (Y);	
				\draw [->,thick] (Y) -- (Z);		
				
				\draw [->,thick] (L1) -- (Y);
				
				\draw [->,thick] (L2) -- (Z);
				\draw [->,thick] (L2) -- (Y);
			\end{tikzpicture}
			\hfill~
			\begin{tikzpicture}
				\draw (-0.5,0.5) node {(iii)};
								
				\draw (0,0) node (X){$X$};
				\draw (1,0) node (Y){$Y$};
				\draw (0.5,1) node (L1){$L_1$};
				\draw (1.5,1) node (L2){$[L_2]$};
				\draw [->,thick] (X) -- (Y);
				
				\draw [->,thick] (L1) -- (X);
				\draw [->,thick] (L1) -- (Y);
				\draw [->,thick] (L2) -- (Y);
			\end{tikzpicture}
			\hfill~
		\end{minipage}
		\ie and we want to compute one of the following:\\
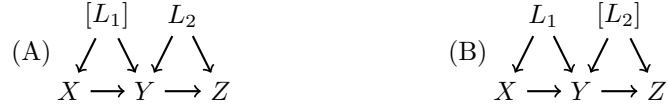
		\begin{minipage}{\textwidth}
			\centering
			\hfill~
			\begin{tikzpicture}
				\draw (-0.5,0.5) node {(A)};
				
				\draw (0,0) node (X){$X$};
				\draw (1,0) node (Y){$Y$};
				\draw (2,0) node (Z){$Z$};
				\draw (0.5,1) node (L1){$[L_1]$};
				\draw (1.5,1) node (L2){$L_2$};
				
				\draw [->,thick] (X) -- (Y);
				\draw [->,thick] (Y) -- (Z);
				
				\draw [->,thick] (L1) -- (X);
				\draw [->,thick] (L1) -- (Y);
				
				\draw [->,thick] (L2) -- (Y);
				\draw [->,thick] (L2) -- (Z);
			\end{tikzpicture}
			\hfill~
			\begin{tikzpicture}				
				\draw (-0.5,0.5) node {(B)};
				
				\draw (0,0) node (X){$X$};
				\draw (1,0) node (Y){$Y$};
				\draw (2,0) node (Z){$Z$};
				\draw (0.5,1) node (L1){$L_1$};
				\draw (1.5,1) node (L2){$[L_2]$};
				
				\draw [->,thick] (X) -- (Y);
				\draw [->,thick] (Y) -- (Z);
				
				\draw [->,thick] (L1) -- (X);
				\draw [->,thick] (L1) -- (Y);
				
				\draw [->,thick] (L2) -- (Y);
				\draw [->,thick] (L2) -- (Z);
			\end{tikzpicture}
			\hfill~
		\end{minipage}		
		First, we simply try to apply the gluing operation
		(Lemma \ref{lemma:glue:apdx})
		with (i) plus one of the other structures.
		This works trivially (\eg after marginalizing $L_1$ or $L_2$ in any of
		(ii,ii$'$, iii)	to restore the same missingness hypothesis again,
		but also only reproduces (i) then.
		We could first simplify (i) by Lemma \ref{lemma:simplify}:
		We can drop $\{Z\}$ (still enforcing $L_2$ as latent as above)
		we cannot drop $\{L_2,Z\}$ (the descendants of $L_2$ contain $Y$),
		we can drop $\{Y,Z,L_2\}$, but then it no longer contains
		the (left-hand side) structural c-component $\graph^c$ of (A) as subgraph
		(on $X,Y,[L_1]$, with an outer node instead of $L_2$; this is also
		a hypothesis of gluing),
		and neither of (ii,ii$'$) do either (for (ii)$'$ note
		that sub-graphs must not contain outer nodes of the containing graph
		as inner nodes by definition \ref{def:subgraphs_simp}~%
		\ref{def:subgraphs:nodesets}).
		Neither the simplification nor (iii) even contain $Z$, so the best
		we could possibly reproduce would be (iii) itself.
		
		Next we try to apply the gluing operation with (ii) or (ii)$'$ and (iii),
		but for (ii)$'$ and (iii), we can only glue (i) (we have to marginalize
		to restore same missingness), and for (ii) and (iii),
		we have to marginalize $L_2$ (can at most glue (B)),
		but neither (ii) nor (iii) contain the right-hand side structural
		c-component.
		
		It seems on first sight, that one could
		disintegrate \eg all observed nodes
		other than $L_1$ (similar for (ii), (ii)$'$)
		and directly "imputate" the result in (i),
		but note that in (i), $Z$ provides information about $L_2$,
		which restricts (via knowledge of $Y$) the plausible origin
		in the (multivariate) $L_1$--$L_2$--argument-plane of $Y$,
		thus it restricts $L_1$
		-- information that cannot possibly be extracted from (iii)
		where nothing is known about where in the argument plane of
		$Y$ (other than the "level-lines" of the preimage of $Y$)
		we are.
		Another indication of this problem is $L_1 \notindependent Z|Y$,
		because the path $L_1 \rightarrow Y \leftrightarrow Z$ is open;
		after realizing the kernel by Lemma \ref{lemma:obs_world_existence}
		the independence-statement makes sense, the question we are asking
		is however not really one of dependence vs.\ independence.
		
		Similarly, it seems on first sight
		plausible to obtain $B$ from (i) and (iii)
		by disintegrating the left-hand side ($X$ and $Y$),
		with the idea that using Lemma \ref{def:disint_product}
		the kernel $\mu_B$ of (B) can be reconstructed from
		$\marginalize{\mu_B}{Z}=\mu'''$ and a disintegration,
		however this disintegration has an argument $L_1$
		which for similar reasons as above cannot be
		reconstructed from (i).
	\end{example}
	
	The problem that we could not find a more general "revealing"
	operation does of cause not prove that there is none.
	However, the situation is clearly much more complicated than
	it may initially seem.
	If we inspect the proofs of the other operations,
	then they are rather simple consequences of the properties of
	"atoms", and the problem encountered here does not seem to
	be amenable to this approach directly either, at least not
	beyond the statement of Cor.\ \ref{cor:revealing:apdx}.

	\section{Details on Models and Observations}
	\label{apdx:models_obs_worlds}
	
	This section provides proofs for results (like the existence of
	model-realizations as "observable worlds") in the main text,
	plus some simple properties of observable and realized worlds
	that will be helpful for filling in the formal details
	when identifying model-properties from observations (§\ref{apdx:id_from_embeddings}).
	
	\begin{rmk}[Subspace-Relations]
		In Def.\ \ref{def:mechanism} and \ref{def:structured_kernel},
		we ask for certain subspace-relations between
		$\topSpace{X}_i$ for different $i\in I$.
		Formally, we should point out that
		notation \ref{notation:index_set_and_spaces}
		should include a suitable notion of intersections;
		for example we could fix a shared topological value space $\bar{\topSpace{X}}$
		and continuous inclusions
		$j_i : \topSpace{X}_i \hookrightarrow \bar{\topSpace{X}}$,
		then intersect in $\bar{\topSpace{X}}$.
		In practice it is usually clear what it means that
		the output of $X$ is in the domain of $Y_x$ so
		$X\otimes Y_x$ makes sense.
	\end{rmk}
	
	\subsection{Results in the Main Text}
	
	First note, that by acyclicity and finite past assumptions,
	causally ordered index sets are essentially $I=\mathbb{N}$ (or finite).
	
	\begin{rmk}[Simplified Form of Index-Sets]
		\label{rmk:index_sets_simplify}
		By acyclicity (Ass.\ \ref{ass:acyclic}), a causal order $\pi_I$
		exists, thus \asswlog $I \subset \mathbb{Z}$.
		By finite past (Ass.\ \ref{ass:finite_past}),
		a $\pi_I$-first element exists and actually \asswlog
		$I = \mathbb{N}$ or $I=\{1,\ldots, N\}$ (if $I$ is finite)
		as ordered set (obviously the countable $I$ is
		has this property as a set anyway).
	\end{rmk}
	
	The shallow distribution was defined as follows:	
	
	\begin{CopyLemma}{def:shallow_distr}{Shallow Distribution}
		Given acyclicity (Ass.\ \ref{ass:acyclic})
		and finite past (Ass.\ \ref{ass:finite_past}),
		then there exists a unique
		probability measure $P_\theta$ on
		$\prod_{i\in I} \topSpace{X}_i$
		with the following property:
		For all finite $I' \subset I$, the marginalization
		to $\topSpace{X}_{I'}$ satisfies
		for all tuples of measurable $B_i\in\sigmaAlgebraBorelIdx{X}{i}$ for $i\in I'$
		the equation
		\begin{equation*}
			\shallowDistr
			\Big(
			\prod_{i'\in I'} B_{i'}
			\Big)
			=
			\bigKernelCompound{
				\otimes^{\pi_I}_{i'\in \Anc_I(I')} f_{J(i')} 
			}
			\Big(
			\big(\prod_{i'\in I'} B_{i'}\big)
			\times				
			\big(\prod_{i\in \Anc_I(I')\setminus I'} \topSpace{X}_i \big)
			\Big)
			\txt.
		\end{equation*}
	\end{CopyLemma}
	
	\begin{proof}
		We use the notation from Rmk.\ \ref{rmk:index_sets_simplify}.
		First, we construct $\shallowDistr^n$ for $n\in\mathbb{N}=I$
		inductively over $n$
		as $\shallowDistr^n := \otimes_{i=1}^{n} f_{J(i)}$
		(cf.\ §\ref{apdx:sparsity} for technical details on extending domains
		from parents to ancestors).
		This is well-defined for any finite $n$ and satisfies
		the claimed equation by construction.
		
		It remains to show, that this is actually well-defined
		and characterizing,
		\ie that there is a unique probability-measure with this property.	
		This follows from the definition of products of $\sigma$-algebras,
		for a standard-result see \eg \citep[Prop.\ 2.2 (p.\,25)]{kallenberg1997foundations}.	
		Let $B \in \sigmaAlgebraBorelRaw{\prod_i \topSpace{X}_i}$ be arbitrary.
		By definition of the product sigma-algebra,
		a basis of measurable sets in
		$\sigmaAlgebraBorelRaw{\prod_i \topSpace{X}_i}$
		is given by "rectangles"
		$(\prod_{i\in I_B} B_i)
		\times (\prod_{i\in I\setminus I_B} \topSpace{X}_i)$
		with a finite $I_B\subset I$ and $B_i \in \sigmaAlgebraBorelIdx{X}{i}$.
		It is enough to fix $P_\theta$ on a basis,
		which the property does.
	\end{proof}

	We start with the existence of realizations by random elements.
	The non-trivial result at the core of the argument is provided by
	\citep[Lemma 2.22 (p.\,34)]{kallenberg1997foundations},
	a standard result in probability-theory.
	
	\begin{CopyDef}{def:obs_world}{Observable World}
		An observable world realizing the model $\mathcal{M}$ is
		a family of random variables
		$\{\anyVar_i: \Omega \rightarrow \topSpace{X}_i\}_{i\in I}$,
		together with measurable maps
		$f_i:\mathcal{X}_{\Pa_I(i)}\times[0,1] \rightarrow \topSpace{X}_i$
		and jointly independent random variables $\eta_i\sim U([0,1])$ called noises,
		such that
		\begin{equation*}
			\anyVar_i = f_i(\anyVar_{\Pa_I(i)},\eta_i)
			\quad\text{and}\quad
			P(\{\anyVar_i\}_{i\in I}) = \shallowDistr
			\txt.
		\end{equation*}
		Mandated by the second equality, we will often simply write
		$\shallowDistr(\{\anyVar_i\}_{i\in I})$.
	\end{CopyDef}
	
	\begin{CopyLemma}{lemma:obs_world_existence}{Observable World Existence}
		Given acyclicity (Ass.\ \ref{ass:acyclic})
		and finite past (Ass.\ \ref{ass:finite_past}),
		an observable world exists,
		it is unique up to equality in distribution
		and it is locally Markov, \ie
		$\anyVar_i \independent \anyVar_K | \Pa_I(i)$ for any
		$K\subset I$ with $K\cap\Dec_I(i)=\emptyset$.
	\end{CopyLemma}
	\begin{proof}
		\textbf{Existence:}
		We use notation from Rmk.\ \ref{rmk:index_sets_simplify}.
		
		Even though $I$ may be infinite,
		by definition of products of sigma-algebras,
		only finitely many factors of a measurable set
		in $\prod_{i\in I} \topSpace{X}_i$ are not the full space,
		thus in order to proof $P(\{\anyVar_i\}_{i\in I}) = \shallowDistr$
		it suffices to proof this statement for all finite subsets $I'\subset I$
		(see \eg \citep[Prop.\ 2.2 (p.\,25)]{kallenberg1997foundations}).
		By finite past (Ass.\ \ref{ass:finite_past})		
		we can replace $I'$ by the (also finite) $I'' := \Anc_I(I')$,
		so \asswlog $I' = \Anc_I(I')=\{1,\ldots, m\} \subset\mathbb{N}=I$.
		
		Our construction will be inductive over $n\in I$,
		\ie we show: $\forall n\in I$ $\exists \anyVar_1, \ldots, \anyVar_n$
		together with $f_i, \eta_i$ such that
		(i) $\forall i\leq n: \anyVar_i = f_i(\anyVar_{\Pa_I(i)},\eta_i)$
		and (ii)
		$P(\anyVar_1, \ldots, \anyVar_n) = \shallowDistr(\anyVar_1, \ldots, \anyVar_n)$.
		
		Inductive start ($n=0$): There is nothing to show;
		the case $n=0$ suffices as inductive start, see
		use of the inductive hypothesis in the inductive step.
					
		Inductive step ($n \mapsto n+1\geq 1$):
		We are given $\anyVar_1, \ldots, \anyVar_n$ with the required
		properties by inductive hypothesis.
		Since $\pi_I$ is a causal order, $\Pa_I(n+1) \subset \{1, \ldots, n\}$.
		If $n+1=1$, then $\Pa_I(n+1)=\emptyset$ and even though the
		inductive hypothesis -- in this case constructed as the inductive start
		-- is trivial and does not give us any random variables
		the construction below works.
		By \citep[Lemma 2.22 (p.\,34)]{kallenberg1997foundations},
		there is a measurable map
		$f_{n+1}:\topSpace{X}_{\Pa_I(n+1)} \times [0,1]
		\rightarrow \topSpace{X}_{n+1}$ and a random variable
		$\eta_{n+1} \sim \UniformDist([0,1])$ with
		$\eta_{n+1} \independent (\eta_1, \ldots, \eta_n)$,
		such that $\forall \pa_{n+1}\in\topSpace{X}_{\Pa_I(n+1)}$
		we have $f_{n+1}(\pa_{n+1}, \eta_{n+1})\sim f_{J(n+1)}(\pa_{n+1})$;
		we write $f_{J(i)}$ to unambiguously distinguish the
		mechanism at $i$ (Def.\ \ref{def:model}) from the
		mapping $f_i$ constructed in this proof.
		Define $\anyVar_{n+1} := f_{n+1}(\anyVar_{\Pa_I(n+1)}, \eta_{n+1})$.
		Property (i) is thus immediately satisfied.
		
		For (ii), by definition of $\shallowDistr$,
		$\shallowDistr(\anyVar_1, \ldots, \anyVar_{n+1})
		=\shallowDistr(\anyVar_1, \ldots, \anyVar_{n})
		\otimes (f_{J(n+1)})_{\Pa_I(n+1)}$.
		Since $f_{n+1}(\pa_{n+1}, \eta_{n+1})\sim f_{J(n+1)}(\pa_{n+1})$
		(see above) $P(\anyVar_{n+1}|\anyVar_{\Pa_I(n+1)}=\pa_{n+1})
		= f_{J(n+1)}(\pa_{n+1})$. 
		Thus for all measurable $B=B_{\leq n} \times B_{n+1}$
		\begin{align*}
			&P(\anyVar_1, \ldots, \anyVar_{n+1} \in B)
			\\
			&\quad=
			\int
			P(dv_{n+1}|\anyVar_{\Pa_I(n+1)}=\pa_{n+1})
			P(dv_1, \ldots, dv_n)
			1_B(v_1, \ldots, v_n, v_{n+1})
			d\vec{v}
			\\
			&\quad=
			\int
			f_{J(n+1)}(\pa_{n+1},dv_{n+1})
			\shallowDistr(dv_1, \ldots, dv_n)
			1_B(v_1, \ldots, v_n, v_{n+1})
			d\vec{v}
			\\
			&\quad=
			\bigKernelCompound{
			\shallowDistr(\anyVar_1, \ldots, \anyVar_{n})
			\otimes
			(f_{J(n+1)})_{\Pa_I(n+1)}
			}(B)
			\txt.
		\end{align*}
		Using that products form a basis of the
		$\sigma$-algebra $\sigmaAlgebraBorelRaw{\topSpace{X}_1
			\times\ldots\times\topSpace{X}_{n+1}}$, this shows claim (ii).
		
		Note that the model (Def.\ \ref{def:model}) by definition
		ensures $\topSpace{X}_{\Pa_J(j)}\subset\topSpace{X}_{\Pa_J}$
		and $\topSpace{X}_J \subset \topSpace{X}_j$,
		so $f_{n+1}$ is defined on all values the parents might take
		and produces values in $\topSpace{X}_{n+1}$.

		\textbf{Properties:}

		Uniqueness up to equality in distribution is clear by definition
		(as any observed world has distribution $\shallowDistr$).
		
		For the local Markov property, 
		first note that by
		construction
		$\anyVar_{n+1} := f_{n+1}(\anyVar_{\Pa_I(n+1)},$ $\eta_{n+1})$,
		so it inductively follows immediately that
		$\anyVar_{n+1}: \Omega\rightarrow\topSpace{X}_{n+1}$
		factors through $\eta_{\Anc_I(n+1)}$, \ie
		there is a measurable mapping $g_{n+1}: [0,1]^{|\Anc_I(n+1)|} \rightarrow \topSpace{X}_{n+1}$ with
		$\anyVar_{n+1}(\omega)$ $= g_{n+1}(\eta_{\Anc_I(n+1)})$.
		Also by 
		$\anyVar_{n+1} := f_{n+1}(\anyVar_{\Pa_I(n+1)}$, $\eta_{n+1})$
		we have $P(\anyVar_{n+1}|\Pa_I(n+1)=\pa_{n+1}) = f_{n+1}(\pa_{n+1}, \eta_{n+1})$.
		
		So for $K\subset I$ with
		$K\cap \Dec_I(n+1)=\emptyset$,
		$P(\anyVar_{n+1},\anyVar_K|\Pa_I(n+1)=\pa_{n+1})
		=P\big(f_{n+1}(\pa_{n+1}, \eta_{n+1}), g_K(\eta_{\Anc_I(K)})\big)$.
		Using $\eta_{n+1} \independent (\eta_1, \ldots, \eta_n)$
		(see above) repeatedly $\eta_{n+1}\independent \eta_{\Anc_I(K)}$;
		by $K\cap \Dec_I(n+1)=\emptyset$ we know $n+1\notin \Anc_I(K)$.
		Measurable transformations of independent variables are independent
		variables, thus the claim follows.
	\end{proof}
	
	\subsection{Properties of Probabilistic Worlds}
	\label{apdx:realized_world_properties}
	
	The following property is essentially a direct consequence of
	definitions, but it will be helpful for making
	the connection to identifiability in §\ref{apdx:id_from_embeddings}
	clearer. A similar result connecting structured queries to the
	realized world will be given in Lemma \ref{lemma:query_from_world}.
	
	\begin{lemma}[Properties]\label{lemma:obs_world_properties}
		Given an observable world $\{\anyVar_i\}_{i\in I}$
		and a model-aligned structural graph $\graphStructural$
		with $\nodesInner=\nodes \subset I$
		and $\edgesStructural = \{(p,c)\in\nodes\times\nodesInner|p\in\Pa_I(c)\}$.
		\begin{enumerate}[label=(\alph*)]
			\item\label{lemma:obs_world_properties:single_level}
				\emph{Single-Level:}
				If $\Anc_I(\nodes)\subset \nodes$
				and $\forall n\in\nodes:\mu^n = f_{J(n)}$,
				then $\forall L\subset \nodes$:			
				\begin{equation*}
					\shallowDistr(\anyVar_{\nodes\setminus L})
					\halfquad=\halfquad
					\mu(\graphStructural, L)
					\txt.
				\end{equation*}
			\item\label{lemma:obs_world_properties:multi_level}
				\emph{Multi-Level:}
				Given $N_0 \subset \nodes$
				with 
				$\Pa_I(\nodes\setminus N_0)\subset \nodes$ and
				$\forall n\in\nodes\setminus N_0:\mu^n = f_{J(n)}$
				and $\forall n\in N_0$: $\mu^n = \delta[x'_n]$,
				then
				there is a structural graph $\graphStructural'$
				with
				$\nodes \subset \nodes'$,
				$\forall n\in\nodes\setminus\Anc_I(N_0):$ $(\mu')^n = \mu^n$,
				$\forall n\in\nodes\setminus\Anc_I(N_0):$ $\Pa'(n) = \Pa_I(n)$,				
				such that
				\begin{equation*}
					\shallowDistr(\anyVar_{\nodes\setminus L}|
					\anyVar_{N_0}=x')
					\halfquad=\halfquad
					\mu(\graphStructural', L')[x']
					\txt,
				\end{equation*}
				where $\mu(\graphStructural', L')[x']$ is a functional
				of $x'$ through alignment  to $\mu^n = \delta[x'_n]$
				at $n\in N_0$.
				Further, given $Y \subset \nodes$
				with $Y\cap \Anc_I(N_0) = \emptyset$ and such that
				it satisfies the graphical property
				(a "latent path" here means a directed path with all nodes except
				for its last one in $L'$,
				\ie $\gamma\cap L'=\gamma\setminus\{b\}$):
				\\
				(i)	Given a latent path $\gamma$ in $\graphStructural'$
				starting at $a\in Y\setminus N_0'$ to a node $b\notin L'$,
				then $b\in Y\setminus L'$.
				
				Then
				$\graphStructural'$ can be chosen to also satisfy the graphical
				property (i') with $a\in N_0$ allowed:\\ 
				(i') Given a latent path $\gamma$ in $\graphStructural'$
				starting at $a\in Y$ to a node $b\notin L'$,
				then $b\in Y\setminus L'$.
				
				\emph{Remark:}
				This graphical property will make more sense after reading
				the proof of Thm.\ \ref{thm:extract_from_backdoor_complete}
				in §\ref{apdx:id_from_embeddings}.
				Essentially (i) appears (also named (i)) in step 1 of that proof,
				for step 3 (the multi-level case) it has to be modified to (i'),
				the statement above essentially says this can always be achieved
				by proper choice of $\graphStructural'$.
		\end{enumerate}
	\end{lemma}
	\begin{proof}
		Using the notation $I=\mathbb{N}$ from Rmk.\ \ref{rmk:index_sets_simplify},
		let $m := \max(\nodes)$.
		By definition
		\begin{equation*}
			\shallowDistr(\anyVar_1, \ldots, \anyVar_m)
			=
			\otimes_{n=1}^m
			f_m
			\txt.
		\end{equation*}
		\textbf{Part \ref{lemma:obs_world_properties:single_level}:}
		In the single-level case, by $\Anc_I(\nodes)\subset \nodes$ we can define
		$B=\{1,\ldots,m\}\setminus \nodes$ with $B=\Dec_{\graph}(B)$ and
		$\graphStructural$ is a simplification (Lemma \ref{lemma:simplify})
		of a similar (using the same edges as the $I$-graph, aligned to $f_i$)
		$\graphStructural^+$ with $\nodes^+ = \{1,\ldots,m\}$.
		By Lemma \ref{lemma:simplify} also
		$\mu(\graphStructural, L)=\mu(\graphStructural^+, L\cup B)$.
		By definition,
		\begin{equation*}
			\mu(\graphStructural^+, L\cup B)
			=
			\marginalize{\otimes_{n=1}^m f_m}{L\cup B}
			=
			\marginalize{\shallowDistr(\anyVar_1, \ldots, \anyVar_m)}{L\cup B}
			=
			\shallowDistr(\anyVar_{\nodes\setminus L})\txt.
		\end{equation*}
		
		\textbf{Part (b)}:
		For the multi-level case, first define $A:= \Anc_I(N_0)$.
		By exploiting sparsity (Lemma \ref{lemma:causal_transformations})
		we can \asswlog arrange nodes (replace $\pi$) such that $A=\{1, \ldots, m_0\}$.
		The conditional distribution
		$\shallowDistr(\anyVar_{\{1,\ldots,m\}\setminus A}|
		\anyVar_A=a)$
		has a regular version (by Ass.\ \ref{ass:standard_borel},
		\citep[Thm.\ 5.3 (p.\,84)]{kallenberg1997foundations})
		and satisfies
		$\shallowDistr(\anyVar_{N_0})
		\otimes\shallowDistr(\anyVar_{\{1,\ldots,m\}\setminus N_0}|
		\anyVar_{N_0}=x')=\shallowDistr(\anyVar_1,\ldots,\anyVar_m)$;
		so by uniqueness of disintegrations (Lemma \ref{def:disint_product}),
		this is the
		same as disintegrating
		$\shallowDistr(\anyVar_{\{1,\ldots,m\}}) = \otimes_{n=1}^m f_n$
		(Def.\ \ref{def:shallow_distr}) by $N_0$ as:
		\begin{equation*}
			\shallowDistr(\anyVar_{\{1,\ldots,m\}\setminus N_0}|
			\anyVar_{N_0}=x')
			\halfquad=\halfquad
			\disint_{N_0}(\otimes_{n=1}^m f_n)_{x'}
			\txt.
		\end{equation*}
		This disintegration can be computed (by uniqueness,
		Lemma \ref{def:disint_product})
		inductively over $a_1 < a_2 < \ldots \in N_0$ moving $a_1$ (then $a_2$ etc.)
		to the first (then second etc.) position in $\otimes_{n=1}^m f_n$
		by anti-causal disintegrations (Def.\ \ref{def:kernel_anti_causal_disint});
		this replaces terms to the left of $a_1$ producing
		\begin{equation*}
			\otimes_{n=1}^{m} f_n
			\halfquad=\halfquad	
			(f_{a_1} \circ \kernelCompound{\otimes_{n<a_1} f_n})
			\otimes (\otimes_{n<a_1} f_n | f_{a_1}) \otimes \ldots
			\txt.
		\end{equation*}
		Before we continue with $a_2$, we note that this only modifies
		terms $\leq a_1$, and $a_1\in N_0\subset A$, thus $a_1\leq m_0$.
		In particular, there are kernels $W^n$ (we can move $a_1$ one position
		at a time to avoid confusion) such that
		(up to transposition, Def.\ \ref{def:transposition}, denoted "$\approx$"
		following notation \ref{notation:uniqueness_transposition})
		\begin{equation*}
			\otimes_{n=1}^{m} f_n
			\halfquad\approx\halfquad
			\kernelCompound{\otimes_{n\in N_0} W^n}
			\otimes
			\kernelCompound{\otimes_{n\in A\setminus N_0} W^n}
			\otimes
			\kernelCompound{\otimes_{n>m_0} f_n}\txt.
		\end{equation*}
		The kernels $W^n$ may have additional arguments\Slash{}parents
		among their ancestors (Lemma \ref{lemma:causal_transformations}~%
		\ref{lemma:causal_transformations:disint_args}); call the new parent-sets $\Pa'(n)$.
		The disintegration of $N_0$ is now (by uniqueness)
		a simple disintegration from the left:
		\begin{equation*}
			\disint_{N_0}(\otimes_{n=1}^m f_n)
			\halfquad=\halfquad			
			\kernelCompound{\otimes_{n\in A\setminus N_0} W^n}
			\otimes
			\kernelCompound{\otimes_{n>m_0} f_n}\txt.
		\end{equation*}
		We can already define a graph $\tilde\graphStructural$
		on $\tilde\nodes=\tilde\nodes_{\nInner}=\{1, \ldots, m\}$, with edges 
		prescribed by $\Pa_I$ (for $n>m_0$) and $\Pa'$ (for $n\leq m_0$),
		and aligned to $W^n$ for $n\leq m_0$,
		to $\delta[x'_n]$ for $n\in N_0$ and
		to $f_n$ for $n\in\{n>m_0\} = \nodes\setminus A$.
		Using $\delta[X=x_0] \otimes \mu_x = \mu_{x=x_0}$
		we readily get the first part of our result
		(using again $B=\{1,\ldots,m\}\setminus (\nodes\cup A)$, as in the proof of part (a);
		note that $\Anc_I(\nodes\cup A)\subset (\nodes\cup A)$ by closedness
		under parents, cf.\ hypothesis):
		\begin{equation*}
			\mu(\tilde\graphStructural, L \cup B)_{n_0=x'}
			=
			\shallowDistr(\anyVar_{\{1,\ldots,m\}\setminus N_0}|
			\anyVar_{N_0}=x')\txt.
		\end{equation*}
		
		However, we cannot yet satisfy the graphical property:
		the nodes in $N_0$ now seem to \emph{confound}
		nodes in $\nodes\setminus A$ with nodes in $A$.
		(In the proof of Thm.\ \ref{thm:extract_from_backdoor_complete},
		we need a sub-graph of the $I$-graph that is a structural
		c-component, if $N_0$ [which will be pinned latents]
		confounds $Y$ [nodes of the subgraph] with other ancestors,
		then those ancestors are also in the c-component.
		The graphical property is a formal means of ensuring that
		ancestors are "disconnected" from the c-component $Y$ after the shuffling
		procedure above.)
		The solution is surprisingly simple. In this case we actually
		\emph{do} know something about the internal structure of
		the confounder. For \emph{singular} (taking only a single value)
		kernels there is no actual
		information in the value taken so they right-distribute:
		\begin{align*}
			\tag{$*$}
			\kernelCompoundConfounded{Z_x \otimes Y_{z,x}}_{[x]} \circ \delta[X=x_0]
			\halfquad&=\halfquad
			Z_{x=x_0} \otimes Y_{z,x=x_0}
			\\&=\halfquad
			(Z_x \circ \delta[X=x_0])
			\otimes
			(Y_{z,[x]} \circ \delta[X=x_0])_z
		\end{align*}
		\Ie, while generally, copying a node in a structural graph
		(adding a node with the same alignment and distributing children across both copies)
		changes the represented structured kernel,
		we \emph{can} copy nodes aligned to singular kernels.
		In particular, we can modify $\tilde\graphStructural$ into
		$\graphStructural'$ with additional nodes
		$\nodes_{\txt{copy}}$ (added to $L$) to ensure the graphical criterion
		as follows:
		Define
		\begin{equation*}
			\nodes_{\txt{copy}}
			\halfquad:=\halfquad
			\Pa_I(Y)\cap N_0
			\txt,
		\end{equation*}
		and $N_0' := N_0 \sqcup \nodes_{\txt{copy}}$ as the set adding a
		disjoint copy
		(we will write $a\in N_0$ and $a'\in \nodes_{\txt{copy}}$
		if the distinction is relevant).
		Define $\graphStructural'$ on nodes
		$\nodes' = \nodesInner' = \nodes \sqcup \nodes_{\txt{copy}}$
		and $L' = L \sqcup \nodes_{\txt{copy}}$,
		with edges as in $\tilde\graphStructural$ except for edges out of $N_0$:
		Here, keep edges out of $\Pa_I(Y)\cap N_0$ into $Y$ and
		out of $N_0 \setminus \Pa_I(Y)$ as is,
		edges out of $p\in\Pa_I(Y)\cap N_0$ to $w\in \nodes\setminus Y$
		are removed and instead added as $p' \rightarrow w$, \ie starting from
		$\nodes_{\txt{copy}}$.
		By equation $(*)$,
		we still have
		\begin{equation*}
			\mu(\graphStructural', L' \cup B)_{n_0=x'}
			=
			\shallowDistr(\anyVar_{\{1,\ldots,m\}\setminus N_0}|
			\anyVar_{N_0}=x')\txt.
		\end{equation*}
		Given $Y \subset \nodes$
		with $Y\cap \Anc_I(N_0) = \emptyset$ and such that
		(i) given a latent path $\gamma$ in $\graphStructural'$
		starting at $a\in Y\setminus N_0'$ to a node $b\notin L'$,
		then $b\in Y\setminus L'$.
		
		Let $\gamma$ be an arbitrary latent path in $\graphStructural'$
		starting at $a\in Y$ to a node $b\notin L'$.
		If $a \notin N_0 \sqcup \nodes_{\txt{copy}}$,
		then by (i), 
		$b\in Y\setminus L'$.
		
		If $a \in (N_0 \sqcup \nodes_{\txt{copy}})\cap Y = N_0$,
		then the child $c$ of $a$ which is the first node on $\gamma$
		after $a$ is, by construction of $\graphStructural'$, in $Y$
		(edges out of $N_0$ into $Y$ were kept in the construction of $\graphStructural'$,
		while edges out of $N_0$ into $\nodes\setminus Y$ were moved to $\nodes_{\txt{copy}}$).
		If $c$ is still in $N_0$, repeat the last argument, otherwise
		apply (i) to the restriction of $\gamma$ starting at $c$,
		to obtain $b \in Y\setminus L'$.
	\end{proof}

	\section{Details on the Extraction of Structured Kernels}
	\label{apdx:extraction_from_data}

	This section provides the technical details about
	the connection from observations to the identification of
	model properties. Before proving statements from the main
	text, we give definitions for the multi-level case in
	parallel to the main text.
	The reader not (currently) interested in the
	multi-level case may skip the first subsection,
	and only needs to know that in the single-level case
	"$\graph=\graphObs=\graphUnderlying$", $\nodesInner=\nodesInnerProper$,
	$\nodesOuter=\nodesOuterProper$ (while $\nodesOuterFixed=\emptyset$).
	The phrasing "$\mu$ is $\delta_i$-identifiable" is equivalent to
	"$\mu$ is identifiable" in the single-level case.

	\subsection{Multi-Level Modifications}

	We need some smaller modifications to families of embeddings and
	to data-sets\Slash{}identifiability-definitions for the multi-level case.

	\paragraph{Families of Embeddings:}

	We formulate the machinery of embeddings on local graphs,
	which are closely related to structural graphs (and structured kernels);
	in the single-level case, $\nodesOuterFixed=\emptyset$
	and $\graph$, $\graphUnderlying$, $\graphObs$
	(see below) agree.
	Starting from a slightly different kind of object
	allows to keep both aspects of the formalism
	reasonably simple.
	
	\begin{definition}[Local Graph]
		\label{def:local_graph_new}
		A local graph $\graph$ is a finite set of nodes
		which is a disjoint union
		of proper, \textbf{pinned (new in multi-level)} and external nodes
		\begin{equation*}
				\nodes
				\halfquad=\halfquad
				\nodesInnerProper
				\halfquad\dot\cup\halfquad
				\nodesOuterFixed
				\halfquad\dot\cup\halfquad
				\nodesOuterProper
			\end{equation*}
		together with a parent structure
		$\mathcal{P}\subset\nodes\times\nodesInnerProper$,
		a set of directed edges (which we will call proper edges)
		ending at a node in $\nodesInnerProper$.
		
		\emph{Alignment:}
		There is an underlying structural graph
		$\graphStructural^0(\graph)$,
		with $\nodes^0:=\nodes$, $\nodesInner^0 := \nodesInnerProper$
		and $\edgesStructural^0 := \parentStructure$.		
		A model-aligned local graph is a local graph together with
		indices $i^x\in I$ for $x\in\nodesOuterFixed$ \textbf{(new in multi-level)} and
		kernels $\mu^n$ plus argument-assignments $\pa^n$
		for $n\in\nodesInnerProper$ such that $\graphStructural^0$
		is model-aligned.
	\end{definition}
	
	The embeddings from Def.\ \ref{def:local_graph_embedding} are slightly modified,
	the only difference is that model-alignment now
	includes additionally the statement "$\forall x\in\nodesOuterFixed$: $i^x=\psi(n)$."
	\begin{definition}[Local Graph Embedding]\label{def:local_graph_embedding:apdx}
		An embedding of
		a model aligned local graph $\graph$, denoted by
		$\psi:\graph \hookrightarrow I$, is an injective mapping
		$\psi:\nodes \hookrightarrow I$
		such that the following conditions are satisfied
		(see Fig.\ \ref{fig:embeddings}):
		\begin{enumerate}[label=(\roman*)]
			\item\label{def:local_graph_embedding:inner_parents_incl:apdx}
			Proper-Node $I$-Parents are Included:
			If $n\in \nodesInnerProper$, then
			$\forall i' \in \Pa_I(\psi(n))$, $\exists n' \in \nodes$
			such that $i' = \psi(n')$.
			\item\label{def:local_graph_embedding:proper_edges_a:apdx}
			Proper Edges equal $I$-Graph Edges:
			For $n' \in \nodesInnerProper$,
			there is a proper edge $n \rightarrow n'$
			in $\graph$ if and only if
			$\psi(n) \in \Pa_I(\psi(n'))$.
			\item\label{def:local_graph_embedding:applicable:apdx}
			Model Alignment:
			$\forall n\in\nodesInnerProper$:
			$\mu^n = f_{J(\psi(n))}$ is a model mechanism (Def.\ \ref{def:model})
			and \textbf{new in multi-level:}
			$\forall x\in\nodesOuterFixed$:
			$i^x=\psi(n)$.
		\end{enumerate}
	\end{definition}

	Families of embeddings can be defined
	as before (Def.\ \ref{def:local_graph_embedding_family}),
	only the notation for $\nodesInner \mapsto \nodesInnerProper$
	changes:
	
	\begin{definition}[Families of Embeddings]
		\label{def:local_graph_embedding_family:apdx}
		A family of local graph embeddings
		\begin{equation*}
			(J_0, \{\psi_j\}_{j\in J_0}, \graph, H, y_0)
		\end{equation*}
		is a collection of local graph embeddings
		$\psi_j: \graph\hookrightarrow I$
		of a single model-aligned local graph $\graph$.
		Further there is a fixed element $y_0 \in \nodes$,
		which will be called the anchor\allowbreak{}\mbox{(-}node),
		a symmetry $H\subset G$ and a range of applicability $J_0\subset I$.
		For $n\in\nodes$ denote:
		\begin{equation*}
			\psi_*(n): J_0 \rightarrow I,
			\halfquad
			j\mapsto \psi_j(n)
			\txt.
		\end{equation*}
		Finally, we will
		require the following conditions to be satisfied:
		\begin{enumerate}[label=(\Roman*)]
			\item\label{def:local_graph_embedding_family:trivial_anchor:apdx}
			Trivial on Anchor:
			$\forall j\in J_0$: $\psi_j(y_0) = j$.
			\item\label{def:local_graph_embedding_family:rigidity:apdx}
			Rigidity:
			$\forall n \in \nodes$:
			$\psi_*(n)$
			is $H$-equivariant.
			\item\label{def:local_graph_embedding_family:freeness:apdx}
			Freeness:
			$\forall n\in\nodesInnerProper$ \textbf{(multi-level: inner
			nodes become proper nodes)}:
			$\psi_*(n)$
			is injective.
		\end{enumerate}
	\end{definition}	
	
	While our models allow only a fixed and finite number of parents
	in the $I$-graph, nodes can have an arbitrary number of children.
	This allows for example for multi-level statics.
	
	\begin{definition}[Singular Nodes]\label{def:pinned_nodes}
		\textbf{This has no analogue in the single-level case:}
		Given a family $\{\psi_j\}_{j\in J_0}$ of
		local graph embeddings, we call a node $x\in \nodes$
		singular embedded at $i\in I$ if $\forall j\in J_0$: $\psi_j(x) = i$.
		\\
		\emph{Remark:}
		For $|J_0|>1$, proper nodes cannot be singular by freeness.
		For pinned alignment
		Def.\ \ref{def:local_graph_embedding}%
		\ref{def:local_graph_embedding:applicable},
		all pinned nodes $\nodesOuterFixed$
		must be singular embedded.
		Conversely singular embedded nodes $x\in\nodesOuterProper$ can always be
		pinned at $i^x=\psi_j(x)$, \ie there is a valid family of embeddings
		with $x\in\nodesOuterFixed$ (this can render the family
		backdoor-free, cf.\ \ref{def:backdoor_free})
	\end{definition}

	For attached structure (decorated graphs) not much new happens
	for the multi-level case. Some general remarks on attached structure
	can be found in the next subsection.
	Some smaller remarks are in place however and
	the definition of ancestral structure should use the
	nomenclature of external and proper (rather than inner and outer) nodes:
	\begin{definition}[Ancestral Structure]
		\label{def:ancestral_structure:apdx}
		Given a local graph $\graph$,
		an ancestral structure $\ancestralStructure$ on $\graph$
		is a set of directed edges $\rightsquigarrow$ (ancestral edges)
		each starting at a \textbf{proper} node and ending at an
		\textbf{external} node,
		such that $\graph$ with proper and ancestral edges is acyclic.
	\end{definition}
	\begin{rmk}[Validity of Ancestral Sets in the Multi-Level Case]
		In the multi-level case,
		a valid ancestral structure remains essentially unchanged
		(compared to the single-level case), so it
		tracks only ancestral paths
		to $\nodesOuterProper$, not to nodes $\nodesOuterFixed$.
		This will suffice for asymptotic results as
		there are only finitely many $j\in J_0$
		with ancestral paths to $\nodesOuterFixed$ (see proof
		of Thm.\ \ref{thm:extract_from_backdoor_complete}, step 3,
		choice of $J_0^{\txt{invalid}}$).
	\end{rmk}
	
	Finally, also in backdoor-freeness (Def.\ \ref{def:backdoor_free}), replace
	$\nodesOuter\mapsto\nodesOuterFixed$.
	This means $\nodesOuterFixed \cap L \neq \emptyset$ is explicitly allowed,
	and is indeed typically (excluding trivial cases) the case.
	
	\begin{definition}[Backdoor-Freeness]\label{def:backdoor_free:apdx}
		Given a decorated family of embeddings $(\{\psi_j\}_{j\in J_0}$,
		$L$, $\ancestralStructure)$, we call an ancestral edge
		$l \rightsquigarrow x \in \ancestralStructure$
		a backdoor if it starts at $l\in L$.
		We call the family backdoor-free, if
		$L\cap\nodesOuterProper=\emptyset$
		\textbf{(multi-level: $\nodesOuter$ is replaced by $\nodesOuterProper$)}
		and there are no backdoors.
	\end{definition}

	\paragraph{Data-Sets:}
	
	We slightly extend the definition of data-sets (Def.\ \ref{def:data_set}),
	note that a data-set in the sense of the main text
	is a dataset with $\randomVar{X}'=\emptyset$.
	Similar to backdoor-freeness, "pinned nodes" (here $\randomVar{X}'$)
	are typically \emph{not} observed; they are however $j$-independent.
	
	\begin{definition}[Data-Set]\label{def:data_set:apdx}
		A data-set is a tuple
		$\mathcal{D} = ((\randomVar{X}_j$, $\randomVar{Y}_j)_{j\in J_0}$,
				$\randomVar{X}')$,
		where $\randomVar{X}_j$, $\randomVar{Y}_j$ and
			$\randomVar{X}'$ \textbf{(multi-level: the inclusion of a
				$\randomVar{X}'$ is new)}
		are tuples
		of $\anyVar_i$, \ie there exist $i_{X,j}^{(1)}, \ldots,
		i_{X,j}^{(n)} \in I$,
		such that $\randomVar{X}_j = (\anyVar_{i_{X,j}^{(1)}},
		\ldots, \anyVar_{i_{X,j}^{(n)}})$ and analogously for $\randomVar{Y}_j$
			and $\randomVar{X}'$.
		
		We call the data-set valid if $X$ and $Y$
		are observed that is
		\begin{equation*}
			J_0^{\txt{obs}}(N)
			\halfquad:=\halfquad
			\{\halfquad
			j\in J_0
			\halfquad|\halfquad
			\forall m:
			i_{X,j}^{(m)} \in \viewport(N),
			\forall m':
			i_{Y,j}^{(m')} \in \viewport(N)
			\halfquad\}
		\end{equation*}
		satisfies $|J_0^{\txt{obs}}(N)|\rightarrow\infty$ as $N\rightarrow\infty$
		and is non-degenerate
		$j\neq j' \Rightarrow \forall m: i_{Y,j}^{(m)}\neq i_{Y,j'}^{(m)}$
		(this last condition is imposed only on $Y$, not on $X$;
		see Rmk.\ \ref{rmk:kernels_uniqueness}).
	\end{definition}

	As before decorated families of embeddings
	induce data-sets (note that even though $\randomVar{X}'$ is not
	observed, its $j$-independent existence is an important condition).
	
	\begin{example}[Data-Set of Decorated Families]
		\label{example:data_set_of_embedding:apdx}
		Given a decorated family of embeddings
		$(\{\psi_j\}_{j\in J_0},L, \ancestralStructure)$ with
		$L\cap\nodesOuterProper=\emptyset$,
		there is a valid data-set
		$((\randomVar{X}_j, \randomVar{Y}_j)_{j\in J_0}$
				, $\randomVar{X'})$
		(Def.\ \ref{def:data_set})
		defined for $j\in J_0$ as the tuples
		\begin{align*}
			\randomVar{X}_j
				\halfquad&=\halfquad
				(\anyVar_{\psi_j(n)})_{n\in\nodesOuterProper
								\cup (\nodesOuterFixed\setminus L)}
				\halfquad\txt{\textbf{(multi-level: add obs.\ pinned $n$)}}
				\\
			\randomVar{Y}_j
			\halfquad&=\halfquad
			(\anyVar_{\psi_j(n)})_{n\in\nodesInnerProper\setminus L}\\
						\randomVar{X}'
						\halfquad&=\halfquad
						(\anyVar_{i^n})_{n\in\nodesOuterFixed\cap L}
			\halfquad\txt{\textbf{(multi-level only)}}
			\txt.
		\end{align*}
	\end{example}

	\paragraph{Identifiability:}
	
	Direct identifiability (Def.\ \ref{def:identification_direct_mt})
	must respect the shared across $j$ value of $\randomVar{X}'$:
	
	\begin{definition}[Direct Identifiability]\label{def:identification_direct:apdx}
		Given a valid data-set (Def.\ \ref{def:data_set})
				$\mathcal{D} = ((\randomVar{X}_j, \randomVar{Y}_j)_{j\in J_0},$
				$\randomVar{X}')$
		and kernels $\{X_j\}_{j\in J_0}$, $Y_{x}[x']$ (a functional of $x'$),
		such that for each $j$ individually 
		the shallow distribution $\shallowDistr$
		(Def.\ \ref{def:obs_world}) satisfies,
		\begin{equation*}
				\forall j\in J_0:\halfquad
				\shallowDistr(\randomVar{X}_j, \randomVar{Y}_j|\randomVar{X}'=x')
				= X_j \otimes Y_{x}[x']
				\halfquad\txt{\textbf{(multi-level: $x'$ is new)}}
				\txt,
			\end{equation*}
		then we call $Y_{x}[x']$ directly identifiable. 
		Formally we also consider the
		known a priori $\knownFunction$ (Def.\ \ref{def:model})
		and the immediately observed $\anyVar_i$ for $i\in\viewportEventual$
		(thus also $\delta[\anyVar_i=x'_i]$; \textbf{multi-level only}) directly identifiable.
		
		\emph{Remark:}
		We usually will not know or observe the value $x'$,
		we only know that $x'$ does \emph{not} depend on $j$.
	\end{definition}
	
	For the multi-level case, we use the following notation.
	
	\begin{notation}[Pinned Kernels]
		\label{def:delta_i_identified}
		Tracking in computations which arguments of a structured kernel
		correspond to pinned nodes would be arduous.
		Instead we adjoin generic kernels
		$\pinnedKernels = \{ \delta_i \}_{i\in I}$
		with zero parents
		(similar to how, for example, polynomial rings adjoin a
		generic variable to an algebraic structure)
		and consider
		$\structureKernels=\modelKernels\cup\knownFunction\cup\pinnedKernels$-%
		structured kernels.
		We will then call an expression $\mu[\delta_{i^1},\ldots,\delta_{i^m}]$,
		$\delta_i$-identifiable, if there is an identifiable kernel $\nu$
		such that (writing $\delta[\anyVar_i]$ for the random measure
		singular at the value taken by $\anyVar_i$):
		\begin{equation*}
			\mu[\delta[\anyVar_{i^1}], \ldots, \delta[\anyVar_{i^m}]]
			\overset{\txt{a.\,s.}}{=} \nu
			\halfquad\txt{(almost surely w.\,r.\,t.\ $\realizedDistr$).}
		\end{equation*}
	\end{notation}
	\begin{example}[Learning Pinned Kernels]
		\label{example:learn_pinned}
		The intuitive idea underlying notation \ref{def:delta_i_identified}
		is that, in each realization,
		hidden contexts (like slope in example \ref{example:multi_level_illustrate})
		do take a value $C=c$. We do not know this value,
		but given observations of many children we gain partial information:
		\Eg if $Y_c = \mathcal{N}(c^2, 1)$ (in SCM notation this could be
		$Y:= C^2 + \eta_Y$), from many observations of $Y$ we gain information
		about $c^2$.		
		So from observations, we know $\realizedDistr(C \in \{\pm\sqrt{E[Y]}\})=1$.
		Indeed, we gain precisely the information about $C$
		needed to compute $Y$;
		we can learn $Y_c$, up to sign of $c$,
		as $Y'=\mathcal{N}(c^2, 1)$.
		Taken together:
		$\realizedDistr(Y_C=Y')=1$,
		where the capital index $Y_C$ means we plug in
		the random-variable $C$ (or equivalently:
		compose $Y_c$ to $Y_C = Y_{c} \circ\delta[C]$).
	\end{example}
	
		We have to associate structural graphs to local graphs
		(the embedded objects) for identification statements.
		
		\begin{definition}[Associated Structural Graphs]
			\label{def:associated_structural_graph}
			We are given
			a model-aligned local graph $(\graph, L)$,
			with $L\cap \nodesOuterProper = \emptyset$.
			Def.\ \ref{def:local_graph_new} introduced
			its underlying structural model $(\graphUnderlying, L)$
			with $\nodes^\gUnderlying:=\nodes$,
			$\nodesInner^\gUnderlying := \nodesInnerProper$
			and $\edgesStructural^\gUnderlying := \parentStructure$.
			
			We additionally define		
			$\graphObs(\graph)$
			with nodes $\nodes^\gObs = \nodes$ split as
			\begin{equation*}
					\nodesInner^\gObs
					:= \nodesInnerProper \cup (\nodesOuterFixed \cap L)
					\halfquad\txt,\qquad
					\nodesOuter^\gObs
					:= \nodesOuterProper \cup (\nodesOuterFixed \setminus L)
				\end{equation*}
			and the edges $\edgesStructural^\gObs = \parentStructure$.
			Then $\graphObs$ is model-aligned with $\mu^n$ for $n\in \nodesInnerProper$
			and $\mu^x:=\delta_{i^x}\in\pinnedKernels$ for $x\in \nodesOuterFixed\cap L$.
		\end{definition}
		
		Thm.\ \ref{thm:extract_from_backdoor_complete} will be formulated
		as $\delta_i$-identifying (notation \ref{def:delta_i_identified})
		$\mu(\graphObs,L)$, the formal details are given below
		at \ref{thm:extract_multi_level}.
	
		\paragraph{Knowledge Sets:}
		
		We modify Def.\ \ref{def:knowledge_set_extracted} to include
		immediately observed variables as (basic) knowledge:

		\begin{definition}[Extracted Knowledge Set]
			\label{def:knowledge_set_extracted:apdx}
			Let $\mathcal{B}$ be the set of all backdoor-free
			families of embeddings.	
			For $(\psi,L,\ancestralStructure)\in\mathcal{B}$
			define a structured model
			$k(\psi,L,\ancestralStructure) := (\graphObs, L)$.
			By Thm.\ \ref{thm:extract_from_backdoor_complete}
			$\mu(\graphObs,L)$ is $\delta_i$-identifiable
			(we simply say $k$ is identifiable).
			
			For $\tilde{f} \in \knownFunction$, define
			a structured model $k(\tilde{f})=(\graphStructural(\tilde{f}),L=\emptyset)$,
			where $\graphStructural(\tilde{f})$ is a graph with a single inner node $y$
			aligned to $\mu^y=\tilde{f}$, $\kappa$ (the number of parents) outer nodes
			$\nodesOuter = \{x_1, \ldots, x_\kappa\}$
			and a proper edge from each $x_k \rightarrow y$.
			By definition (Def.\ \ref{def:model} and \ref{def:identification_direct_mt}),
			elements of $\knownFunction$ are considered known,
			thus $k(\tilde{f})$ is identifiable for all $\tilde{f} \in \knownFunction$.
			\textbf{Multi-Level only:}
			Finally, what is immediately observed is known
			(Def.\ \ref{def:identification_direct_mt}), so for $i\in \viewportEventual$,
			define $k(i)$ as the structured graph with a single inner node
			aligned to $\delta[\anyVar_i]$.

			Define the, thus (element-wise)
			$\delta_i$-identifiable \textbf{($\delta_i$ added in multi-level)},
			basic knowledge set as
			\textbf{(multi-level: add last union)}
			\begin{equation*}
					\knowledgeSet_{\txt{basic}}
					\halfquad=\halfquad
					\bigcup_{(\psi,L,\ancestralStructure)\in\mathcal{B}}
					k(\psi,L,\ancestralStructure)
					\halfquad\cup\halfquad
					\bigcup_{\tilde{f}\in\knownFunction}
					k(\tilde{f})
					\halfquad\cup\halfquad
					\bigcup_{i\in\viewportEventual}
					k(i)
					\txt.
				\end{equation*}
		\end{definition}
	
	\subsection{Embedding-Families and Attached Structure}
	
	We add some remarks and simple properties,
	beyond what was stated in the main text,
	that are useful for later reference
	and for the practical interpretation of these definitions.
	
	\begin{rmk}[Change of Anchor]\label{rmk:change_of_anchor}
		Given a family of embeddings
		$\{\psi_j\}_{j\in J_0}$ anchored on $x\in\nodes$
		and a fixed choice $y\in\nodes$ with $\psi_*(y)$ injective
		(\eg
		any proper node by freeness
		\ref{def:local_graph_embedding_family}~%
		\ref{def:local_graph_embedding_family:freeness}).
		Define $J_0' := \img(\psi_*(y)) \subset I$, clearly when
		restricting the target to the image
		$\psi_*(y):J_0 \twoheadrightarrow J_0'$ is also surjective, thus
		$\psi_*(y)$ is bijective (to its image $J_0'$);
		we will write $\psi_y^{-1}(j') := (\psi_*(y))^{-1}(j')$ for its inverse
		at $j'\in J_0'$.
		Define a family of local graph embeddings anchored on $y$
		indexed by $J_0'$ as
		\begin{equation*}
			\psi'_j : \nodes \hookrightarrow I,
			\quad
			n \mapsto \psi_{\psi_y^{-1}(j)}(n)\txt.
		\end{equation*}
		Note that given an anchor (here $x\in\nodes$),
		even though freeness \ref{def:local_graph_embedding_family:freeness}
		is given only at proper nodes, the mapping $\psi_*(x)$ at the original
		anchor is of course also injective (by
		\ref{def:local_graph_embedding_family:trivial_anchor}).
	\end{rmk}
	
	\begin{rmk}[Relevant Latent Subsets]\label{rmk:relevant_latent}
		In practice, for the finite sample case,
		the definition of $\minLatentSets_\psi$ can easily
		be adapted as follows:
		A $L\subset\nodes$ is \emph{relevant}, if
		there is no $L' \subsetneq L$ such that (at least) the same
		number of data-points
		is available $|J_0^{\txt{obs}}(N)| \leq |(J_0^{\txt{obs}})'(N)|$
		(cf.\ Def.\ \ref{def:observed})
		and the number $|J_0^{\txt{obs}}(N)|$ of available data-points
		is larger some (possibly $\psi$-dependent;
		estimator and precision-target specific) minimum number
		(if precision-measures are available for a specific estimator,
		those can of course used directly instead of a minimum number of
		samples, this may involve expensive evaluations on
		many potentially relevant branches however).
		It may however make sense to track not only minimal latent subsets,
		as trading additional latents for additional samples might be favorable.
	\end{rmk}
	
	It will be helpful to have a notion of sub-graphs and sub-families
	(in loose analogy to Def.\ \ref{def:subgraphs_simp}).
	\begin{definition}[Local Subgraphs]
		\label{def:subgraphs_local}
		A local subgraph $\graph^A\leq\graph^B$
		of a local graph $\graph^B$
		is a local graph $\graph^A$ such that:
		\begin{enumerate}[label=(\roman*)]
			\item\label{def:subgraphs_local:nodesets}
			\emph{Node Sets:}
			$\nodes^A\subset\nodes^B$,
			$\nodesInnerProper^A\subset\nodesInnerProper^B$
			and $\nodesOuterFixed^A \subset \nodesOuterFixed^B$.
			\item\label{def:subgraphs_local:edges}
			\emph{Edge Sets:}
			$\parentStructure^B \cap (\nodes^A\times\nodesInnerProper^A)$,
			\ie edges are exactly those in $\graphStructural^B$
			from nodes in $\graphStructural^A$
			to inner nodes of $\graphStructural^A$.
			\item\label{def:subgraphs_local:inner_parents}
			\emph{Proper Parents:}
			$\forall n\in\nodesInnerProper^A$:
			if $p\in\Pa_{\graphStructural^B}(n)$,
			then $p\in\nodes^A$.
			\item\label{def:subgraphs_local:alignment}
			\emph{Alignment:} For model-aligned graphs
			$\forall n\in\nodesInnerProper^A$
			(inner nodes of $\graphUnderlying(\graph^A)$):
			$(\mu^A)^n = (\mu^B)^n$.
		\end{enumerate}
		A sub-family of embeddings $\{\psi^A_j\}_{J_0^A}\leq \{\psi^B_j\}_{J_0^B}$
		of a family of embeddings $\{\psi^B_j\}_{J_0^B}$
		is a family of embeddings $\{\psi^A_j\}_{J_0^A}$,
		such that
		\begin{enumerate}[label=(\roman*)]
			\item
				\emph{Graphs:}
				$\graph^A \leq \graph^B$.
			\item
				\emph{Applicability:}
				$\{\psi^A_j\}_{J_0^A}$ and $\{\psi^B_j\}_{J_0^B}$
				have the same anchor $n_0\in\nodes^A$,
				$J_0^B \subset J_0^A$ (note that larger graphs
				have \emph{worse} symmetry, cf.\ Rmk.\ \ref{rmk:families_and_symmetries})
				and $\forall j\in J_0^B: \psi_j^A = \psi_j^B|_{\nodes^A}$.
		\end{enumerate}		
		A decorated sub-family $(\{\psi^A_j\}_{J_0^A}, L^A, \ancestralStructure^A)
		\leq (\{\psi^B_j\}_{J_0^B}, L^B, \ancestralStructure^B)$
		of a decorated family $(\{\psi^B_j\}_{J_0^B},$ $L^B,$ $\ancestralStructure^B)$
		is a decorated family $(\{\psi^A_j\}_{J_0^A},$ $L^A,$ $\ancestralStructure^A)$,
		such that
		\begin{enumerate}[label=(\Roman*)]
			\item
			\emph{Underlying Families:}
			$\{\psi^A_j\}_{J_0^A}\leq \{\psi^B_j\}_{J_0^B}$.
			\item
			\emph{Latents:}
			$L^A \subset L^B \cap \nodes^A$.
			\item 
			\emph{Ancestral Structure:}
			$\ancestralStructure^A \subset
			P_{B\setminus A}$, where $P_{B\setminus A}$ contains
			exactly those pairs $(y,x)\in\nodesInnerProper^A \times \nodesOuterProper^A$,
			such that there is directed path in $\graph^B$ exclusively along edges
			in $\ancestralStructure^B\cup(\parentStructure^B\setminus\parentStructure^A)$
			from $y$ to $x$.
		\end{enumerate}	
	\end{definition}
	
	Larger graphs have larger latent sets in the following sense,
	which will make algorithmic search simpler:
	\begin{lemma}[Latent Monotonicity]\label{lemma:monotonicity_of_L}
		Given a sub-family of embeddings
		$\{\psi^A_j\}_{J_0^A} \leq \{\psi^B_j\}_{J_0^B}$,
		then
		\begin{equation*}
			\forall L^B \in \minLatentSets^B:
			\halfquad
			\exists L^A \in \minLatentSets^A:
			L^A \subset L^B\cap\nodes^A
			\txt.
		\end{equation*}
	\end{lemma}
	\begin{proof}
		Let $L^B\in \minLatentSets^B$
		be arbitrary.
		We first show: $O^A := \nodes^A \setminus (L^B\cap\nodes^A)$
		is $\psi^A$-observed.
		Since $L^B \minLatentSets^B$,
		$O^B=\nodes^B\setminus L^B$ is $\psi^B$-observed,
		thus (by definition, \ref{def:observed})
		$|J^{B,\text{obs}}_0(N,O^B)| \rightarrow \infty$ for $N\rightarrow\infty$,
		where $J^{B,\text{obs}}_0(N,O^B) \subset J_0^B \subset J_0^A$
		by definition.
		Since $O^A = \nodes^A \setminus (L^B\cap\nodes^A) \subset O^B$,
		we have $J^{B,\text{obs}}_0(N,O^B) \subset J^{A,\text{obs}}_0(N,O^A)$.
		So also $|J^{A,\text{obs}}_0(N,O^A)| \rightarrow \infty$,
		and by definition \ref{def:observed}
		$O^A$ is $\psi^A$-observed.
		
		Finally,
		$\exists L^A \in \minLatentSets^A$ with $L^A\subset L^B\cap\nodes^A$,
		by definition of $\minLatentSets^A$ and
		$O^A$ being observed.
	\end{proof}
	
	\begin{rmk}[Data Availability]\label{rmk:data_availability}
		Given a family of embeddings $\{\psi_j\}_{j\in J_0}$
		with a minimal latent-set $L$ then (by definition)
		$O=\nodes\setminus L$ is observed.
		By Def.\ \ref{def:observed},
		thus
		$|J^{\text{obs}}_0(N,O)|\rightarrow \infty$ for
		$N\rightarrow\infty$, where
		$J^{\text{obs}}_0(N,O) =
		\{j \in J_0| \psi_j(O)\subset \viewport(N)\}$.
		So there is, asymptotically, infinite data available;
		this is implicit in the definition of decorated families of embeddings.
		Note, that freeness at inner nodes (Def.\ %
		\ref{def:local_graph_embedding_family}%
		\ref{def:local_graph_embedding_family:freeness}),
		together with $|\nodesInner|<\infty$,
		protects inner nodes sufficiently from degeneracies in this data-set;
		see also Lemma \ref{lemma:valid_datasets}.
	\end{rmk}
	
	\begin{rmk}[Relevant Sparse Families]\label{rmk:relevant_sparse}
		In analogy to Rmk.\ \ref{rmk:relevant_latent},
		we call a decorated family of embeddings $(\{\psi_j\}_{j\in J_0}$, $L,
		\ancestralStructure)$
		weakly maximal if for any $(\{\psi'_j\}_{j\in J'_0}$, $L',
		\ancestralStructure')$
		that satisfies all requirements for
		$(\{\psi'_j\}_{j\in J'_0}$, $L',
		\ancestralStructure')\leq (\{\psi_j\}_{j\in J_0}$, $L,
		\ancestralStructure)$
		other than $J'_0\subset J_0$, also $J'_0 \subset J_0$ holds.
	\end{rmk}
	
	Larger graphs have (essentially) larger ancestral structure,
	thus together with Lemma \ref{lemma:monotonicity_of_L}
	all attached structure is monotonic.
	
	\begin{lemma}[Ancestral\Slash{}Decorated Monotonicity]
		\label{lemma:monotonicity_of_Anc}
		Given a family of embeddings
		$\{\psi^A_j\}_{J_0^A}$
		and a decorated family $(\{\psi^B_j\}_{J_0^B}, L^B, \ancestralStructure^B)$
		with
		$\{\psi^A_j\}_{J_0^A} \leq \{\psi^B_j\}_{J_0^B}$,
		then
		\begin{equation*}
			\exists L^A \in \minLatentSets^A
			\exists \txt{ valid } \ancestralStructure^A
			\txt{ such that }
			(\{\psi^A_j\}_{J_0^A}, L^A, \ancestralStructure^A)
			\leq
			(\{\psi^B_j\}_{J_0^B}, L^B, \ancestralStructure^B)
			\txt.
		\end{equation*}
	\end{lemma}
	\begin{proof}
		A suitable $L^A \subset L^B\cap\nodes^A$ exists
		by Lemma \ref{lemma:monotonicity_of_L}.
		We may shrink $J_0^A$ as long as $J_0^A \subset J_0^B$,
		thus it is enough to check $j\in J_0^B$ for ancestral paths.
		Let $\gamma$ be an ancestral path, \ie a directed path in the $I$-graph from
		$\psi^A_j(y)$ to $\psi^B_j(x)$ (see Def.\ \ref{def:ancestral_structure})
		for $y\in\nodesInnerProper^A$ to $x\in\nodesOuterProper^A$
		along $I$-graph edges covered by edges in $\graph^A$.
		Then either $\gamma$ is covered by edges in $\graph^B$
		(and thus $(y,x)\in P_{B\setminus A}$) or
		$\gamma$ leaves $\graph^B$.
		If $\gamma$ leaves $\graph^B$, since $\parentStructure^B$ contains
		only edges pointing to $\nodesInnerProper^B$, the last covered edge
		points at $\psi^B_j(y')$ for $y'\in\nodesInnerProper^B$.
		To end up at $x\in \nodes^A$, $\gamma$ must at some
		point (possibly at $x$) reenter $\nodes^B$,
		by Def.\ \ref{def:local_graph_embedding}~%
		\ref{def:local_graph_embedding:inner_parents_incl}
		and \ref{def:local_graph_embedding:proper_edges_a}
		this can only happen at $\psi^B_j(x')$ for $x'\in\nodesOuterFixed^B$.
		The restriction of $\gamma$ to $\psi_j^B(y') \rightsquigarrow \psi_j^B(x')$
		is an ancestral path, so by validity of $\ancestralStructure^B$
		we must have $(y',x')\in\ancestralStructure^B$.
		Doing this whenever $\gamma$ leaves $\graph^B$ again,
		we get a path in
		$\ancestralStructure^B\cup(\parentStructure^B\setminus\parentStructure^A)$
		(there are no edges in $\parentStructure^A$ as $\gamma$ is outside
		of $\graph^A$ as ancestral path for $\graph^A$).
		This shows $(y,x)\in P_{B\setminus A}$.
	\end{proof}

	\subsection{Identifiability}
	\label{apdx:identifiablity}
	
	We extend on the remarks made in §\ref{sec:identifiability}
	about meaningful definitions of identifiability.
	
	\paragraph{Statistical Estimation:}
	We intentionally separate the question of
	having data for a problem (our notion of identifiability)
	from the question of the existence of statistical estimators
	and notions of estimator consistency.
	
	Already (absolute) density-estimation is usually an ill-posed problem
	\citep[§2.5 (p.\,36ff)]{VapnikEstimation},
	also problems like the hardness of conditional independence testing arise from
	conditional density estimation \citep[p.\,1517]{shah2020hardness}.
	Many approaches to such problems, restoring meaningful notions
	of consistency and of the actually estimable, are known.
	Influential and informative approaches include
	VC-theory \citep{VapnikEstimation} or Skorohod-embeddings
	\citep[§12 (p.\,220ff)]{kallenberg1997foundations}
	(which includes for example functional CLT formulations
	like Donsker's theorem, cf.\ \citep[Thm.\,12.9]{kallenberg1997foundations}).
	
	While learning conditional distributions is not generally possible
	in a naive sense, these (or often also much simpler) ideas
	still allow to learn meaningful and practically useful results.
	For example, in practice one is often interested in causal
	effects (expectation-values), optimal choices of interventional parameters
	or approximations by specific
	function-classes.
	These problems often \emph{can} be resolved in a satisfactory way,
	and provide statistical estimators that provide
	interpretable results given enough (and suitably curated) data.
	
	We focus on the aspect of "suitably curating" data for sub-problems,
	without regard to what kind of statistical estimator is used in the end.
	That is, \emph{we focus on the causal and symmetry aspects of the problem}.
	The kernels "directly identified" in our notion contain all information
	about that particular (combination of) model-mechanisms, queries etc.,
	thus any estimator learning some property of that kernel should
	be applicable to the curated data-set if the actual model-property in 
	question is compatible with estimator-assumptions (like linearity etc.).

	\paragraph{Stationarity:}
	
	The random walk example \ref{example:random_walk}
	may raise, for the time-series case,
	the question about connections to stationarity.
		
	In some cases, it may be possible to restore an IID-like notion
	of identifiability via a stationarity (or similar) condition.
	But is such a condition necessary, or even helpful?
	Fundamentally, the properties of the model and its predictions
	(like effects) do not seem to have anything to do with stationarity,
	so unsurprisingly it turns out to be unnecessary in a systematic
	approach via symmetry (see main text).
	Concerning the actual helpfulness of such assumptions,
	there are two cases that should be distinguished:
	Stationarity, as an a-priori assumption about the form of
	a \emph{specific} model\allowbreak{}\mbox{(-}realization),
	can be helpful, see §\ref{sec:queries_beyond_basic}, §\ref{sec:knowledge_closures}.
	On the other hand stationarity as a means of structuring the
	problem, besides conceptual doubts (see above), does not
	seem to solve, but rather to move and obfuscate the difficulty.

	\subsection{Identification from Families of Embeddings}
	\label{apdx:id_from_embeddings}
	
	The main technical content of §\ref{sec:extraction_from_data}
	--
	besides finding workable definitions for what embeddings
	should look like and how additional structure
	like latent-sets or ancestral information can be encoded
	--
	is in the actual identification of structured kernels.	
	
	We first need a simple result about valid data-sets	
	provided by families of embeddings.
	
	\begin{lemma}[Valid Data-Sets]
		\label{lemma:valid_datasets}
		Given a decorated family of embeddings
		$(\{\psi_j\}_{j\in J_0},L)$ and
		$N_Y \subset \nodesInnerProper \setminus L$,
		$N_X \subset \nodes\setminus (N_Y \cup L)$,
		and $N_0 \subset \nodesOuterFixed \cap L$
		then
		the data-set
		$\mathcal{D}=((\randomVar{X}_j, \randomVar{Y}_j)_{j\in J_0}$,
		$\randomVar{X}')$
		defined for $j\in J_0$ by the tuples
		\begin{align*}
			\randomVar{X}_j
			\halfquad&=\halfquad
			(\anyVar_{\psi_j(n)})_{n\in N_X}\\
			\randomVar{Y}_j
			\halfquad&=\halfquad
			(\anyVar_{\psi_j(n)})_{n\in N_Y}\\
			\randomVar{X}'
			\halfquad&=\halfquad
			(\anyVar_{i^n})_{n \in N_0}
			\txt.
		\end{align*}
		is valid (Def.\ \ref{def:data_set}).
	\end{lemma}
	\begin{rmk}
		\label{rmk:datasets_and_support_multi_level}
		Our definition of valid datasets (Def.\ \ref{def:data_set},
		Def.\ \ref{def:data_set:apdx})
		does not require any guarantee on the non-degeneracy of $\randomVar{X}_j$.
		In particular, it is always legal to move nodes from $\randomVar{X}'$
		to $\randomVar{X}_j$ if they are not hidden (or do not
		overlap $\graphObs$ in the proof of Thm.\
		\ref{thm:extract_from_backdoor_complete} below).
		This choice should be understood in the context
		of Ass.\ \ref{ass:support} and Rmk.\ \ref{rmk:kernels_uniqueness}:
		If $\randomVar{X}_j$ is degenerate (the same random
		variable $\anyVar_i$ appears asymptotically infinitely often),
		then the directly identified kernel (Def.\ \ref{def:identification_direct_mt})
		is learned on an observational support which has
		singular points of non-zero probability.
		This is no different from other problems related to observational
		support. If \emph{all} $j$ lead to the same value,
		there is no observational support elsewhere and
		Ass.\ \ref{ass:support} severely restricts allowed conclusions about
		transfer.
		But similar problems can occur for non-$X$-degenerate data-sets
		as well. It seems conceptually more meaningful to
		treat this problem as an issue of observational support,
		rather than as a problem of multi-level statistics, see also
		§\ref{apdx:iid:mz_transport}.
	\end{rmk}
	\begin{proof}
		By definition of decorated families,
		$L$ is a minimal latent subset, thus
		$O:=\nodes \setminus L$ is observed (Def.\ \ref{def:observed}),
		\ie
		\begin{align*}
			&J^{\text{obs}}_0(N,O)
			\halfquad:=\halfquad
			\{j \in J_0| \psi_j(O)\subset \viewport(N)\}
			\txt,\\
			&\txt{satisfies }
			|J^{\text{obs}}_0(N,O)|
			\rightarrow \infty
			\quad\txt{for } N\rightarrow\infty
			\txt.
		\end{align*}
		The sets $J_0^{\txt{obs}}(N)$ defined in Def.\ \ref{def:data_set} as
		\begin{equation*}
			J_0^{\txt{obs}}(N)
			\halfquad:=\halfquad
			\{\halfquad
			j\in J_0
			\halfquad|\halfquad
			\forall m:
			i_{X,j}^{(m)} \in \viewport(N),
			\forall m':
			i_{Y,j}^{(m')} \in \viewport(N)
			\halfquad\}
			\txt,
		\end{equation*}
		by $N_X\cup N_Y \subset O$, satisfy
		$J^{\text{obs}}_0(N,O) \subset J_0^{\txt{obs}}(N)$,
		in particular $|J_0^{\txt{obs}}(N)|\rightarrow\infty$
		for $N\rightarrow \infty$.
		The non-degeneracy of elements in $N_Y$,		
		$j\neq j' \Rightarrow \forall m: i_{Y,j}^{(m)}\neq i_{Y,j'}^{(m)}$,
		is satisfied by Def.\ \ref{def:local_graph_embedding_family}~%
		\ref{def:local_graph_embedding_family:freeness} and
		$N_Y\subset \nodesInnerProper$.
	\end{proof}
	
	Next, we show that backdoor free embedded families
	allow for the identification of structured kernels.
	The intuition is that $\psi_j(\graphObs)$ behaves like a
	structural c-component of the $I$-graph.
	We need to find suitable sub-structures that
	are both $j$-independent and have suitable data-sets
	attached to them, so that they become directly identifiable.
	
	\begin{CopyThm}{thm:extract_from_backdoor_complete}{}
		\label{thm:extract_multi_level}
		Given a backdoor free family of embeddings 
		$(\{\psi_j\}_{j\in J_0},L,\ancestralStructure)$
		(Def.\ \ref{def:backdoor_free:apdx}),
		then 
		\begin{equation*}
			\mu(\graphObs(\graph),L)
			\txt{ (Def.\ \ref{def:model_alignment_and_kernel}) is $\delta_i$-identifiable (Notation \ref{def:delta_i_identified}).}
		\end{equation*}				
	\end{CopyThm}
	\begin{proof}
		We first focus on the single-level case, the multi-level case is
		then studied in step 3.
		Let $j\in J_0$ be arbitrary.
		We build a model-aligned structural graph 
		$\graphStructural^j$ with only
		inner nodes. To this end, set $\nodes^j=\nodesInner^j=\Anc_I(\psi_j(\nodes))$
		and define an edge-set $\edgesStructural^j
		=\{(p,c)\in\nodes^j\times\nodesInner^j| p\in\Pa_I(c)\}$
		the same as in the $I$-graph.
		Align (inner) nodes to the model-mechanisms,
		\ie $\mu^i := f_i$ for $i\in\nodesInner^j$.
		Then define $L^j := L \cup (\nodes^j \setminus \nodes)$.
		
		With $\psi_j: \nodes \hookrightarrow I$ having
		$\img(\psi_j) \subset \nodes^j$ by construction
		of $\nodes^j$, we can identify nodes in $\graphObs$
		with their images in $\nodes^j$ (by injectivity of $\psi_j$,
		Def.\ \ref{def:local_graph_embedding}) making
		it a sub-graph $\graphObs \leq \graphStructural^j$
		(the alignment of nodes agrees by
		Def.\ \ref{def:local_graph_embedding}~%
		\ref{def:local_graph_embedding:applicable}).
		
		\textbf{Step 1:} We show, by backdoor freeness, $\graphObs \leq \graphStructural^j$
		is a union of structural c-components
		(\ie every structural c-component $\graphStructural^c\leq\graphStructural^j$
		touching, in the sense of non-empty overlap of inner nodes,
		$\graphObs\leq \graphStructural^j$ is also a structural c-component of $\graphObs$
		and every structural c-component of $\graphStructural^c\leq\graphObs$ is also a
		structural c-component of $\graphStructural^j$).
		
		We write $n \approx_j n'$ for $n \approx_{L} n'$ in $(\graphStructural^j,L^j)$
		and $n \approx_\gObs n'$ for $n \approx_L n'$ in $(\graphObs,L)$.
		
		We have to show: 
		(a) Given $n,n'\in\nodesInner^\gObs$ with $n \approx_\gObs n'$,
		then $n \approx_j n'$ and
		(b) given $n\in\nodesInner^\gObs$, $n'\in\nodesInner^j$
		with $n \approx_j n'$,
		then $n' \in \nodesInner^\gObs$ and 
		$n \approx_\gObs n'$.
		
		(a)
		Let $n,n'\in\nodesInner^\gObs$ with $n \approx_\gObs n'$ be arbitrary.
		By definition of $\approx_\gObs$, there are
		$y,w\in\nodesInner^\gObs\setminus L$ (possibly equal to $n$ or $n'$),
		latent (all nodes on $\gamma$
		other than $y$ are in $L$) paths $\gamma_y: n \rightsquigarrow y$
		and $\gamma_w: n' \rightsquigarrow w$ in $\graphObs$ (possibly trivial),
		an element $l\in L$ and latent
		paths $\gamma'_y: l \rightsquigarrow y$
		and $\gamma'_w: l \rightsquigarrow w$ in $\graphObs$.
		
		By $L\subset L^j$ non-endpoint nodes on these paths
		are also latent in $\graphStructural^j$.
		Edges on these paths are in $\edgesStructural^j$ by
		Def.\ \ref{def:local_graph_embedding}~%
		\ref{def:local_graph_embedding:proper_edges_a}
		and the definition of $\edgesStructural^j$ as
		$I$-graph edges.
		Thus these paths are latent paths in $\graphStructural^j$ and
		thus $n \approx_j n'$.
		
		(b)
		Let $n\in\nodesInner^\gObs$, $n'\in\nodesInner^j$
		with $n \approx_j n'$ be arbitrary.
		By definition of $\approx_j$, there are
		$y,w\in\nodesInner^j\setminus L^j$ (possibly equal to $n$ or $n'$),
		latent paths $\gamma_y: n \rightsquigarrow y$
		and $\gamma_w: n' \rightsquigarrow w$ in $\graphStructural^j$
		(possibly trivial),
		an element $l\in L^j$ and latent
		paths $\gamma'_y: l \rightsquigarrow y$
		and $\gamma'_w: l \rightsquigarrow w$ in $\graphStructural^j$.
		
		We show: (i)
		Given a latent path $\gamma$ in $\graphStructural^j$
		starting at $a\in\nodesInner^\gObs$ to a node $b\notin L^j$,
		then $b\in\nodesInner^\gObs\setminus L$.
		
		Proof of (i):
		Recall that $L^j = L \cup (\nodes^j \setminus \nodes^\gObs)$,
		thus $b \in\nodes^\gObs\setminus L$.
		By contradiction, assume it were $b\in\nodesOuter^\gObs$.
		With $a\in \nodesInner^\gObs$, we have $a\neq b$
		and $\gamma$ is non-trivial (contains more than one node).
		By definition (Def.\ \ref{def:structural_graph}),
		there are no (proper) edges in $\graphObs$
		ending at $b$, thus the last node before $b$ on $\gamma$
		is not in $\nodes^\gObs$ (using Def.\ \ref{def:local_graph_embedding}~%
		\ref{def:local_graph_embedding:proper_edges_a}).
		The starting point $a$ of $\gamma_y$ is in $\nodes^\gObs$,
		thus there is a last (along  $\gamma_y$) node
		$l^\gObs\in\nodes^\gObs \cap L^j=\nodes^\gObs\cap L$.
		Restrict $\gamma$ to $l^\gObs\rightsquigarrow b$.
		This is a path in the $I$-graph (by construction of $\graphStructural^j$)
		not through $\graphObs$.
		By validity of the ancestral structure $\ancestralStructure$
		(Def.\ \ref{def:ancestral_structure}),
		there is $l^\gObs \rightsquigarrow b$ in $\ancestralStructure$.
		This is a contradiction to backdoor freeness of
		$(\{\psi_j\}_{j\in J_0}, L, \ancestralStructure)$.
		
		Next we show: (ii) Given a latent path
		$\gamma$ in $\graphStructural^j$ ending at $b\in\nodesInner^\gObs\setminus L$,
		then $\gamma$ is a latent path in $\graphObs$.
		
		Proof of (ii):
		By contradiction. If $\gamma$ is not a latent path in $\graphObs$,
		then by $L\subset L^j$, there is a last (along $\gamma$)
		node $a$ with $a\in L^j\setminus\nodesInner^\gObs$.
		Clearly $a\neq b$, because $b\in\nodesInner^\gObs$.
		
		By $\nodesOuter^\gObs = \nodesOuterProper \cup (\nodesOuterFixed \setminus L)$
		by definition of $\graphObs$ and
		$\nodesOuterProper \cap L = \emptyset$ by backdoor freeness,
		we have $\nodesOuter^\gObs \cap L = \emptyset$.
		By definition of $L^j$ also $\nodesOuter^\gObs \cap L^j = \emptyset$.
		In particular $(*)$ $a\notin \nodes^\gObs$.
		
		But $a\in\Pa_I(b')$ where $b'$ is the node directly after $a$
		(on $\gamma$). Note that $b'\in\nodesInner^\gObs$,
		because $a$ by definition was the \emph{last} node (along $\gamma$)
		not in $\nodesInner^\gObs$.
		By Def.\ \ref{def:local_graph_embedding}~%
		\ref{def:local_graph_embedding:inner_parents_incl}
		also $a\in\nodes^\gObs$, contradicting the previous result $(*)$.
		
		Finalizing the proof of (b):
		By (i) applied to $\gamma_y$:
		$y\in \nodesInner^\gObs\setminus L$.
		Thus by (ii), $\gamma_y$ and $\gamma_y'$ are latent paths in $\graphObs$.
		By (i), applied to $\gamma_w$,
		$w\in \nodesInner^\gObs\setminus L$.
		Thus by (ii), $\gamma_w$ and $\gamma_w'$ are latent paths in $\graphObs$.
		In particular $n'\in \nodesInner^\gObs$.
		Finally the latent paths $\gamma_y$, $\gamma_y'$
		and $\gamma_w$, $\gamma_w'$ in $\graphObs$ show
		$n \approx_\gObs n'$.
		
		\textbf{Step 2}:
		We show the single-level ($\nodesOuterFixed=\emptyset$) case.
		
		For $n\in \nodesInner^\gObs\setminus L$,
		contained in the structural c-component $\graphStructural^c \leq
		\graphObs, \graphStructural^j$ (by step 1),
		we have $A^n(\graphObs,L) = A^n(\graphStructural^c,L^c)
		=A^n(\graphStructural^j,L^j)$ by applying Lemma
		\ref{lemma:properties_c_components}~%
		\ref{lemma:properties_c_components:atoms} twice.
		Crucially $A^n(\graphObs,L)$ (and thus $A^n(\graphStructural^j,L^j)$)
		does not depend on $j$.
		Further, by Lemma \ref{lemma:atom_properties}~%
		\ref{lemma:atom_properties:product},
		$\mu(\graphStructural^j,L^j) = \prod_{n'\in\nodesInner^j\setminus L} A^{n'}(\graphStructural^j, L^j)$ (note that only terms $n'$ also in
		$\nodesInner^\gObs\setminus L$ do not depend on $j$), thus
		for  $n\in \nodesInner^\gObs\setminus L$:
		\begin{equation*}
			\marginalize{\mu(\graphStructural^j,L^j)}{n'>n}
			=
			\BigKernelCompound{\prod_{n'\in\nodesInner^j\setminus L^j, n'<n}
				A^{n'}(\graphStructural^j, L^j)}
				\otimes
				A^n(\graphObs,L)\txt.
		\end{equation*}
		This is already in the correct form for the right-hand-side in 
		Def.\ \ref{def:identification_direct_mt}.
		For $n\in \nodesInner^\gObs\setminus L$,
		by Lemma \ref{lemma:valid_datasets},
		there is further a valid data-set $\mathcal{D}_n$ for
		$N_X^n=\{n'\in\nodesInner^j\setminus L^j|n'<n\}$
		(by $\nodesInner^j\setminus L^j=\nodes^\gObs\setminus L$
		Lemma \ref{lemma:valid_datasets} applies),
		$N_Y^n = \{n\}$ and $N_0^n = \emptyset$.
		
		Finally, by Lemma \ref{lemma:obs_world_properties}~%
		\ref{lemma:obs_world_properties:single_level}
		\begin{equation*}
			\shallowDistr(\randomVar{X}_j,\randomVar{Y}_j)
			= \marginalize{\mu(\graphStructural^j,L^j)}{n'>n}\txt,
		\end{equation*}
		so Def.\ \ref{def:identification_direct_mt} applies
		and $A^n(\graphObs,L)$ is directly identifiable.		
		By Lemma \ref{lemma:atom_properties}~%
		\ref{lemma:atom_properties:product},
		$\mu(\graphObs,L)$ is identifiable (Def.\ \ref{def:identification}).
		
		\textbf{Step 3}:
		We show the multi-level case.
				
		Define the structural graph $\graphStructural^j$
		as before.
		Applying Lemma \ref{lemma:obs_world_properties}~%
		\ref{lemma:obs_world_properties:multi_level}
		(cf.\ below for satisfiability of the hypothesis)
		with $N_0 = \psi_j(\nodesOuterFixed)$
		provides a $(\graphStructural^j)'$;
		the statement of the lemma ensures
		$\graphObs[\delta_i=\delta[\mathcal{D}_i]]\leq (\graphStructural^j)'$.
		Note that, if $\nodesInner^\gObs \cap \Anc_I(\nodesOuterFixed) = \emptyset$,
		its graphical property (see Lemma \ref{lemma:obs_world_properties}~%
		\ref{lemma:obs_world_properties:multi_level}) applies
		for $Y=\nodesInner^\gObs$ by result (i) of step 1 above
		and then replaces (i) in the remainder of step 1.
		Then step 1 produces the same result as before.
		
		The condition $\nodesInner^\gObs \cap \Anc_I(\nodesOuterFixed) = \emptyset$
		in the hypothesis of Lemma \ref{lemma:obs_world_properties}~%
		\ref{lemma:obs_world_properties:multi_level}
		may not hold true for all $j\in J_0$. It fails however only for
		finitely many, and thus we can replace
		$J_0$ by a (still infinite) $J_0' := J_0 \setminus J_0^{\txt{invalid}}$.
		Removing a finite ($N$-independent)
		subset does not affect asymptotic properties: if the used
		data-set $\mathcal{D}$ is valid (using an infinite subset of $J_0$),
		then the removal of any intersection with $J_0^{\txt{invalid}}$ does not
		change this validity (infinitude of the subset of $J_0$).
		The subset we want to remove is:
		\begin{equation*}
			J_0^{\txt{invalid}}
			\halfquad:=\halfquad
			\big\{\halfquad
			j\in J_0
			\halfquad\big|\halfquad	
			\psi_j(\nodesInnerProper) \cap \Anc_I( \psi_j(\nodesOuterFixed) )
			\neq \emptyset
			\halfquad\big\}
			\txt.
		\end{equation*}
		This set is finite:
		By the finite past assumption (Ass.\ \ref{ass:finite_past}),
		$\Anc(\psi_j(\nodesOuterFixed))$ is finite.
		For any $n\in\nodesInnerProper$, by freeness
		Def.\ \ref{def:local_graph_embedding_family}%
		\ref{def:local_graph_embedding_family:freeness},
		$\psi_*(n)$ is injective, thus
		$\psi_*(n)^{-1}(\Anc(\psi_j(\nodesOuterFixed))) \allowbreak\subset J_0$
		is finite.
		The number of inner nodes $|\nodesInnerProper|<\infty$ is also finite
		(by Def.\ \ref{def:local_graph_new})
		thus $\cup_{n\in\nodesInnerProper}\psi_*(n)^{-1}(\Anc(\psi_j(\nodesOuterFixed)))$
		is finite.
		Finally we show
		\begin{equation*}
			J_0^{\txt{invalid}} \subset
			\cup_{n\in\nodesInnerProper}\psi_*(n)^{-1}(\Anc(\psi_j(\nodesOuterFixed)))
		\end{equation*}
		(thus $J_0^{\txt{invalid}}$ is finite as subset of a finite set):
		Let $j\in J_0^{\txt{invalid}}$ be arbitrary.
		By construction, 
		$\psi_j(\nodesInnerProper) \cap \Anc_I(\psi_j(\nodesOuterFixed))\neq \emptyset$,
		\ie
		$\exists n\in \nodes \nodesInnerProper: \psi_j(n)\in\Anc_I(\psi_j(\nodesOuterFixed))$.
		Therefore $j\in \psi_*(n)^{-1}(\Anc(\psi_j(\nodesOuterFixed)))$, in particular
		\begin{equation*}
			j\in \cup_{n\in\nodesInnerProper}\psi_*(n)^{-1}(\Anc(\psi_j(\nodesOuterFixed)))
			\txt.
		\end{equation*}

		Returning to step 3, with step 1 established, it remains to modify step 2.
		$\mu(\graphObs,L)[\delta_i]$ now
		is a functional of $\delta_i$ for $i\in \psi_j(\nodesOuterFixed)$,
		but
		step 2 works essentially as before,
		Lemma \ref{lemma:valid_datasets} is applied
		with $N_X^n=\{n'\in\nodesInner'\setminus (L' \cup N_0)|n'<n\}$
		$N_Y^n = \{n\}$ and $N_0^n = N_0$
		to obtain a data-set $\mathcal{D}_n$
		on which Lemma \ref{lemma:obs_world_properties}~%
		\ref{lemma:obs_world_properties:multi_level} yields
		\begin{equation*}
			\shallowDistr(\randomVar{X}_j,\randomVar{Y}_j|\randomVar{X}'=x')
			= \marginalize{\mu(\graphStructural',L')[x']}{n'>n}\txt,
		\end{equation*}
		thus Def.\ \ref{def:identification_direct_mt} applies
		and $A^n(\graphObs[\delta_{N_0}[x']]$
		is directly identifiable, where $\randomVar{X}' = x'$
		is the value taken by nodes in (the $j$-independent) $N_0$
		in the data-set $\mathcal{D}_n\subset\mathcal{D}(\omega)$,
		\ie $x' = \anyVar_{\psi_j(\nodesOuterFixed)}(\omega)$.

		The claim of the Lemma is that $\mu(\graphObs,L)$ is
		$\delta_i$-identifiable, \ie it remains to show, that
		\begin{equation*}
			\mu(\graphObs,L)[\delta[\anyVar_{\psi_j(\nodesOuterFixed)}(\omega)]]
			\overset{\txt{a.\,s.}}{=}
				\mu(\graphObs, L)[\delta[\anyVar_{\psi_j(\nodesOuterFixed)}]]
				\txt.
		\end{equation*}
		The left-hand-side is a constant (not random), thus we have to show
		$\realizedDistr$ a.\,s.\ the right-hand-side takes this value.
		Let $\omega' \in\Omega$ with $\mathcal{D}(\omega') = \mathcal{D}(\omega)$
		be arbitrary
		(by definition $\realizedDistr(\cdot)=
		\shallowDistr(\cdot|\{\anyVar_i\}_{i\in\viewportEventual}
		=\mathcal{D}(\omega))$, so we know
		$\mathcal{D}(\omega') = \mathcal{D}(\omega)$ almost surely).
		By the above result 
		$\mu(\graphObs,L)[\delta[\anyVar_{\psi_j(\nodesOuterFixed)}(\omega')]]$
		is (uniquely) identifiable from $\mathcal{D}(\omega') = \mathcal{D}(\omega)$,
		thus
		\begin{equation*}
			\mu(\graphObs,L)[\delta[\anyVar_{\psi_j(\nodesOuterFixed)}(\omega')]]
			=
			\mu(\graphObs,L)[\delta[\anyVar_{\psi_j(\nodesOuterFixed)}(\omega)]]
			\txt.
		\end{equation*}
	\end{proof}
	
	\begin{rmk}[Practical Computation of Extraction]
		\label{rmk:practical_computation_of_extraction}
		Inspecting step 2 of the proof of
		Thm.\ \ref{thm:extract_multi_level},
		it becomes evident that we first have to learn
		($j$-independent) atoms $A^n$ for $n\in\nodesInner\setminus L$.
		This can be achieved in practice as described in
		Rmk.\ \ref{rmk:compute_atoms_practice}:
		We know the relevant arguments and data-set, thus
		we simply apply our favorite estimator.
		Knowing all the atoms of $(\graphObs, L)$,
		by Lemma \ref{lemma:atom_properties}~%
		\ref{lemma:atom_properties:product},
		we can compute $\mu(\graphObs,L)$ as their product.
	\end{rmk}

	We have thus seen, that backdoor-free families of embeddings
	allow for the identification of associated kernels from data.
	On the other hand, backdoor-free families of embeddings
	are a rather generic means of curating data-sets with invariant properties,
	thus it seems plausible that Conj.\ \ref{conjecture:extraction}
	might be true.
	
	\begin{rmk}
		It seems that internal structure (like linearity)
		can be employed post hoc as a computational means
		(see §\ref{sec:knowledge_closures})
		without the need to modify this claim.
		It is however conceivable that using internal structure might
		allow to effectively improve symmetry locally, in which case
		the above conjecture \ref{conjecture:extraction} could hold at most in the generic case.
	\end{rmk}

	\subsection{Algorithmic Construction of Backdoor-Free Families}
	\label{apdx:algo_bd_free}

	Algorithm \ref{algo:fcs_r} shows how to
	extend a given decorated family $(\psi, L, \ancestralStructure)$ to
	obtain (non-unique) c-connected backdoor-free
	families of embeddings.
	\begin{definition}[C-Connected Backdoor-Free Families]
		\label{def:minimal_backdoor_complete}
		We call a backdoor-free family of embeddings
		$(\{\psi_j\}_{j\in J_0},L,\ancestralStructure)$
		c-connected, if $\graphObs(\graph,L)$ is c-connected.
	\end{definition}
	It is enough to find all c-connected backdoor-free
	families:
	\begin{lemma}\label{lemma:c_conn_extr_enough}
		Let $\mathcal{B}_{\txt{C}}$ be the set of all
		c-connected backdoor-free families of embeddings
		and $k'=(\graphStructural',L') \in
			\bigcup_{(\psi,L,\ancestralStructure)\in\mathcal{B}}
			k(\psi,L,\ancestralStructure)$
		(cf.\ Def.\ \ref{def:knowledge_set_extracted}),
		then $\exists (\psi^1,L^1,\ancestralStructure^1),\ldots,
		(\psi^m,L^m,\ancestralStructure^m) \in\mathcal{B}_{\txt{C}}$
		such that
		by repeated application of Lemma \ref{lemma:glue:mt}
		$\mu(\graphStructural',L')$ can be regularly computed from
		$\mu(\graphObs_1,L^1), \ldots, \mu(\graphObs_m,L^m)$.
		In particular by Thm.\ \ref{thm:extract_from_backdoor_complete},
		$k'\in\knowledgeSetBasic$
		is identifiable from only c-connected backdoor-free families.
	\end{lemma}
	\begin{proof}
		Let $\graphStructural^c\leq \graphStructural'=\graphObs(\graph',L')$
		be a structural c-component.
		Define a local graph
		$\graph''$ with inner nodes
		$\nodesInner^c\setminus	(\nodesOuterFixed'\cap L)$,
		pinned nodes $\nodesInner^c \cap (\nodesOuterFixed'\cap L)$
		and external nodes $\nodesOuter^c$.
		Then the restrictions of $\psi'_j|_{\nodes''}$ to $\graph''$
		are again a family of embeddings.
		The choices $L''=L'\cap \nodes''$ and
		$\ancestralStructure'' = \ancestralStructure'
		\cap (\nodesInnerProper''\times\nodesOuterFixed'')$
		are such that $\nodes''\setminus L''$
		is observed and $\ancestralStructure''$ is valid,
		but they may not be minimal anymore.
		If a smaller (minimal) $L'''\subset L''$ exists,
		repeat the previous proof-steps with (potentially
		smaller) c-components of $((\graphObs)'', L''')$;
		note: smaller $\ancestralStructure'''$ will result in
		$\graphObs$ having more simplifications (Lemma \ref{lemma:simplify}),
		which is needed to get smaller pieces,
		but is not needed to recover the glued $\graphStructural'$.
		
		Lemma \ref{lemma:glue:mt} automatically applies with the c-components
		(which each contain themselves).
		If the above steps were repeated with smaller $L'''$,
		the gluing target (larger graph) is the one with $L=\cup_c L'''_c$
		(where $L'''_c$ is $L'''$ of c-component $c$)
		a finial marginalization will connected this result
		to the original $\graphStructural'$.
	\end{proof}

	\begin{algorithm}
		\renewcommand{\thealgorithm}{ExtractCS-R}
		\caption{extract\_cstructures\_relative}
		\label{algo:fcs_r}
		\textbf{Input:} A decorated family of embeddings $F^B=(\{\psi_j^B\}_{j\in J_0^B},L^B,\ancestralStructure^B)$.\\
		\textbf{Output:} The set of backdoor-free families of embeddings		
		$F=(\{\psi_j\}_{j\in J_0},L,\ancestralStructure)$
		c-connected and with $F^B \leq F$ (Def.\ \ref{def:subgraphs_local}).
		\begin{enumerate}[label=\arabic*)]
			\item\label{algo:find_min_backdoor_complete:min_L}
			\emph{Attach monotonic structure:}\\
			$L \forkAssign$ \texttt{minimal\_hidden}$(\psi, L)$.
			\Comment{No-op in first iteration.}\\
			$\ancestralStructure \forkAssign$
			\texttt{minimal\_ancestral$(\psi, L, \ancestralStructure)$}.
			\Comment{Relative choice by Lemma \ref{lemma:monotonicity_of_Anc}.}
			
			\item\label{algo:find_min_backdoor_complete:abs_pa}			
			\emph{Absorb hidden parent-chains:}\\
			$\psi \forkAssign$ \texttt{absorb\_parents}
			$(\psi, L\cap(\nodesOuter\setminus\nodesOuterFixed))$.\\
			\textbf{Repeat} (1--2) until
			\begin{equation*}
				\tag{$*$}
				L\cap(\nodesOuter\setminus\nodesOuterFixed)=\emptyset
				\txt.
			\end{equation*}
			
			\item\label{algo:find_min_backdoor_complete:abs_ch}
			\emph{Absorb starting nodes of backdoor paths:}\\
			$\psi, w \forkAssign$ \texttt{absorb\_child\_ancestral}%
			$(\psi, L, \ancestralStructure)$.\\
			\textbf{Repeat} (1--3) until
			\begin{equation*}
				\tag{$**$}
				(\{\psi_j\}_{j\in J_0},L,\ancestralStructure) \txt{ is backdoor-free.}
			\end{equation*}
			
			\item\label{algo:find_min_backdoor_complete:yield}		
			\emph{Collect branch into output set:}\\
			\textbf{Yield} $(\{\psi_j\}_{j\in J_0},L,\ancestralStructure)$.
		\end{enumerate}
	\hrule
	\vspace*{0.3em}
	\emph{Notation:}
	For set-returning sub-algorithms,
	we write $\forkAssign$ short for iterate over all elements of.
	Yield returns to iterating those other elements, appending
	the yielded result to the output set.
	Details on subroutines are given in \ref{apdx:algorithm_details}.
	\end{algorithm}
	
	\begin{algorithm}
	\renewcommand{\thealgorithm}{ExtractCS}
	\caption{extract\_cstructures}
	\label{algo:fcs}
	\textbf{Input:} Implicitly the model and viewport.\\
	\textbf{Output:} The set of c-connected backdoor-free families of embeddings		
	$F=(\{\psi_j\}_{j\in J_0},L,\ancestralStructure)$.
	\begin{enumerate}[label=\arabic*)]
		\item 
		\emph{Initialize:}\\
		$k := 0$. $R_0 := \emptyset$.\\
		\textbf{For each $J\in\mathcal{J}(M)$:}\\
		\hspace*{2em}
		Add \texttt{FCS-R}( \texttt{direct\_embbed}($J$) ) to $R_0$.
		\item\label{algo:fcs:step}
		\emph{Step:}\\
		$R_{k+1}=R_k$. $k := k+1$.\\
		\textbf{For each}
		$F \in R_k$ and 
		\textbf{For each}
		$\chi \in \RelevantCh(F)$:\\
		\hspace*{2em}
		Add \texttt{FCS-R}( \texttt{absorb\_child}$(F, \chi)$ )
		to $R_{k}$.\\
		\textbf{Repeat} \ref{algo:fcs:step}
		\textbf{until} $R_k = R_{k-1}$ \textbf{or} limit $k_{\txt{max}}$ reached.
		\item 
		\emph{Finalize:}\\
		\textbf{Return} $R_k$.
	\end{enumerate}
	\hrule
	\vspace*{0.3em}
	\emph{Remark:}
	We separate Algo.\ \ref{algo:fcs}
	from Algo.\ \ref{algo:fcs_r} because, for most practical
	applications, for example for use as input to
	Algo.\ \ref{algo:svs},,
	searching all relevant children in step \ref{algo:fcs:step}
	is not necessary, rather the \texttt{match} (and even \texttt{gluable})
	filters of \ref{algo:svs} can be integrated directly with Algo.\ \ref{algo:fcs_r}
	to only enumerate useful families of embeddings.
	\texttt{direct\_embbed}($J$) returns the
	embedding of Lemma \ref{lemma:direct_embedding}.
	\end{algorithm}
	
	In the algorithms \ref{algo:fcs_r} and \ref{algo:fcs} we sketch the general logic
	for the systematic construction of c-connected backdoor-free families of embeddings
	(absorbing operations are given in §\ref{apdx:algorithm_details}).
	Important to note are the many non-unique partial constructions,
	discernible from the $\forkAssign$ operators in \ref{algo:fcs_r}
	and for-each statements in \ref{algo:fcs}. These constructions
	in the IID case (§\ref{apdx:iid})
	are usually unique (Lemma \ref{lemma:iid_psi_classification}).
	Both algorithms are intentionally written such that ignoring
	these non-uniqueness challenges,
	they are easy to read and make sense for the IID case.
	Indeed the construction of c-components -- leaving non-uniqueness aside
	(each mechanism can occur locally in many different c-components)
	-- in the IID case is very similar.
	
	There are however two notable differences:
	First, in enumerating all (or all matching etc.) c-components \ref{algo:fcs}
	always starts from embedding a single variable (in \texttt{fcs-r}$($
	\texttt{direct\_embed}$(J))$),
	then adds children (ancestral or relevant) relative to this embedded
	variable. There is a simple reason for this approach: Embedding disconnected
	graphs does not usually make sense in our formalism (any form of "rigidity",
	equivariance relative to an anchor in another component, must be
	introduced in some ad-hoc way).
	
	Second, and maybe more importantly, we cannot, in general,
	start from maximal size c-components and go to c-subgraphs
	as we can in the IID case. The "largest" local c-structure
	is also the local c-structure with worst symmetry
	(Rmk.\ \ref{rmk:families_and_symmetries}),
	it may therefore happen that smaller c-structures (proper sub-graphs
	of this largest one) are observed (have smaller than the trivial intersection
	minimal latent set $L$), while the largest one is not.
	Thus it may happen that these smaller sub-graphs lead to better
	identification results.
	
	It seems usually not advisable to eagerly discover all elements of
	$\mathcal{B}_{\txt{C}}$. Rather, given a query, only parents and children
	matching the structure of the query have to be inspected (see Algo.\ %
	\ref{algo:svs}). These matching results can be discovered lazily
	(see §\ref{apdx:algorithm_details}).

	\section{Details on Query Identification}
	\label{apdx:queries}

	The formal claims of §\ref{sec:queries}
	fall in two main categories:
	The representation of basic queries by structured queries
	and the relation of underlying basic queries, in particular
	their target in the realized distribution,
	to the material of §\ref{sec:structured_kernels}
	and §\ref{sec:extraction_from_data}.
	These aspects are discussed separately in the two subsections 
	§\ref{apdx:query_structural_representation} and
	§\ref{apdx:query_identification}.
	As in §\ref{apdx:extraction_from_data}, we start
	by stating some small modifications that allow to account for the multi-context
	case.
	Finally, we briefly outline how our work fits into a larger
	formal context and can be extended in future work,
	for example to include instrumental variable arguments
	or stationarity considerations for time-series.

	\subsection{Multi-Level Modifications}
	
	Proofs are given in the sub-sections below, here we only
	provide the correct statements to be proved for the multi-level case.	
	In the multi-level case, we need to identify
	"contexts", \ie variables $\anyVar_i$ shared by the observations and a query.
	In example \ref{example:multi_level_illustrate},
	a query about an intervention in a \emph{specific} context (\eg river-site),
	has relevant contexts (\eg slope) as auxiliary index (see below).
	A query about an unknown (generic) context on the other hand
	instead contains a copy (a separate node sharing the mechanism
	but not observed descendants)
	of this context (\ie it is a Bayesian estimate integrating out
	the prior over, for example, slopes; see §\ref{apdx:multi_level_examples}).
	
	\begin{definition}[Auxiliaries]\label{def:aux}
		\textbf{This has no single-level analogue:}
		Given an embedded local graph $\psi:\graph\hookrightarrow I$,
		we call $\IAux(\psi) := \Anc_I(\psi(\nodes))\cap \Anc_I(\viewportEventual)$
		$\psi$-auxiliary indices.
	\end{definition}
	
	While it is not immediately obviously at a first glance,
	the basic query implicitly accounts for multi-level structure:
	this is because it is formulated on $\realizedDistr$, not
	on $\shallowDistr$. The hierarchical structure emerges in
	embeddings and structural queries from these basic queries,
	despite the absence of an explicit hierarchy in either
	the basic query or the model. 
	Compared to Def.\ \ref{def:structured_query},
	\emph{structured} queries in
	the multi-level case have to additionally account for
	auxiliary nodes:
	
	\begin{definition}[Structured Query]
		\label{def:structured_query:apdx}
		A structured query is an embedded local graph
		$\psi:\graph\hookrightarrow I$ together with a set of interest $Y \subset \nodesInnerProper$,
		such that $\ancestralStructure=\emptyset$ is valid
		(Def.\ \ref{def:ancestral_structure}) and
		\textbf{(multi-level case only:)}
		$\psi(\nodesOuterFixed)=\IAux(\psi)\cap\psi(\nodes\setminus\nodesOuterProper)$.
		
		\emph{Underlying Query:}
		Given a structured query $(\psi, Y)$,
		there is an underlying basic query $q(\psi,Y):=(\tilde{Y}=\psi(Y),
		\tilde{X}=\psi(\nodesOuterProper), \theta)$.
		We call $(\psi, Y)$ identifiable, if the underlying
		basic query is identifiable.
	\end{definition}
	
	Before connecting basic queries
	to structured queries in the multi-level case, we
	first need to make a small modification to associated structural graphs
	(for reasons explained below in Rmk.\ \ref{rmk:associated_structural_graph_query}).
	In §\ref{apdx:extraction_from_data}, for the multi-level
	extraction, we used an associated structural graph $\graphObs$
	(Def.\ \ref{def:associated_structural_graph}), the
	next definition is almost identical, except for the
	treatment it gives to observed pinned nodes:
	
	\begin{definition}[Associated Structural Query-Graphs]
		\label{def:associated_structural_graph_query}
		We are given
		a model-aligned local graph $(\graph, L)$,
		with $L\cap \nodesOuterProper = \emptyset$.
		
		We define a model-aligned structural graph	
		$\graphQuery(\graph)$
		with nodes $\nodes^\gQuery = \nodes$ split as
		\begin{equation*}
			\nodesInner^\gQuery
			:= \nodesInnerProper \cup \nodesOuterFixed
			\halfquad\txt,\qquad
			\nodesOuter^\gQuery
			:= \nodesOuterProper
		\end{equation*}
		and the edges $\edgesStructural^\gQuery = \parentStructure$.
		Then $\graphQuery$ is model-aligned with $\mu^n$ for $n\in \nodesInnerProper$
		and $\mu^x:=\delta_{i^x}\in\pinnedKernels$ for $x\in \nodesOuterFixed$.
	\end{definition}
	\begin{rmk}[Comparison of Associated Structural Graphs]
		\label{rmk:associated_structural_graph_query}
		Extraction and queries take a different perspective on
		\emph{observed} pinned nodes.
		
		\emph{Extraction:}
		We want to learn the most informative object from data.
		Given an observed pinned node, we can treat it as
		an outer node with the understanding that we
		know the associated $\delta_i=\delta[\anyVar_i]$
		since ($\anyVar_i$ is observed), thus we can plug in
		the value of $\anyVar_i$ in case we want to recover the
		context-specific result, or we may use the structured
		kernel like any other kernel. We will not have actually
		useful observational support (Rmk.\ \ref{rmk:kernels_uniqueness})
		for doing much else
		(other than plugging in the actual value of $\anyVar_i$) with it,
		but conceptually, since we are working under Ass.\ \ref{ass:support},
		it seems more logical to retain the general object
		rather than suddenly mixing support-arguments into our
		otherwise support-agnostic treatment.
		
		\emph{Query:}
		We have to produce a specific target kernel (Def.\ \ref{def:query}),
		and must fill in all information required to compute that target.
		The value of an observed $\anyVar_i$ is fixed in $\realizedDistr$,
		thus if it affects the query (is in its structural query),
		the information about the (known) value of $\anyVar_i$
		must be used, and the structural graph of the query cannot
		leave "open" additional external nodes. Thus
		and inner node must be aligned to $\delta_i$.
		
		\emph{Identification:}
		Note that for identification, we typically have to 
		glue (Lemma \ref{lemma:glue:mt}) the structural
		query graph from smaller (extracted, §\ref{sec:extraction_from_data})
		pieces. These pieces can contain a composition
		with $\delta_i$ (if they are extracted with pinned nodes,
		cf.\ Thm.\ \ref{thm:extract_multi_level}),
		but it is also legal to extract a (standard, no $\delta_i$)
		kernel with the pinned query-node as outer-node
		and glue $\delta_i$ (as a structured kernel whose
		graph has a single node,
		cf.\ $k(i)$ in Def.\ \ref{def:knowledge_set_extracted:apdx}),
		which is (directly) identifiable if $i\in\viewportEventual$,
		\ie if $\anyVar_i$ is observed
		(cf.\ Def.\ \ref{def:identification_direct:apdx}).
	\end{rmk}
	
	This more sophisticated structured query-formulation can
	capture additional basic queries, \ie the hypothesis
	in Lemma \ref{lemma:assoc_structured_query} can be weakened
	(we still exclude "queries causing observations"
	by enforcing $\tilde{Y} \cap \Anc_I(\viewportEventual) = \emptyset$,
	cf.\ Rmk.\ \ref{rmk:obs_before_query:apdx},
	but this is a much weaker condition than disconnecting $\tilde{Y}$ from
	$\viewportEventual$ in the $I$-graph):
	\begin{lemma}[Associated Structured Query]
		\label{lemma:assoc_structured_query:apdx:state}		
		Given an abstract query \query{},
		then 
		there is a structured query $(\psi(q),\tilde{Y})$
		with underlying query $q$,
		if and only if
		$\tilde{Y} \cap \Anc_I(\viewportEventual) = \emptyset$
		\textbf{(replacing the stronger disconnectivity-hypothesis
			of the single-level case)}
		and $\tilde{X}$ cannot be bypassed in the sense of:
		For $y\in\tilde{Y}$
		and $w \in \Anc_I(\tilde{X})\setminus(\tilde{X}\cup\psi(\nodesOuterFixed))$,
		there
		is no directed path $\gamma: w \rightsquigarrow y$ in the $I$-graph
		with $\gamma\cap(\tilde{X}\cup\psi(\nodesOuterFixed)) = \emptyset$
		\textbf{(multi-level: $\gamma$ now also does not pass through the
			previously empty $\psi(\nodesOuterFixed)$)}.
	\end{lemma}
	\begin{rmk}[Assumptions for Structured Representation]
		\label{rmk:obs_before_query:apdx} 
		The condition
		$\tilde{Y} \cap \Anc_I(\viewportEventual) = \emptyset$
		intuitively speaking means the query does not cause the observations.
		For "interventional" formulations (like the do-calculus)
		this is always true;
		one plausible use-case violating this assumption would be missing-value imputation,
		the problem we avoid by this simplification is to give a
		formal interpretation of what "imputation" should mean
		precisely (for example, should it use downstream observations
		to narrow down imputed distributions?).
	\end{rmk}
	
	Finally, in Lemma Thm.\ \ref{thm:query_id}
	the hypothesis has to be (and can be) weakened to $\delta_i$-identifiable inputs:
	\begin{lemma}
		\label{lemma:query_id:apdx}
		Given a structured query
		$(\psi,Y)$,
		using
		$L:=(\nodesInnerProper\setminus Y)\cup\nodesOuterFixed$
		\textbf{(multi-level: includes $\nodesOuterFixed$)},
		if there is a regular functional $F$ computing
		$\mu(\graphQuery(\graph),L)=F[\mu_1,\ldots,\mu_n]$
		with $\mu_1,\ldots,\mu_n$ all $\delta_i$-identifiable
		\textbf{(multi-level: replace identifiable by $\delta_i$-identifiable)},
		then the underlying query
		$q(\psi,Y)$ is identifiable.
	\end{lemma}

	\subsection{Structured Representation of Basic Queries}
	\label{apdx:query_structural_representation}

	The representation of basic queries by structured queries
	in the main text (Lemma \ref{lemma:assoc_structured_query})
	was phrased as a single-level simplification of
	(and immediately follows from) the following
	
	\begin{lemma}[Associated Structured Query]
		\label{lemma:assoc_structured_query:apdx}
		Given an abstract query \query{},
		then 
		there is a structured query $(\psi(q),\tilde{Y})$
		with underlying query $q$,
		if and only if
		$\tilde{Y} \cap \Anc_I(\viewportEventual) = \emptyset$
		and $\tilde{X}$ cannot be bypassed in the sense of:
		For $y\in\tilde{Y}$
		and $w \in \Anc_I(\tilde{X})\setminus(\tilde{X}\cup\psi(\nodesOuterFixed))$,
		there
		is no directed path $\gamma: w \rightsquigarrow y$ in the $I$-graph
		with $\gamma\cap(\tilde{X}\cup\psi(\nodesOuterFixed)) = \emptyset$.
	\end{lemma}
	\begin{proof}
		By slight abuse of notation, we will
		use injectivity of $\psi$ (Def.\ \ref{def:local_graph_embedding:apdx})
		to identify $\nodes$ with $\psi(\nodes)$.
		
		"$\Leftarrow$":
		We start by constructing a local graph $\graph$.
		First, the node-set is built inductively as a subset of $I$
		as follows:
		We start from $\nodesInnerProper^0 := \nodes^0 := \tilde{Y}$.
		Then
		\begin{align*}
			\nodes^{k+1} &:= \nodes^k \cup \Pa_I(\nodesInnerProper^k) \txt,
			\\
			\nodesOuterProper^{k+1} &:= \nodes^{k+1} \cap \tilde{X} \txt,
			\\
			\nodesOuterFixed^{k+1} &:=
			(\nodes^{k+1}\setminus\nodesOuterProper^{k+1})
			\cap \Anc_I(\viewportEventual) \txt,
			\\
			\nodesInnerProper^{k+1} &:=
			\nodes^{k+1}\setminus(\nodesOuterProper^{k+1}\cup\nodesOuterFixed^{k+1})
			\txt.
		\end{align*}
		By finiteness of $\tilde{Y}$ (Def.\ \ref{def:query})
		and finite past (Ass.\ \ref{ass:finite_past}),
		$|\Anc_I(\tilde{Y})|<\infty$, so
		we may define the finite sets $\nodes := \nodes^\infty := \cup_k \nodes^k$ etc.,
		and obtain a graph $\graph$ by adding edges
		from the $I$-graph, $\parentStructure :=
		\{(p,c)\in\nodes \times\nodesInnerProper| p\in\Pa_I(c)\}$.
		Then
		define $\psi: \nodes \hookrightarrow I, n \mapsto n$ as the (restricted)
		identity mapping,
		align to model-mechanisms $\mu^i = f_{J(i)}$
		and pin $x\in\nodesOuterFixed$ at $i^x = x$,
		finally define the node-set of interest as $Y:=\tilde{Y}$.
		Then by construction the underlying query is $q$.
		
		We have to check that
		(a) $Y\subset \nodesInnerProper$,
		(b) $\ancestralStructure=\emptyset$ is valid,		
		(c) $\psi(\nodesOuterFixed)=\IAux(\psi)\cap\psi(\nodes)$
		and that
		(d) $\psi: \graph \hookrightarrow I$ is indeed an embedded
		local graph (satisfies Def.\ \ref{def:local_graph_embedding}).
		
		\textbf{Part (a):}
		By Def.\ \ref{def:query} $\tilde{X}\cap\tilde{Y} = \emptyset$,
		thus $\nodesOuterProper\cap Y = \emptyset$.
		By hypothesis,
		$\tilde{Y} \cap \Anc_I(\viewportEventual) = \emptyset$
		thus $\nodesOuterFixed\cap Y = \emptyset$.
		So, indeed, $Y\subset \nodesInnerProper$.
		
		\textbf{Part (b):}
		We have to check Def.\ \ref{def:ancestral_structure}.
		Let $w\in\nodesInnerProper$ and $x\in\nodesOuterProper$
		be arbitrary.
		By contradiction. Assume there were a directed path
		$\gamma: w \rightsquigarrow x$
		in the $I$-graph not through $\graph$.
		In particular
		$w\in\Anc_I(\tilde{X})\setminus(\tilde{X}\cup\psi(\nodesOuterFixed))$.
		By $w\in\nodesInnerProper$ and by construction of $\nodesInnerProper$
		there is path $\eta:w\rightsquigarrow y$
		to $y\in\tilde{Y}$ with $\eta\subset\nodesInnerProper$
		(the chain of parents of $y$ that inductively were included
		in subsequent $\nodesInnerProper^k$).
		But $\eta$ also is a path as described in the statement
		of the lemma (it starts at 
		$w\in\Anc_I(\tilde{X})\setminus(\tilde{X}\cup\psi(\nodesOuterFixed))$
		and ends at $y\in \tilde{Y}$), thus
		by hypothesis $\eta\cap(\tilde{X}\cup\psi(\nodesOuterFixed))\neq\emptyset$,
		which contradicts $\eta\subset\nodesInnerProper$.
		
		\textbf{Part (c):}
		By
		construction $\nodesOuterFixed = (\nodes\setminus\nodesOuterProper) \cap \Anc_I(\viewportEventual)$
		and thus $\nodesOuterFixed = \IAux(\psi)\cap (\nodes\setminus\nodesOuterProper)$
		by definition (Def.\ \ref{def:aux}).		
		
		\textbf{Part (d):}
		First note that $\graph$ is indeed a local graph,
		\ie edges point only at proper nodes by construction.
		For Def.\ \ref{def:local_graph_embedding}, $\psi$ is
		the identity mapping, thus injective.
		(i) $\Pa_I(\nodesInnerProper)\subset\nodes$ by construction.
		(ii) Edges equal $I$-graph edges by construction.
		(iii) Model-alignment is to $f_{J(i)}$ and to $i^x=x=\psi(x)$
		by construction.
		
		"$\Rightarrow$":
		\textbf{Part 1 ($\psi(Y) \cap \Anc_I(\viewportEventual) = \emptyset$):}
		$Y\subset\nodesInnerProper$, but
		$\Anc_I(\viewportEventual)\cap\psi(\nodes) \subset \IAux(\psi)\cap\psi(\nodes)$
		(by Def.\ \ref{def:aux})
		and
		$\IAux(\psi)\cap\psi(\nodes)\subset \psi(\nodesOuterProper\cup\nodesOuterFixed)$
		(by Def.\ \ref{def:structured_query}),
		so by $\nodesInnerProper \cap (\nodesOuterProper\cup\nodesOuterFixed)
		=\emptyset$ (by Def.\ \ref{def:local_graph_new})
		$\psi(Y) \cap \Anc_I(\viewportEventual) = \emptyset$.
		
		\textbf{Part 2 (no bypassing of $\tilde{X}$):}
		By contradiction.
		Assume there were
		$y\in\tilde{Y}$
		and $w \in \Anc_I(\tilde{X})\setminus(\tilde{X}\cup\psi(\nodesOuterFixed))$, and
		a directed path $\gamma: w \rightsquigarrow y$ in the $I$-graph
		with $\gamma\cap(\tilde{X}\cup\psi(\nodesOuterFixed)) = \emptyset$.
		With proper parents being included
		Def.\ \ref{def:local_graph_embedding:apdx}~%
		\ref{def:local_graph_embedding:inner_parents_incl:apdx},
		and $\gamma\cap(\tilde{X}\cup\psi(\nodesOuterFixed)) = \emptyset$,
		all nodes on $\gamma$ (going backwards from $y$ step by step)
		must be in $\tilde{Y}$, in particular
		$w\in\tilde{Y}$,
		contradicting
		validity of $\mathcal{A}$.
	\end{proof}

	\subsection{Identification of Queries}
	\label{apdx:query_identification}
	
	Finally, for the identification of structured queries, we
	connect their structured kernels to the realized world,
	then confirm that regular computation works
	on $\delta_i$-identifiable objects as is,
	and then assemble both to a proof of Thm.\ \ref{thm:query_id}.
	
	Having fixed a structured graph associated
	to the structured query (Def.\ \ref{def:associated_structural_graph_query},
	note that in the single-level
	case $\graphQuery=\graph$), we start by associating
	it to the realized world distribution; this is actually already the
	main technical step of the query-formalism.
	
	\begin{lemma}
		\label{lemma:query_from_world}
		Given a structured query
		$(\graph,Y)$, then,
		using
		$L:=(\nodesInnerProper\setminus Y)\cup\nodesOuterFixed$,
		\begin{equation*}
			\realizedDistr(\anyVar_{\psi(Y)}|
			\anyVar_{\psi(\nodesOuterFixed)},
			\anyVar_{\psi(\nodesOuterProper)}=\tilde{x})
			\halfquad=\halfquad
			\big(\mu(\graphQuery,L)[\delta_i=\delta[\anyVar_i]]\big)_{x=\tilde{x}}\txt.
		\end{equation*}
		If $\exists \nu$ a probability kernel with
		$\mu(\graphQuery,L)[\delta_i=\delta[\anyVar_i]]=\nu$
		almost surely $\realizedDistr$,
		then
		\begin{equation*}
			\realizedDistr(\anyVar_{\psi(Y)}|
			\anyVar_{\psi(\nodesOuterProper)}=x)
			\halfquad=\halfquad
			\nu_x\txt.
		\end{equation*}
		\emph{Remark:}
		The last equation, by our notation of equality of kernels
		(cf.\ discussion in and around Rmk.\ \ref{rmk:kernels_uniqueness})
		is an equality almost everywhere, so the "almost surely" part did not
		magically vanish, it was just absorbed into our
		notation of kernel equality\Slash{}uniqueness.
	\end{lemma}
	
	\begin{proof}
		We will write $\tilde{Y}:=\psi(Y)$ and $\tilde{X}:=\psi(\nodesOuterProper)$.
		We start by using an embedded local graph analogue
		of simplification (Lemma \ref{lemma:simplify})
		which allows us to assume \asswlog
		$\psi(\nodes) \subset \Anc_I(\psi(Y))$.
		Applying Lemma \ref{lemma:simplify} with the unique maximal $B$
		for $L=(\nodesInnerProper\setminus Y)\cup\nodesOuterFixed$
		(see statement of the Lemma)
		we get $(\graphQuery)'$ with $\nodesInner'\subset \Anc_{\graph}(Y)$,
		and $\mu(\graphQuery,L)=\mu((\graphQuery)', L')$.
		The restriction of $\psi$ to $(\graphQuery)'$ is still a
		structured query.
		The hypothesis of the lemma is thus still satisfied, and
		the claimed result is the same, thus we may proof the Lemma for this
		new structured query instead. This new structured query
		satisfies
		$\psi(\nodes) \subset \Anc_I(\psi(Y))$.
		
		By finite past (Ass.\ \ref{ass:finite_past}),
		$\Anc_I(\psi(Y))$ is finite.
		We use Rmk.\ \ref{rmk:index_sets_simplify}
		to write $I=\mathbb{N}$ as totally ordered set
		and we may assume \asswlog that $\Anc_I(\psi(Y))=\{1,\ldots,m\}$.
		Next, we show that by $\ancestralStructure=\emptyset$ being valid
		the query-graph is not "interlaced" with any of its ancestors,
		that is we can further assume \asswlog
		$\psi(\nodesInnerProper)=\{m_0, m_0+1, \ldots, m-1, m\}$, \ie
		we can put $\psi(\nodesInnerProper)$ a the $\pi$-end of $\Anc_I(\psi(Y))$.
		To see this, note that by validity of $\ancestralStructure=\emptyset$
		and inclusion of proper parents Def.\ \ref{def:local_graph_embedding}~%
		\ref{def:local_graph_embedding:inner_parents_incl:apdx},
		nodes in $\Anc_I(\psi(\nodesInnerProper))=\Anc_I(\psi(\nodes))$
		cannot have arguments in any $\psi(\nodesInnerProper)$,
		so by sparsity (Lemma \ref{lemma:causal_transformations})
		we can move them all to the left of $\psi(\nodesInnerProper)$.
		
		By definition of $\shallowDistr$ (Def.\ \ref{def:shallow_distr}):
		\begin{equation*}
			\shallowDistr(\anyVar_{\Anc_I(\tilde{Y})})
			=
			\kernelCompound{\otimes_{k=1}^{m_0-1} f_k}
			\otimes
			\kernelCompound{\otimes_{k=m_0}^{m} f_k}
		\end{equation*}
		By uniqueness of disintegrations, with
		$\shallowDistr(\anyVar_{m_0,\ldots,m}|\anyVar_{1,\ldots,m_0-1})$
		satisfying the characterizing property of Lemma \ref{def:disint_product},
		and disintegrating $m_0-1$ times from the left:
		\begin{equation*}
		\shallowDistr(\anyVar_{m_0,\ldots,m}|\anyVar_{1,\ldots,m_0-1})
		=
		\kernelCompound{\otimes_{k=m_0}^{m} f_k}
		\txt.
		\end{equation*}
		Combining this with the definition of structured kernels
		(and model-alignment of $\psi$, Def.\ \ref{def:local_graph_embedding}~%
		\ref{def:local_graph_embedding:applicable})
		using the \emph{extraction}-associated graph
		$\graphObs(\graph, L\setminus\nodesOuterFixed)$, which has
		outer nodes $\nodesOuterProper \cup \nodesOuterFixed$
		(denoted with value $x'=(\tilde{x},d)$, where
		$\tilde{x}$ is associated to $\nodesOuterProper$ and
		thus $\tilde{X}$, while $d$ the pinned part related
		to the dataset),
		we thus already find
		\begin{align*}
			\mu
			:=\mu(\graphObs(\graph, L\setminus\nodesOuterFixed)_{x'=(\tilde{x},d)}
			\halfquad&=\halfquad
			\marginalize{
				\kernelCompound{\otimes_{k=m_0}^{m} f_k}
			}{L}
			\\&=\halfquad
			\shallowDistr(\anyVar_{\tilde{Y}}|\anyVar_{1,\ldots,m_0-1}=(\tilde{x},d))
			\txt.
		\end{align*}	
		Importantly, if we treat the model-aligned multi-level
		kernel $\mu(\graphQuery,L)[\delta_i]$ as a functional
		of its place-holder kernels $\delta_i$, then
		with all proper nodes aligned the same in 
		$\mu(\graphQuery,L)$ and $\mu(\graphObs(\graph, L\setminus\nodesOuterFixed)$,
		\begin{equation*}
			\tag{$*$}
			\mu:=
			\mu(\graphObs(\graph, L\setminus\nodesOuterFixed),L)_{x'=(\tilde{x},d)}
			\halfquad=\halfquad
			\big(\mu(\graphQuery,L)[\delta_i=\delta[d_i]]
			\big)_{x'=(\tilde{x})}\txt,
		\end{equation*}
		and it suffices to compute $\mu$,
		for $d_i=\anyVar_i$ (\ie plug the \emph{random variable} $\anyVar_i$
		into the argument $d_i$; this is the right-hand side
		of the claim of the Lemma).
		
		For structured queries $\psi(Y) \cap \Anc_I(\viewportEventual) = \emptyset$
		(Lemma \ref{lemma:assoc_structured_query:apdx} is an if and only if statement).		
		So by applying the Markov-property
		(Lemma \ref{lemma:obs_world_existence}):
		\begin{align*}
			\shallowDistr(\anyVar_{\tilde{Y}}|\anyVar_{1,\ldots,m_0-1})
			\halfquad=\halfquad
			\shallowDistr(\anyVar_{\tilde{Y}}|
			\anyVar_{1,\ldots,m_0-1},\anyVar_{\viewportEventual})
			\txt.
		\end{align*}
		On the other hand $\mu(\graphObs,L\setminus\nodesOuterFixed)$ has arguments
		only in $\psi(\nodesOuterProper\cup\nodesOuterFixed)$;
		by Def.\ \ref{def:local_graph_embedding}~%
		\ref{def:local_graph_embedding:inner_parents_incl}
		and \ref{def:local_graph_embedding:proper_edges_a} this is also
		true for $\otimes_{k=m_0}^{m} f_k$ in the definition of $\shallowDistr$.
		Thus by successive application of the Markov-property
		(Lemma \ref{lemma:obs_world_existence})
		\begin{align*}
			\shallowDistr(\anyVar_{\tilde{Y}}|\anyVar_{1,\ldots,m_0-1},\anyVar_{\viewportEventual})
			\halfquad=\halfquad
			\shallowDistr(\anyVar_{\tilde{Y}}|
			\anyVar_{\tilde{X}\cup\psi(\nodesOuterFixed)},\anyVar_{\viewportEventual})
			\txt.
		\end{align*}
		Putting the previous results together:
		\begin{equation*}
			\mu
			\halfquad=\halfquad
			\shallowDistr(\anyVar_{\tilde{Y}}|
			\anyVar_{\tilde{X}\cup\psi(\nodesOuterFixed)}=(x,d),
			\anyVar_{\viewportEventual}=\mathcal{D})
			\txt.
		\end{equation*}
		The right-hand-side, by Def.\ \ref{def:realized_world},
		is $\realizedDistr(\anyVar_{\tilde{Y}}|
		\anyVar_{\tilde{X}}=x)$, thus, plugging in $(*)$,
		we have proved the first claim of the lemma.
		(Note that, writing $P(Y|X=x)|_{x=X}$ for
		the random measure obtained by plugging the random
		variable $X$ into a regular version of $P(Y|X=x)$,
		we have $P(Y|X=x)|_{x=X}=P(Y|X)$;
		cf.\ \eg \citep[Thm.\ 5.3 (p.\,84)]{kallenberg1997foundations}).
		
		For the second part,
		we have
		$\realizedDistr(\anyVar_{\psi(Y)}|
		\anyVar_{\psi(\nodesOuterFixed)},
		\anyVar_{\psi(\nodesOuterProper)}=x)=\nu_x$
		almost surely $\realizedDistr$,
		where the right-hand side is a probability kernel
		(not a random measure or kernel) by hypothesis.
		The left-hand side is (see above,		
		cf.\ \eg \citep[Thm.\ 5.3 (p.\,84)]{kallenberg1997foundations})
		\begin{align*}
			&\realizedDistr(\anyVar_{\psi(Y)}|
			\anyVar_{\psi(\nodesOuterFixed)},
			\anyVar_{\psi(\nodesOuterProper)}=x)\\
			&\quad=\realizedDistr(\anyVar_{\psi(Y)}|
			\anyVar_{\psi(\nodesOuterFixed)}=d,
			\anyVar_{\psi(\nodesOuterProper)}=x)%
			|_{d=\anyVar_{\psi(\nodesOuterFixed)}(\omega)}
			\txt.
		\end{align*}
		In particular
		$\realizedDistr(\anyVar_{\psi(Y)}|
		\anyVar_{\psi(\nodesOuterFixed)}=t,
		\anyVar_{\psi(\nodesOuterProper)}=x)=\nu_x$		
		for $\realizedDistr(\anyVar_{\psi(\nodesOuterFixed)})$
		almost all $t$.
		Plugging this into
		\begin{align*}
			&\realizedDistr(\anyVar_{\psi(Y)}|
			\anyVar_{\psi(\nodesOuterProper)}=x)\\
			&=\int
			\realizedDistr(\anyVar_{\psi(Y)}|
			\anyVar_{\psi(\nodesOuterFixed)}=t,
			\anyVar_{\psi(\nodesOuterProper)}=x)
			\realizedDistr(\anyVar_{\psi(\nodesOuterFixed)}\in\dt)\\
			&=\nu_x
			\int 
			\realizedDistr(\anyVar_{\psi(\nodesOuterFixed)}\in\dt)\\
			&=\nu_x
			\txt,
		\end{align*}
		we obtain the claimed result;
		we denoted the integral of a measurable $f(t)$ in the variable
		$t$ over the measure $\realizedDistr(\anyVar_{\psi(\nodesOuterFixed)})$
		as $\int f(t) \realizedDistr(\anyVar_{\psi(\nodesOuterFixed)}\in\dt)$.
	\end{proof}

	\begin{lemma}
		\label{lemma:regular_computation_with_delta_i_id}
		Given a regular functional computing
		$\mu[\delta_i]=F[\mu_1,\ldots,\mu_n]$
		with $\mu_1,\ldots,\mu_n$ all $\delta_i$-identifiable,
		then $\mu$ is $\delta_i$-identifiable.
		Given representing
		$\nu_1, \ldots, \nu_n$ with $\nu_i = \mu_i$ almost surely
		$\realizedDistr$,
		then $\nu = F[\nu_1,\ldots,\nu_n]$ represents $\mu$
		almost surely $\realizedDistr$.
	\end{lemma}
	\begin{proof}
		By definition (\ref{def:delta_i_identified}),
		there are $\nu_1, \ldots, \nu_n$ with $\nu_i = \mu_i$ almost surely
		$\realizedDistr$ after plugging in $\delta_i(\omega)=\delta[\anyVar_i(\omega)]$.
		Let $\omega\in\Omega$ be arbitrary.
		With probability $1$, the finitely many
		$\nu_i = \mu_i(\omega)$ simultaneously.
		Thus $\mu(\omega)=F[\mu_1(\omega),\ldots,\mu_n(\omega)]=F[\nu_1,\ldots,\nu_n]$
		with probability $1$.
		In particular $\exists \nu := F[\nu_1,\ldots,\nu_n]$
		with $\mu = \nu$ almost surely
		$\realizedDistr$ after plugging in $\delta_i(\omega)=\delta[\anyVar_i(\omega)]$
		and $\mu$ is $\delta_i$-identified
		by definition (\ref{def:delta_i_identified}).
	\end{proof}

	\begin{CopyThm}{thm:query_id}{}
		Given a structured query
		$(\psi,Y)$, using
		$L:=(\nodesInnerProper\setminus Y)\cup\nodesOuterFixed$,
		if there is a regular functional computing
		$\mu(\graphQuery(\graph),L)=F[\mu_1,\ldots,\mu_n]$
		with $\mu_1,\ldots,\mu_n$ all $\delta_i$-identifiable,
		then $\mu(\graphQuery,L)$ and the underlying query
		$q(\psi,Y)$ are identifiable.
	\end{CopyThm}
	\begin{proof}
		We have to identify
		the target
		$\realizedDistr(\tilde{Y}|\tilde{X}=x)$
		of the underlying query.
		We apply Lemma \ref{lemma:query_from_world},
		using Lemma \ref{lemma:regular_computation_with_delta_i_id}
		to apply the second part,
		to obtain
		\begin{equation*}
			\realizedDistr(\anyVar_{\psi(Y)}|
			\anyVar_{\psi(\nodesOuterProper)}=x)
			\halfquad=\halfquad
			\nu_x\txt,
		\end{equation*}
		where $\nu = F[\nu_1,\ldots,\nu_n]$
		with $\nu_1,\ldots,\nu_n$ identifiable representatives
		of $\mu_1,\ldots,\mu_n$ (see Lemma
		\ref{lemma:regular_computation_with_delta_i_id}).
		In particular $\nu$ is identifiable,
		and hence so is
		$\realizedDistr(\tilde{Y}|\tilde{X}=x)$
		(using that by definition of the underlying query,
		$\tilde{X}=\psi(\nodesOuterProper)$, $\tilde{Y}=\psi(Y)$).
	\end{proof}

	\subsection{Beyond Basic Queries} 
	\label{sec:queries_beyond_basic}
	
	Basic queries contain conventional do- or soft-interventions
	as special cases, but are substantially more flexible and
	by their simple form easily extensible.
	
	\paragraph{Meta Queries:}
	
	Our approach is intentionally kept simple, describing only
	causal relations, no other preconceptions about data.
	By "meta-queries" we mean queries that stack other kinds of reasoning
	on top of our queries.
	
	This includes for example compiling the language of IID do-interventions
	to basic queries as described in example \ref{example:iid_do_intervention},
	or accounting for internal structure of singular interventions
	to reproduce "rule 2" conditioned interventions (§\ref{apdx:iid_conditional})
	and thus the results of \citep{shpitser2006identification1}.
	
	More generally, assumptions like stationarity of time-series
	can be leveraged: To ensure Ass.\ \ref{ass:finite_past} (finite past) on
	a \emph{stationary} $k$-Markov
	time-series (with not too complicated missingness structure),
	estimate the stationary distribution jointly over $k$ time-steps from data
	and prefix the observed time-interval in the $I$-graph
	by a "known" (in $\knownFunction=\theta$) $k$-step block (without ancestors),
	then in identification-formulas of queries, plug in the estimated
	stationary distribution for these "known" kernels (where they appear).
	
	Similarly, if a query $\tilde{Y}_x$ can also be formulated in form of a Markov-kernel
	(\ie $\tilde{X}$ is a "past" and $\tilde{Y}$ a "present" in the sense
	that $\tilde{Y}_{y_0} \otimes \tilde{Y}_{y_1} \otimes \tilde{Y}_{y_2}
	\otimes \ldots$ is a valid expression), one may post hoc want to know
	what stationary distributions of this process (with Markov-kernel $\tilde{Y}_x$)
	exist. \Ie it is (from a formal perspective)
	readily possible to ask for example for the new stationary distribution(s)
	\emph{after an intervention}; our formalism provides
	a kernel suitable to assess this fixed-point question.
	
	Consider also the following example (which also illustrates the
	relevance of the set $\tilde{X}$ in basic queries):
	\begin{example}[Query Argument]\label{example:query_argument}
		Consider a simple time-homogeneous (\ie there is only one model-variable\Slash{}kernel $f_X$), Markov (single time-lag)
		univariate time-series
		$\ldots\rightarrow X_{t-1} \rightarrow X_t
		\rightarrow X_{t+1} \rightarrow \ldots$ without latents
		(\ie $\viewport(N)$ contains $X_1, \ldots, X_N$).
		For time-series, different meaningful notions of interventions compete.
		For example we could replace $X_t$ by a do intervention $X_{\PearlDo}$
		(without parents)
		for a single $t$. At least if $\tilde{Y}$ contains only
		elements $X_{t'}$ with $t'>t$, we can obviously
		extract the kernel $X_{x^{\txt{past}}}$ -- note that it does at this point not 
		matter in our formalism if the time-series is stationary (up to support problems),
		\ie neither the joint "$P(X_{t-1},X_t)$" nor the marginals "$P(X_t)$" must
		be defined, cf.\ example \ref{example:random_walk},
		Rmk.\ \ref{rmk:kernels_uniqueness}
		-- then compute any finite future sub-sequence
		(thus elements of $\tilde{Y}$ jointly) as
		\begin{equation*}
			X_{\PearlDo} \otimes X_{x^{\txt{past}}} \otimes X_{x^{\txt{past}}}
			\otimes \ldots X_{x^{\txt{past}}}\txt.
		\end{equation*}
		Now, what about soft interventions? Or even intervention-free queries
		(we have extracted $X_{x^{\txt{past}}}$, but how do we ask for it as a query)?
		If $\tilde{f}_{x}$ replacing $X_t$ has a single parent (its immediate past),
		then even asking
		for future values of $\tilde{Y}$ makes no sense if the time-series is
		not stationary (indeed this is also \emph{not} a finite
		past query, Ass.\ 
		\ref{ass:finite_past}).
		The definition of query we employ allows in such cases to 
		ask the (actually meaningful) question about
		\begin{equation*}
			X'_{t-1} \otimes \tilde{f}_x \otimes X_{x^{\txt{past}}} \otimes X_{x^{\txt{past}}}
			\otimes \ldots X_{x^{\txt{past}}}
			\quad\txt{with}\quad
			\tilde{X} = \{ X'_{t-1} \}
			\txt.
		\end{equation*}
		We can reason about the effect of this soft intervention on
		future values $\tilde{Y}$, if we can reason about
		$X'_{t-1}\in\knownFunction$.
		For example if we believe our intervention happens in a
		stationary time-series,
		if our observations are also stationary simply intervene
		to the distribution over observations;
		the last point need not be true however:
		we could extract $X_{x^{\txt{past}}}$ also from
		a non-stationary time-series (\eg alternating $X_{x^{\txt{past}}}$
		and another kernel two-periodically).
		In this case, while we \emph{can} reason about
		$X'_{t-1}\in\knownFunction$, we may not actually care about
		any specific choice of $X'_{t-1}$, and in this case
		can use $\tilde{X}=\{X'_{t-1}\}$ instead to query for
		a "partial" result $X_{x^{\txt{past}}}$.
		Then $X'_{t-1}$ should be an invariant
		state of this (partial) result $X_{x^{\txt{past}}}$,
		in particular it can (at least in principle) be estimated from
		this first query. In a second query (the one from before),
		we can then use this estimated $X'_{t-1}$
		as an intervention.
	\end{example}
	
	\paragraph{Nested Queries:}
	
	We focused on queries (and graphical operations) with structured kernels
	as results.
	The authors are not aware of queries outside of this category
	having been considered in the literature, however, such queries
	exist, can be identifiable, and might make sense.
	It should, with the same graphical operations plus "disintegration nesting"
	(align nodes to disintegrations of other queries),
	be possible to approach problems like this:
	\begin{example}
		\label{example:nested_query}
		We are given system $A$ which describes
		an ecosystem of micro-organisms in a Petri-dish
		(with multiple causally interrelated biological variables).
		System $B$ is an experimental protocol that leads to
		samples (full ecosystems/Petri-dishes) being discarded
		or kept\Slash{}multiplied for future test-generations
		based on some of the measured properties of $A$.
		If we want to describe the long term result of repeatedly
		generating time-series from $A$, then selecting to initial states from 
		$B$ and again propagating by $A$ many times,
		our query has disintegrated (from selection) sub-systems of
		$A$ as sub-graphs (or "sub-structured" nodes).
	\end{example}

	\subsection{Examples for Multi-Level Queries}
	\label{apdx:multi_level_examples}
	
	A multi-level description allows for very fine-grained
	query specification.
	We start from a modification of example \ref{example:multi_level_illustrate}.
	
	\begin{figure}[ht]
		\colorlet{dblue}{blue!50!black}
\colorlet{dgreen}{green!50!black}
\colorlet{dyellow}{yellow!50!black}
\colorlet{dpurple}{purple!50!black}
\colorlet{dorange}{orange!50!black}
\colorlet{c4}{dblue!20!white}
\colorlet{c3}{dgreen!20!white}
\colorlet{c1}{dyellow!20!white}
\colorlet{c2}{dpurple!20!white}
\colorlet{c5}{dorange!20!white}
\colorlet{dgray}{gray!50!black}
\colorlet{c0}{dgray!20!white}
\colorlet{b0}{c0!40!white}
\colorlet{b5}{c5!40!white}
\colorlet{b3}{c3!40!white}
\colorlet{b4}{c4!40!white}

\pgfdeclarelayer{bg}
\pgfsetlayers{bg,main}

\begin{minipage}{\textwidth}
	\centering

\begin{tikzpicture}[xscale=1.25]
	
	\draw (0,0) node (C1)[rectangle, fill=c0, rounded corners, inner sep=0.4em,
	draw=dgray, thick, dashed]
	{$C_{\txt{sediment}}$};
	
	\draw (2,0) node (C2)[rectangle, fill=c0, rounded corners, inner sep=0.4em,
	draw=dgray, thick, dashed]
	{$C_{\txt{geom}}$};
	
	
	\draw [->, very thick] (C1) -- (C2);
	
	\foreach \s in {1,...,4}{		
		\draw (-1+1.5*\s, 1.5)
			node (p\s) [outer sep=0, inner sep=0] {\scalebox{0.4}[0.4]{
		\begin{tikzpicture}[xscale=1.25]			
			\draw (0,1) node (X\s)[rectangle, fill=c3, rounded corners, inner sep=0.4em,
			draw=dgreen,	thick]
			{$X$};
			\draw (0,0) node (Y\s)[rectangle, fill=c3, rounded corners, inner sep=0.4em,
			draw=dgreen,	thick]
			{$Y$};
			
			\draw [->, very thick] (X\s) -- (Y\s);	
		\end{tikzpicture}
		}};	
	}
	\draw[->, thick] (C2.north) -- ([yshift=-0.1em] p1.south);
	\draw[->, thick] (C2.north) -- ([yshift=-0.1em] p2.south);
	\draw[->, thick] (C2.north) -- ([yshift=-0.1em] p3.south);
	\draw[->, thick] (C2.north) -- ([yshift=-0.1em] p4.south);
	\draw ([xshift=1cm] p4.south east) node {$\cdots$};

	\draw (1,-1) node (Xa)[rectangle, fill=c2, inner sep=0.4em,
		draw=dpurple,	thick]	{$X'$};
	\draw (2,-1) node (Ya)[circle, fill=c4, inner sep=0.2em,	draw=dblue,	thick]
		{$Y_1$};
	\draw[->, thick] (Xa) -- (Ya);
	\draw[->, thick] (C2) -- (Ya);

	\draw (0,-2) node (C2b)[rectangle, fill=c0, rounded corners, inner sep=0.4em,
	draw=dgray, thick, dashed]
	{$C_{\txt{geom}}$};
	\draw (-1,-3) node (Xb)[rectangle, fill=c2, inner sep=0.4em,
	draw=dpurple,	thick]	{$X'$};
	\draw (0,-3) node (Yb)[circle, fill=c4, inner sep=0.2em,	draw=dblue,	thick]
	{$Y_2$};
	
	\draw[->, thick] (Xb) -- (Yb);
	\draw[->, thick] (C2b) -- (Yb);
	\draw[->, thick] (C1) -- (C2b);

	\draw (4,-2) node (C1c)[rectangle, fill=c0, rounded corners, inner sep=0.4em,
	draw=dgray, thick, dashed]
	{$C_{\txt{sediment}}$};
	\draw (6,-2) node (C2c)[rectangle, fill=c0, rounded corners, inner sep=0.4em,
	draw=dgray, thick, dashed]
	{$C_{\txt{geom}}$};
	\draw (5,-3) node (Xc)[rectangle, fill=c2, inner sep=0.4em,
	draw=dpurple,	thick]	{$X'$};
	\draw (6,-3) node (Yc)[circle, fill=c4, inner sep=0.2em,	draw=dblue,	thick]
	{$Y_3$};
	
	\draw[->, thick] (Xc) -- (Yc);
	\draw[->, thick] (C2c) -- (Yc);
	\draw[->, thick] (C1c) -- (C2c);
\end{tikzpicture}

\end{minipage}
		\caption{Different multi-level queries for the same system, see
		example \ref{example:multi_level_queries}.
		Green (solid outline, rounded corner rectangles)
		are observations a single site (other observed sites are not shown).}
		\label{fig:multi_level_queries}
	\end{figure}
	
	\begin{example}\label{example:multi_level_queries}
		We observe many river sites (water-levels, throughput, \dots),
		each over many time-steps.
		Each river-site has a context $C_{\txt{sediment}}$ "sediment type".
		This context is extremely complex (there are a lot of aspects
		to the actual sediment that affect to form of a resulting
		river-bed geometry and other relevant parameters)
		and we thus consider $C_{\txt{sediment}}$ unobserved.
		However, $C_{\txt{sediment}}$ will not change over time (at least not
		within realistic time-frames of observation).
		There is a distribution\Slash{}prior $P_{\txt{sediment}}$ 
		that describes how likely different types of sediment among
		the river-sites we observe are; we do not know $P_{\txt{sediment}}$,
		but we fix notation here for later reference and interpretation.
		
		Further we assume there is a context $C_{\txt{geom}}$ describing
		the "geometry" of a riverbed.
		This context is extremely complex (there are a lot of aspects
		to the actual geometry that are relevant to water-flow,
		like the distribution over cross-sections encountered
		but also the relation between cross-sections along the river,
		roughness measures and many more)
		and we thus consider $C_{\txt{geom}}$ unobserved.
		However, $C_{\txt{geom}}$ will not change over time (at least not
		within realistic time-frames of observation).
		There is a distribution\Slash{}prior $P_{\txt{geom}|\txt{sediment}}$,
		over geometries encountered given a sediment-type; this is
		a probability-kernel from sediment-types to geometries.
		Again, we do not know this kernel.
		
		For each river-site, there are \emph{time-dependent}
		variables like water-throughput $X$ and water-level $Y$.
		For the sake of argument, assume that water-throughput
		$X$ is driven externally (throughput must match influx over
		time-scales where rising water-levels and "water-capacity" at
		the site are irrelevant; for simplicity let's assume
		samples far enough apart in time, that they are effectively IID
		for each site)
		and that water-levels $Y$ depends on $X$ and $C_{\txt{geom}}$ only,
		in particular the kernel\Slash{}mechanism $Y_{x,g}$ is assumed
		to be the same across all river-sites (site-dependent
		aspects are captured by $C_{\txt{geom}}$ in good enough approximation);
		$X$ may vary by site (this will not matter for the purpose
		of this example).
		
		We ask a question about $Y$, when changing $X$ (for example by
		diverting part of the water-influx to a channel or large-scale reservoir).
		Indeed we ask not one question, but rather three questions
		that are similar, but practically important to distinguish.
		\begin{enumerate}[label=\arabic*)]
			\item 
				Enforcing influx\Slash{}throughput $X'$
				at a particular site $s$,
				how will $Y$ change at that site $s$?
			\item 
				Enforcing influx\Slash{}throughput $X'$
				at a particular previously unobserved site $s'$
				that is geographically very close to an observed
				site $s$ and can plausibly be assumed to share
				geological sediment properties $C_{\txt{sediment}}$
				with site $s$, how will $Y$ change\Slash{}behave
				at that site $s'$?
			\item 
				Enforcing influx\Slash{}throughput $X'$
				at a particular previously unobserved site $s''$,
				how will $Y$ change\Slash{}behave
				at that site $s''$?
		\end{enumerate}
		The actual question is always the same,
		so we should be more precise as to what we mean by
		"how will $Y$ change at that site $s$?".
		Indeed we will typically want to know,
		what is the "most narrow" (smallest noise)
		prediction we can make.
		There are both a finite-sample aspect
		and an aspect of the precise formulation of the question
		(and infinite-sample answer) to this question.
		Figure \ref{fig:multi_level_queries} illustrates these
		three different queries for $Y_1,Y_2, Y_3$ (in this order);
		the mechanism at each $Y_i$ is $Y_{x,g}$ as in the observations,
		the additional index is just for book-keeping.
		
		Note, that these queries are identifiable under different conditions:
		For $Y_1$, we need many time-points in context $s$ to
		learn the structural query where $C_{\txt{geom}}$ is pinned;
		for $Y_2$, we either need to observe $C_{\txt{geom}}$
		(as variable, with distribution $C_{\txt{geom}} \circ C_{\txt{sediment}}$)
		across many sites or we need to
		observe many sites (not just $s$, $s'$) with "shared"
		$C_{\txt{sediment}}$ (while not shown in Fig.\ \ref{fig:multi_level_queries}
		this could be read off the $I$-graph, where a single $C_{\txt{sediment}}$
		occurrence\Slash{}instance may have many $C_{\txt{geom}}$-children)
		to learn
		$Y_{x,[g]} \circ C_{\txt{geom}}$ with $C_{\txt{sediment}}$ pinned
		(the argument of $C_{\txt{geom}}$);
		for $Y_3$ we need to learn
		$Y_{x,[g]} \circ C_{\txt{geom}} \circ C_{\txt{sediment}}$,
		which requires many sites
		(note that our "freeness" condition means we use only one sample
		for $Y$ for each site, which is likely not a good strategy in the
		finite sample case; however one has to ensure to suitably reweigh
		different sites if they have different numbers of time-points).
		
		In principle, staying more local is more precise:
		For example the variance of $Y_1$ is essentially the variance of
		noise at $Y$, while the variance of $Y_2$ (or $Y_3$)
		additionally captures variance from $P_{\txt{geom}|\txt{sediment}}$
		(or from both $P_{\txt{sediment}}$ and $P_{\txt{geom}|\txt{sediment}}$);
		\ie priors are relevant to the interpretation of resulting distributions.
		However, as pointed out above, conditions for identification are different.
		Further, in the finite-sample case, the actual variance of an estimator
		depends on both the variance of the true predicted distribution
		but also on finite-sample error -- and the number of available
		samples favors the opposite order ($Y_3$ can use, somewhat over-simplified,
		more data-points than $Y_2$ than $Y_1$),
		so these non-trivial trade-offs have to be considered in practice.
		
		Finally the mathematically inclined reader may have noticed an interesting
		detail: If we observe multiple $Y$ (time-points)
		at the new site, then they are IID for query $Y_1$,
		but only exchangeable (conditionally IID given $C_{\txt{geom}}$
		for that particular site) for queries $Y_2$, $Y_3$.
		\Ie $\realizedDistr((Y_2)_{t_0}, \ldots, (Y_2)_{t_n})$ is \emph{not}
		IID (but exchangeable), while
		$\realizedDistr((Y_1)_{t_0}, \ldots, (Y_1)_{t_n})$
		is IID. 
		This subtlety is correctly captured by our query and identification
		formalism:
		The query (for $(Y_2)_{t_0}, \ldots, (Y_2)_{t_n}$) in this case
		has a single $C_{\txt{geom}}$ node $i_c\in I$ (for the new site $s'$)
		which is a child of a $C_{\txt{sediment}}$ node shared with $s$
		(but $s$ has a different $C_{\txt{geom}}$ node).
		This single node $i_c$ in turn has $n+1$ children
		$(Y_2)_{t_0}, \ldots, (Y_2)_{t_n}$.
		Above we sketched two possible identification-strategies for
		$Y_2$: (1) "observe" (somehow) $C_{\txt{geom}}$
		to learn $\bullet \rightarrow C_{\txt{geom}}$
		(where $\bullet$ is a pinned node for the shared $C_{\txt{sediment}}$)
		and $\circ \rightarrow Y \leftarrow \circ$ (with two outer external nodes
		for $X$ and the now assumed observed $C_{\txt{geom}}$).
		The query (see above) can immediately be glued from these:
		In practice this means, to jointly
		draw from $\realizedDistr((Y_2)_{t_0}, \ldots, (Y_2)_{t_n})$,
		first draw from $\bullet \rightarrow C_{\txt{geom}}$ once,
		then plug the result into $\circ \rightarrow Y \leftarrow \circ$ for
		all $n+1$ observations. 
		To draw again from
		$\realizedDistr((Y_2)_{t_0}, \ldots, (Y_2)_{t_n})$, draw a new
		value from $\bullet \rightarrow C_{\txt{geom}}$ (once) and so on.
		(2) given many sites sharing $C_{\txt{sediment}}$,
		learn $\bullet \rightarrow [C_{\txt{geom}}] \rightarrow Y \leftarrow \circ$
		(where $C_{\txt{geom}}$ is hidden, indicated by square-brackets),
		for a single $Y$, plug in $X'$.
		However, this second strategy does \emph{not} apply for the joint
		$\realizedDistr((Y_2)_{t_0}, \ldots, (Y_2)_{t_n})$,
		because the single node for $C_{\txt{geom}}$ in the query would
		enforce an overlap which is not of the form allowed by
		Lemma \ref{lemma:glue:mt} (gluing;
		formally we work on $\graphObs$ with the pinned node
		aligned to some $\delta_i$, however the query has a c-component
		containing the hidden $C_{\txt{geom}}$ and its children
		$(Y_2)_{t_0}, \ldots, (Y_2)_{t_n}$, none of
		our learned structures contains this full c-component as subgraph).
		Indeed from
		$\bullet \rightarrow [C_{\txt{geom}}] \rightarrow Y \leftarrow \circ$
		we only know how to draw $Y$ once before redrawing $C_{\txt{geom}}$.
		
		Small issues like the last one are often difficult to even spot.
		In our graphical operations and query perspective, they
		become surprisingly evident.
	\end{example}
	
	As this example shows, there are subtly different,
	but meaningful, queries that are formally distinguished by
	our formalism.
	These queries differ in how the query-part of the $I$-graph is attached to the
	observational part of the $I$-graph, in particular this distinction cannot
	be made in the single-level case (where query and observational
	part of the $I$-graph are disconnected,
	cf.\ Lemma \ref{lemma:assoc_structured_query}).

	\section{The IID Case}\label{apdx:iid}
	
	We compare our results to well-established standard results for the IID case.
	This enables us to validate that answers obtained from our formalism
	are indeed sound. Further, while we do not show that our algorithms
	are complete in the general case, they
	turn out to be complete in IID settings.

	\subsection{Standard Results for the Single-Context Case}
	\label{apdx:iid:standard_results}
	
	Some of the assumptions implicit in this notation become clearer in
	the next subsection §\ref{apdx:translate_IID_setup} where they
	appear as explicit constraints on the form of the model (in the sense of
	Def.\ \ref{def:model}).
	We start from the standard formulation via structural causal models (SCM) \citep{PearlBook,Elements}:
	
	\begin{notation}[Variables]
		\label{notation:iid:vars}
		For some finite index set $\IVars\times\ISample$
		(typically $\ISample=\mathbb{N}$),
		$\IVars$ is a finite set,
		fix a set of endogenous random variables
		$\{\iidVarAny_{v,s}\}_{(v,s)\in \IVars\times\ISample}$,
		taking values in $\val{X}_v$, that is
		measurable mappings $\iidVarAny_{v,s}: \Omega\rightarrow\val{X}_v$.
		We assume these are IID (in $\ISample$), \ie $\forall s,s'\in\ISample$:
		$(*)$ $(\iidVarAny_{v,s})_{v\in\IVars} \overset{\txt{d}}{=}
		(\iidVarAny_{v,s'})_{v\in\IVars}$ are (jointly) equal in distribution
		and for different $s\in\ISample$ are (jointly) independent,
		\ie
		\begin{align*}
			P(\{\iidVarAny_{v,s}\}_{(v,s)\in \IVars\times\ISample})
			\halfquad&=\halfquad
			\prod_{s\in\ISample}
			P(\{\iidVarAny_{v,s}\}_{v\in \IVars})
			\\
			&\overset{(*)}{=}\halfquad
			\prod_{s\in\ISample}
			P(\{\iidVarAny_{v}\}_{v\in \IVars})
			\txt,
		\end{align*}
		by slight abuse of notation (we just drop the index $s$ mandated by $(*)$);
		as is common in the literature we will (justified by IIDness) talk about
		variables $\{\iidVarAny_{v}\}_{v\in \IVars}$ (omitting an index $s\in\ISample$).
		We will sometimes use indices $v\in \IVars$ and their
		associated random variable $\iidVarAny_v$ interchangeably,
		this should not lead to confusion (it is equivalent to saying
		the $\IVars$, which as a finite set is defined essentially by the number
		of its elements \emph{is} the set containing as elements these
		random variables).
	\end{notation}
	\begin{definition}[SCM]
		\label{def:iid:scm}
		An SCM $M$ consists of the following information:
		
		Fix an index-set $\IVars\times\ISample$
		and an observed subset $\mathcal{O} \subseteq \IVars$
		with complement $L := \IVars\setminus\mathcal{O}$.
		Using notation~\ref{notation:iid:vars}, there are (IID)
		variables $\{\iidVarAny_{v}\}_{v\in \IVars}$
		and (IID) exogenous noises (hidden) $\{\eta_v\}_{v\in \IVars}$,
		taking values in $\val{N}_v$.
		Noises are jointly independent also within a sample, \ie
		we additionally (to IID) have
		$P(\{\eta_v\}_{v\in \IVars})=\prod_{v\in\IVars}P(\eta_v)$.
		
		Further, for each $v \in \IVars$ there is a set of parents
		$\Pa_v \subset \IVars \setminus \{v\}$ and a mechanism
		(a measurable mapping)
		$f_v : \val{X}_{\Pa_v} \times \val{N}_v \rightarrow \val{X}_v$
		such that the endogenous variables satisfy the structural equations
		\begin{equation*}
			\iidVarAny_v := f_v(\iidVarAny_{\Pa_i}, \eta_i)
		\end{equation*}
		relative to their parents and noise-term.
		Parent-sets are assumed to satisfy a suitable minimality condition,
		\eg \citep[Def.\ 2.6]{BongersCyclic}, that is $f_v$ are not trivial in
		any argument.
		
		We assume SCMs are uniquely solvable (for example acyclic),
		meaning the distribution of endogenous variables
		is fixed uniquely by noise-terms and mechanisms.
		
		We assume (in the IID case, this is essentially
		w.\,l.\,o.\,g., cf.\ Rmk.\ \ref{rmk:latent_proj})
		each $l\in L$ has no parents $\Pa_l=\emptyset$
		and at most two children (usually we start from
		$M$ where there are always exactly two children, but after a hard-intervention,
		see below, there might be less children), and its
		children are observed.
	\end{definition}
	\begin{rmk}[Latent Projections]
		\label{rmk:latent_proj}
		The form of latent-structure in assumed in the previous
		definition is standard in the literature, see \eg
		\citep[Def.\ 2.6.1 (p.\,52)]{PearlBook},
		because it simplifies graphical representation and
		logic, and a model can also be transformed into this
		form, see \eg \citep[Thm.\ 2.6.2 (p.\,52)]{PearlBook}.
		This is however primarily a consequence of the strong restrictions
		IIDness puts on symmetries of the model.
		In general, there does not seem to be any simple "standard form"
		in a comparable sense.
		Note that our approach makes latents always explicit, also
		in IID models, for example in Fig.\ \ref{fig:intro}
		it is very clear that $(Z_2)_{x,[l_2]} \circ L_2$
		(the second blue-box labeled "learn") does depend on $L_2$
		and can only be transferred to a query, where the latent
		has the \emph{same distribution}. For transfer between contexts,
		it is important to be transparent about such requirements.
	\end{rmk}
	\begin{definition}[Interventions]
		\label{def:iid:interventions}
		Interventions (or "actions") are defined as follows:
		Given an SCM $M$, an intervened on $\tilde{X}\subset\mathcal{O}$
		SCM $M'$ is an SCM
		with the same exogenous noises, and parent-sets $\Pa'_v$, mechanisms $g_v$
		such that
		for $v\in\IVars\setminus\tilde{X}:$ $\Pa'_v=\Pa_v$ and $g_v=f_v$.
		An intervention is called hard if
		for $v\in\tilde{X}:$ $\Pa'_v=\emptyset$.
		A hard intervention is a do-intervention,
		if for $v\in\tilde{X}:$ $g_v=x_v$ is a constant.
		We denote a do-intervention on $\tilde{X}$ by Pearl's do-operator
		$\PearlDo(\tilde{X}=x)$, for example
		the interventional distribution
		of $\tilde{Y}\subset\mathcal{O}$
		conditioned on ("in a given context" in \citep{shpitser2006identification1})
		$\tilde{Z}\subset\mathcal{O}$
		\begin{equation*}
			P(\tilde{Y}|\PearlDo(\tilde{X}=x), \tilde{Z})
			\halfquad:=\halfquad
			P(\{\iidVarAny'_v\}_{v\in\tilde{Y}}|\{\iidVarAny'_v\}_{v\in\tilde{Z}})
			\txt,
		\end{equation*} 
		where $\iidVarAny'_v$ are the endogenous variables of the model $M'$
		(since we assumed SCMs are uniquely solvable, this is well-defined).
		We also write $M^{\PearlDo(\tilde{X}=x)}$ or simply
		$M^{\PearlDo}$ (if clear from context) in this case.
		We will assume that the intervened value is in the value-space
		$x\in\val{X}_{\tilde{X}}$ (\ie we do not intervene
		\eg categorical variables to arbitrary real values).
	\end{definition}
	\begin{definition}[Causal Graph]
		The causal graph including latents $\graphClassical^L$ of a SCM $M$ is
		defined as the directed graph with nodes
		$\nodesClassical^L = \IVars$ and a directed edge $x\rightarrow y$
		in $\edgesClassical^L$
		if $x\in\Pa_y$.
		The corresponding graph of the do-interventional model $M^{\PearlDo}$
		is denoted $\graphClassicalDo^L$, and can also be obtained graphically
		by "amputating" all incoming edges to nodes in $\tilde{X}$.
		
		Define $\graphClassical$ as the graph with nodes
		$\nodesClassical = \mathcal{O}$ and a directed edge
		$x\rightarrow y$ in $\edgesClassical$
		if $x\in\Pa_y$
		and a bi-directed edge $x\leftrightarrow y$ in $\edgesClassical$
		if $L\cap\Pa_x\cap\Pa_y\neq\emptyset$.
		
		Causal graphs including latents $\graphClassical^L$ will typically be DAGs
		(directed \emph{acyclic} graphs) below, in this case
		the graph $\graphClassical^L$ is called Markovian
		and $\graphClassical$ is called semi-Markovian.
	\end{definition}
	
	For this section we adapt:
	\begin{definition}[Causal Effect Identifiability]
		\label{def:iid:identify}
		Quoting \citep[Def.\ 2 (p\,3)]{shpitser2006identification1}
		with only notation adapted and comments in square brackets added:
		
		The causal effect of an action $\PearlDo(\tilde{X})=x$
		on a set of variables $\tilde{Y}$
		in a given context $\tilde{Z}=z$ such that
		$\tilde{X}, \tilde{Y}, \tilde{Z}$ are disjoint
		[subsets of $\mathcal{O}$]
		is said to be identifiable from $P$
		[$P$ is the joint distribution $P(\{\iidVarAny_v\}_{v\in\mathcal{O}})$]
		in $\graphClassical$, if
		$P(\tilde{Y}|\PearlDo(\tilde{X}=x), \tilde{Z})$
		is (uniquely) computable from $P$ in any causal model which
		induces $\graphClassical$.
	\end{definition}
	\begin{rmk}
		\label{rmk:IID:computable}
		Here computable (seems to) mean "uniquely determined",
		\ie given the causal graph $\graphClassical$ and any
		observed distribution $P$, then
		any two SCMs that produce the same $P$
		and have causal graph $\graphClassical$ will induce the same interventional
		distribution.
		
		Here "any SCM" puts us in the generic case, \ie internal structure
		(like linearity) cannot be used (as there could still be a non-linear
		SCM with the same graph and observed distribution).
		
		Non-identifiability (and thus completeness) is then
		proved by counter-example (\ie given some property
		of $\graphClassical$, show there are two SCMs with the same
		observational distribution and different interventional distribution).
		
		This last point is relevant:
		We are talking about SCMs compatible with $\graphClassical$
		and producing the same, but \emph{any} observed distribution,
		not \emph{the} observed distribution.
	\end{rmk}

	\citep{shpitser2006identification1,shpitser2006identification2} implicitly%
	\footnote{Unfortunately neither of both references seem to define this notion,
		however in \citep{shpitser2006identification2}, in the example given in
		the paragraph immediately above Thm.\,4 (p.\,1223) it becomes evident
		that this is what their definitions assume.}
	use the following definition of "subgraphs"
	\begin{definition}
		A edge-subgraph 
		$\graphClassical' \approxsubset \graphClassical$ is a
		graph on the same nodes $\nodesClassical'=\nodesClassical$
		but with
		a \emph{subset} of edges $\edgesClassical'\subset\edgesClassical$.
		We call a edge-subgraph of a (node-)subgraph a sparse subgraph.
	\end{definition}
	\begin{rmk}
		We use, throughout this paper,
		the convention that ancestors $\Anc(X)$ of $X$ include $X$ as
		an element $X\in\Anc(X)$.
		The definitions below and in \citep{shpitser2006identification1,
			shpitser2006identification2} also follow this convention.
		However some care has to be taken, as in places
		(for example the definition of a root set
		\citep[p.\,1220, left column]{shpitser2006identification2}
		and \citep[p.\,2, right column]{shpitser2006identification1})
		the given references deviate from this convention.
	\end{rmk}
	
	\begin{definition}[C-Component]
		\citep[Def.\,3 (p.\,1221)]{shpitser2006identification2},
		see also \citep{tian2002general}:
		Given a semi-Markovian graph $\graphClassical$
		such that a subset of its bi-directed arcs forms a spanning
		tree over all vertices in $\graphClassical$.
		Then $\graphClassical$ is called a c-component.
		\\
		\emph{Remark:} This is analogous (and the inspiration for)
		our definition of c-connectivity (on $\graphClassical^L$),
		\ie here a "c-component"
		is a graph with a single c-component in our notation.
		Where this would lead to confusion, we clarify what is meant;
		usually the wording "$\graphClassical$ is a c-component"
		is used, which seems rather unambiguous in either nomenclature.
	\end{definition}
	
	\begin{definition}[C-Tree]
		\citep[Def.\,4 (p.\,1221)]{shpitser2006identification2}:
		Given a semi-Markovian graph $\graphClassical$
		such that $\graphClassical$ is a c-component,
		all observable nodes have at most one child,
		and there is a node $Y\in\nodesClassical$ such that
		$\Anc_{\graphClassical}(Y)=\nodesClassical$.
		Then $\graphClassical$ is a $Y$-rooted c-tree.
		\\
		\emph{Remark:}
		Usually we are only interested in whether there is a
		edge-subgraph that is a c-tree.
	\end{definition}
	\begin{definition}[C-Forest]
		\label{def:iid:c_forest}
		\citep[Def.\,5 (p.\,1222)]{shpitser2006identification2}:
		Given a semi-Markovian graph $\graphClassical$
		and $R\subset\nodesClassical$ (called the root set) such that
		$\Anc_{\graphClassical}(R)=\nodesClassical$.
		Then $\graphClassical$ is a $R$-rooted c-forest,
		if $\graphClassical$ is a c-component,
		and all observable nodes have at most one child.
		\\
		\emph{Remark:}
		Usually we are only interested in whether there is a
		edge-subgraph that is a c-forest.
	\end{definition}
	
	\begin{definition}[Hedge]
		\label{def:iid:hedge}
		\citep[Def.\,6 (p.\,1223)]{shpitser2006identification2}:
		Given a semi-Markovian graph $\graphClassical$.
		Let $X, Y \subset\nodesClassical$ with $X\cap Y=\emptyset$.
		Let $F,F'$ be $R$-rooted c-forests
		[both are sparse subgraphs of $\graphClassical$] such that
		(for their node-sets which are subsets of the node-set of $\graphClassical$)
		$F\cap X\neq \emptyset$,
		$F'\cap X=\emptyset$, $F'$ is a sparse subgraph of $F$,
		and $R\subset \Anc_{\graphClassicalDo}(Y)$.
		Then $F$ and $F'$ form a hedge for $(X, Y)$ in $\graphClassical$.
	\end{definition}

	\begin{lemma}[Hedge Criterion]
		\label{lemma:iid:hedge_criterion}
		We use the formulation from
		\citep[Prop.\ 1 (p.\,7)]{shpitser2023does},
		where the authors correct a small erratum in
		their statement of \citep[Cor.\ 3 (p.\,1225)]{shpitser2006identification2}:
		
		An interventional distribution $P(\tilde{Y}|\PearlDo(\tilde{X}=x))$,
		is identifiable (Def.\ \ref{def:iid:identify}), if and only if
		there is no hedge for $(\tilde{X},\tilde{Y})$ in $\graphClassical$.
	\end{lemma}
	
	\begin{lemma}[Completeness of ID-Algorithm]
		\label{lemma:iid:id_algo_complete}
		\citep[Thm.\,8 (p.\,1226)]{shpitser2006identification2}:
		If the id-algorithm (see \citep{tian2002general,shpitser2006identification2})
		fails to identify 
		$P(\tilde{Y}|\PearlDo(\tilde{X}=x))$,
		then it outputs a hedge $(F,F')$ for $(\tilde{X},\tilde{Y})$.
		
		As a corollary by using the hedge-criterion
		\citep[Cor.\,5 (p.\,1226)]{shpitser2006identification2}
		thus: The id-algorithm is complete.
	\end{lemma}

	\subsection{Translating the Standard IID-Setup}
	\label{apdx:translate_IID_setup}

	All IID models can be described with the following strongly restricted form 
	of symmetry and viewport.
	See also examples \ref{example:model_iid}, \ref{example:model_iid_distr}
	and \ref{example:iid_do_intervention}.
	
	\begin{definition}[IID-Observations with Uniform Missingness]
		\label{def:iid:model}
		We call a model IID (see Rmk.\ \ref{rmk:exchangable}),
		if it can be written in the following form:
		$I=\IDataset\times\IVars\times\ISample$, where $|\ISample|=\infty$,
		while $|\IVars|<\infty$ and $|\IDataset|<\infty$
		with symmetry-group
		$G=\mathfrak{S}_{\IDataset\times\ISample}$
		(the permutations\Slash{}symmetric group acting trivially
		on the middle factor $\IVars$ of $I$ and mixing datasets and samples),
		and variables $J_{(C,v)} = C \times \{v\} \times \ISample$
		with symmetry $H_v = \mathfrak{S}_{C\times\ISample}$
		for some subset $C\subset \IDataset$ (but disjoint, \ie every
		$c$ is in exactly one $C$ of some $J_{(C,v)}$, denoted $C(c)$) and
		$\mathfrak{S}_C\subset \mathfrak{S}_{\IDataset}$ embedded as subgroup
		by extending as the identity on elements outside of $C$.
		We assume $\IDataset=\IDataset^{\txt{obs}}\sqcup\IDataset^{\txt{query}}$
		is a disjoint union of observed data-sets (see viewport below)
		and pseudo-datasets $\IDataset^{\txt{query}}$ to encode queries
		(in the single context case, both are one-element sets:
		$\IDataset^{\txt{obs}}=\{\alpha\}$ and $\IDataset^{\txt{query}}=\{\beta\}$).
		By slight abuse of notation we will pretend that
		$\ISample$ depends on the "context" $c\in\IDataset$
		(as only the viewport, see below, really matters this is purely notational).
		
		There is a subset $\mathcal{O}\subset\IVars$
		with complement $L=\IVars\setminus\mathcal{O}$ and the viewport $\viewport(N)$
		is such that $\forall (c,v,s) \in \IDataset^{\txt{obs}}\times\mathcal{O}\times\ISample$
		$\exists N_0 \in\mathcal{N}$ such that $\forall N\geq N_0$:
		$(c,v,s)\in \viewport(N)$,
		while for $\forall (c,v,s) \in
		(\IDataset^{\txt{obs}}\times L\times\ISample)
		\cup (\IDataset^{\txt{query}}\times \IVars\times\ISample)$
		(\ie for the complement in $I$)
		$\forall N\in\mathcal{N}$: $(c,v,s) \notin \viewport(N)$.\\
		In words: All variables in $\mathcal{O}$ are asymptotically observed for 
		infinitely many samples in all observed contexts $\IDataset^{\txt{obs}}$.
		A typical viewport of this form, for $\ISample=\mathbb{N}$, is given by
		$\mathcal{N} = |\IDataset^{\txt{obs}}\times\mathcal{O}|\mathbb{N}$
		(we observe samples of size $M:=|\IDataset^{\txt{obs}}\times\mathcal{O}|$,
		thus the viewport $\viewport(N)$ is defined for $N$ that are multiples
		of this number),
		and for $N=n M\in\mathcal{N}$
		observed indices are $\viewport(N)=\IDataset^{\txt{obs}}\times\mathcal{O}\times
		\{1,\ldots,n\}$.
		
		Note that by equivariance of parent-sets,
		the $I$-graph decomposes into disjoint (no mutual edges)
		subgraphs on nodes $\{c\}\times\IVars\times\{s\}$,
		and further by $\mathfrak{S}_{\ISample}\subset H_J$
		(and acting as permutations on $\ISample$),
		the equivariant parent sets cannot change between samples
		$s$ (only between contexts $c$; here lower indices $c$ on
		graphs will denote contexts, elsewhere in this paper upper
		indices $c$ on graphs denote structural c-components,
		this should not lead to confusion), thus
		there are local graphs $\graph_c$ (with only proper nodes)
		and families of embeddings (anchored at any $v_0\in\IVars$)
		$\{\psi^c_{j=(c,v_0,s)}\}_{s\in\ISample}$
		embedding $\graph_c$ with image $\{c\}\times\IVars\times\{s\}$
		into the $I$-graph; we can identify $\graph_c$ with
		$\graphStructural_c=\graphObs(\graph_c)$
		(there are no external or pinned nodes;
		see also Def.\ \ref{def:translate_model}).
	\end{definition}
	\begin{rmk}[Exchangeability]\label{rmk:exchangable}
		Exchangeability refers to the property of a
		distribution being invariant under permutations.
		We explicitly keep track of multiple datasets via
		$I = \IDataset \times \IVars \times \ISample$
		and given the data-set index
		we are essentially in a simple case of
		de Finetti's theorem (see \eg \citep{kallenberg2005probabilistic}):
		If the samples are exchangeable, than they are conditionally
		IID.
		In our case: conditional on the data-set index.
		Thus even though we only assume a symmetry,
		by explicitly tracking data-set association
		the models defined above are effectively IID (not just exchangeable)
		and the naming-convention makes sense.
	\end{rmk}
	
	Model mechanisms in the sense of Def.\ \ref{def:mechanism}
	under these symmetry-restrictions are then
	push-forwards of noise-distributions.
	
	\begin{lemmaDef}[Translating Mechanisms]
		\label{lemma:iid:translate_mech}
		Given a measurable map $f : \val{X} \times \val{N} \rightarrow \val{Y}$
		and a probability measure $P_\eta$ on $\val{N}$,
		define a probability-kernel $\mu(f,P_\eta)$ from $\val{X}$
		to $\val{Y}$ given by $x\mapsto f(x,\cdot)_*P_\eta$ (see Rmk.\ \ref{rmk:kernels_measurevalued})
		or equivalently for $x\in\val{X}$ and $B\in\sigmaAlgebraBorelRaw{\val{Y}}$:
		\begin{equation*}
			\mu(f,P_\eta)(x,B)
			\halfquad:=\halfquad
			\int 1_B\big( f(x, u) \big) P_\eta(\du)\txt.
		\end{equation*}
		If a random variable $\randomVar{Y}$ satisfies
		$P(\randomVar{Y}|\randomVar{X}=x)=f(x,\eta)$ for a random
		variable $\eta\independent \randomVar{X}$ with law $P(\eta)=P_\eta$, then
		$P(\randomVar{X},\randomVar{Y}) = P_\eta \otimes \mu(f,P_\eta)$.
		
		In particular for an SCM $M$, there are probability-kernels
		$\mu_v:=\mu(f_v,P(\eta_v))$ for $v\in\IVars$ such that aligning $\graphStructural=\graphClassical^L$
		(with only inner nodes) to these $\mu_v$ has a structured kernel
		$\mu(\graphStructural,L,\pi)=P(\{\iidVarAny_v\}_{v\in\mathcal{O}})$.
		
		\emph{Remark:} Indeed $\mu(f_y,P(\eta_y))_{\pa_y}
		= P(\randomVar{Y}|\Pa_y=\pa_y)$,
		however it seems conceptually clearer to define it directly relative
		mechanisms and noise-distributions, rather than through intermediate
		realizations as random variables.
	\end{lemmaDef}
	\begin{proof}
		See examples in §\ref{apdx:ptheo_basics}.
	\end{proof}
	
	Commonly, \eg in \citep{PearlBook, Elements,
		tian2002general,shpitser2006identification2,
		shpitser2006identification1, bareinboim2012TransportCompleteness},
	the following assumption is implicit in the assumption of "generic" mechanisms:
	
	\begin{assumption}[Unique Mechanisms]
		\label{ass:iid:unique_mech}
		There are no two different variables with the same mechanism:
		Given an SCM $M$, then
		$\forall v\neq v' \in \IVars: \mu(f_v,P(\eta_v))\neq \mu(f_{v'},P(\eta_{v'}))$.
	\end{assumption}
	
	We are now able to "translate" standard SCMs
	to (restricted symmetry) models in the sense of Def.\ \ref{def:model}.
	
	\begin{lemmaDef}[Translated Model]
		\label{def:translate_model}
		Given an SCM $M$ with variables $\IVars$ and
		observed variables $\mathcal{O}\subset\IVars$,
		and an intervention $P(\tilde{Y}|\PearlDo(\tilde{X}=x))$,
		there is an IID model $\mathcal{M}(M)$ (Def.\ \ref{def:model})
		given by $\IDataset^{\txt{obs}} = \{\alpha\}$,
		$\IDataset^{\txt{query}}=\{\beta\}$,
		variables as required for being an IID model (Def.\ \ref{def:iid:model}),
		\ie $J_{(C,v)} = C \times \{v\} \times \ISample$ with
		symmetry
		$H_v = \mathfrak{S}_{C\times\ISample}$.
		Here $C$ is $C=\{\alpha,\beta\}$ if $v\notin \tilde{X}$,
		if $v\in \tilde{X}$ then there are two variables
		$J_{(\{\alpha\},v)}$ and $J_{(\{\beta\},v)}$ associated to $v$.
		
		If $\alpha \in C$, then the
		mechanisms are $f_{J_{(C,v)}} = \mu(f_v,P(\eta_v))$ (Def.\ \ref{lemma:iid:translate_mech})
		and (equivariant) parent-mappings
		at $j=(v,s)\in J_v$ take values
		$\PaIdx{k}_{J_v}(j) = (p_k,s)$, where $p_k$ is the
		$k$-th (in order of arguments of $f_v$) parent of $v$,
		\ie as a set $\Pa_{J_v}(j=(v,s))=\Pa_v\times \{s\}$.
		Value- and parent-value-spaces are the value-spaces $\val{X}_v$
		of the SCM.
		
		If $C=\{\beta\}$ (and thus $v\in\tilde{X}$),
		then mechanisms are \emph{known} $f^{\PearlDo}_v\in\knownKernels$.
		Indeed for a do-intervention, they are singular at value
		$x_v$ (the $v$-component of the tuple $x$ to which $\tilde{X}$ was intervened)
		without parents. We will default to the value-space
		$\val{X}_v$ (assuming $x_v\in\val{X}_v$ as required
		in Def.\ \ref{def:iid:interventions}).
		
		The viewport is as required by Def.\ \ref{def:iid:model},
		typically 
		$\viewport(N)=\IDataset^{\txt{obs}}\times\mathcal{O}\times
			\{1,\ldots,n\}$.
		
		This model $\mathcal{M}(M)$ has the following properties.
		The $I$-graph of $\mathcal{M}(M)$ is the disjoint union of
		its restrictions to
		$\{\alpha\}\times\IVars\times\{s\}$
		which are of the form $\graphClassical^L$ (the causal graph including latents
		of $M$) and
		its restriction to
		$\{\beta\}\times\IVars\times\{*\}$
		which is of the form $\graphClassicalDo^L$,
		\ie the families of embeddings of the last paragraph
		of Def.\ \ref{def:iid:model} embed
		$\graphStructural_\alpha=\graph_\alpha = \graphClassical^L$
		and $\graphStructural_\beta=\graph_\beta = \graphClassicalDo^L$ as
		local graphs with only proper nodes and aligned to $\mu(f_v,P(\eta_v))$.
		
		As an immediate consequence		
		$\forall s\in\ISample$:
		\begin{align*}
			P=P(\{\iidVarAny_v\}_{v\in\mathcal{O}})
			\halfquad&=\halfquad
			\shallowDistr(\{\anyVar_i\}_{i\in\{\alpha\}\times\mathcal{O}\times\{s\}})
			\\			
			P(\tilde{Y}|\PearlDo(\tilde{X}=x))
			\halfquad&=\halfquad
			\shallowDistr(\{\anyVar_i\}_{i\in\{\beta\}\times\tilde{Y}\times\{*\}})
			\txt,
		\end{align*}
		where $*$ is the single (relevant) "sample" of the query.
	\end{lemmaDef}
	\begin{proof}
		The claims about the $I$-graph are true by construction.
		By Lemma \ref{lemma:iid:translate_mech}
		thus for each $s\in\ISample$:
		$\mu(\graph_\alpha,L,\pi) = P(\{\iidVarAny_v\}_{v\in\IVars})=P$.
		On the other hand, by
		$\shallowDistr(\{\anyVar_i\}_{i\in\{\alpha\}\times\mathcal{O}\times\{s\}})
		= \mu(\graph_\alpha,L,\pi)$
		by Def.\ \ref{def:shallow_distr}.
		
		Similarly 
		$\mu(\graph_\beta,L,\pi) = P(\{\iidVarAny^{\PearlDo}_v\}_{v\in\IVars})=
		P(\tilde{Y}|\PearlDo(\tilde{X}=x))$
		and $\shallowDistr(\{\anyVar_i\}_{i\in\{\beta\}\times\tilde{Y}\times\{*\}})
		=\mu(\graph_\beta,L,\pi)$.
	\end{proof}
	
	The commonly considered do-interventions
	are a special case of structured queries.
	We already built them into the translated model
	as a context $\beta$, and the last result,
	$P(\tilde{Y}|\PearlDo(\tilde{X}=x))=
	\shallowDistr(\{\anyVar_i\}_{i\in\{\beta\}\times\tilde{Y}\times\{*\}})$
	suggests a rather obvious choice of query (Def.\ \ref{def:query}):
		
	\begin{definition}[Translated Query]
		\label{def:iid:query_graph}
		In $\mathcal{M}(M)$, construct a query as follows:
		$q=(Y_q, X_q, \theta)$ with $X_q=\emptyset$
		and $Y_q=\{\beta\}\times\tilde{Y}\times\{*\}$ and $\theta=x=(x_k)_{k\in\tilde{X}}$
		fully parametrizes $\knownFunction=\{\delta[x_k]\}_{k\in\tilde{X}}$.
		The structured query associated by Lemma \ref{lemma:assoc_structured_query}
		is $\graphClassicalDo^L|_{\AncDo(\tilde{Y})}$
		(cf.\ proof of Lemma \ref{lemma:assoc_structured_query}
		or its explanation in the main text),
		using the identification $\graphStructural_\beta = \graphClassicalDo^L$
		and restricting this graph to
		$\AncDo(\tilde{Y}):=\Anc_{\graphClassicalDo^L}(\tilde{Y})$.
		Note that $\IAux=\emptyset$ (the view-port and the query are
		in disconnected graph-components of the $I$-graph); thus (by Markov-property,
		Lemma \ref{lemma:obs_world_existence}),
		the underlying query is
		\begin{equation*}
			\realizedDistr(\{\anyVar_i\}_{i\in\{\beta\}\times\tilde{Y}\times\{*\}})
			=
			\shallowDistr(\{\anyVar_i\}_{i\in\{\beta\}\times\tilde{Y}\times\{*\}})
			=			
			P(\tilde{Y}|\PearlDo(\tilde{X}=x))
			\txt.
		\end{equation*}
	\end{definition}
		
	In the IID-case, no difficulties arise from non-uniqueness (of minimal
	latent-sets, ancestral structures or absorbed parents\Slash{}children),
	which allows for a very simple classification of families of embeddings:
	
	\begin{lemma}[IID Classification of Embeddings]\label{lemma:iid_psi_classification}
		For an IID-model with
		$\graphClassical^L = (\nodesClassical^L,\edgesClassical^L)$,
		decorated families of embeddings with maximal $J_0$ 
		are of the form: $J_0 = C\times\ISample$, where
		$\emptyset\neq C\subset \IDataset^{\txt{obs}}$
		are the relevant contexts
		(for the present sub-section $C=\{\alpha\}$ is the only valid choice),
		$\nodes\subset\nodesClassical^L=\IVars$, $\nodesOuterFixed=\emptyset$ and
		$\psi_{j=(c,s)}(n) = (c,n,s)\in \IDataset\times\IVars\times\ISample = I$.
		There is always a unique minimal latent subset given by $\nodes^L\cap L$.
		Edges are fixed by the subset of $\edgesClassical^L$
		(by Def.\ \ref{def:local_graph_embedding}).
		\\
		\emph{Remark:} Some care should be taken, because ancestral structure
		is part of decorated families of embeddings (but hidden in this notation),
		in a multi-context setup,
		there might be different choices for the relevant contexts
		$C\subset \IDataset$ differing only by validity of respective ancestral structures.
	\end{lemma}
	\begin{proof}
		As has been noted in Def.\ \ref{def:iid:model},
		there are families of embeddings of $\graph_c$
		to $\{c\}\times\IVars\times\{s\}$ indexed by $s\in J_0=\ISample$.
		We can also embed subgraphs $\graph'_c$ of these $\graph_c$.
		The index set $J'_0$ can be potentially enlarged to
		$J_0^C =C\times\{v_0\}\times\ISample$, where $C$ is the largest
		subset of $\IDataset^{\txt{obs}}$ such that
		model-alignment is still valid.
		However by Ass.\ \ref{ass:iid:unique_mech}
		(unique mechanisms) each mechanism $\mu_v$ appears
		exactly once for each $(c,s)\in\IDataset^{\txt{obs}}\times\ISample$
		so these are indeed the largest possible $J_0$.
		
		By the form of the viewport, hiddenness of $i=(c,v,s)$ depends
		only on $v\in L$ and all contexts are observed (asymptotically)
		infinitely often, so the claims minimal latent sets are immediate.		
	\end{proof}
	\begin{notation}
		The classification-lemma \ref{lemma:iid_psi_classification}
		justifies denoting IID embeddings $\psi$ as
		a triple $(C,\nodesInnerProper,\nodes)$ with $C\subset\IDataset$,
		$\nodesInnerProper\subset \nodes \subset \IVars$.
		In the single context case, there is only one observed context, so
		$C=\{\alpha\}$,
		thus until §\ref{apdx:translate_mz_transport} we will
		simply write $\psi=(\nodesInnerProper,\nodes)$ for statements about extraction
		(the set of minimal latents is implicitly given by $L\cap\nodes$).
		Indeed, as long as contexts individually have infinitely many
		samples (asymptotically), the size of $C$ is actually irrelevant (as long
		as it is non-empty).
	\end{notation}
	\begin{lemma}[Direct Identification]
		\label{lemma:iid:direct_id}
		In the IID case,
		the claim of Thm.\ \ref{thm:extract_from_backdoor_complete}
		(extraction from backdoor-free families of embeddings)
		holds true with the notion of "direct identifiability" from
		Def.\ \ref{def:identification_direct_mt} replaced by
		"knowledge of the observational distribution $P=\shallowDistr(\{\anyVar_i\}_{i\in\{\alpha\}\times\mathcal{O}\times\{s\}})$
		and its conditional distributions".
		
		\emph{Remark:}
		This results will become much clearer in
		the light of the analysis of the ED-phase in
		§\ref{apdx:iid:single_context_empty_cond}.
	\end{lemma}
	\begin{proof}
		The proof
		of Thm.\ \ref{thm:extract_from_backdoor_complete}
		only uses data-sets induced via Lemma \ref{lemma:valid_datasets},
		which are of the form of marginalizations of the observational distribution.
		Direct identification than claims that if
		for each sample	$P(X_s,Y_s)= \mu_s \otimes \nu_x$
		for $s$-independent $\nu$ (the disintegration
		$\nu_x = P(Y|X=x)$ always independent of $s$ by IIDness),
		$\nu$ is "directly identifiable",
		and $\nu_x = P(Y|X=x)$ is, indeed, a conditional
		distribution of $P$.
	\end{proof}
	\begin{cor}[Identification]
		Regular functionals are given uniquely
		by their argument,
		thus
		identifiability in the sense of Def.\ \ref{def:identification}
		implies in the IID case
		identifiability in the sense of Def.\ \ref{def:iid:identify}.
		
		If we can show that our algorithm identifies a query in the sense
		of Def.\ \ref{def:identification},
		then this implies
		our algorithm identifies this query
		in the sense of Def.\ \ref{def:iid:identify}.
		
		Finally, identifying the structured query
		$\graphClassicalDo^L|_{\AncDo(\tilde{Y})}$
		identifies $P(\tilde{Y}|\PearlDo(\tilde{X}=x))$
		(cf.\ \ref{def:iid:query_graph})
		in the sense
		of Def.\ \ref{def:iid:identify}.
	\end{cor}
	\begin{proof}
		Regular functionals are composed of finitely many operations,
		which individual produce a well-defined (unique) result
		(§\ref{apdx:ptheo_basics}).
		Thus regular functionals compute a unique result.
		The other claims follow directly from the definitions
	\end{proof}

	\subsection{Single-Context with Empty Conditioning Set}
	\label{apdx:iid:single_context_empty_cond}
	
	We first analyze what results are produced by the
	extract--decompose (ED)-phase(s) in Algo.\ \ref{algo:id},
	steps \ref{algo:id:E} and \ref{algo:id:D}.
	Then we separately analyze the assembly (A)-phase in Algo.\ \ref{algo:id},
	step \ref{algo:id:A}.
	
	\paragraph{ED-Phase:}
	We discus the "ED-phase" combined, as it will turn out that in the IID-case
	the extraction (E)-phase can be replaced by the claim
	"$\mu(\graphClassical^L,L)=P$ (the observational distribution)
	is identifiable" (which also makes it clearer why
	Lemma \ref{lemma:iid:direct_id} was essentially trivial to prove).
	To illustrate this, we recall the (not finite-context IID) example
	\ref{example:extraction_frontdoor}
	where the E-phase is clearly necessary.
	\begin{example}[Non-Trivial Extraction]
		\label{example:non_trivial_extraction}
		Consider a setup (inspired by JPCMCI \citep{JPCMCI}) given as follows:
		There are time-series for time-points $T$ at river-sites $S$ for
		multiple variables. Some of the mechanisms vary between
		sites $s\in S$, others over time $t\in T$ or are shared.
		We assume there are many sites and time-points,
		\ie the relevant asymptotic limit has both
		$|S|\rightarrow\infty$ and $|T|\rightarrow\infty$ as
		$N\rightarrow\infty$.
		Lets consider a very simple such model of the form
		$X \rightarrow Z_s \rightarrow Y_t$,
		where $X$ is the same mechanism for all $s\in S$ and $t\in T$,
		$Z_s$ is the same for all $t\in T$, but depends on $s\in S$,
		while $Y_t$ is the same for all $s\in S$, but depends on $t\in T$.
		Most frameworks, like JCI \citep{JCI}, mz-transport
		\citep{bareinboim2012TransportCompleteness} or JPCMCI
		\citep{JPCMCI} would inject
		"context-variables" for $S$ and $T$.
		
		We want to predict an intervention on site $s$ and at time $t$ on
		for $\PearlDo(X=x)$. Conventional wisdom tells us this effect is
		identifiable and given simply by
		$P(Y|X=x, S=s, T=t)$. This is correct in IID-terms.
		It is also, unfortunately, not helpful:
		We have exactly one data-point at $S=s$ and $T=t$.
		So what could we do?
		The E-phase of our approach would actually search for smallest
		(best symmetry) backdoor-free graphs.
		It would find the following:
		$\circ \rightarrow Z_s$ is identifiable on the
		(asymptotically infinite) data-set $T\times \{s\} \times \{(x,z)\}$.
		$\circ \rightarrow Y_t$ is identifiable on the
		(asymptotically infinite) data-set $\{t\}\times S \times \{(z,y)\}$.
		Then the A-phase simply glues these partial results.		
	\end{example}
	
	In this example the E-phase does essentially two things:
	(1) It validates that only local graphs with sufficient available data
	are taken into account. Thus it avoids, even if $Z$ is unobserved,
	dubious identification claims and rather fails gracefully.
	(2) It enables the A-phase to find the
	(actually useful) "frontdoor-like" identification
	strategy, even though it would conventionally be entirely unclear
	what makes this strategy better than any other plausible one.
	
	The point here is, that smaller local graphs have better symmetry,
	thus larger available data-sets.
	This is not just
	a finite sample issue, as the example above illustrates, taking
	the wrong local graphs, we can end up with finite data-sets (or even literally a
	single data-point) to learn from, which also breaks
	asymptotic identification results; this problem only gets more severe
	once latents are added.
	As already remarked \ref{rmk:families_and_symmetries},
	the symmetry of a local graph containing two subgraphs $A$ and $B$
	of symmetry $H_A$ and $H_B$ has symmetry $H_A\cap H_B$.
	In some cases, we can simply intersect all
	symmetries of all variables $J$ to get $\cap_{J\in\mathcal{J}} H_J$.
	We then get a "fundamental-domain" local graph containing "everything"
	about the model. In the IID case, all mechanisms
	are at least $\mathfrak{S}_{\ISample}$-invariant,
	and essentially the fundamental domain would be $\graphClassical$
	(or in the finite-context\Slash{}mz-transport case:
	a disjoint union of causal and interventional graphs,
	one for each context and experiment,
	cf.\ Fig.\ \ref{fig:mz_transport}).
	Another example: A time-series with $m$- and $k$-periodic mechanisms
	has fundamental domain $\tau$-periodic, where $\tau$ is
	the smallest common multiple of $m$ and $k$ (translation
	symmetries here are $m\mathbb{Z}$ and $k\mathbb{Z}$ with
	$m\mathbb{Z}\cap k\mathbb{Z}= \tau\mathbb{Z}$).
	As the example above shows, generally the fundamental domain
	may not be identifiable (here:
	$(\mathfrak{S}_{S}\times\{0\})\cap(\{e\}\times\mathbb{Z}_{T})=\{(e,0)\}$
	is the trivial sub-group).
	This is the deeper reason why
	the E-phase is necessary in general, but could be
	skipped for the IID-case.
	
	The result of the combined ED-phase is the same
	as applying only the D-phase (Algo.\ \ref{algo:decomp})
	to $\mu(\graphClassical^L,L)$; see also Fig.\ \ref{fig:graphical_ops}
	for an illustration.
	
	\begin{lemma}[ED-Phase]
		\label{lemma:iid:ed_phase}
		The ED-phase (the first two steps of Algo.\ \ref{algo:id})
		produces after a finite number of steps
		the finite set of all
		c-connected structured kernels (and associated identification strategies)
		which can be obtained from $\mu(\graphClassical^L,L)$
		by repeated decomposition into c-components
		(Lemma \ref{lemma:c_components:mt}) and simplification
		(Lemma \ref{lemma:simplify}); this includes
		intermediate c-connected results (cf.\ Algo.\ \ref{algo:decomp}).
		
		\emph{Notation:} We call c-connected sub-graphs obtained by repeated
		decomposition into c-components
		(Lemma \ref{lemma:c_components:mt}) and simplification
		(Lemma \ref{lemma:simplify}), c-fragments.
		The claim above then reads: The ED-phase outputs all
		c-fragments of $\mu(\graphClassical^L, L)$.
	\end{lemma}
	\begin{proof}
		The E-phase produces all minimal (c-connected)
		backdoor-free families, which are, by classification
		of embeddings (Lemma \ref{lemma:iid_psi_classification}),
		are all c-subgraphs of $\graphClassical^L$.
		\Ie the E-phase
		produces the decomposition of $\graphClassical^L$
		into c-components
		(Lemma \ref{lemma:c_components:mt}) and simplifications
		(Lemma \ref{lemma:simplify}).
		The D-phase then further iterates these steps.
		
		Since all non-trivial decomposition-steps (c-components or simplifications)
		reduce the size of the graph
		(the produced pieces are smaller than the original one)
		and Algo.\ \ref{algo:decomp} yields on trivial ones,
		the sizes of graphs in step $m$ is at most
		$|\IVars| - m$ nodes (because $\graphClassical^L$ has $|\IVars|$ nodes),
		and the algorithm terminates after a finite number of iterations.
		There is a finite number of forks
		(sub-iterations, \eg c-components, simplifications)
		in each iteration (there are only finitely many c-components and
		subsets of $L$ in finite graphs).
	\end{proof}

	\paragraph{A-Phase:}
	We finally discuss the assembly (A)-phase in Algo.\ \ref{algo:id},
	step \ref{algo:id:A}. We start by relating c-fragments
	(cf.\ \ref{lemma:iid:ed_phase}) to c-forests (Def.\ \ref{def:iid:c_forest},
	cf.\ \citep{shpitser2006identification2}).
	
	\begin{lemma}[C-Fragments and C-Forests]
		\label{lemma:iid:cfragments_to_forests}
		Some properties of c-structures are:
		\begin{enumerate}[label=(\alph*)]
			\item\label{lemma:iid:cfragments_to_forests:a}
			Given a c-fragment $\graphStructural^c\leq(\graphStructural, L)$,
			then there is a c-forest in $\graphStructural$ 
			on $\nodesInner^c\setminus L$
			rooted on
			\begin{equation*}
				R
				\halfquad:=\halfquad
				\{\halfquad
				y \in \nodesInner^c\setminus L
				\halfquad|\halfquad
				\Ch_{\graphStructural}(y)\cap\nodesInner^c=\emptyset
				\halfquad\}\txt.
			\end{equation*}
			\item\label{lemma:iid:cfragments_to_forests:forest_do}
			A c-forest $F$ in $\graphClassicalDo$
			is a c-forest in $\graphClassical$.
			\item\label{lemma:iid:cfragments_to_forests:overlap_x}
			Given a $R$-rooted c-forest $F'$ in $\graphClassical$ such that
			there is no c-fragment $\graphStructural^c\leq(\graphClassical^L, L)$,
			with $\nodesInner^c\setminus L=F'$, then there is a $R$-rooted
			c-forest $F\supsetneq F'$ in $\graphClassical$.
		\end{enumerate}
	\end{lemma}
	\begin{proof}
		\textbf{Part \ref{lemma:iid:cfragments_to_forests:a}:}		
		We first show: C-fragments contain an edge-subgraph that is a forest.
		By definition of $R$, every node $y\in\nodesInner^c\setminus (L\cup R)$
		has a child $c\in\nodesInner^c$.
		Using the special form of hidden structure
		assumed in the IID-case (cf.\ Def.\ \ref{def:iid:scm},
		Rmk.\ \ref{rmk:latent_proj}), hidden nodes have no parents,
		thus $c\in \nodesInner^c\setminus L$.
		Pick one of these children and call it $c(y)$.
		For $y\in R$.
		Then define a forest $F$ rooted on $R$
		with nodes $\nodesInner^c$
		and an edges $y\rightarrow c(y)$ for each
		$y\in\nodesInner^c\setminus (L\cup R)$.
		C-fragments are c-connected (a c-component in the sense
		of \citep{tian2002general,shpitser2006identification2}) by definition,
		thus this is a c-forest.

		\textbf{Part \ref{lemma:iid:cfragments_to_forests:forest_do}:}
		A c-forest in $\graphClassicalDo$ is
		an edge-subgraph $F\subset \graphClassicalDo$ with a single c-component.
		Also $\graphClassicalDo\subset\graphClassical$ is an edge-subgraph,
		so $F \subset \graphClassical$ is an edge-subgraph.
		Since the nodes of $F$ are in a single c-component
		of $\graphClassicalDo$, there is a path of
		bi-directed arrows between any two nodes, but again
		by $\graphClassicalDo\subset\graphClassical$ being an edge-subgraph,
		this path of bi-directed arrows is also in $\graphClassical$.
		
		\textbf{Part \ref{lemma:iid:cfragments_to_forests:overlap_x}:}
		By definition the c-forest $F'$ is a single c-component,
		\ie all its nodes are in the same c-component
		$F'\subset C\subset\nodesClassical$
		of $\graphClassical$.
		Call the associated structural c-component
		(with $\nodesInner^c\setminus L = C$)
		$\graphStructural^c\leq\graphClassical^L$.
		
		This $\graphStructural^c$ need not be rooted on $R$.
		So we define a simplification $\graphStructural'$
		of $\graphStructural^c$ as follows:
		Define $R_0 := R$, then
		\begin{equation*}
			R_{k+1} := R_k \cup
			\big(\nodesInner^c \cap \Pa_{\graphClassical}(R_k)\big)
			\txt.
		\end{equation*}
		After a finite number of steps $R_{n+1} = R_n$.
		Finally define $B:=\nodesInner^c \setminus (R_n \cup L)$
		and add latents that are not a parent
		(or ancestor by form of hidden structure,
		cf.\ Def.\ \ref{def:iid:scm}, Rmk.\ \ref{rmk:latent_proj})
		of any node in $\nodesInner^c$ not in $B$ to get
		$B^L := B \cup (L\setminus \Pa_{\graphClassical^L}
		(\nodesInner^c \setminus B))$.
		By construction $B=\Dec_{\graphStructural^c}(B)$ and
		there is a $B$-simplification $\graphStructural'$
		of $\graphStructural^c$.
		
		As a simplification of the structural c-component $\graphStructural^c$,
		$\graphStructural'$ is a c-fragment of $\graphClassical^L$.
		Thus
		by hypothesis $\nodesInner' \setminus L \neq F'$.
		Note that initially $F\subset\nodesInner^c$, and the simplification
		cannot remove nodes in $F$ (because they are ancestors of the
		observable $R$), thus $\nodesInner' \setminus L \subsetneq F'$.
		By \ref{lemma:iid:cfragments_to_forests:a}, there
		is a c-forest $\tilde{F}$ on the nodes $\nodesInner' \setminus L$.
		Note that $\tilde{F}$ is rooted on $R$,
		by construction of $R_n$.
		Define $F$ on the nodes $\nodesInner' \setminus L$
		containing $F'$ as a subgraph as follows:
		For nodes $y\in F'\subset \tilde{F}$ insert all edges (there is actually
		at most one) out of $y$ in $F'$, for nodes $y\in \tilde{F}\setminus F$,
		insert all edges (there is actually exactly one) out of $y$ in $\tilde{F}$.
	\end{proof}

	\begin{lemma}[Identification of Do-Interventions]
		\label{lemma:iid:id_of_do_interventions}
		Applying algorithm \ref{algo:id}
		to the structured query (Def.\ \ref{def:iid:query_graph})
		of the intervention $P(\randomVar{Y}|\PearlDo(\randomVar{X})=x)$
		returns after a finite number of steps.
		It outputs
		a non-empty set (thus at least one sound identification strategy)
		if and only if there is a no hedge $(F,F')$ for
		$P(\tilde{Y}|\PearlDo(\tilde{X}=x))$.
	\end{lemma}
	\begin{proof}
		We equivalently proof: The output of Algo.\ \ref{algo:id}
		is empty $\Leftrightarrow$ there is a hedge.
		Recall that $\graphQuery = \graphClassicalDo^L|_{\AncDo(Y)}$
		(Def.\ \ref{def:iid:query_graph}).
		Further, we note that by construction structural c-components
		$\graphStructural^c\leq \graphClassicalDo^L|_{\AncDo(Y)}$
		of $\graphClassicalDo^L|_{\AncDo(Y)}$ are either such that
		$\nodesInner^c \cap \tilde{X} = \emptyset$
		or $\nodesInner^c=\{x\}$ for a single node $x\in\tilde{X}$.
		This is because hard-interventions do not have parents,
		thus also no hidden parents.
		
		"$\Rightarrow$":
		Given, the output of Algo.\ \ref{algo:id}
		is empty.
		
		\textbf{Step 1:} $\exists \graphStructural^c \leq \graphClassicalDo^L|_{\AncDo(Y)}$,
		a structural c-component of $\graphClassicalDo^L|_{\AncDo(Y)}$,
		where $\graphStructural^c$ is not one of the single-node c-components
		with inner node in $\tilde{X}$,
		such that $\graphStructural^c \leq \graphClassical^L$
		is \emph{not} a c-fragment
		of $\graphClassical^L$.
		
		\emph{Proof of Step 1:}
		By contradiction. Let
		$\graphStructural^c \leq \graphClassicalDo^L|_{\AncDo(Y)}$
		be an arbitrary structural c-component
		of $\graphClassicalDo^L|_{\AncDo(Y)}$.
		If $\graphStructural^c$ is one of the single-node c-components
		with inner node in $\tilde{X}$, then $\graphStructural^c$ is identified
		by definition of $\knownFunction$
		and added to the output of the ED-phase (in step \ref{algo:id:A} of
		Algo.\ \ref{algo:id}).
		Otherwise,
		$\graphStructural^c \leq \graphClassical^L$
		is a c-fragment
		of $\graphClassical^L$ (this is the negation of the claim considered in the
		proof by contradiction).
		Thus by Lemma \ref{lemma:iid:ed_phase},
		$\graphStructural^c$ is identifiable and in the output of
		the ED-phase.		
		By repeated application of gluing (Lemma \ref{lemma:glue:mt}),
		thus $\graphClassicalDo^L|_{\Anc(Y)}$ is identifiable.
		The A-phase (Algo.\ \ref{algo:svs}) searches all such gluings
		on the output of the ED-phase, thus it finds (and returns)
		this identification-strategy.
		This is a contradiction to the output of Algo.\ \ref{algo:id}
		being empty.
		
		\textbf{Step 2:} There is a hedge for $P(\randomVar{Y}|\PearlDo(\randomVar{X})=x)$ in $\graphClassical^L$.
		
		\emph{Proof of Step 2:}
		Let $\graphStructural^c \leq \graphClassicalDo^L|_{\AncDo(Y)}$
		be the structural c-component from step 1.
		By applying
		Lemma \ref{lemma:iid:cfragments_to_forests}~%
		\ref{lemma:iid:cfragments_to_forests:a}
		(with $\graphStructural=\graphClassicalDo^L$ and using
		that c-components are trivially c-fragments),
		then \ref{lemma:iid:cfragments_to_forests}~%
		\ref{lemma:iid:cfragments_to_forests:forest_do}
		there is a c-forest $F'$ in $\graphClassical$ on nodes
		$\nodesInner^c\setminus L$
		rooted on (some) $R\subset\Anc_{\graphClassical}(Y)$
		(by $\graphStructural^c \leq \graphClassicalDo^L|_{\Anc(Y)}$
		we have $\nodesInner^c \subset \AncDo(Y)\subset
		\Anc_{\graphClassical^L}(Y)$, with
		$\Anc_{\graphClassical}(Y)= \Anc_{\graphClassical^L}(Y)\setminus L$,
		thus $\nodesInner^c\setminus L\subset \Anc_{\graphClassical}(Y)$).
		Since $\graphStructural^c$ is not one of the single-node c-components
		with inner node in $\tilde{X}$,
		this c-forest satisfies $\nodesInner^c\cap\tilde{X}=\emptyset$.
		Also by step 1,
		$\graphStructural^c \leq \graphClassical^L$
		is \emph{not} a c-fragment.
		Thus, by Lemma \ref{lemma:iid:cfragments_to_forests}~%
		\ref{lemma:iid:cfragments_to_forests:overlap_x},
		there is a c-forest $F\supsetneq F'$ in $\graphClassical^L$ rooted on $R$.
		Let $x\in F\setminus F'$ be arbitrary.
		By $F$ being a c-forest
		there is $y\in F'$, $x\leftrightarrow y$ in $\graphClassical$
		(\ie $\exists l\in\graphClassical^L$: $x\leftarrow l\rightarrow y$
		in $\graphClassical^L$).
		Since $x\notin \nodesInner^c=F'$,
		one of the two edges $x\leftarrow l\rightarrow y$
		is not in $\graphClassicalDo^L$
		(otherwise $x$ would be in the c-component $\nodesInner^c$).
		The only edges in $\graphClassical^L$ not in $\graphClassicalDo^L$,
		by definition,
		are edges into $\tilde{X}$. Since $y\in F'$ and $F'\cap \tilde{X}=\emptyset$,
		thus $x \in \tilde{X}$. In particular
		$F\cap \tilde{X}\neq\emptyset$.
		Thus $(F,F')$ is a hedge for $P(\randomVar{Y}|\PearlDo(\randomVar{X})=x)$ in $\graphClassical^L$.
		
		"$\Leftarrow$":
		This direction follows indirectly by soundness and
		the hedge-criterion \citep{shpitser2006identification2} (Lemma
		\ref{lemma:iid:hedge_criterion}).
		Alternatively, note that the proof of the other direction
		can (essentially) be traced backward, very briefly:
		Since $F'$ is c-connected and $F'\cap\tilde{X}=\emptyset$,
		there is a c-component $\graphStructural^c$
		in $\graphClassicalDo$ containing
		$F'$ in its inner nodes.
		By $F\cap\tilde{X}\neq\emptyset$ and $F$ being rooted on $R$
		there is a chain of children in $F$ from $x\in\tilde{X}\cap F$
		to $F'$; or from the perspective of its endpoint $y$ in $F'$:
		there is a chain of (observable) parents in $F$ to $x$.
		Parents of observables cannot be simplified (Lemma \ref{lemma:simplify}) away, so they must remain in any c-fragment of $\graphClassical$.
			
		\textbf{Finite number of steps:}
		The query-graph is finite, and gluing grows
		the intermediate result by at least one node.
		Further there are finitely many elements in $\knowledgeSet$ to 
		choose from (for gluing), thus Algo.\ \ref{algo:svs} terminates
		after a finite number of steps. For the other algorithms,
		see Lemma \ref{lemma:iid:ed_phase}.
	\end{proof}
	\begin{cor}[IID Completeness, single context, empty conditioning set]
		\label{cor:iid:complete_single_context_uncond}
		By the hedge-criterion (Lemma \ref{lemma:iid:hedge_criterion};
		cf.\ \citep{shpitser2006identification2,shpitser2023does}),
		algorithm \ref{algo:id}
		is complete in the (single context) IID case
		for do-interventions with empty conditioning set.
		\\
		\emph{Remark:}
		This result is analogous to the completeness of the id-algorithm
		(Lemma.\ \ref{lemma:iid:id_algo_complete}).
	\end{cor}

	\subsection{Single-Context with Non-Empty Conditioning Set}	
	\label{apdx:iid_conditional}

	In the conditional case, the reader may notice that seemingly
	our criterion and the one given by \citep{shpitser2006identification1} disagree.
	We start by an example, and clarifying explanation on this problem.
	Afterwards we briefly return to conditional effects more generally.

	Our formalism allows for two different approaches to
	conditional effects $\realizedDistr(Y|X=x)$.
	The first one is by directly setting $\tilde{X}=X$ in 
	the basic query. However, it is of course also possible to
	consider $\realizedDistr(Y\cup X)$, \ie setting $\tilde{Y}=Y\cup X$,
	then disintegrate\Slash{}condition for $\realizedDistr(Y|X=x)$ afterwards;
	generally we can split $X$, \ie we only need $Y\subset\tilde{Y}$,
	$\tilde{X}\subset X$ (these two conditions ensure
	information in the query is retained)
	and $X\subset\tilde{X}\cup\tilde{Y}$ (this condition ensures
	the conditional distribution can be obtained).
	These two viewpoints are in very close analogy to
	using rule 2 of the do-calculus to move conditions into the intervened
	set, which is described in detail in \citep[§5]{shpitser2006identification1}.
	However some subtle differences arise:
	\begin{example}
		\label{example:rule2}
		Given observations and a simple intervention $\randomVar{X}_{\PearlDo}\in\knownFunction$.\\
		\begin{minipage}{0.2\textwidth}
			\centering
			\begin{tikzpicture}
				\draw (0,1.5) node[align=center] {Observations:};
				
				\draw (-1,0) node (X) {$\randomVar{X}$};
				\draw (1,0) node (Y) {$\randomVar{Y}$};
				\draw (0,1) node (M) {$\randomVar{M}$};
				
				\draw[->] (X) -- (M);
				\draw[->] (M) -- (Y);
				\draw[->] (X) -- (Y);
				\draw[<->] (X) edge[bend left] (M);
			\end{tikzpicture}
		\end{minipage}
		\begin{minipage}{0.2\textwidth}
			\centering
			\begin{tikzpicture}
				\draw (0,1.5) node[align=center] {Query:};
				
				\draw (-1,0) node (X) {$\randomVar{X}_{\PearlDo}$};
				\draw (1,0) node (Y) {$\randomVar{Y}$};
				\draw (0,1) node (M) {$\randomVar{M}$};
				\draw (-1,1) node (L) {$[\randomVar{L}]$};
				
				\draw[->] (X) -- (M);
				\draw[->] (M) -- (Y);
				\draw[->] (X) -- (Y);
				\draw[->] (L) -- (M);
			\end{tikzpicture}
		\end{minipage}
		\begin{minipage}{0.55\textwidth}
			We want to know $P(\randomVar{Y}|\PearlDo(\randomVar{X}=x),\randomVar{M}=m)$.
			\begin{itemize}
				\item
				According to \citep{shpitser2006identification1} this query is identifiable.
				\item
				According to our logic it is not.
			\end{itemize}
		\end{minipage}\\
		As it turns out, both assessments are correct. This happens, because we do \emph{not} assume
		any (known\Slash{}exploitable) internal structure on mechanisms (observed or intervened)
		while \citep{shpitser2006identification1} implicitly does through rule 2 of the do-calculus.
	\end{example}
	
	We first illustrate the actual problem: Rule 2 of the do-calculus allows the
	replacement 
	\begin{equation*}
		P(\randomVar{Y}|\PearlDo(\randomVar{X}=x),\randomVar{M}=m)
		=
		P(\randomVar{Y}|\PearlDo(\randomVar{X}=x,\randomVar{M}=m))
		\txt.
	\end{equation*}
	This replacement is correct, \emph{if $P(\randomVar{X}_{\PearlDo})$ is singular}
	(only takes a single value),
	which rule 2 assumes, but we do not (this is an internal structure of the intervened
	mechanism $\randomVar{X}_{\PearlDo}$).
	If the intervention were for example chosen as $\randomVar{X}_{\PearlDo}\sim\mathcal{N}(0,1)$,
	then conditioning on $M$ generally will introduce selection-bias on $\randomVar{X}_{\PearlDo}$.
	For estimation of
	$\shallowDistr(\randomVar{Y}|\randomVar{M}=m)$ (the intervention $\randomVar{X}_{\PearlDo}$
	is now encapsulated as part of $\theta$), we have to estimate the joint distribution
	$\shallowDistr(\randomVar{Y},\randomVar{M})$ first --
	which is not identifiable from our observations (due to the confounder between $X$ and $M$).
	The reader may convince themselves, that this assessment is correct,
	on basis of the confounder preventing an unbiased estimation of the selection-bias
	$P(\randomVar{X}_{\PearlDo}|\randomVar{M})$.
	
	So in general, the conditional effect discussed above is not identifiable, but
	under the constraint of $\randomVar{X}_{\PearlDo}$ taking only a single value it is.
	Both conclusions in example \ref{example:rule2} make sense,
	depending on the question we initially intended to ask.
	The evident next question to ask is: Can we, in our formulation, recover the
	results of do-calculus? \Ie can we draw conclusions under assumptions on the internal
	structure of $\randomVar{X}_{\PearlDo}$?
	As it turns out, there are actually two different ways to achieve this:	
	\begin{enumerate}[label=(\alph*)]
		\item
		Conventionally, expressions like $P(\randomVar{Y}|\PearlDo(\randomVar{X}=x),\randomVar{X}=x)$ are not considered.
		If the intervention on $\randomVar{X}$ is singular, there is no further information provided by its value.
		However, our formalism allows to ask for
		$\shallowDistr(\randomVar{Y}|\randomVar{M}=m,\randomVar{X}_{\PearlDo}=x)$,
		and this \emph{does} make sense if $\randomVar{X}_{\PearlDo}$ is not singular.
		Interestingly this query is identifiable (because $Y_{x,m}$ is identifiable),
		so one can directly recover the notion of conditional queries
		as interpreted by the do-calculus (and \citep{shpitser2006identification1})
		by \emph{conditioning on intervened variables} (including them in $\tilde{X}$
		for the basic query).
		
		This approach does not require any additional theory and works as is.
		It does, however, not provide any insights about other internal structure
		affecting selection-bias.
		\item
		Similar to, for example, instrumental variables (example \ref{example:instrumental_variables}),
		we can extend our notion of regularity\Slash{}computation:
		The rule we need is that for some kernels in
		$\mathcal{F}_{\txt{special}}\subset\knownFunction\cup\modelKernels$
		(containing for example all singular interventions in queries or experimental data)
		selection-bias is trivial, \eg allowing for a regular functional
		$(\mu|\nu) \mapsto \mu$ if $\mu\in\mathcal{F}_{\txt{special}}$.
		Such additional classification of kernels (similar to instrumental variables)
		can be useful more generally. For example if $\randomVar{X}$ is multi-dimensional
		or categorical of the form $\randomVar{X}\in\{1,\ldots,m\}\times \{1,\ldots,n\}$
		and $\randomVar{M}$ depending only on one coordinate-projection $\pi_1$,
		$\randomVar{Y}$ on the other coordinate-projection $\pi_2$, then if $P(\randomVar{X})$ is a product (\ie the components
		relevant to $\randomVar{M}$ and $\randomVar{Y}$ are independent
		$\randomVar{X}_1 \independent \randomVar{X}_2$),
		the query of example \ref{example:rule2} is identifiable.
		
		This approach is more complicated to implement, as it requires a custom
		notion of regularity, but if one is interested in complex problems
		with categorical or multi-dimensional variables, it may provide substantial
		insights.
	\end{enumerate}
	
	In conclusion, the confusion seems to arise from the do-calculus conventionally
	being thought of as a set of \emph{causal} inference rules. However, through rule 2,
	it allows additionally for conclusions based on the internal structure of interventions
	(and unrelated to causality). To reproduce the same results from a causal rule-set,
	one will have to add such an assumption about internal structure separately.
	This is possible for our approach in different ways, simply by conditioning on interventions,
	or by clarifying the rules of allowed computations for interventions known to satisfy
	suitable constraints.

	\paragraph{General Approach:}
	
	Generally, smaller $\tilde{Y}$ makes queries easier to identify.
	So we should (similar to \citep[§5]{shpitser2006identification1})
	move "conditions" into $\tilde{X}$ if possible;
	\ie if compatible with the "no bypassing" hypothesis in
	Lemma \ref{lemma:assoc_structured_query}.
	In the IID case there is a unique maximal choice of $\tilde{X}$,
	cf.\ \citep[Cor.\ 1 (p.\,5)]{shpitser2006identification1}.	
	If we follow (a) outlined above,
	then (an analogue of) \citep[Thm. 6 (p.\,5)]{shpitser2006identification1}
	applies to show that our approach is complete for IID do-interventions
	also in the conditional case.
	
	For the example \ref{example:rule2} above,
	$\tilde{X}=\{M\}$ is not initially allowed, as there is a path from $X$ to $Y$
	bypassing $M$ (the direct edge $X\rightarrow Y$)
	so we have to first identify $\tilde{Y}=\{X_{\PearlDo},M,Y\}$
	which we cannot: The only information we have on $M$ is the extracted
	$X \rightarrow M$ plus the latent confounder $L$. With $\Dec(Y)=\{M,Y\}$
	in this c-connected extracted graph there is no simplification still containing
	information on $M$. To glue the query, which includes
	$X_{\PearlDo}$ in the place of $X$, we thus cannot use our knowledge on $M$.
	If we follow approach (a) above, then $\tilde{X}=\{X_{\PearlDo}, M\}$,
	which is allowed. The only information needed to glue this query,
	implicitly assuming singularity of $X_{\PearlDo}$,
	is its single inner node $\circ \rightarrow Y \leftarrow \circ$,
	which we can extract because $Y$ and both its parents are observed.

	\subsection{Standard Results for the Multi-Context, Experimental Case}
	\label{apdx:mz_transport}
	
	In a seminal series
	of papers
	\citep{Bareinboim2013TransportAlgo, bareinboim2012TransportCompleteness,
	Bareinboim2016TransportOverview}
	(see also \citep{pearl2022external} for an overview)
	Bareinboim and Pearl studied the use of "selection-variables"
	for combining observations and experiments from multiple
	related contexts.
	This idea was also employed and further developed
	(\eg under the name "context-variables") for causal discovery
	for example in \citep{CD-NOD,JCI,JPCMCI}
	and has formed the technical and conceptual basis for
	most available results concerning causality in multi-context
	systems.
	
	Many of the problems encountered and insights made possible by these
	approaches provided substantial motivation for the present paper.
	It hence seems natural to assess how our formalism
	describes their "mz-transportability" setups.
	We primarily discuss the topic along the lines of
	\citep{bareinboim2012TransportCompleteness} in this subsection,
	but reproduce its contents in language consistent with
	§\ref{apdx:iid:standard_results}.
	
	The setup from before is modified two-fold:
	Instead of one system, multiple systems (with known relations)
	are available, with further potentially different experiments
	(observations from intervened, by do-interventions, models).
	We start with two contexts $\pi$ and $\pi^*$ (and add experiments later).
	The	reader may want to have a look at Fig.\ \ref{fig:mz_transport},
	which illustrates many of the encountered concepts on an example.
	Results below are quoted from \citep{bareinboim2012TransportCompleteness}
	as indicated; some comments were added in square-brackets for clarity and context.
	
	\begin{definition}[Selection Diagrams]
		\citep[Def.\ 1 (p.\,4)]{bareinboim2012TransportCompleteness}:
		Let $(M,M^*)$ be a pair of SCMs
		relative to domains $(\pi,\pi^*)$
		[domains are "labels" for the contexts]
		sharing a causal graph $\graphClassical$.
		$(M,M^*)$ is said to induce a selection diagram
		$D$ if $D$ is constructed as follows:
		$D$ has nodes $\nodesClassical \cup S$, where
		$S$ contains a "selection variable" $S_v$ for each
		$v\in\mathcal{O}$ with $f_v\neq f_v^*$ or $P(\eta_v)\neq
		P(\eta_v^*)$,
		every edge in $\graphClassical$ is also an edge in $D$;
		$D$ contains an extra edge $S_v\rightarrow v$
		for each $S_v\in S$.
	\end{definition}
	\begin{rmk}[Changing Latents]
		This definition is not compatible with changing latent variables,
		but can immediately be interpreted in that case by
		adding a selection variable $S_v$ also if
		a latent parent $l\in\Pa_v\cap L$ has
		changing mechanism $f_l\neq f_l^*$
		or changing noise $P(\eta_l)\neq
		P(\eta_l^*)$.
	\end{rmk}

	\begin{definition}[mz-Transportability]
		\label{def:iid:mztransport}
		\citep[Def.\ 2 (p.\,4)]{bareinboim2012TransportCompleteness}:
		Let $\mathEulervmBold{D}=\{D^{(1)}, \ldots, D^{(n)}\}$
		be a collection of selection diagrams relative to source domains
		$\Pi=\{\pi_1,\ldots,\pi_n\}$ and a target domain $\pi^*$, respectively
		[\ie there are SCMs $M^*$ and $M_1, \ldots, M_n$
		and $D^{(k)}$ is a selection-diagram for
		$(M_k,M^*)$, in particular they all share a causal graph
		$\graphClassical$],
		and $Z_k$ (and $Z^*$)
		[subsets of $\IVars$]
		be the variables in which experiments can be conducted in domain $\pi_k$
		(and $\pi^*$).
		Let $(P^k, I_z^k)$ be the pair of observational
		and interventional distributions of $\pi_k$ [on subsets of $Z_k$;
		\ie $P^k=P(\{\iidVarAny^k_v\}_{v\in\mathcal{O}})$,
		where $\iidVarAny^k_v$ are the endogenous variables of $M_k$,
		and $I^k_z=\cup_{Z'\subset Z_k}
		P(\{\iidVarAny^k_v\}_{v\in\mathcal{O}}|\PearlDo(Z'=z'))$,
		where the union-symbol '$\cup$' means the collection of these
		distributions is considered],
		and in an analogous manner, $(P^*, P^*_z)$ be the observational and
		interventional distributions of $\pi^*$.
		
		[\emph{Identifiability:}]
		The causal effect $R=P^*(Y|\PearlDo(X=x))$ is said to be mz-transportable,
		from $\Pi$ to $\pi^*$ in $\mathEulervmBold{D}$
		if $P^*(Y|\PearlDo(X=x))$ is uniquely computable from
		$\cup_k (P^k,P^k_z) \cup (P^*,P^*_z)$ in any model
		[any collection of SCMs $M^*, M_1,\ldots,M_n$] that induces
		$\mathEulervmBold{D}$.
		
		\emph{Remark:}
		This notion of identifiability used is conceptual very similar
		to Def.\ \ref{def:iid:identify}.
	\end{definition}

	\begin{definition}[mz*-shedge]
		\label{def:iid:shedge}
		\citep[Def.\ 5 (p.\,6)]{bareinboim2012TransportCompleteness}:
		Let $\mathEulervmBold{D}=\{D^{(1)}, \ldots, D^{(n)}\}$
		be a collection of selection diagrams relative to source domains
		$\Pi=\{\pi_1,\ldots,\pi_n\}$ and a target domain $\pi^*$, respectively,
		$S_k$ represents the collection of $S$-variables in the selection-diagram
		$D^{(i)}$ and let $D^{(*)}=\graphClassical^*$
		be the causal diagram [graph] of $\pi^*$.
		Let $(P^k, I_z^k)$ be the pairs of observational
		and interventional [on subsets of $Z_k$] distributions of $\pi_k$ [as before],
		similar for $(P^*, I_z^*)$.
		
		Consider an $R$-rooted hedge $(F,F')$ (Def.\ \ref{def:iid:hedge}).
		We say that the induced collection of pairs of $R$-rooted c-forests
		over each diagram $((F_1,F'_1), \ldots, (F_n,F'_n))$,
		is an mz-shedge for $P^*(Y|\PearlDo(X=x))$ relative to 
		experiments $(I^*_z, I_z^1, \ldots, I_z^n)$
		if they are all hedges and one of the following conditions hold for each domain
		$\pi_k$ ($k\in\{*,1,\ldots,n\}$):
		\begin{enumerate}[label=\arabic*.]
			\item 
			There exists at least one variable $s\in S_k$ pointing
			to the induced diagram $F'_k$, or
			\item 
			$(F_k\setminus F'_k)\cap Z_k=\emptyset$, or
			\item 
			The collection of pairs
			of c-forests induced over diagrams,
			$((F_1,F'_1), \ldots, (F_k\setminus Z_k^*, F'_k),$ $\ldots,$
			$(F_n,F'_n))$,
			is also an mz-shedge relative to
			$(I^*_z, I_z^1, \ldots, I_{z\setminus z^*_k}, \ldots, I_z^n)$,
			where $Z_k^* = (F_k\setminus F'_k)\cap Z_k$.
		\end{enumerate}
		Furthermore, we call mz*-shedge the mz-shedge in which there exist[s]
		one directed path from $R\setminus (R\cap\Dec_F(X))$
		to $R\cap\Dec_F(X)$ not passing through $X$.
	\end{definition}
	
	\begin{rmk}
		Unfortunately the quoted definition is in places slightly unclear.
		The appendix to \citep{bareinboim2012TransportCompleteness} containing
		proof-details also seems to not be available online (anymore?).
		
		\textbf{(a) The path from $R\setminus (R\cap\Dec_F(X))$
			to $R\cap\Dec_F(X)$:}
		For the distinction between mz*-shedges and mz-shedges
		\citep[p.\,6f]{bareinboim2012TransportCompleteness} give the following
		example (lhs):\\
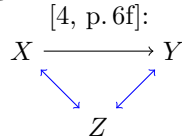
		\begin{minipage}{0.4\textwidth}
			\centering
			\citep[p.\,6f]{bareinboim2012TransportCompleteness}:\\
			\begin{tikzpicture}
				\draw (0,0) node (X) {$X$};
				\draw (2,0) node (Y) {$Y$};
				\draw (1,-1) node (Z) {$Z$};
				
				\draw[->] (X) -- (Y);
				\draw[<->,blue] (X) -- (Z);
				\draw[<->,blue] (Z) -- (Y);
			\end{tikzpicture}
		\end{minipage}
		\begin{minipage}{0.58\textwidth}
			\centering
			Extended Example:\\
			
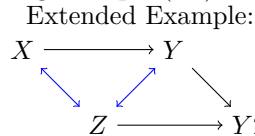
\begin{tikzpicture}
				\draw (0,0) node (X) {$X$};
				\draw (2,0) node (Y) {$Y$};
				\draw (1,-1) node (Z) {$Z$};
				\draw (3,-1) node (Y2) {$Y2$};
				
				\draw[->] (X) -- (Y);
				\draw[->] (Z) -- (Y2);
				\draw[->] (Y) -- (Y2);
				\draw[<->,blue] (X) -- (Z);
				\draw[<->,blue] (Z) -- (Y);
			\end{tikzpicture}
		\end{minipage}\\
		As \citep[p.\,6f]{bareinboim2012TransportCompleteness} correctly
		point out $P(Y|\PearlDo(X=x))$ is identifiable,
		while $P(Y,Z|\PearlDo(X=x))$ is not.
		The requirement for the existence of a path
		(last sentence of the definition) is motivated as follows:
		$F=\{X,Y,Z\}$ and $F'=\{Z,Y\}$ are c-forests rooted
		at $R=\{Z,Y\}$ and the directed path
		"from $R\setminus (R\cap\Dec_F(X))=\{Z\}$
		to $R\cap\Dec_F(X)=\{Y\}$ not passing through $X$" does not exist.
		However, the definition of a hedge already asks
		for $R \subset \AncDo(\tilde{Y})$
		(where $\tilde{Y}$ is the set of target-nodes).
		In particular $(F,F')$ is a hedge for $\tilde{Y}=\{Y,Z\}$
		(consistent with $P(Y,Z|\PearlDo(X=x))$ being \emph{not}
		identifiable), but it is \emph{not} a hedge
		for $\tilde{Y}=\{Y\}$, because $Z$ is not an $\graphClassicalDo$-ancestor
		of $Y$. Indeed for $\tilde{Y}=\{Y\}$ no hedge exists (removing $Y$ from
		$F,F'$ would break them being c-forests).
		
		Note, that directed paths "not passing through $X$" are the same
		as paths in $\graphClassicalDo$ (if the path does not start at a node in $X$).
		So the definition seems to attempt to move the ancestral requirement from
		$\graphClassicalDo$-ancestors of $Y$ "up" to
		$\graphClassicalDo$-ancestors of $(R\cap\Dec_F(X))$
		(there is also a small additional technical\Slash{}notational
		problem here: if one of the sets $R\setminus (R\cap\Dec_F(X))$
		or $R\cap\Dec_F(X)$ is empty -- which occurs especially for simple examples
		quite often -- there obviously never exists such a path; the intended
		statements seems to be one about ancestral relationships in $\graphClassicalDo$).
		However, the weaker requirement of
		$R \subset \AncDo(\tilde{Y})$
		already seems to prevent identification. Consider the second, extended,
		example above: If we want to identify (from a single context without experiments)
		$P(Y_2|\PearlDo(X=x))$,	then we notice that $(F,F')$ is a hedge
		as $R\subset \AncDo(\{Y_2\})$.
		Yet there still exists no directed path
		from $R\setminus (R\cap\Dec_F(X))=\{Z\}$
		to $R\cap\Dec_F(X)=\{Y\}$.
		So for theorems 2 and 3 in \citep[p.\,7]{bareinboim2012TransportCompleteness}
		to hold, one cannot insist on the existence of such a path
		(otherwise the extended example above is a simple counter-example).
		
		\textbf{(b) The "inducedness" and point 3:}
		The authors say a mz-shedge is "induced" from a pair $(F,F')$.
		Since all contexts in $\Pi$ share a graph this makes sense.
		However, point 3 then modifies a single $F_i$
		which clearly breaks the "inducedness" property, so clearly
		the result cannot "also [be] a mz-shedge".
		It appears that the statement intended by the authors was
		something like: "$((F_1,F'_1), \ldots, (F_n,F'_n))$ has the
		mz-shedging property, if 1--3 hold" (where in 3 it is only required
		that the replacement "also has the mz-shedging property");
		this can be written inductively in a non-cyclic form.
		Then an mz-shedge is a pair $(F,F')$ such that the induced
		$((F_1,F'_1), \ldots, (F_n,F'_n))$ has the shedging-property.
	\end{rmk}
	
	Our attempt at a slightly clearer phrasing is the following:
	
	\begin{definition}[mz-shedge revisited]		
		\label{def:iid:shedge_revisited}
		Interpreting \citep[Def.\ 5 (p.\,6)]{bareinboim2012TransportCompleteness}:
		Let $\mathEulervmBold{D}=\{D^{(1)}, \ldots,$ $D^{(n)}\}$
		be a collection of selection diagrams relative to source domains
		$\Pi=\{\pi_1,\ldots,\pi_n\}$ and a target domain $\pi^*$, respectively,
		$S_k$ represents the collection of $S$-variables in the selection-diagram
		$D^{(i)}$ and let $D^{(*)}=\graphClassical^*$
		be the causal diagram [graph] of $\pi^*$.
		Let $(P^k, I_z^k)$ be the pairs of observational
		and interventional [on subsets of $Z_k$] distributions of $\pi_k$ [as before],
		similar for $(P^*, I_z^*)$.
		
		We say that a collection of pairs of $R$-rooted c-forests,
		one for each diagram, $((F_1,F'_1), \ldots,$ $(F_n,F'_n))$,
		has the mz-shedging$_0$ property for $P^*(Y|\PearlDo(X=x))$ relative to 
		experiments $(I^*_z, I_z^1, \ldots, I_z^n)$
		if they are all hedges (Def.\ \ref{def:iid:hedge})
		and one of the following conditions holds for each domain
		$\pi_k$ ($k\in\{*,1,\ldots,n\}$):
		\begin{enumerate}[label=\arabic*.]
			\item 
			There exists at least one variable $s\in S_k$ pointing
			to the induced diagram $F'_k$, or
			\item 
			$(F_k\setminus F'_k)\cap Z_k=\emptyset$.
		\end{enumerate}
		Inductively, 
		$((F_1,F'_1), \ldots, (F_n,F'_n))$,
		has the mz-shedging$_{m+1}$ property for $P^*(Y|\PearlDo(X=x))$ relative to 
		experiments $(I^*_z, I_z^1, \ldots, I_z^n)$
		if they are all hedges (Def.\ \ref{def:iid:hedge})
		and one of the following conditions holds for each domain
		$\pi_k$ ($k\in\{*,1,\ldots,n\}$):
		\begin{enumerate}[label=\arabic*.]
			\item 
				There exists at least one variable $s\in S_k$ pointing
				to the induced diagram $F'_k$, or
			\item 
				The collection of pairs
				of c-forests induced over diagrams,
				$((F_1,F'_1), \ldots,$
				$(F_k\setminus Z_k^*, F'_k),$ $\ldots,$
				$(F_n,F'_n))$,
				is has the mz-shedging$_{m}$ property relative to
				$(I^*_z, I_z^1, \ldots, I_{z\setminus z^*_k}, \ldots, I_z^n)$,
				where $Z_k^* = (F_k\setminus F'_k)\cap Z_k$.
		\end{enumerate}
		$((F_1,F'_1), \ldots, (F_n,F'_n))$,
		has the mz-shedging property for $P^*(Y|\PearlDo(X=x))$ relative to 
		experiments $(I^*_z, I_z^1, \ldots, I_z^n)$
		if $\exists m\in\mathbb{N}_0$ such that it has the
		mz-shedging$_{m}$ property for $P^*(Y|\PearlDo(X=x))$ relative to 
		experiments $(I^*_z, I_z^1, \ldots, I_z^n)$.
		
		Consider an $R$-rooted hedge $(F,F')$ (Def.\ \ref{def:iid:hedge}).
		$(F,F')$ is a mz-shedge,
		for $P^*(Y|\PearlDo(X=x))$ relative to 
		experiments $(I^*_z, I_z^1, \ldots, I_z^n)$,
		if the induced 
		$((F_1,F'_1), \ldots, (F_n,F'_n))$
		has the mz-shedging property
		for $P^*(Y|\PearlDo(X=x))$ relative to 
		experiments $(I^*_z, I_z^1, \ldots, I_z^n)$.
	\end{definition}
	
	\begin{lemma}[Shedge-Criterion]
		\label{def:iid:sheding_criterion}
		Interpreting \citep[Thm.\ 3 (p.\,7)]{bareinboim2012TransportCompleteness}:
		Let $\mathEulervmBold{D}=\{D^{(1)}, \ldots,$ $D^{(n)}\}$
		be a collection of selection diagrams relative to source domains
		$\Pi=\{\pi_1,\ldots,\pi_n\}$ and a target domain $\pi^*$, respectively,
		and $\{I_z^i\}$ for $i=\{*,1,\ldots,n\}$ defined appropriately.
		If there is a mz*-shedge for the effect $R=P^*(Y|\PearlDo(X=x))$
		relative to experiments $(I_z^*,I_z^1,\ldots,I_z^n)$ in $\mathEulervmBold{D}$,
		$R$ is not mz-transportable (Def.\ \ref{def:iid:mztransport})
		from $\Pi$ to $\pi^*$ in $\mathEulervmBold{D}$.
	\end{lemma}
	
	\begin{lemma}[Completeness of TR\textsuperscript{mz}-Algorithm]
		\citep[Thm.\,4 (p.\,8)]{bareinboim2012TransportCompleteness}:
		If the TR\textsuperscript{mz}-algorithm
		(\citep[Fig.\ 3 (p.\,7)]{bareinboim2012TransportCompleteness})
		fails to identify 
		$P(\tilde{Y}|\PearlDo(\tilde{X}=x))$,
		then it outputs a graph-pair that
		contains as edge-subgraphs c-forests that span a mz*-shedge for
		$P(\tilde{Y}|\PearlDo(\tilde{X}=x))$.
		[To our understanding, the same problem as with 
		\citep[Cor.\ 3 (p.\,1225)]{shpitser2006identification2} appears,
		see \citep[Prop.\ 1 (p.\,7)]{shpitser2023does} and Lemma \ref{lemma:iid:hedge_criterion}].
		
		As a corollary by using the shedge-criterion thus
		\citep[Thm.\,5 (p.\,8)]{bareinboim2012TransportCompleteness}:
		The TR\textsuperscript{mz}-algorithm is complete.
	\end{lemma}
	
	\subsection{Translating the Transportability Setup}
	\label{apdx:translate_mz_transport}
	
	It is possible to translate the mz-transport setup
	in two equivalent ways: as different systems sharing
	properties or as copies of the same system with selection-variables
	providing a multi-level structure (where the upper level
	is trivial\Slash{}consists exclusively of nodes without
	edges connecting them).
	
	\paragraph{Viewpoint -- Different Systems Sharing Properties:}
	
	As long as every context $c\in\IDataset^{\txt{obs}}$
	has an infinite number of observations,
	it is only relevant what properties are shared with the query
	(which is always in $\pi^*$).
	For clarity we nevertheless give the full model (and modified
	Def.\ \ref{def:iid:model}):
	\begin{definition}[IID-Model, Shared Properties]
		We call a model IID (shared properties),
		if it can be written in the following form:
		$I=\IDataset\times\IVars\times\ISample$, where $|\ISample|=\infty$,
		with symmetry-group $G=\mathfrak{S}_{\IDataset\times\ISample}$ (acting trivially
		on the middle factor $\IVars$ of $I$),
		and variables $J_{(C,v)} = C \times \{v\} \times \ISample$
		(where $C\subset \IDataset$)
		with symmetry $H_{C,v} = \mathfrak{S}_{C\times\ISample}$, where
		$\mathfrak{S}_{C} \subset \mathfrak{S}_{\IDataset}$ is embedded
		as the permutations on elements of $C$ extended by the identity on other elements
		of $\IDataset$.
	\end{definition}
	
	\begin{lemmaDef}
		Given a mz-transport setup $M^{\txt{mz}}$
		Def.\ \ref{def:iid:mztransport},
		and an intervention $P(\tilde{Y}|\PearlDo(\tilde{X}=x))$,
		then there is an 
		IID (shared properties) model (Def.\ \ref{def:model})
		$\mathcal{M}(M^{\txt{mz}})$
		with $|\IDataset^{\txt{obs}}|=\sum_{k\in\{*,1,\ldots,n\}} (|I_z^k| + 1)$,
		call its single elements $(\pi_k, Z')$ (where $Z'$ is the actual intervention,
		cf. Def.\ \ref{def:iid:mztransport}),
		$|\IDataset^{\txt{query}}|=1$, call its single element $\beta$,
		such that
		$\forall c=(\pi_k, Z')\in\IDataset^{\txt{obs}}, s\in\ISample$
		\begin{align*}
			P=P(\{\iidVarAny^{\pi_k}_{v}\}_{v\in\mathcal{O}}\setminus Z'
			|\PearlDo(Z'=z'))
			\halfquad&=\halfquad
			\shallowDistr(\{\anyVar_i\}_{i\in\{c\}\times(\mathcal{O}\setminus Z')\times\{s\}})
			\\			
			P(\tilde{Y}|\PearlDo(\tilde{X}=x))
			\halfquad&=\halfquad
			\shallowDistr(\{\anyVar_i\}_{i\in\{\beta\}\times\tilde{Y}\times\{*\}})
			\txt,
		\end{align*}
		where $\iidVarAny^{\pi_k}_{v}$ are the variables of SCM $M_k$
		(corresponding to context $\pi_k$)
		and $*$ is the single (relevant) "sample" of the query.
		This model is such that for $i=(\beta,\tilde{x},*)$ with $\tilde{x}\in\tilde{X}$
		and for $i=((\pi_k,Z'),\tilde{z}', s)$ with $\tilde{z}'\in Z'$
		the corresponding mechanism
		$f_i\in\knownFunction$ is known: It is singular at value $x$ or $z'$
		respectively.
		The $I$-graph of $\mathcal{M}(M^{\txt{mz}})$ is the disjoint union of
		its restrictions to
		$\{c\}\times\IVars\times\{s\}$
		which are, for $c=(\pi_k,Z')$ of the form $\graphClassical^L$ 
		(the shared causal graph including latents) if $Z'=\emptyset$ and
		$\graphClassicalDo^{L,Z'}$ for the intervention on $Z'$
		and its restriction to
		$\{\beta\}\times\IVars\times\{*\}$
		which is of the form $\graphClassicalDo^L$
		(for the intervention on $\tilde{X}$).
	\end{lemmaDef}
	
	The query graph (embedded into $\{\beta\}\times\IVars\times\{*\}$, which
	is as before) can be defined as in the single-context case
	(Def.\ \ref{def:iid:query_graph}).
	Similarly the classification of embeddings and
	results on definitions of identification carry over.

	\paragraph{Viewpoint -- Multi Level Structure:}
	
	This is the viewpoint taken for example by
	\citep{JCI} (even though the connection to multi-level statistics
	\citep{Gelman2006} is not made explicit).
	\begin{definition}[IID-Model, Shared Properties]
		We call a model IID (shared model),
		if it can be written in the form Def.\ \ref{def:iid:model}
		after extending the system as follows:
		Let $S:=\IVars$ be a disjoint copy of $\IVars$.
		To avoid \textbf{issues with observational support}
		(Rmk.\ \ref{rmk:kernels_uniqueness})
		we follow the philosophy of \citep{JPCMCI}:
		We think of "context" as a hidden property of the system,
		and the (know) context $\pi_k$ (the dataset of origin)
		as a proxy for this property (this also clarifies how to define
		context-\Slash{}selection-variables as random elements):
		They are not realized per-sample but rather (at most) per data-set,
		more precisely:
		Add for all $v\in\IVars$ a variable $s^0_v$ to $S$,
		and add (the single node in the $I$-graph) $s^0_v$ as a parent
		to (the infinitely many nodes in the $I$-graph) of
		$((\pi^*,Z'), v, s)$ for all $Z'\subset Z_*$ and for all $s\in\ISample$
		an of $(\beta, v, *)$ (making it auxiliary for the query-graph).
		If $D_k$ contains
		a selection-variable $S_v$, then add a index $s^k_v$ to $S$,
		and add (the single node in the $I$-graph) $s^k_v$ as a parent
		to (the infinitely many nodes in the $I$-graph) of
		$((\pi_k,Z'), v, s)$ for all $Z'\subset Z_k$ and for all $s\in\ISample$.
		If $D_k$ does \emph{not} contain
		a selection-variable $S_v$, instead add
		$s^0_v$ as a parent of these variables.
		Finally add $S$ to $I$, \ie replace $I$ by $I' = I \sqcup S$
		(a disjoint union).
		Further we assume the viewport is such that $S\cap\viewport(N)=\emptyset$ for
		all $N$.
	\end{definition}

	\subsection{Multi-Context and Transportability}
	\label{apdx:iid:mz_transport}
	
	\begin{figure}[ht]
\colorlet{context}{black!20!orange}
\begin{tabular}[c]{m{0.14\textwidth}||m{0.24\textwidth}|m{0.24\textwidth}|m{0.24\textwidth}}
	\toprule
	& $\pi^*$ (target) & $\pi^1$ & $\pi^2$\\\midrule
	Selection-\newline
	Diagram
	$D_i$
	&
	&
	\begin{tikzpicture}
		\draw (0,0) node(X) {$X$};
		\draw (-0.5,1) node(Z1) {$Z_1$};		
		\draw (1,0) node(Z2) [circle, fill=lightgray, inner sep=0.1em] {$Z_2$};	
		\draw (2,0) node(Y) {$Y$};
		\draw (-0.7,0.25) node (S1) [fill,rectangle] {};
		\draw (0.8,-0.75) node (S2) [fill,rectangle] {};
		
		\draw[->] (Z1) -- (X);		
		\draw[->] (X) -- (Z2);		
		\draw[->] (Z2) -- (Y);
		\draw[->] (S1) -- (Z1);
		\draw[->] (S2) -- (Z2);
		
		\draw[<->,blue] (Z1) edge[bend left] (X);
		\draw[<->,blue] (Z1) edge[bend left] (Z2);
		\draw[<->,blue] (Z1) edge[bend left] (Y);
	\end{tikzpicture}
	&
	\begin{tikzpicture}
		\draw (0,0) node(X) {$X$};
		\draw (-0.5,1) node(Z1) [circle, fill=lightgray, inner sep=0.1em] {$Z_1$};		
		\draw (1,0) node(Z2) {$Z_2$};	
		\draw (2,0) node(Y) {$Y$};
		\draw (1.8,-0.75) node (S3) [fill,rectangle] {};
		
		\draw[->] (Z1) -- (X);		
		\draw[->] (X) -- (Z2);		
		\draw[->] (Z2) -- (Y);
		\draw[->] (S3) -- (Y);
		
		\draw[<->,blue] (Z1) edge[bend left] (X);
		\draw[<->,blue] (Z1) edge[bend left] (Z2);
		\draw[<->,blue] (Z1) edge[bend left] (Y);
	\end{tikzpicture}
	\\\midrule
	Observations
	&
	\begin{tikzpicture}
		\draw (0,0) node(X) {$X$};
		\draw (-0.5,1) node(Z1) {$Z_1$};		
		\draw (1,0) node(Z2) {$Z_2$};	
		\draw (2,0) node(Y) {$Y$};
		
		\draw[->] (Z1) -- (X);		
		\draw[->] (X) -- (Z2);		
		\draw[->] (Z2) -- (Y);
		
		\draw[<->,blue] (Z1) edge[bend left] (X);
		\draw[<->,blue] (Z1) edge[bend left] (Z2);
		\draw[<->,blue] (Z1) edge[bend left] (Y);
	\end{tikzpicture}
	&
	\begin{tikzpicture}
		\draw (0,0) node(X) {$X$};
		\draw (-0.5,1) node(Z1) {\textcolor{context}{$Z_1'$}};		
		\draw (1,0) node(Z2) {\textcolor{context}{$Z_2'$}};
		\draw (2,0) node(Y) {$Y$};
		
		\draw[->] (Z1) -- (X);		
		\draw[->] (X) -- (Z2);		
		\draw[->] (Z2) -- (Y);
		
		\draw[<->,blue] (Z1) edge[bend left] (X);
		\draw[<->,blue] (Z1) edge[bend left] (Z2);
		\draw[<->,blue] (Z1) edge[bend left] (Y);
	\end{tikzpicture}
	&	
	\begin{tikzpicture}
		\draw (0,0) node(X) {$X$};
		\draw (-0.5,1) node(Z1) {$Z_1$};		
		\draw (1,0) node(Z2) {$Z_2$};
		\draw (2,0) node(Y) {\textcolor{context}{$Y'$}};
		
		\draw[->] (Z1) -- (X);		
		\draw[->] (X) -- (Z2);		
		\draw[->] (Z2) -- (Y);
		
		\draw[<->,blue] (Z1) edge[bend left] (X);
		\draw[<->,blue] (Z1) edge[bend left] (Z2);
		\draw[<->,blue] (Z1) edge[bend left] (Y);
	\end{tikzpicture}
	\\
	Intervention
	&
	\begin{tikzpicture}
		\draw (0.5, 1.5) node {Query:};
		
		\draw (0,0) node(X) {\textcolor{context}{$X^{\txt{do}}$}};
		\draw (1,0) node(Z2) {$Z_2$};	
		\draw (2,0) node(Y) {$Y$};
		
		\draw (1.3,0.75) node(L2) {\textcolor{blue}{$L_2$}};
		\draw (2.3,0.75) node(L3) {\textcolor{blue}{$L_3$}};
		
		\draw[->] (X) -- (Z2);		
		\draw[->] (Z2) -- (Y);
		
		\draw[->,blue] (L2) -- (Z2);
		\draw[->,blue] (L3) -- (Y);
	\end{tikzpicture}
	&
	\begin{tikzpicture}
		\draw (0.5, 1.5) node {Experiment:};
		
		\draw (0,0) node(X) {$X$};
		\draw (-0.5,1) node(Z1) {\textcolor{context}{$Z_1'$}};		
		\draw (1,0) node(Z2) {\textcolor{context}{$Z_2^{\txt{do}}$}};
		\draw (2,0) node(Y) {$Y$};
		\draw (0.8,0.5) node(L) {\textcolor{blue}{$L_2$}};
		
		\draw[->] (Z1) -- (X);			
		\draw[->] (Z2) -- (Y);
		
		\draw[<->,blue] (Z1) edge[bend left] (X);
		\draw[->,blue] (L) -- (Z1);
		\draw[<->,blue] (Z1) edge[bend left] (Y);
	\end{tikzpicture}
	&	
	\begin{tikzpicture}
		\draw (0.5, 1.5) node {Experiment:};
		
		\draw (0,0) node(X) {$X$};
		\draw (-0.5,1) node(Z1) {\textcolor{context}{$Z_1^{\txt{do}}$}};		
		\draw (1,0) node(Z2) {$Z_2$};
		\draw (2,0) node(Y) {\textcolor{context}{$Y'$}};
		
		\draw (0.3,0.75) node(L1) {\textcolor{blue}{$L_1$}};
		\draw (1.3,0.75) node(L2) {\textcolor{blue}{$L_2$}};
		\draw (2.3,0.75) node(L3) {\textcolor{blue}{$L_3$}};
		
		\draw[->] (Z1) -- (X);	
		\draw[->] (X) -- (Z2);				
		\draw[->] (Z2) -- (Y);
		
		\draw[->,blue] (L1) -- (X);
		\draw[->,blue] (L2) -- (Z2);
		\draw[->,blue] (L3) -- (Y);
	\end{tikzpicture}
	\\\midrule
	Extract
	&
		$X^{\txt{do}} \in \mathcal{F}^{\txt{known}}$,\newline
		glue query from\newline
		known and extracted.
	&
	\hspace*{1.5em}
	\begin{tikzpicture}
		\draw (1,0) node(Z2) {$\circ$};	
		\draw (2,0) node(Y) {$Y$};
		
		\draw (2.3,0.75) node(L3) {\textcolor{blue}{$L_3$}};
		
		\draw[->] (Z2) -- (Y);
		
		\draw[->,blue] (L3) -- (Y);
	\end{tikzpicture}
	&	
	\hspace*{1.5em}
	\begin{tikzpicture}
		\draw (0,0) node(X) {$\circ$};
		\draw (1,0) node(Z2) {$Z_2$};	
		
		\draw (1.3,0.75) node(L2) {\textcolor{blue}{$L_2$}};
		
		\draw[->] (X) -- (Z2);
		
		\draw[->,blue] (L2) -- (Z2);
	\end{tikzpicture}
	\\\bottomrule
\end{tabular}
		\caption{Example from \citep[Fig. 1 c,\,d (p.\,3)]{bareinboim2012TransportCompleteness} (first row),
		filled squares are the selection-variables for each selection-diagram,
		experiments are possible in $Z_2$ (for $\pi^1$) and in $Z_1$ (for $\pi^2$)
		respectively (gray circle nodes).
		Below (rows 2 and 3), are observations, observed interventions (experiments)
		and the query-graph. This is the data-available and the problem to solve (query).
		Finally the last line shows relevant extracted c-subgraphs (in this case
		both non-trivially extracted kernels are obtained from the respective experiment).
		Mechanisms with a prime (dark orange) differ compared to $\pi^*$.}
		\label{fig:mz_transport}
	\end{figure}

	\paragraph{Viewpoint -- Different Systems Sharing Properties:}
	
	\begin{lemma}[ED-Phase]
		\label{lemma:IID:ed_phase_multi}
		The ED-phase (the first two steps of Algo.\ \ref{algo:id})
		outputs, after a finite number of steps,
		all c-fragments (cf.\ Lemma \ref{lemma:iid:ed_phase})
		in one of: $\graphClassical^{\pi_k}$
		(the context $\pi_k$ affects which kernels each node is aligned to) or
		$\graphClassicalDo^{\pi_k,Z'}$ for $Z'\subset Z_k$.
		together with
		identification formulas.
	\end{lemma}
	\begin{proof}
		As for Lemma \ref{lemma:iid:ed_phase}, note that
		different contexts correspond to disconnected
		graph-components in the $I$-graph.
	\end{proof}
	
	\begin{lemma}[Identification of Do-Interventions]
		Applying algorithm \ref{algo:id}
		to a query $P(\randomVar{Y}|\PearlDo(\randomVar{X})=x)$
		returns after a finite number of steps.
		It outputs
		a non-empty set (thus at least one sound identification strategy)
		if and only if there is a no mz-shedge $(F,F')$ for
		$P(\tilde{Y}|\PearlDo(\tilde{X}=x))$.
	\end{lemma}
	\begin{proof}[Proof Sketch]
		Formal details about c-forests etc.\ are similar
		as in the proof of Lemma \ref{lemma:iid:id_of_do_interventions}.
		Additionally we use the following information:
		
		We can use an extracted c-fragment for the mz-transport query,
		if its nodes are aligned to the same kernels as in
		the query-graph (as in $\pi^*$), this explains the
		condition "1." in the mz*-shedge definition (Def.\ \ref{def:iid:shedge}).
		
		We extract c-forests of $\graphClassicalDo^{\pi_k,Z'}$,
		which have no edges into $Z'$ (the targets of interventions
		in experimental data-sets), in particular nodes in $Z'$
		are their own c-components (aligned to $\knownFunction$ however),
		but also they thus cannot be part of any other c-fragment,
		which explains the appearance of conditions "2." and "3."
		in the mz*-shedge definition (Def.\ \ref{def:iid:shedge}).
	\end{proof}
	
	\paragraph{Viewpoint -- Multi Level Structure:}
	
	This time, we will extract many results with pinned nodes:
	the hidden (potentially very complex) system properties encoded
	by selection-variables are realized once for one or multiple systems,
	potentially including $\pi^*$ and thus our query.
	
	\begin{definition}[Multi-Level Query-Graph]
		The query-graph from
		Def.\ \ref{def:iid:query_graph}
		needs some small modification.
		We add, for each node in $v\in\nodesInnerProper^{\txt{query}}$,
		a pinned\Slash{}auxiliary node
		aligned to the index $i=s_v^0\in S\subset I$.
		The underlying query thus contains a kernel $\delta_{s_v^0}$
		in this place and we ask for the target
		$\realizedDistr(\tilde{Y})$
		to be $\delta_i$-identified.
	\end{definition}
	
	Condition "1." in the mz*-shedge definition (Def.\ \ref{def:iid:shedge})
	now arises from agreement (or disagreement) of aligned indices at pinned nodes.
	In gluing $\graphObs(\graphStructural^{\txt{query}})$
	this structural graph has aligned kernels $\delta_{s_v^0}$ in $\pi^*$
	and $\delta_{s_v^0}$ or $\delta_{s_v^k}$ in $\pi^k$
	which must be matched (as before), see Algo.\ \ref{algo:svs}.
	
	As before, experimental knowledge is the reason for the appearance of
	"2." and "3.".

	\subsection{Mediation and Counterfactuals}
	\label{apdx:iid:mediation}

	\begin{figure}[ht]
		\begin{minipage}{0.3\textwidth}
			\centering
			Observations $\graphClassical$:
			\begin{tikzpicture}
				\draw (0,0) node (X) {$X$};
				\draw (1,1) node (M) {$M$};
				\draw (2,0) node (Y) {$Y$};
				
				\draw[->] (X) -- (M);
				\draw[->] (M) -- (Y);
				\draw[->] (X) -- (Y);
			\end{tikzpicture}
		\end{minipage}
		\begin{minipage}{0.3\textwidth}
			\centering
			Query-Graph:
			\begin{tikzpicture}
				\draw (0,0) node (X) {$X_{\PearlDo}^B$};
				\draw (-1,1) node (X2) {$X_{\PearlDo}^A$};
				\draw (1,1) node (M) {$M$};
				\draw (2,0) node (Y) {$Y$};
				
				\draw[->] (X2) -- (M);
				\draw[->] (M) -- (Y);
				\draw[->] (X) -- (Y);
			\end{tikzpicture}
		\end{minipage}
		\caption{Query graph yielding the mediation-formula.
		}\label{fig:iid:mediation}
	\end{figure}
	\begin{lemma}
		On IID-observations for $\graphClassical$ in Fig.\ \ref{fig:iid:mediation},
		extraction yields kernels $X=P(X)$, $M_x=P(M|X=x)$ and $Y_{x,m}=P(Y|X=x,M=m)$ plus
		$X_{\PearlDo}^A,X_{\PearlDo}^B\in\knownFunction$, singular
		at $\tilde{x}^a$ and $\tilde{x}^b$.
		The structured query in Fig.\ \ref{fig:iid:mediation}
		with nodes of interest $\{Y\}$
		is identifiable (and identified by Algo.\ \ref{algo:id}) as
		for all $B\in \sigmaAlgebraBorel{\val{Y}}$:
		\begin{align*}
			&Y_{[m,x^b]} \circ \kernelCompound{
				X_{\PearlDo}^B
				\otimes
				(M_{[x^a]}\circ X_{\PearlDo}^A)
			}(B)\\
			&=
			\int X_{\PearlDo}^B(dx^b)
			\int X_{\PearlDo}^A(dx^a)
			\int M_{x^a}(\dm)
			\int Y_{x^b, m}(\dy)
			1_B(y)\\
			&=
			\int M_{\tilde{x}^a}(\dm)
			\int Y_{\tilde{x}^b, m}(\dy)
			1_B(y)\\
			&=
			\int P(M=m|X=\tilde{x}^a)
			\left(
			\int P(Y=y|X=\tilde{x}^b, M=m)
			1_B(y)
			\dy \right)
			\dm
			\txt.
		\end{align*}
		Taking the expectation of $Y$ we get
		\begin{align*}
			&E\Big[Y_{[m,x^b]} \circ \kernelCompound{
				X_{\PearlDo}^B
				\otimes
				(M_{[x^a]}\circ X_{\PearlDo}^A)
			}\Big]\\
			&=
			\int P(M=m|X=\tilde{x}^a)
			E[Y=y|X=\tilde{x}^b, M=m]
			\dm
			\txt.
		\end{align*}
		\Ie our algorithm finds the mediation-formula for natural direct effects
		(see \eg \citep[Eq.\,8, (p.\,6)]{Pearl2001}; usually natural direct effects are
		defined as a difference of expectations, but taking the expectation in
		$Y$ computes both terms of the difference immediately).
		With additional observed confounders $Z_1, \ldots$, 
		extraction yields kernels with more arguments,
		the same confounders appear in the query, but are not in 
		the variables of interest, thus get integrated out
		as in \citep[Eq.\,8, (p.\,6)]{Pearl2001}.
	\end{lemma}
	
	Considering that we never introduced counterfactuals,
	it may seem quite surprising, that we can
	describe natural direct effects, which conventionally
	are formalized via counterfactuals.
	However, maybe more accurately, the idea behind
	natural direct effects can be directly articulated
	with our more flexible query-language, thus counterfactuals
	are not actually needed to phrase the questions we wanted to ask
	in the first place.
	Many issues arising for counter-factuals,
	see \eg \citep{robins2010alternative},
	seem related to formulating practical questions
	via counterfactuals, rather than from the original questions.
	
	From the formal perspective,
	counterfactuals are usually seen as fundamentally
	different from interventional queries, \citet{PearlBook}
	even elevates them to a separate rank on his causal ladder.
	This special place is usually justified by the observation
	that there are counterfactuals that cannot be identified
	even given full experimental knowledge
	\citep[Lemma 2 (p\,5)]{avin2005identifiability}.
	However such counter-examples usually rely
	on SCM-mechanisms $f_i(\pa_i,\eta_i)$ which are not
	injective in their respective noise-term argument.
	This type of counter-example seems impossible
	with mechanisms formulated as probability-kernels,
	as pushing forward along non-injective mappings removes
	precisely the ambiguity in the model that allows
	examples like
	\citep[Lemma 2 (p\,5)]{avin2005identifiability}.
	This raises the question if this kind of counterfactual
	is meaningful in the first place,
	which in turn puts the special role of counterfactuals in doubt.

	\subsection{Conclusion}
	
	We focused our comparison to literature on the IID case and do-interventions
	(both in queries and experimental data).
	The primary reason for this decision is, that most results in the literature
	are in these setups or build on these setups.
	The relation to do-calculus
	and selection-\Slash{}context-variables
	seems of central importance to understand -- and validate -- our
	approach in the light of well-understood problem-statements.
	While we can generalize some of the results immediately
	within the space of questions available in standard causal language
	(for example to soft-interventions; such individual problems
	have often also been studied, see \eg \citep{correa2020general}),
	the purpose of our formalism is not primarily on
	the extension to a single setup, rather it is
	in \emph{unifying} many apparently completely different
	causal questions (like mediation, transportability, etc.)
	into a simple shared language
	and to expand the scope of accessible data-structures and
	accessible queries -- not just algorithmically, but also
	in terms of what assumptions about data and what queries can be
	formulated and be made formally available in the first place.

	\section{Future Work}
	\label{apdx:future_work}
	
	Our approach is very flexible by design,
	therefore it is impossible to fully explore all its consequences
	in a single paper.
	We summarize some insights and observations relevant to future work
	here.
	This should also clarify why the scope of the present
	paper was chosen as presented.

	\subsection{Statistical Considerations}
	\label{apdx:statistica_considerations}
	
	Our approach will usually extract \emph{small}
	compounds of mechanisms, then compose (possibly many)
	such smaller pieces (via Lemma \ref{lemma:glue:mt}).
	A main reason for doing so, is the high symmetry of small compounds
	(Rmk.\ \ref{rmk:families_and_symmetries}),
	besides improving asymptotic results, these can,
	in the finite data case, be learned on a larger part of the data.
	They might (having better symmetry) also be more robust under structure
	miss-specification.
	Further, the results of \citep{guo2022efficient}
	for linear models suggest that at least for simple (low parameter-count
	per learned kernel) estimators such a local approach via many
	small pieces can actually be statistically very efficient.
	
	It is also clear however, that for very complex estimators,
	or very large queries, statistical properties can also be adversely
	affected by this approach.
	The extraction phase, and especially the query-construction phase,
	should therefore potentially take statistical finite-sample properties
	into account.
	
	Our identification results are phrased relative to probability kernels,
	similar to how classical identification results are often phrased relative
	to conditional distributions \citep{tian2002general,shpitser2006identification1,
		Perkovic2015,OSets}.
	Depending on the actual quantity of interest, one
	will want to choose one's estimator wisely, cf.\ \ref{apdx:identifiablity}.
	It is in practice of course quite complex to do so,
	as suitable choices of estimators have to be
	propagated through regular computations in a
	correct, but also statistically efficient way.
	An instructive example for constructing non-trivial
	estimators for such a compound problem can be found \eg in
	\citep{margueritte2026learning}.

	\subsection{Selection Bias and Correlated Missingness}
	\label{apdx:selection_bias}
	
	Correlated missingness refers to statistical dependence
	between what is observed and the values taken by the observed
	(see \eg \citep{rubin1976inference} and the many subsequent
	works inspired by it).
	We could model this by allowing the viewport $\viewport$ itself
	to be random, for example by adding a (always observed)
	binary variable $\randomVar{B}_i$ for each $i\in I$
	such that $\anyVar_i$ is observed if and only if $\randomVar{B}_i=1$.
	
	This problem is closely related with selection bias:
	Selecting for $\anyVar_i\in A_i$ essentially amounts
	to $\randomVar{B}_i = 1$ if $\anyVar_i \in A_i$ and
	$\randomVar{B}_i = 0$ otherwise. $\randomVar{B}_i$
	has in this case a single parent ($\anyVar_i$), and
	selection-bias is modeled as a special case of correlated missingness.
	
	Allowing for correlated missingness in this sense (\ie generally\Slash{}without
	any further restriction) is evidently
	quite non-trivial:
	Already finding a suitable notion of what an asymptotic
	limit is supposed to be is not obvious anymore.
	Providing a comprehensive treatment is beyond the scope of
	this paper.
	
	Special cases on the other hand should be rather
	accessible. Nested queries, §\ref{sec:queries_beyond_basic},
	already allow for phrasing questions for systems under selection-bias.
	In the extraction phase, selection bias from a descendant
	can be for example included by "growing" a backdoor-free graph
	to envelop the selecting descendant,
	a disintegration by the selected nodes
	of the joint distribution of the enlarged
	$\graphObs$ should be identifiable from data.
	In this case a computation in the sense of §\ref{sec:structured_kernels}	
	can no longer be complete (as nested graphs produced by the extraction-phase
	must be considered for regular computations).

	\subsection{Structure Discovery}
	\label{apdx:structure_discovery}
	
	This paper describes models and how to use them for
	the identification of queries from data.
	The separate topic of discovery of such structure
	is left to future work.
	We want to briefly remark, that the model-structure as
	described is not explicitly bounded in its complexity,
	meaning it can become more complicated with more data.
	In practice one does of course require some
	assumptions about the prior distribution of the
	model-structure that might appear.
	As a special case, the prior over model-structures
	can be chosen to allow only IID-models (with fixed
	product structure on $I$) or only time-series etc.,
	in these cases conventional causal discovery ideas apply.
	If the result of query-estimation does not depend
	on the choice of representative in a discovered Markov
	equivalence-class, it is identifiable.
	Note that Ass.\ \ref{ass:standard_borel}
	(value-spaces are standard Borel)
	is closed under finite products, so finite
	unorientable sub-graphs may be aggregated 
	into a meta-mechanism and collapsed to a single node in the $I$-graph.
	
	Generally, it is well-understood that employing
	some notion of invariant mechanism can substantially
	aid structure discovery \citep{CD-NOD,JCI,JPCMCI,peters2016causal}.
	Considering that in our approach it is formally very clear
	what this means, holds promise to better understand
	the discovery of causal structures from data.

	\subsection{Algorithmic Approach}
	\label{apdx:algorithm_details}
	
	We briefly explain what "absorbing" nodes in the main text
	should mean formally. Then we provide some
	clarification about some of the algorithms.
	The algorithmic approach is, as pointed out repeatedly in the main text,
	considered future work. We provide sketches for algorithms,
	because these seem quite helpful to better understand
	the concepts of the main text.

	\paragraph{Constructing Families:}	
	Families of embeddings can be constructed from
	individual mechanisms by successive "absorption"
	of neighbors. This sub-section briefly explains what this means concretely.
	
	\begin{lemma}[Direct Embedding of Model Mechanisms]\label{lemma:direct_embedding}
		Given a mechanism $(f_J, J, H_J, \Pa_J)$ (Def.\ \ref{def:mechanism})
		of the model $M$ (Def. \ref{def:model}),
		define the local graph $\graph$ consisting of a single inner node $y$
		and external nodes
		$x_1$, \dots, $x_{\kappa}$ (where $\kappa$ is the number of arguments
		of the mechanisms kernel)
		with exactly one proper edge from each external to the single inner node.
		Model-alignment is to
		$\mu^y := f_J$ (and assigning arguments by $\pa^y(i) = x_i$ in 
		the structured kernel, cf.\ \ref{def:structured_kernel}).
		This local graph is embedded by the mappings
		$\psi_j(y)=j$ and $\psi_j(x_k)=\PaIdx{k}_J(j)$.
		The collection $\{\psi_j\}_{j\in J}$ is a family
		of local graph embeddings anchored at $y$.
	\end{lemma}
	\begin{proof}
		We unravel definitions carefully.
		First $\graph$ is an model-aligned local graph
		(Def.\ \ref{def:local_graph_new}) by construction.
		Next we check,
		that $\psi_j$ is a local graph embedding
		(Def.\ \ref{def:local_graph_embedding}):
		
		Injectivity of $\psi_j$:
		By Def.\ \ref{def:mechanism},
		$\forall k: \PaIdx{k}(j)\neq j$ and
		$k\neq k'$ $\Rightarrow$ $\forall j\in J$:
		$\PaIdx{k}(j) \neq \PaIdx{k'}(j)$.
		Thus $\psi_j$ is injective.
		
		Proper node $I$-parents included
		\ref{def:local_graph_embedding:inner_parents_incl}:
		The only proper node is $y$, by definition of the $I$-graph
		(Def.\ \ref{def:Igraph}),
		$\psi_j(y)$ has parents $\Pa_I(\psi_j(y))=
		\{i\in I|\exists k:i = \PaIdx{k}_{J(\psi_j(y))}(\psi_j(y))\}$.
		By construction, $\psi_j(y)=j\in J$ the (region of applicability of)
		the mechanism we are working on, so by
		definition $J(\psi_j(y)) = J$,
		also by $\psi_j(y)=j$, thus $\PaIdx{k}_{J(\psi_j(y))}(\psi_j(y))\}
		=\PaIdx{k}_J(j)=\psi_j(x_k)$.
		In particular $i\in \Pa_I(\psi_j(y))$
		$\Rightarrow$ $\exists k$: $i=\psi_j(x_k)$.
		
		Proper edges equal $I$-graph edges
		\ref{def:local_graph_embedding:proper_edges_a}:
		The only proper node is $y$, by the previous point (and injectivity),
		there are exactly $\kappa$ parents $\psi_j(x_k)$
		in the $I$-graph, these are also precisely
		the proper parents in the local graph.
		
		Model Alignment \ref{def:local_graph_embedding:applicable}:
		The only proper node $y$ is aligned to $f_{J(\psi_j(y))}=f_J$
		(see above) by construction.
		There are no pinned nodes.

		We check that
		the collection $\{\psi_j\}_{j\in J}$ is a family
		of local graph embeddings (Def.\ \ref{def:local_graph_embedding_family}).
		
		Trivial on Anchor \ref{def:local_graph_embedding_family:trivial_anchor}:
		$\psi_j(y)=j$ by construction.
		
		Rigidity \ref{def:local_graph_embedding_family:rigidity}:
		Using the notation $\psi_*(n): J_0 \rightarrow I, j\mapsto \psi_j(n)$.
		Let $h\in H=H_J$ (by construction),
		and $j' = h \cdot j$ with both $j,j' \in J$,
		we have to show $\psi_{h\cdot j}(n) = h\cdot \psi_j(n)$
		(cf.\ Def.\ \ref{def:symmetry}, which explicitly
		defined equivariance restricted to subsets of $I$ in this way).
		Case 1 ($n=y$):
		$\psi_{h\cdot j}(y)=h\cdot j = \psi_{h\cdot j}(y)$
		by construction.
		Case 2 ($n=x_k$):
		$\psi_{h\cdot j}(x_k) =\PaIdx{k}_J(h\cdot j)
		=h\cdot \PaIdx{k}_J(j) = h \cdot \psi_{j}(x_k)$
		by $H_J$-equivariance of $\PaIdx{k}_J$
		(Def.\ \ref{def:mechanism}).
		
		Freeness \ref{def:local_graph_embedding_family:freeness}:
		For the (only) proper node $y$, we have to show
		$\psi_*(y)$ (see last point) is injective.
		This is the identity mapping $\id_J:J\rightarrow J, j\mapsto j$
		and clearly injective.
	\end{proof}

	\begin{definition}[Absorbed Parents]\label{def:absorb_parents}
		Given a family of local graph embeddings
		$\{\psi_j\}_{j\in J_0}$
		and an external node
		$x\in\nodesOuterProper$,
		we say a family of local graph embeddings
		$\{\psi'_j\}_{j\in J'_0}$, on $J'_0\subset J_0$,
		is an $x$-parent extension of		
		$\{\psi_j\}_{j\in J_0}$,
		if exactly one of
		\begin{equation*}
			\begin{cases}
				\nodesInner' = \nodesInner \cup \{x\}
				&\txt{ and }
				\nodesOuterFixed' = \nodesOuterFixed
				\txt{, or}\\
				\nodesInner = \nodesInner
				&\txt{ and }
				\nodesOuterFixed' = \nodesOuterFixed \cup \{x\}
				\txt.
			\end{cases}
		\end{equation*}
		A unique minimal set of outer nodes, and edges,
		are determined by Def.\ \ref{def:local_graph_embedding_family}.
		We assume $J'_0$ is chosen maximal with the above properties.
		We write for the set of $x$-parent extensions
		\begin{equation*}
			\absorbPa(\psi, x)
			\halfquad:=\halfquad
			\big\{\halfquad
			\{\psi'_j\}_{j\in J'_0}
			\halfquad\big|\halfquad
			\{\psi'_j\}_{j\in J'_0} \txt{ is an $x$-parent extension of	}
			\{\psi_j\}_{j\in J_0}
			\halfquad\big\}
			\txt.
		\end{equation*}
		We say an $x$-parent extension introduces the mechanism 
		$J\in\mathcal{J}$
		(Def.\ \ref{def:model}, Def.\ \ref{def:mechanism}),
		if $x\in\nodesInner'$ and $\mu^x = f_J$. 
		For a mechanism $J\in\mathcal{J}$,
		we write $\absorbPa(\psi, x, J)$ for the subset of $x$-parent 
		extensions
		introducing $J$ and $\absorbPa(\psi, x, *)$
		for the subset of $x$-parent extensions
		with $x\in\nodesOuterProper'$.
	\end{definition}
	
	\begin{rmk}[Freeness and Non-Uniqueness]
		\label{rmk:absorb_freeness}
		The original family $\{\psi_j\}_{j\in J_0}$,
		on the external node $x$, may be degenerate.
		In order for $\{\psi'_j\}_{j\in J'_0}$,
		to be free (Def.\ \ref{def:local_graph_embedding_family}~%
		\ref{def:local_graph_embedding_family:freeness}),
		if $x$ is absorbed as proper node,
		it might be necessary to shrink $J_0$,
		and the resulting $J'_0$ need not be unique.
		In practice the best choice for $J_0'$ will also depend on
		latent- and ancestral-structure,
		see §\ref{apdx:practical_considerations_algos}.
	\end{rmk}
	
	\begin{definition}[Absorbed Children]\label{def:absorb_children}
		Given a family of local graph embeddings
		$\{\psi_j\}_{j\in J_0}$
		and an inner node
		$l\in\nodesInner$,
		we say a family of local graph embeddings
		$\{\psi'_j\}_{j\in J'_0}$, on $J'_0\subset J_0$,
		is an $l$-child extension of		
		$\{\psi_j\}_{j\in J_0}$,
		if $\nodesInner' = \nodesInner \cup \{w\}$
		for $l \in \Pa_{\graph'}(w)$ and
		$\nodesOuterFixed = \nodesOuterFixed$.
		A unique minimal set of outer nodes, and edges,
		are determined by Def.\ \ref{def:local_graph_embedding_family}.
		We assume $J'_0$ is chosen maximal with the above properties.
		We write for the set of $l$-child extensions
		\begin{equation*}
			\absorbCh(\psi, l)
			\halfquad:=\halfquad
			\big\{\halfquad
			\{\psi'_j\}_{j\in J'_0}
			\halfquad\big|\halfquad
			\{\psi'_j\}_{j\in J'_0} \txt{ is an $x$-parent extension of	}
			\{\psi_j\}_{j\in J_0}
			\halfquad\big\}
			\txt.
		\end{equation*}
		We say an $l$-child extension introduces $(J,k)$
		if $\mu^w = f_J$ and $l$ is the $k$th parent of $w$.
		We call the $l$-child extension ancestral, if
		every valid ancestral structure $\ancestralStructure$
		for $\psi'$ contains $w \rightsquigarrow x$ to
		some $x\in\nodesOuterProper$.
		
		Finally, we call a tuple $\chi=(l,J,k)$ a relevant child,
		if an $l$-child extension introducing $(J,k)$
		with $J_0\neq\emptyset$ exists. We write
		$\RelevantCh(\psi)$ for the set of relevant children of $\psi$.
	\end{definition}

	\paragraph{Constructing C-Structures:}
	\label{apdx:construction_c_structures}	
	We briefly provide some context on how theses
	definitions tie into the
	algorithms \ref{algo:fcs}, \ref{algo:fcs_r}.

	\begin{definition}[Sub-Algorithms of Algo.\ \ref{algo:fcs}]
		\label{def:alg_bd_complete_sub}
		The following sub-routines are involved:
		\begin{itemize}
			\item
			\texttt{minimal\_hidden}$(\psi,L)$ for a family of embeddings
			and $L\subset \nodes$
			returns the set of minimal latent subsets containing $L$, \ie
			$\{ L'\subset\minLatentSets_\psi | L \subset L' \}$.
			\\
			\emph{Remark:}
			By Lemma \ref{lemma:monotonicity_of_L}
			this recursive relative search finds all $L$.
			\item
			\texttt{minimal\_ancestral\_structure}$(\psi,L,\ancestralStructure)$
			for a family of embeddings		
			and $L\in\minLatentSets_\psi$
			returns the set of minimal valid ancestral 
			structures subsets containing $\ancestralStructure$
			(in the sense of Def.\ \ref{def:subgraphs_local}).		
			\emph{Remark:}
			By Lemma \ref{lemma:monotonicity_of_Anc}
			this recursive relative search finds all $\ancestralStructure$.
			\item 
			\texttt{absorb\_parents}$(\psi, P)$ for a family of embeddings
			and $P\subset \nodesOuter$
			returns the set of local graph embeddings
			$\absorbPa(\psi, P)$ (absorbing parents one-by-one
			in any order).
			\item
			\texttt{absorb\_children\_ancestral}($\psi$,$C$)
			for a family of embeddings
			and $C\subset \nodes$
			returns the set of tuples $(\psi',J,k)$ of
			ancestral local graph embeddings $\psi'$ introducing $(J,k)$,
			\ie	$\absorbCh(\psi, C)$ such that
			no additional ancestral child can be absorbed
			(\ie continue until no more ancestral absorptions are
			possible, discard "incomplete" intermediate results).
			
			\item
			\texttt{absorb\_child}($\psi$,$\chi$)
			for a family of embeddings
			and a relevant child $\chi$ (cf.\ Def.\ \ref{def:absorb_children})
			returns the set of tuples $(\psi',J,k)$ of
			ancestral local graph embeddings $\psi'$ introducing $\chi$.
		\end{itemize}
	\end{definition}

	\paragraph{Practical Considerations:}
	\label{apdx:practical_considerations_algos}
	From a practical perspective, the main challenge
	to an algorithmic realization seems to be in
	issues of non-uniqueness and ensued problems in
	structuring intermediate results.
	
	It is also clear, that in very simple cases,
	like IID or multi-context\Slash{}transport IID
	cases an algorithmic realization,
	in light of classifications of embedded families
	(Lemma \ref{lemma:iid_psi_classification}),
	should be rather easy to attain.
	It may be helpful to restrict the setup initially.
	
	There are some obvious practical modifications that should be made:
	Instead of infinitude of data-sets, it only matters if they
	are large enough to give an estimate suitable for the problem at hand.
	Also, for most cases Algo.\ \ref{algo:id} should probably
	not be evaluated eagerly, however the "matching"
	filters of Algo.\ \ref{algo:svs} can effectively be propagated
	to absorption-operations, so these can be evaluated lazily with
	constraints making them much more attainable in practice.
	
	It also seems that instead of inflating absorption-results
	to ensure freeness (Rmk.\ \ref{rmk:absorb_freeness})
	eagerly, it may be more practical to instead
	track possibly degenerate data-sets, then account for degeneracy
	in the very end. This not only simplifies absorption-operations
	in practice, but likely estimators can, in the finite sample case,
	benefit from additional data,
	even if it is to some degree degenerate
	(enforcing freeness corresponds to what in the missing data-literature
	seems to usually be called "deletion").

	\section{Mathematical Basics and Notation}
	\label{apdx:ptheo_basics}
	
	This section revisits standard concepts from measure theory and
	probability theory, fixing notation for the remainder of the paper.
	We start by some standard results from measure-theory,
	in particular disintegrations. The focus is on writing
	statements as found in literature, under Ass.\ \ref{ass:standard_borel}
	(spaces are standard Borel) in a simplified language
	suitable for our purposes.
	Next we recall probability kernels and some commonly employed operations.
	
	Afterwards results presented in §\ref{apdx:kernels_order_independent_ops}
	and §\ref{apdx:sparsity} are still elementary, but tailored
	more specifically to the needs of the present discussion,
	which makes them hard to find in this form in the literature
	(we provide proofs for these results).
	Finally, even more customized, and to our knowledge novel in this form,
	results are presented separately in §\ref{apdx:structured_kernels}
	where kernels are systematically "structured" using certain graphs.

	\subsection{Measure-Theory}
	\label{apdx:measure_theory}

	We recall some standard definitions and translate results
	from literature,
	specifically from \citep{kallenberg1997foundations, chang1997conditioning},
	into a more uniform language -- we will only consider standard
	Borel spaces (Ass.\ \ref{ass:standard_borel}; see below)
	which will allow for a number of simplifications.

	A topological space $(X,\mathcal{O})$ has an induced Borel $\sigma$-algebra
	$\mathcal{B}=\sigma(\mathcal{O})$; the pair $(X,\mathcal{B})$ is then
	a measurable space, a measure on this measurable space is called
	a Borel measure on $X$.
	A measurable space $(X,\mathcal{B})$ is called Borel space
	if it is Borel-isomorphic
	to $([0,1], \mathcal{B}_{[0,1]})$, where $\mathcal{B}_{[0,1]}$ is the
	Borel $\sigma$-algebra induced by the standard (subspace in $\mathbb{R}$)
	topology on $[0,1]$.
	Here Borel-isomorphic means, there is a Borel-isomorphism;
	a Borel-morphism is a measurable map,
	an isomorphism is a morphism $f:X\rightarrow Y$ for which an inverse morphism $g$
	-- \ie satisfying $fg=\id_Y$ and $gf=\id_X$ -- exists.
	Sometimes the term "Borel space" is defined in other ways,
	we follow \citep[p.\,20]{kallenberg1997foundations}) here, as we only
	need this property for results from this reference.
	
	A Polish space $X$ is a separable (there exists a dense countable subset)
	topological space $(X,\mathcal{O})$
	such that there exists a metric $d$ inducing the topology $\mathcal{O}$
	which makes $(X,d)$ into a complete (Cauchy-sequences converge)
	metric space.
	Polish spaces, endowed with their Borel $\sigma$-algebra are Borel spaces
	\citep[Thm.\,A1.6]{kallenberg1997foundations}.
	
	Given a metric space $X$, a Borel measure $\lambda$ on $X$
	is called a Radon measure \citep[p.\,293]{chang1997conditioning}, if
	for all compact sets $K$, $\lambda K< \infty$
	(compact subsets of Hausdorff-spaces -- and metrizable spaces are Hausdorff -- are closed, thus in the Borel $\sigma$-algebra),
	and $\forall B\in\mathcal{B}:\lambda B = \sup_{K\subset B} \lambda K$
	(the supremum is over all compact $K\subset B$).
	A finite Borel-measure (thus in particular any probability measure)
	on a Polish space is Radon (see \eg
	\citep[p.\,293]{chang1997conditioning} and references therein).
	
	A measure $\lambda$ on $(X,\mathcal{F})$ is called $\sigma$-finite
	if there is a countable union $X=\cup_i X_i$ of $X_i \in \mathcal{F}$
	such that $\forall i: \lambda X_i <\infty$. In particular every finite measure,
	and thus every probability-measure, is also $\sigma$-finite.
	
	Finally note, that since metric spaces are Hausdorff and in 
	Hausdorff-spaces points (subsets $\{t\}$ containing a single point)
	are closed, they are in the Borel $\sigma$-algebra $\{t\}\in\mathcal{B}$.
	Further separable \emph{metric} spaces are second countable,
	\ie there is a countable basis of the topology $\mathcal{O}$, thus
	of the induced Borel $\sigma$-algebra.
	In particular both is the
	case for Polish spaces with their Borel $\sigma$-algebra
	(\ie for standard Borel spaces): individual points $\{t\}$ are measurable
	and the $\sigma$-algebra is countably generated.
	
	We will later specialize to $P=X\times Y$ and $\pi=\pi_X: P \rightarrow X$
	the coordinate-projection, but first note the following general result.
	\begin{definition}[Disintegration]\label{def:disint_gen}
		See \eg \citep[Def.\ 1 (p.\,292)]{chang1997conditioning}:
		Let $\pi:(P,\mathcal{A})\rightarrow (X,\mathcal{B})$ be a measurable map,
		$\lambda$ a $\sigma$-finite measure on $\mathcal{A}$,
		$\mu$ a $\sigma$-finite measure on $\mathcal{B}$,
		then:\\
		We say $\lambda$ has a $(\pi,\mu)$-disintegration $\{\lambda_x\}_{x\in X}$, if
		\begin{enumerate}[label=(\roman*)]
			\item
			$\lambda_x$ is a sigma-finite measure on $\mathcal{A}$,
			concentrated on $\{\pi(p)=x\}$
			(\ie $\lambda_x \{\pi(p)\neq x\}=0$ for $\mu$-almost all $t$).
			\item 
			$\forall$ non-negative measurable $f:P \rightarrow \Reals_{\geq 0}$:
			$x \mapsto \int \lambda_x(\intdp) f(p)$ is measurable.
			\item 
			$\forall$ non-negative measurable $f:P \rightarrow \Reals_{\geq 0}$:
			$\int \lambda(\intdp) f(p) = \int \mu(\dx) \int \lambda_x(\intdp) f(p)$.
		\end{enumerate}
	\end{definition}
	\begin{thm}[Existence and Uniqueness of Disintegrations]\label{thm:disint_ex}
		See \eg \citep[Thm.\ 1 (p.\, 293)]{chang1997conditioning}:
		Let $\lambda$ be a $\sigma$-finite Radon measure on a metric space $P$
		(with its Borel $\sigma$-algebra),
		and let $\pi:P\rightarrow (X,\mathcal{B})$ be measurable.
		Let $\mu$ be a $\sigma$-finite measure on $\mathcal{B}$
		that dominates $\pi_*\lambda$
		(\ie $\pi_*\lambda B=0 \Rightarrow \mu B=0$).
		If $\mathcal{B}$ is countably generated and contains
		all singleton sets, $\forall x\in X: \{x\}\in\mathcal{B}$,
		then $\lambda$ has a $(\pi,\mu)$-disintegration.
		
		Uniqueness: If $\{\lambda_x\}$, $\{\lambda'_x\}$ are two
		$(\pi,\mu)$-disintegrations of $\lambda$, then
		$\mu\{ x\in X | \lambda_x \neq \lambda_x' \}=0$,
		\ie the disintegration is $\mu$-almost everywhere unique.
	\end{thm}
	\begin{cor}
		\label{cor:disint_measure}
		Let $X$ and $Y$ be standard Borel spaces
		and $\lambda$ a probability-measure on $X\times Y$.
		Denote by $\pi_X : X \times Y \rightarrow X$ the projection
		to the first coordinate.
		Set $\mu = (\pi_X)_* \lambda$.
		Then there is a $\mu$-almost everywhere
		unique probability-kernel $\nu$ from $X$ to $Y$,
		such that $\lambda = \mu \otimes \nu$
		(cf.\ Def.\ \ref{def:kernel_product}).
	\end{cor}
	\begin{proof}
		Both $\lambda$ and $\mu$ are probability measures.
		Apply Thm.\ \ref{thm:disint_ex} to get
		$\lambda_x$ with properties (i--iii) of Def.\ \ref{def:disint_gen}.
		Define for $B\in\mathcal{B}(Y)$
		\begin{equation*}
			\nu_x(B) := \lambda_x(\{x\}\times B)
			\txt.
		\end{equation*}
		For all $B' \in \mathcal{B}(X\times Y)$
		and for all $x\in X$ by property (i):
		\begin{align*}
			\tag{$*$}
			\lambda_x(B') &= \lambda_x(B'\cap\{x\}\times Y) +
			\lambda_x(B'\cap(X-\{x\})\times Y)\\
			&\overset{\txt{(i)}}{=}
			\lambda_x(B'\cap\{x\}\times Y)
			\quad\txt{$\mu$-almost everywhere}\\
			&=\nu_x(\pi_Y(B'\cap\{x\}\times Y))
		\end{align*}
		We first use this relation to immediately conclude from
		properties (i) and (ii) (which holds in particular for indicator functions)
		that $\nu_x$ is a probability kernel (for normalization,
		see \citep[Thm.\ 2\,(iii) (p.\,294)]{chang1997conditioning}).
		Second, we use this relation to rewrite
		property (iii) with
		\begin{equation*}
			f(p)=1_{B'}(p)
			=1_{\pi_X(B')}(\pi_X(p))
			1_{\pi_Y(B'\cap\{\pi_X(p)\}\times Y)}(\pi_Y(p))
		\end{equation*}
		as follows:
		\begin{align*}
			\lambda(B') &=
			\int \lambda(\intdp) 1_{B'}(p)\\
			&\overset{\txt{(iii)}}{=}
			\int \mu(\dx) \int \lambda_x(\intdp) 1_{B'}(p)\\
			&=\int \mu(\dx) 1_{\pi_X(B')}(x) \int \lambda_x(\intdp)
			1_{\pi_Y(B'\cap\{\pi_X(p)\}\times Y)}(\pi_Y(p))\\
			&\quad\txt{($(*)$, plus $\mu=(\pi_X)_*\lambda$, so definitionally
				$\mu(\dx)=\lambda(\pi_X^{-1}(\dx))$)}\\
			&=\int \mu(\dx) 1_{\pi_X(B')}(x)
			\lambda(B'\cap\{\pi_X(p)\}\times Y) \\
			&=\int \mu(\dx) 1_{\pi_X(B')}(x)
			\nu_x(\pi_Y(B'\cap\{\pi_X(p)\}\times Y)) \\
			&=\int \mu(\dx) 1_{\pi_X(B')}(x) \int \nu_x(\dy)
			1_{\pi_Y(B'\cap\{\pi_X(p)\}\times Y)}(y)\\
			&=\int \mu(\dx) 1_{\pi_X(B')}(x) \int \nu_x(\dy)
			1_{B'}(p)
			\quad= (\mu \otimes \nu)(B')
			\txt.
		\end{align*}
		This demonstrates existence.
		For uniqueness, assume $\nu_x$ and $\nu'_x$ are two probability-kernels
		with $\lambda = \mu \otimes \nu = \mu \otimes \nu'$. 
		Then both $\lambda_x(B') := \nu_x(\pi_Y(B' \cap \{x\}\times Y))$
		and $\lambda'_x(B') := \nu'_x(\pi_Y(B' \cap \{x\}\times Y))$
		are disintegrations (by reverting the arguments given above;
		the only one that is not evidently an equivalence is
		that our definition of kernel requires measurability only on indicator
		functions, but general integrals are limits over such indicator-function
		integrals, thus not more general in this case, \ie Def.\ \ref{def:disint_gen}ii
		follows from Def.\ \ref{def:kernels}i)
		and thus $\mu\{ x\in X | \lambda_x \neq \lambda_x' \}=0$
		by the uniqueness-part of Thm.\ \ref{thm:disint_ex}.
		In particular, by equation $(*)$ also $\nu_x = \nu'_x$,
		$\mu$-almost everywhere.
	\end{proof}
	\begin{cor}\label{cor:disintegrate_apdx}
		Let $S$, $X$ and $Y$ be standard Borel spaces
		and $\lambda_s$ a probability-kernel from $S$ to $X\times Y$.
		Denote by $\pi_X : X \times Y \rightarrow X$ the projection
		to the first coordinate.
		Set $\mu_s = (\pi_X)_* \lambda_s$.
		Fix a supporting probability-measure $\vartheta$
		on $S$.
		Then there is a
		probability-kernel $\nu_{s,x}$ from $S\times X$ to $Y$,
		such that $\lambda_s = (\mu \otimes \nu)_s$
		and this kernel is $\vartheta\otimes\mu$-almost everywhere
		unique in the sense that for
		and any two such $\nu$ and $\nu'$:
		\begin{equation*}
			(\vartheta \otimes \mu)
			\{ (s,x) | \nu_{s,x} \neq \nu'_{s,x} \}
			=0
		\end{equation*}
	\end{cor}
	\begin{proof}
		Apply Cor.\ \ref{cor:disint_measure} to $\vartheta \otimes \lambda$
		with $\pi_1 : (S\times X) \times Y \rightarrow S\times X$.
	\end{proof}

	\subsection{Kernels and Basic Properties}
	\label{apdx:kernels}
	
	There are multiple equivalent ways to define kernels, we
	use the following definition; the preliminary
	definition form the main text can be made mathematically
	rigorous and then agrees with this notion (Rmk.\ \ref{rmk:kernels_measurevalued}).

	\begin{definition}[Kernels]\label{def:kernels}
		See \eg \citep[p.\,19]{kallenberg1997foundations}:
		Given measurable spaces $\measurableSpace{S}$ and $\measurableSpace{T}$,
		a mapping $\mu: \topSpace{S} \times \sigmaAlgebra{T} \rightarrow \bar\Reals_+$
		where $\bar\Reals_+ = \Reals_{\geq 0} \cup \{\infty\}$ are the
		non-negative real numbers plus an element at infinity,
		then $\mu$ is called a (probability) kernel from $\measurableSpace{S}$
		to $\measurableSpace{T}$ if
		\begin{enumerate}[label=(\roman*)]
			\item
			$\forall B\in\sigmaAlgebra{T}$:
			\begin{equation*}
				\mu_s B : \topSpace{S} \rightarrow \bar\Reals_+
				\txt,\quad
				s\mapsto \mu(s,B)
			\end{equation*}
			is $\sigmaAlgebra{S}$-measurable.
			\item 
			$\forall s\in \topSpace{S}$:
			\begin{equation*}
				\mu_s : \sigmaAlgebra{T} \rightarrow \bar\Reals_+
				\txt,\quad
				B \mapsto \mu(s,B)
			\end{equation*}
			is a (probability) measure.
		\end{enumerate}		
	\end{definition}
	\begin{rmk}\label{rmk:kernels_measurevalued}
		Given a measurable space $\measurableSpace{S}$,
		on the set $\mathcal{M}(\sigmaAlgebra{S})$ of the $\sigma$-finite measures
		on $\sigmaAlgebra{S}$,
		one defines the $\sigma$-algebra induced by mappings of the form
		$\pi_B: \mathcal{M}(\sigmaAlgebra{S}) \rightarrow \bar\Reals_+,
		\mu \mapsto \mu B$ (\ie the smallest $\sigma$-algebra for which
		all these projects are measurable w.\,r.\,t.\ the Borel
		$\sigma$-algebra on the reals; see \eg \citep[p.\,18]{kallenberg1997foundations}).
		In particular the probability measures
		$\mathcal{P}(\sigmaAlgebra{S}) = \pi_{\topSpace{S}}^{-1}(\{1\})$
		are a measurable subset
		$\mathcal{P}(\sigmaAlgebra{S}) \subset \mathcal{M}(\sigmaAlgebra{S})$,
		thus possess an induced $\sigma$-algebra.
		Then $\mu$ is a probability-kernel
		from $\measurableSpace{S}$ to $\measurableSpace{T}$
		according to Def.\ \ref{def:kernels},
		if and only if
		$\topSpace{S} \rightarrow \mathcal{P}(\sigmaAlgebra{T}), s \mapsto \mu_s$
		is measurable \citep[Lemma 1.37 (p.\,19)]{kallenberg1997foundations}.
		\Ie the intuition that probability kernels are measure-valued mappings
		is formally correct in this sense and
		the preliminary Def.\ \ref{def:kernels_prelim} can be made
		precise by fixing this $\sigma$-algebra.
	\end{rmk}
	\begin{example}\label{example:kernels}
		Given a probability-kernel $\mu$
		from $\measurableSpace{S}$ to $\measurableSpace{T}$,
		then there is a mapping $f:\topSpace{S}\times [0,1]\rightarrow \topSpace{T}$
		and a uniform random variable $\eta \sim U([0,1])$
		such that $\forall s\in \topSpace{S}: f(s, \eta)$ has distribution
		$\mu_s$ \citep[Lemma 2.22 (p.\,34)]{kallenberg1997foundations}.
		The reader accustomed to SCMs in the notation given a noise-term $\eta$
		and a mapping $f$ may thus think about probability-kernels
		as the analogue of the law of a random variable for causal mechanisms.
		The domain $\topSpace{S}$ corresponds to the value-space of the parents.\\
		\emph{Remark:} Given random variables $\vartheta_1, \ldots, \vartheta_n$,
		the new random variable $\eta$ can be chosen independent
		$\eta \independent (\vartheta_1, \ldots, \vartheta_n)$.
		
		Conversely, given a measurable mapping
		$f:\topSpace{S}\times N\rightarrow \topSpace{T}$
		and a random element $\eta \in N$,
		then mapping to the push-forward
		$s \mapsto f(s,\eta)_*P = f(s,\cdot)_*\eta_*P \in \mathcal{P}(T)$
		is measurable, thus defines a probability kernel (see remark
		\ref{rmk:kernels_measurevalued}).
	\end{example}
	\begin{rmk}
		There are in general many choices for $f$ and $\eta$ in the example
		\ref{example:kernels}
		above. This ambiguity goes far beyond simply
		transforming $\eta$ to a universal uniform. Indeed even after
		standardizing on uniform noise, besides flipping it to $1-\eta$, without continuity-requirements	on $f$ we can arbitrarily tear apart and reassemble
		the unit interval (and adapt $f$ accordingly) or choose non-injective $f$
		and so on.
		The corresponding probability-kernel may then be thought of as a
		means of representing the equivalence-class	of indistinguishable pairs $(f,P(\eta))$.
	\end{rmk}
	
	Kernels can easily be arranged together to describe more complex
	joint distributions\Slash{}kernels:
	
	\begin{definition}[Product of Kernels]
		\label{def:kernel_product}
		See \eg \citep[p.\,19]{kallenberg1997foundations}.
		Given measurable spaces $\measurableSpace{S}$, $\measurableSpace{T}$,
		$\measurableSpace{U}$ and
		(probability-)kernels
		$\mu: \topSpace{S}\times\sigmaAlgebra{T} \rightarrow \mathbb{R}_{\geq 0}$
		and
		$\nu: \topSpace{S}\times\topSpace{T}\times\sigmaAlgebra{U} 
		\rightarrow \mathbb{R}_{\geq 0}$,
		define the product
		$\mu \otimes \nu : \topSpace{S}\times(\sigmaAlgebra{T}\otimes\sigmaAlgebra{U})
		\rightarrow \mathbb{R}_{\geq 0}$
		as the (probability-)kernel given by
		\begin{equation*}
			\forall B \in \sigmaAlgebra{T}\otimes\sigmaAlgebra{U}:\quad
			\kernelCompound{\mu \otimes \nu}(s,B)
			=
			\int \mu(s,\dt) \int \nu(s,t,\du) 1_B(t, u)
			\txt.
		\end{equation*}
	\end{definition}	
	\begin{example}\label{example:kernel_product}
		Given
		probability-kernels $\mu$ from $\{*\}$ (the one-element set)
		to $\topSpace{S}$, and $\nu$ from $\{*\}\times \topSpace{S}=\topSpace{S}$
		to $\topSpace{T}$,
		then by example \ref{example:kernels},
		there are mappings $g:\{*\} \rightarrow \topSpace{S}$,
		$f:\topSpace{S}\times [0,1]\rightarrow \topSpace{T}$
		and independent uniform
		random variables $\eta_x, \eta_y \sim U([0,1])$
		such that $\randomVar{X}:=g(\eta_x)$ has distribution $\mu$ and
		$\forall s\in \topSpace{S}: f(s, \eta_y)$ has distribution
		$\nu_s$.
		If we define $\randomVar{Y} := f(\randomVar{X}, \eta_y)$,
		then the pair $(\randomVar{X},\randomVar{Y})$ has joint distribution
		$\mu \otimes \nu$, as can be verified by a simple computation.
		
		The reader accustomed to SCMs in the notation given an noise-terms $\eta$
		and a mappings $f$ may thus think about products of probability-kernels
		as a generalization ($\mu$ need not, in general,
		be defined on the trivial space $\{*\}$)
		of the law of the joint distribution of
		a causal model $\randomVar{X} \rightarrow \randomVar{Y}$.
		Concerning the choice of domain for $\mu$,
		note that in the case of $\nu$ above
		(which does have a non-trivial domain $\topSpace{S}$),
		$\nu$ encodes $P(\randomVar{Y}|\randomVar{X})$, and its domain corresponds to
		the value-space of $\randomVar{X}$.	
		Similarly, the domain of $\mu$ relates
		to conditioning on causal ancestors,
		\emph{if a suitable joint distribution to be conditioned on exists}
		(already for non-stationary time-series, this is typically not the case,
		see example \ref{example:random_walk}).
	\end{example}	
	\begin{lemmaDef}[Properties of Products]\label{lemma:kernel_product_properties}
		Products have the following properties:
		\begin{enumerate}[label=(\roman*)]
			\item
			Associativity:
			\begin{equation*}
				\xi \otimes \kernelCompound{\mu \otimes \nu}
				=
				\kernelCompound{\xi \otimes \mu} \otimes \nu
			\end{equation*}
			\item\label{lemma:kernel_product_properties:marginalize_right}
			Trivial marginalization on the right
			($U\in\sigmaAlgebra{U}$ is the full space):
			\begin{equation*}
				\forall B \in \mathcal{T}:\quad
				\kernelCompound{\mu \otimes \nu} (B \times \topSpace{U})
				=
				\mu(B)
			\end{equation*}
			We will write $\marginalizeRight(\mu\otimes \nu)=\mu$
			for this operator, uniquely computing $\mu$ from $\mu\otimes \nu$.
		\end{enumerate}
	\end{lemmaDef}
	\begin{proof}
		(i)
		Use $\pi_{uw} : \topSpace{T}\times\topSpace{U}\times\topSpace{W}
		\rightarrow \topSpace{U}\times\topSpace{W}$
		and set
		$B_t := \pi_{uw}(B\cap \{t\}\times\topSpace{U}\times\topSpace{W})$,
		then
		by definition $\forall B \in
		\sigmaAlgebra{T}\otimes\sigmaAlgebra{U}\otimes\sigmaAlgebra{W}
		=\sigmaAlgebraBorelRaw{\topSpace{T}\times\topSpace{U}\times\topSpace{W}}$
		(using Ass.\ \ref{ass:standard_borel}):
		\begin{align*}
			\kernelCompound{\xi \otimes \kernelCompound{\mu \otimes \nu}}(s,B)
			&=
			\int \xi(s,\dt) \int \kernelCompound{\mu \otimes \nu}((s, t), \du\dw)
			1_B(t, u, w)\\
			&=
			\int \xi(s,\dt) \int \kernelCompound{\mu \otimes \nu}((s, t), \du\dw)
			1_{B_t}(u, w)\\
			&= \int \xi(s,\dt)
			\kernelCompound{\mu \otimes \nu}((s, t), B_t)\\
			&= \int \xi(s,\dt)
			\int \mu((s,t), \du) \int \nu((s, t, u), \dw) 1_{B_t}(u, w)\\
			&= \int \xi(s,\dt)
			\int \mu((s,t), \du) \int \nu((s, t, u), \dw) 1_{B}(t, u, w)\txt.
		\end{align*}
		A similar argument expands
		$\kernelCompound{\xi \otimes \mu} \otimes \nu$ into the same
		("un-associated") expression.
		
		(ii)
		By Ass.\ \ref{ass:standard_borel},
		$\sigmaAlgebra{T}\otimes\sigmaAlgebra{U}
		=\sigmaAlgebraBorelRaw{\topSpace{T}\times\topSpace{U}}$,
		in particular $\forall B\in\sigmaAlgebra{T}$,
		also $B \times \topSpace{U}\in \sigmaAlgebra{T}\otimes\sigmaAlgebra{U}$
		is measurable (so the expression is well-defined).
		By definition,
		\begin{align*}
			\kernelCompound{\mu \otimes \nu}(s,B\times\topSpace{U})
			&=
			\int \mu(s,\dt) \int \nu(s,t,\du) 1_{B\times\topSpace{U}}(t, u)\\
			&=
			\int \mu(s,\dt) \int \nu(s,t,\du) 1_{B}(t)\\
			&=
			\int \mu(s,\dt)  1_{B}(t) \int \nu(s,t,\du)\\
			&=
			\int \mu(s,\dt) 1_{B}(t) = \mu(s,B)
			\txt.
		\end{align*}
		The last line is obtained by $\nu$ being a probability-kernel
		(normalized to $1$) and definition of the integral.
	\end{proof}

	\begin{definition}[Composition of Kernels]
		\label{def:kernel_composition}
		Given measurable spaces $\measurableSpace{S}$, $\measurableSpace{T}$,
		$\measurableSpace{U}$ and
		(probability-)kernels
		$\mu: \topSpace{S}\times\sigmaAlgebra{T} \rightarrow \mathbb{R}_{\geq 0}$
		and
		$\nu: \topSpace{S}\times\topSpace{T}\times\mathcal{U} 
		\rightarrow \mathbb{R}_{\geq 0}$,
		define the composition
		$\kernelComposition{\nu}{\mu} : \topSpace{S} \times\sigmaAlgebra{U}
		\rightarrow \mathbb{R}_{\geq 0}$
		as the (probability-)kernel given by
		\begin{equation*}
			\forall B \in \sigmaAlgebra{U}:\quad
			(\kernelComposition{\nu}{\mu})(B) :=
			\kernelCompound{\mu \otimes \nu}(\topSpace{T} \times B)
			\text.
		\end{equation*}
	\end{definition}
	\begin{example}\label{example:kernel_composition}
		Extending example \ref{example:kernel_product},
		note that the composition $\kernelComposition{\nu}{\mu}$
		\begin{align*}			
			(\kernelComposition{\nu}{\mu})(B) &=
			\kernelCompound{\mu \otimes \nu}(\topSpace{T} \times B) \\
			&=\int \mu(\dt) \int \nu(t,\ds) 1_B(s)
		\end{align*}
		simply integrates out $\mu$ corresponding to $\randomVar{X}$ in
		example \ref{example:kernel_product}.
		The result is thus simply $P(\randomVar{Y})$.
		Why is this a composition? Since $\randomVar{Y}$ causally depends
		on $\randomVar{X}$, it also "sees" the noise $\eta_X$
		injected at $\randomVar{X}$, and the distribution of $\randomVar{Y}$
		is only fully described from a suitable combination of the \emph{mechanisms}
		$\mu$ at $\randomVar{X}$ and $\nu$ at $\randomVar{Y}$.
	\end{example}
	\begin{lemma}[Properties of Compositions]\label{lemma:kernel_composition_properties}
		Compositions are associative:
		\begin{equation*}
			\kernelComposition{(\kernelComposition{\nu}{\mu})}{\xi}
			=
			\kernelComposition{\nu}{(\kernelComposition{\mu}{\xi})}
		\end{equation*}
		\emph{Warning:}
		The emerging ring-like algebraic structure does not right-distribute,
		generally
		\begin{equation*}
			\kernelComposition{\kernelCompound{\mu \otimes \nu}}{\xi}
			\quad\neq\quad
			(\kernelComposition{\mu}{\xi}) \otimes
			(\kernelComposition{\nu}{\xi})
			\txt{, see Lemma \ref{lemma:kernel_anticausal_properties}ii,
				\ref{lemma:causal_transformations}iii.}
		\end{equation*}
	\end{lemma}
	\begin{proof}
		We again make use of
		$\sigmaAlgebra{T}\otimes\sigmaAlgebra{U}\otimes\sigmaAlgebra{W}
		=\sigmaAlgebraBorelRaw{\topSpace{T}\times\topSpace{U}\times\topSpace{W}}$
		(by Ass.\ \ref{ass:standard_borel}).
		By definition $\forall B\in\sigmaAlgebra{W}$,
		\begin{align*}
			\big(\kernelComposition{(\kernelComposition{\nu}{\mu})}{\xi}\big)(s,B)
			&=			
			\kernelCompound{\xi\otimes(\kernelComposition{\nu}{\mu})}
			(s,\topSpace{T}\times B)\\
			&=
			\int \xi(s,\dt)
			\int (\kernelComposition{\nu}{\mu})((s,t), \dw) 1_{\topSpace{T}\times B}(t,w)\\
			&=
			\int \xi(s,\dt)
			\int (\kernelComposition{\nu}{\mu})((s,t), \dw) 1_{B}(w)\\
			&=
			\int \xi(s,\dt)
			(\kernelComposition{\nu}{\mu})((s,t), B)\\
			&=
			\int \xi(s,\dt)
			\kernelCompound{\mu\otimes\nu}((s,t), \topSpace{U}\times B)\\
			&=
			\int \xi(s,\dt)
			\int \mu((s,t), \du)
			\int \nu((s,t,u), \dw)
			1_{\topSpace{U}\times B}(u,w)\\
			&=
			\int \xi(s,\dt)
			\int \mu((s,t), \du)
			\int \nu((s,t,u), \dw)
			1_{B}(w)
		\end{align*}
		Again, the same "un-associated" form can similarly be obtained
		for $\kernelComposition{\nu}{(\kernelComposition{\mu}{\xi})}$.
	\end{proof}
	
	Knowing $\mu$ and $\nu$, we can obviously obtain $\mu \otimes \nu$,
	interestingly the opposite is also true (beware of the warning
	in the previous lemma \ref{lemma:kernel_composition_properties} however)
	via disintegrations (see \citep{chang1997conditioning} for a very nice overview):
	\begin{lemmaDef}[Disintegration]
		\label{lemma:disintegration}
		Let $\measurableSpaceBorel{S}$, $\measurableSpaceBorel{X}$
		and $\measurableSpaceBorel{Y}$ be standard Borel spaces
		(Ass.\ \ref{ass:standard_borel})
		and $\lambda_s$ a probability-kernel from $\topSpace{S}$
		to $\topSpace{X}\times \topSpace{Y}$.
		Denote by $\pi_{\topSpace{X}} :
		\topSpace{X} \times \topSpace{Y} \rightarrow \topSpace{X}$ the projection
		to the first coordinate and
		set $\mu_s = (\pi_{\topSpace{X}})_* \lambda_s$
		(then $\mu_s(B) = \lambda_s(B\times \topSpace{Y})$, so this marginalizes the
		$\topSpace{Y}$-coordinate away).
		Then there is a
		probability-kernel $\nu$ from $\topSpace{S}\times \topSpace{X}$
		to $\topSpace{Y}$, such that $\lambda = \mu \otimes \nu$.
		This kernel is unique almost everywhere relative to
		(any fixed choice of) a
		supporting probability measure $\vartheta$ on $\topSpace{S}$,
		in the sense that
		for any two such $\nu$ and $\nu'$:
		\begin{equation*}
			\kernelCompound{\vartheta \otimes \mu}
			\{ (s,x) | \nu_{s,x} \neq \nu'_{s,x} \}
			=0
		\end{equation*}
		We will write $\disintLeft(\mu\otimes \nu)=\nu$
		for this operator, uniquely (in the sense defined above;
		see also Rmk.\ \ref{rmk:kernels_uniqueness})
		computing $\nu$ from $\mu\otimes \nu$.
	\end{lemmaDef}
	\begin{proof}
		This was already proofed as Cor.\ \ref{cor:disintegrate_apdx},
		as a consequence of standard results \citep{chang1997conditioning}
		in §\ref{apdx:measure_theory} above.
	\end{proof}
	
	\begin{example}
		\label{example:disint}
		If we think of $\lambda$ as a joint distribution
		$P(\randomVar{X},\randomVar{Y})$, then
		$\mu$ (by marginalization from the right, lemma
		\ref{lemma:kernel_product_properties}ii)
		contains precisely the information about $P(\randomVar{X})$
		(see example \ref{example:kernel_product}).
		By example \ref{example:kernels}
		the constructed kernel $\nu$ (for trivial $S=\{*\}$)
		corresponds to $P(\randomVar{Y}|\randomVar{X})$, so
		disintegration states in this language
		that the (regular version \citep[Thm.\ 5.3 (p.\,84)]{kallenberg1997foundations} of) the conditional distribution
		$P(\randomVar{Y}|\randomVar{X})$
		is unique and thus can be reconstructed from knowledge of the joint distribution;
		see also \citep[Thm.\ 5.4 (p.\,85)]{kallenberg1997foundations}.
	\end{example}
	
	It does not generally make sense to talk about commutativity
	of $\otimes$ as domains of the involved kernels are not generally suitable
	for exchanging the order, we will come back to this in Lemma
	\ref{lemma:causal_transformations}i.
	For now, we make the following observation:
	\begin{definition}[Anti-Causal Disintegration]\label{def:kernel_anti_causal_disint}
		Given are standard Borel spaces
		$\measurableSpaceBorel{S}$, $\measurableSpaceBorel{T}$,
		$\measurableSpaceBorel{U}$ and
		(probability-)kernels
		$\mu: \topSpace{S}\times\sigmaAlgebraBorel{T} \rightarrow \mathbb{R}_{\geq 0}$
		and
		$\nu: \topSpace{S}\times\topSpace{T}\times\sigmaAlgebraBorel{U} 
		\rightarrow \mathbb{R}_{\geq 0}$.
		For $B \subset \topSpace{X}\times \topSpace{Y}$ denote by
		$B^t = \{ (y,x) \in \topSpace{Y}\times \topSpace{X}| (x,y)\in B \}$
		the "transposed" set\footnote{%
			Note that $\sigmaAlgebraBorelRaw{\topSpace{X}\times \topSpace{Y}}
			=\sigmaAlgebraBorel{X} \times \sigmaAlgebraBorel{Y}$,
			so transposing in this sense maps measurable sets to measurable sets,
			\ie it induces a well-defined mapping
			$(\cdot)^t:\sigmaAlgebraBorelRaw{
				\topSpace{X}\times \topSpace{Y}}
			\rightarrow\sigmaAlgebraBorelRaw{
				\topSpace{Y}\times \topSpace{X}}$;
			indeed the element-wise transposition $(x,y)\mapsto(y,x)$ is a Borel
			isomorphism on the Borel product.}.
		For a kernel $\lambda$ from $\topSpace{S}$
		to $\topSpace{X}\times \topSpace{Y}$ denote
		\begin{equation*}
			\forall B \in \sigmaAlgebraBorelRaw{
				\topSpace{Y}\times \topSpace{X} }:\quad
			\lambda_s^t(B) := \lambda_s(B^t)
			\txt.
		\end{equation*}
		If $\lambda = \mu \otimes \nu$, then lemma \ref{lemma:disintegration}
		produces a disintegration of $\lambda^t = \xi \otimes \vartheta$.
		We will denote $\xi:=\marginalizeRight(\lambda^t)$
		by $\kernelComposition{\nu}{\mu}$ (as will be justified
		by Lemma \ref{lemma:kernel_anticausal_properties}i)
		and $\vartheta:=\disintLeft(\lambda^t)$ by $(\mu|\nu)$, thus (definitionally)
		\begin{equation*}
			\kernelCompound{\mu \otimes \nu}^t = (\kernelComposition{\nu}{\mu}) \otimes (\mu|\nu)
			\txt,
		\end{equation*}
		where $\kernelComposition{\nu}{\mu}$ is a probability-kernel
		from $\topSpace{S}$ to $\topSpace{U}$
		and $(\mu|\nu)$ is a probability-kernel from
		$\topSpace{S}\times\topSpace{U}$ to $\topSpace{T}$.
	\end{definition}	
	\begin{example}\label{example:kernel_anticausal_disintegr}
		Given a joint distribution for $(\randomVar{X},\randomVar{Y})$,
		from a causal model $\randomVar{X} \rightarrow \randomVar{Y}$,
		there is a "causal" conditioning
		$P(\randomVar{Y}|\randomVar{X})$
		which is given directly by the causal mechanism $\nu$ \emph{alone}.
		But there is of course also an anti-causal $P(\randomVar{X}|\randomVar{Y})$
		that is uniquely (in a suitable sense) determined by the joint distribution.
		The joint distribution is $P(\randomVar{X},\randomVar{Y})
		= P(\randomVar{X}|\randomVar{Y})P(\randomVar{Y})$ --
		however $P(\randomVar{Y})$ is not trivially given by the mechanism on $\randomVar{Y}$,
		rather by the composition
		$\kernelComposition{\nu}{\mu}$ (cf.\ example \ref{example:kernel_composition}).
		Similarly the anti-causal disintegration is a property of the
		joint distribution and depends on \emph{both} $\mu$ and $\nu$.
	\end{example}
	\begin{lemma}[Properties of Anti-Causal Disintegrations]
		\label{lemma:kernel_anticausal_properties}
		The kernels $\xi$ and $\vartheta =: (\mu|\nu)$ in 
		Def.\ \ref{def:kernel_anti_causal_disint} have the following properties:
		\begin{enumerate}[label=(\roman*)]
			\item 
			$\xi = \kernelComposition{\nu}{\mu}$ (justifying the
			notation used in Def.\ \ref{def:kernel_anti_causal_disint}).
			\item\label{lemma:kernel_anticausal_properties:compose} 
			$\kernelComposition{(\mu|\nu)}{(\kernelComposition{\nu}{\mu})} = \mu$ and
			$(\mu \circ \xi) \otimes (\nu \circ (\xi|\mu))
			= \kernelCompoundConfounded{ \mu \otimes \nu } \circ \xi$.
		\end{enumerate}
	\end{lemma}
	\begin{proof}
		(i)
		This follows immediately from
		Lemma \ref{lemma:kernel_product_properties}\,(ii)
		and the definition of composition (using again Ass.\ \ref{ass:standard_borel}):
		\begin{align*}
			\xi(s,B) &= \kernelCompound{\xi \otimes \vartheta}(s,B\times \topSpace{T})\\
			&= \kernelCompound{\mu \otimes \nu}^t(s,B\times \topSpace{T})\\
			&= \kernelCompound{\mu \otimes \nu}(s, \topSpace{T}\times B)\\
			&= (\kernelComposition{\nu}{\mu})(s,B)\txt.
		\end{align*}
		
		(ii) 
		By definition of composition and the anti-causal disintegration,
		and Lemma \ref{lemma:kernel_product_properties}\,(ii):
		\begin{align*}
			\big((\mu|\nu) \circ (\nu \circ \mu)\big)(s,B)
			&= \kernelCompound{(\nu \circ \mu) \otimes (\mu|\nu)}
			(s,\topSpace{U} \times B)\\
			&= \kernelCompound{\mu \otimes \nu}^t(s,\topSpace{U} \times B)\\
			&= \kernelCompound{\mu \otimes \nu}(s,B \times \topSpace{U})\\
			&= \mu(s,B)
		\end{align*}		
		Given that compositions are defined as marginalized $\otimes$-products,
		it should not surprise us, that also the second statement works
		(by using associativity, Lemma \ref{lemma:kernel_product_properties}):
		\begin{align*}
			(\mu \circ \xi) \otimes (\nu \circ (\xi|\mu))(s, B\times B')
			&=
			(\mu \circ \xi) \otimes \kernelCompound{(\xi|\mu) \otimes \nu}
			(s, B\times \topSpace{T} \times B')\\
			&=
			\kernelCompound{(\mu \circ \xi) \otimes (\xi|\mu)} \otimes \nu
			(s, B\times \topSpace{T} \times B')\\
			&=
			\kernelCompound{\xi \otimes \mu}^t \otimes \nu
			(s, B\times \topSpace{T} \times B')\\
			&=
			\kernelCompound{\xi \otimes \mu} \otimes \nu
			(s, \topSpace{T} \times B \times B')\\
			&=
			\xi \otimes \kernelCompound{ \mu \otimes \nu }
			(s, \topSpace{T} \times B \times B')\\
			&=
			\kernelCompoundConfounded{ \mu \otimes \nu } \circ \xi
			(s, B \times B')	
		\end{align*}
		Note, that by definition
		$(\nu \circ \mu) \otimes (\mu|\nu) = \kernelCompound{\mu \otimes \nu}^t$,
		so this statement also tells us what happens if
		we replace the $\otimes$-product in the anti-causal disintegration
		by a confounded one.
		Recall, that we said earlier, that disintegrating \emph{confounded}
		terms is non-trivial -- it does not simply produce the right-hand side kernel;
		this latest statement also tells us what the
		non-trivial result actually is:
		The disintegration of the right-hand side is of course the same
		as the disintegration of the (equivalent) expression on the left-hand side
		-- thus we have also found, that disintegrating
		$\kernelCompoundConfounded{ \mu \otimes \nu } \circ \xi$
		will yield $(\nu \circ (\xi|\mu))(s,t)$.
		This does make intuitive sense:
		Conditioning on $\mu$ selects values of $\xi$, namely
		$(\xi|\mu)$, and the effect of confounding is
		precisely that conditioning $\nu$ on $\mu$ will see
		this selection of its (other, marginalized) parent $\xi$
		on top of the causal effect.
	\end{proof}

	\subsection{Order Independent Operators}
	\label{apdx:kernels_order_independent_ops}
	
	In general the order of variables in a joint distribution
	does not matter -- yet properties like the factorization
	of a joint distribution into a marginal and a conditional
	remain untouched --
	and similarly we can reorder $\otimes$-products
	in the same way we have already seen for
	anti-causal disintegrations (Def.\ \ref{def:kernel_anti_causal_disint}).
	
	\begin{definition}[Transposition]
		\label{def:transposition}
		Given a $\otimes$-product
		$\mu_1 \otimes \ldots \otimes \mu_n$ with $n$ factors
		from $\topSpace{X}$ to $\topSpace{Y}_1 \times \ldots \times \topSpace{Y}_n$
		and a permutation $T\in \mathfrak{S}_n$,
		there is a Borel-isomorphism
		\begin{equation*}
			T_B
			:
			\topSpace{Y}_{T(1)} \times \ldots \times \topSpace{Y}_{T(n)}
			\rightarrow
			\topSpace{Y}_1 \times \ldots \times \topSpace{Y}_n,
			(y_1, \ldots, y_n)
			\mapsto
			(y_{T^{-1}(1)}, \ldots, y_{T^{-1}(n)})
		\end{equation*}
		and we write
		$(\mu_1 \otimes \ldots \otimes \mu_n)^T(x, B) := (\mu_1 \otimes \ldots \otimes \mu_n)(x, T_B(B))$, which is a kernel
		from $\topSpace{X}$
		to $\topSpace{Y}_{T(1)} \times \ldots \times \topSpace{Y}_{T(n)}$.
	\end{definition}
	
	We defined marginalizations on the right
	Lemma \ref{lemma:kernel_product_properties}\,%
	\ref{lemma:kernel_product_properties:marginalize_right}
	and disintegrations on the left
	Lemma \ref{lemma:disintegration}
	however in the light of anti-causal disintegrations and more generally
	reordering by transposition, it is unclear how exactly
	the order in which multiple marginalizations or disintegrations
	on different orders of the variables are executed matters.
	The present section shows, that there is a suitable order-independent
	notion, that clarifies this behavior.
	
	First note, that marginalization more generally can immediately be
	defined as follows:
	\begin{definition}[Marginalizing Out]
		\label{def:kernel_marginalization}
		Given a product of kernels $\mu_1 \otimes \ldots \otimes \mu_n$
		(with respective target $\topSpace{X}_k$)
		and a subset $L \subset \{ 1, \ldots, n \}$
		we define a marginalization of $L$, denoted
		$\marginalize{\mu_1 \otimes \ldots \otimes \mu_n}{L}$
		on
		$\prod_{k\notin L} B_k \in \sigmaAlgebraBorelRaw{\prod_{k\notin L} \topSpace{X_k}}$ as
		\begin{align*}
			\marginalize{\mu_1 \otimes \ldots \otimes \mu_n}{L}(s,\prod_{k\notin L} B_k)
			\halfquad&:=\halfquad
			\kernelCompound{\mu_1 \otimes \ldots \otimes \mu_n}
			(s,\prod W_k)\\
			\txt{where }
			W_k &=
			\begin{cases}
				\topSpace{X}_k & \txt{ if }k\in L\\
				B_k & \txt{ if }k\notin L
				\txt.
			\end{cases}
		\end{align*}
	\end{definition}
	\begin{lemma}[Properties of Marginalization]
		\label{lemma:kernel_marginalization_safe}
		If $L=L_1\cup L_2$, then
		\begin{align*}
			\marginalize{\mu_1 \otimes \ldots \otimes \mu_n}{L}
			&= \marginalize{
				\marginalize{\mu_1 \otimes \ldots \otimes \mu_n}{L_1}}
			{L_2\setminus L_1}\\
			&= \marginalize{
				\marginalize{\mu_1 \otimes \ldots \otimes \mu_n}{L_2}}
			{L_1\setminus L_2}
			\txt.
		\end{align*}
	\end{lemma}
	\begin{proof}
		This follows immediately from plugging in the definitions.
	\end{proof}

	\begin{lemmaDef}[Joint from Marginal and Conditional]
		\label{def:disint_product}
		Given a product $\mu = \mu_1 \otimes \ldots \otimes \mu_m$
		of kernels
		and a subset $C \subset \{ 1, \ldots, m \}$,
		label its complement as $L:= \{ 1, \ldots, m \} \setminus C$,
		then there is a unique (in the sense of Lemma \ref{lemma:disintegration})
		kernel $\exists!\disint_C(\mu)$, such that
		\begin{equation*}
			\marginalize{\mu}{L}
			\otimes
			\disint_C(\mu)
			=
			\mu^T
			\txt.
		\end{equation*}
	\end{lemmaDef}
	\begin{proof}
		Part \ref{lemma:product_operator_properties:marg_prod} of
		Lemma \ref{lemma:product_operator_properties} is proved independently
		of this result, and will be used here as result $(*)$.
		The (non-cyclic) logical proof order is thus:
		Show Lemma \ref{lemma:product_operator_properties}
		\ref{lemma:product_operator_properties:marg_prod},
		then show the present Lemma \ref{def:disint_product},
		then show Lemma \ref{lemma:product_operator_properties}
		\ref{lemma:product_operator_properties:disint_prod}
		and \ref{lemma:product_operator_properties:disint_marg}.

		Inductively over $|C|$.
		Start of induction ($C=\emptyset$):		
		If $C=\emptyset$, then define $\disint_\emptyset(\mu) = \mu$.
		Further this can of course be written as a product of $m-|C|=m$ terms
		and the last $m$ (\ie all) factors remain the same.
		The identity functional is regular by definition.
		
		Inductive step ($|C|=k\mapsto |C|=k+1$, $k\geq 0$):
		We can write $C=\{ c_1, \ldots, c_{k+1} \}$
		with $c_1 < c_2 < \ldots < c_{k+1}$.
		Define $C' := \{ c_1, \ldots, c_k \}$
		and set $c:= c_{k+1}$.
		By inductive hypothesis, there is
		$\disint_{C'}(\mu)$ with
		$\marginalize{\mu}{L'}\otimes\disint_{C'}(\mu)=\mu$,
		where $L' = L\cup\{c\}$ is the complement of $C'$.
		Further $\disint_{C'}(\mu)=\nu_1 \otimes \ldots \otimes \nu_{c_k - k}
		\otimes \mu_{c_k +1} \otimes \ldots \otimes \mu_m$.
		With $c_{k+1} \geq c_k+1$ the corresponding term $\mu_{c_{k+1}}$
		is in the second part.
		We use anti-causal disintegration to write this as
		\begin{align*}
			\disint_{C'}(\mu)
			\halfquad&=\halfquad
			\mu_{c_{k+1}} \circ
			\kernelCompound{
				\nu_1 \otimes \ldots \otimes \nu_{c_k - k}
				\otimes
				\mu_{c_k +1} \otimes \ldots \otimes \mu_{c_{k+1} -1}
			}\\
			&\otimes
			\kernelCompound{
				(\nu_1 \otimes \ldots \otimes \nu_{c_k - k}
				\otimes
				\mu_{c_k +1} \otimes \ldots \otimes \mu_{c_{k+1} -1}
				|\mu_{c_{k+1}})
				\otimes
				\mu_{c_{k+1} +1} \otimes \ldots \otimes \mu_{m}
			}
			\txt.
		\end{align*}
		This can be elementarily disintegrated
		(from the left; Lemma \ref{lemma:disintegration}),
		\ie by a regular functional we obtain:
		\begin{equation*}
			\disint_{C}(\mu)
			\halfquad:=\halfquad
			(\nu_1 \otimes \ldots \otimes \nu_{c_k - k}
			\otimes
			\mu_{c_k +1} \otimes \ldots \otimes \mu_{c_{k+1} -1}
			|\mu_{c_{k+1}})
			\otimes
			\mu_{c_{k+1} +1} \otimes \ldots \otimes \mu_{m}
			\txt.
		\end{equation*}
		This is regular and has the required product-form, it remains to check,
		that indeed
		$\marginalize{\mu}{L}\otimes\disint_{C}(\mu)=\mu$.
		First, note that
		$\mu_{c_{k+1}} \circ
		\kernelCompound{
			\nu_1 \otimes \ldots \otimes \nu_{c_k - k}
			\otimes
			\mu_{c_k +1} \otimes \ldots \otimes \mu_{c_{k+1} -1}
		} = \marginalize{\disint_{C'}(\mu)}{L}$.
		\begin{align*}
			\mu
			&=
			\marginalize{\mu}{L'}\otimes\disint_{C'}(\mu)
			&\\
			&=
			\marginalize{\mu}{L'}
			\otimes\marginalize{\disint_{C'}(\mu)}{L}
			\otimes\disint_{C}(\mu)
			&\\
			&=
			\marginalize{\marginalize{\mu}{L'} \otimes \disint_{C'}(\mu)}{L}
			\otimes\disint_{C}(\mu)
			&\txt{use $(*)$ with $L^A=\emptyset$, $L^B=L$}
			\\ 
			&=
			\marginalize{\mu}{L}
			\otimes\disint_{C}(\mu)\txt.&
		\end{align*}
		\emph{Remark:}
		In the line applying $(*)$ (cf.\ initial paragraph of this proof),
		we use the same indexing convention explained in detail in
		the proof of Lemma \ref{lemma:product_operator_properties}
		\ref{lemma:product_operator_properties:disint_marg}.
		
		Finally, uniqueness follows from the uniqueness
		in Lemma \ref{lemma:disintegration}
		applied to $\mu^T$
		(a kernel on $\topSpace{X}_{C} \times \topSpace{X}_L$)
		factorizing it into the marginalization
		$\marginalize{\mu}{L}$ on $\topSpace{X}_{C}$
		(marginalizing $L$ leaves its the value-space components of its complement $C$),
		and a unique second factor.
		Any kernel satisfying the characterizing property 
		can take the role of this second factor,
		thus by its uniqueness is always the same as the kernel constructed above.
	\end{proof}

	\begin{lemma}[Properties]
		\label{lemma:product_operator_properties}
		These operators have the following properties:
		\begin{enumerate}[label=(\alph*)]			
			\item\label{lemma:product_operator_properties:marg_prod}
			Given products of kernels $\mu^A = \mu^A_1 \otimes \ldots \otimes \mu^A_m$
			and $\mu^B = \mu^B_{m+1} \otimes \ldots \otimes \mu^B_n$
			of kernels and subsets
			$L^A \subset \{1, \ldots, m\}$,
			$L^B \subset \{m+1, \ldots, n\}$,
			such that arguments of $\mu^B$ are
			not in $L^A$,
			then
			\begin{equation*}
				\marginalize{\mu^A}{L^A}
				\otimes
				\marginalize{\mu^B}{L^B}\\
				\halfquad=\halfquad
				\marginalize{
					\mu^A \otimes \mu^B
				}{L^A \cup L^B}
				\txt.	
			\end{equation*}
			
			\item\label{lemma:product_operator_properties:disint_prod}
			Given products of kernels $\mu^A = \mu^A_1 \otimes \ldots \otimes \mu^A_m$
			and $\mu^B = \mu^B_{m+1} \otimes \ldots \otimes \mu^B_n$
			of kernels and subsets
			$C^A \subset \{1, \ldots, m\}$,
			$C^B \subset \{m+1, \ldots, n\}$,
			such that arguments of $\mu^B$ are in $C^A$
			or shared (with $\mu^A$, \ie not in $\{1,\ldots,m\}\setminus C^A$),
			then
			\begin{equation*}
				\disint_{C^A}(\mu^A)
				\otimes
				\disint_{C^B}(\mu^B)\\
				\halfquad=\halfquad
				\disint_{C^A\cup C^B}(\mu^A \otimes \mu^B)
				\txt.	
			\end{equation*}	
			\item\label{lemma:product_operator_properties:disint_marg}
			Given a product $\mu = \mu_1 \otimes \ldots \otimes \mu_m$
			of kernels
			and subsets $C, L \subset \{1, \ldots, m\}$,
			with $C\cap L= \emptyset$, then
			\begin{equation*}
				\disint_C\big(\marginalize{\mu}{L}\big)
				=
				\marginalizeOp_L\big(\disint_C(\mu)\big)
				\txt.		
			\end{equation*}
		\end{enumerate}
	\end{lemma}
	\begin{proof}
		\textbf{Part \ref{lemma:product_operator_properties:marg_prod}:}
		With arguments of $\mu^B$ in $\{1, \ldots, m\}\setminus L^A$ (by hypothesis),
		the left-hand side is defined.
		By definition, for $B_i\in\sigmaAlgebraBorelIdx{X}{i}$
		with $B_i = \topSpace{X}_i$ for $i\in L^A,L^B$,
		\begin{align*}
			\kernelCompound{
				\marginalize{\mu^A}{L^A}
				\otimes
				\marginalize{\mu^B}{L^B}
			}(\prod_{i=1, i\notin L^A\cup L^B}^n B_i)
			&=
			\kernelCompound{
				\mu^A
				\otimes
				\mu^B
			}(\prod_{i=1}^n B_i)\\
			&=
			\marginalize{\mu^A\otimes\mu^B}{L^A\cup L^B}
			(\prod_{i=1, i\notin L^A\cup L^B}^n B_i)\txt.
		\end{align*}
		These products form a basis of the product sigma-algebra,
		thus this proves the claim.

		\textbf{Part \ref{lemma:product_operator_properties:disint_prod}:}
		Define $L^A:=\{1,\ldots, m\}\setminus C^A$, $L^B:=\{m+1,\ldots,n\}\setminus C^B$.
		The right-hand side is (by definition)
		the unique kernel with the universal property
		$\marginalize{\mu^A\otimes\mu^B}{L^A\cup L^B}\otimes
		\disint_{C^A\cup C^B}(\mu^A \otimes \mu^B)=\mu^A\otimes\mu^B$,
		thus it suffices to show,
		that the left-hand side also has this characterizing property.
		To this end, we compute (using the hypothesis
		that arguments of $\mu^B$ are in $C^A$ in the second and third line):
		\begin{align*}
			&\marginalize{\mu^A\otimes\mu^B}{L^A\cup L^B}
			\otimes
			\kernelCompound{
				\disint_{C^A}(\mu^A)
				\otimes
				\disint_{C^B}(\mu^B)
			}
			\\
			=
			&\marginalize{\mu^A}{L^A}\otimes\marginalize{\mu^B}{L^B}
			\otimes
			\disint_{C^A}(\mu^A)
			\otimes
			\disint_{C^B}(\mu^B)
			&\txt{by \ref{lemma:product_operator_properties:marg_prod}}
			\\
			=
			&
			\kernelCompound{
				\marginalize{\mu^A}{L^A}
				\otimes
				\disint_{C^A}(\mu^A)
			}
			\otimes
			\kernelCompound{
				\marginalize{\mu^B}{L^B}
				\otimes
				\disint_{C^B}(\mu^B)
			}
			&\txt{by sparsity \ref{lemma:causal_transformations}}
			\\
			=
			&\mu^A \otimes \mu^B
			&\txt{by definition}
			\txt.
		\end{align*}
		The claim follows by uniqueness of disintegrations.

		\textbf{Part \ref{lemma:product_operator_properties:disint_marg}:}
		Define $L':= \{1, \ldots, m\} \setminus C$.
		By $L\cap C=\emptyset$, we have $L\subset L'$.
		We will again use the universal property of disintegrations,
		for example
		$\marginalize{\mu}{L'}\otimes\disint_C(\mu)=\mu$.
		To apply previous parts, we first have to clarify our indexing convention,
		as terms $\mu_i$ in $\mu$ appear twice, once in the first term, once in the second term. It does not matter which one of these remains (\ie corresponds to the
		right-hand side), only that we get each factor exactly once and in the right order.
		We retain $\mu^A = \mu_1 \otimes \ldots \otimes \mu_m$ as given,
		and define $\mu^B := \mu_{m+1} \otimes \ldots \otimes \mu_n$, where
		$n=2m$ and $\mu{m+i} := \mu_i$, in particular $\mu=\mu^A=\mu^B$
		(as kernels, ignoring indexing).
		For subsets of $\{1,\ldots, m\}$ we write $L+m$ etc.\ for the
		respective shifted versions.
		
		We use the same trick as in \ref{lemma:product_operator_properties:disint_prod},
		and notice that
		the left-hand side is characterized by the property
		$\marginalize{\marginalize{\mu}{L}}{L'\setminus L}
		\otimes\disint_C
		\marginalize{\mu}{L} = \marginalize{\mu}{L}$.
		Thus we compute:
		\begin{align*}
			&\marginalize{\marginalize{\mu}{L}}{L'\setminus L}
			\otimes\marginalize{\disint_C \mu}{L}
			&\\
			&=\marginalize{\marginalize{\mu^A}{L}}{L'\setminus L}
			\otimes\marginalize{\disint_{C+m} \mu^B}{L+m}
			&\txt{by $\mu=\mu^A=\mu^B$}
			\\
			&=\marginalize{\marginalize{\marginalize{\mu^A}{L}}{L'\setminus L}}{\emptyset}
			\otimes\marginalize{\disint_{C+m} \mu^B}{L+m}
			&
			\\
			&=\marginalize{
				\marginalize{\mu^A}{L'}
				\otimes\disint_{C+m} \mu^B
			}{L+m}
			&\txt{by \ref{lemma:product_operator_properties:marg_prod}
				and \ref{lemma:kernel_marginalization_safe}} 
			\\
			&
			=\marginalize{
				\mu
			}{L}
			&\txt{by $\mu=\mu^A=\mu^B$ and Def.}
		\end{align*}
		The claim follows by uniqueness of disintegrations.
	\end{proof}

	\subsection{Sparsity and Contraction}
	\label{apdx:sparsity}

	In principle the constructions introduced so far
	already describe acyclic causal models, in the sense
	that a suitable joint distribution can be described by a product
	over mechanisms in causal order: Each node can be interpreted to
	depend on all its causal-order predecessors, albeit the dependence on
	predecessors that are not direct parents is trivial.
	However, for causal reasoning, precisely this sparsity in arguments,
	of each node depending \emph{only} on its parents, is of central importance.
	Thus a formal description that "forgets" about this structure may be simple,
	but is unlikely to be useful.
	
	We approach the problem in two steps: First, in the present sub-section,
	we formally define what we mean by "sparsity" in arguments.	
	This notion is unfortunately rather unwieldy to work with,
	so §\ref{apdx:structural_graphs} instead
	encapsulates it in a graphical description.
	For most practical purposes, the reader may want to
	continue with §\ref{apdx:structural_graphs}, as it seems rather
	evident that the graphical definitions given there \emph{can}
	be expressed directly in a formal language.
	Nevertheless,
	in case the reader really wants to know, these
	ideas can be expressed formally as follows (we borrow language from
	tensor-networks, and call the wiring of input-arguments to sources contractions):
	
	\begin{definition}[Contraction]\label{def:contraction}
		Given finite sets $I_S, I_T, I_U, I_W$,
		given a kernel $\mu$ from $\topSpace{S}:=\prod_{i\in I_S}\topSpace{X}_i$ to
		$\topSpace{T}:=\prod_{i\in I_T}\topSpace{X}_i$,
		a kernel $\nu$ from  $\topSpace{U}:=\prod_{i\in I_U}\topSpace{X}_i$ to
		$\topSpace{W}:=\prod_{i\in I_W}\topSpace{X}_i$,
		and a mapping $w: I_U \rightarrow I_S \cup I_T$
		such that $\topSpace{X}_{w(i_u)}\subset \topSpace{X}_{i_u}$,
		define a kernel $\mu \otimes_w \nu$
		from $\topSpace{S}$
		to $\topSpace{T}\times\topSpace{W}$
		by
		\begin{align*}
			&\forall s=(x_i)_{i\in I_S} \in \topSpace{S},\halfquad
			\forall (B_i)_{i\in I_T\cup I_W} \in \sigmaAlgebraBorelRaw{
				\topSpace{T}\times\topSpace{W}}:\\
			&\mu \otimes_w \nu((x_i)_{i\in I_S}, (B_i)_{i\in I_T\cup I_W})\\
			&:=\halfquad
			\int \mu(s, (\dx_i)_{i\in I_T})
			\int \nu((x_{w(i)})_{i\in I_U}, (\dx_i)_{i\in I_W})
			1_B((x_i)_{i\in I_T\cup I_U})
			\txt.
		\end{align*}
		Given a  set $\{\mu^n\}_{n\in\nodes}$ of kernels
		indexed by a finite totally ordered (by $\leq$) set $\nodes$
		(\asswlog $\nodes=\{1,\ldots,n\}$ to simplify notation),
		from $\prod_{i\in I^n}\topSpace{X}_i$
		to $\topSpace{X}_n$ (for finite $I^n$),
		a finite set $I_S$ (and $S:= \prod_{i\in I_S} \topSpace{X}_i$) and
		$w^n: I^n \rightarrow I_S \cup \{ m\in \nodes | m< n\}$,
		such that $\topSpace{X}_{w_n(i)}\subset\topSpace{X}_{i}$
		then there is a kernel
		\begin{equation*}
			\otimes_{n\in\nodes}^{w^*} 
			\mu^n
			\halfquad:=\halfquad
			\tilde{\mu}^1 \otimes_{w^2}
			\mu^2 \otimes_{w^3}
			\ldots \otimes_{w^n}
			\mu^n
			\halfquad
			\txt{ from $S$ to }
			\prod_{n\in\nodes}\topSpace{X}_n
			\txt,
		\end{equation*}
		where $\tilde{\mu}^1(s=(x_i)_{i\in I^1}, B) := \mu^1((x_{w^1(i)})_{i\in I^1}, B)$.
		If one or multiple of the $\mu^k$ have no arguments (are measures\Slash{}kernels
		on the one-element set $\{*\}$), replace "$\otimes_{w^k}$" by "$\otimes$",
		if $\mu^1$ has no arguments, then $\tilde{\mu}^1 = \mu^1$.
	\end{definition}

	\begin{notation}[Causal Wiring]\label{notation:causal_wiring}
		In explicit computations,
		we denote kernels by capital roman letters, \eg $X$.
		If the same letter appears as a lowercase roman letter in an index,
		\eg in $Y_x$, then this argument is set to the integration-variable
		taking the value of their capitalized version, \eg
		a mediator plus direct effect (a) vs.\ a chain (mediator only; b)
		then read:
		\begin{align*}
			&\txt{(a)}\quad &
			\kernelCompound{X \otimes M_x \otimes Y_{x,m}}(B)
			&=
			\int X(\dx) \int M_x(\dm)
			\int Y_{x,m}(\dy)
			1_B(x,m,y)
			\\
			&\txt{(b)}\quad &
			\kernelCompound{X \otimes M_x \otimes Y_m}(B)
			&=
			\int X(\dx) \int M_x(\dm) \int Y_m{\dy}
			1_B(x,m,y)
			\txt.
		\end{align*}
		For readability, we repeat inputs that come from "outside" of a bracket around kernels
		\begin{equation*}
			\kernelCompound{X \otimes M_x} \otimes Y_{x,m}
			=
			X \otimes \kernelCompound{M_x \otimes Y_{x,m}}_x
		\end{equation*}
		and put composed upon arguments in square brackets
		($m\mapsto[m]$)
		\begin{equation*}
			(\kernelComposition{Y_{x,[m]}}{M_x})_x
			\txt.
		\end{equation*}
		
		In light of Lemma \ref{lemma:kernel_anticausal_properties}ii
		(see also the warning on Lemma.\ \ref{lemma:kernel_composition_properties}),
		we will say	for kernels $X_{s,t}$ and $Y_{x,s,u}$
		that "overlap" in at least one composed with argument (here: $s$),
		that their product is confounded (see example \ref{example:kernels_confounding})
		and change the notation of the brackets to
		\begin{equation*}
			\kernelCompoundConfounded{X_{s,t} \otimes Y_{x,s,u}}_{[s],t,u} \circ S
		\end{equation*}
		Note that, while this seems to occur less rarely, similar problems
		arise with "confounded" compositions (which are defined through products),
		and we adopt a similar notation:
		\begin{align*}
			\big( \kernelComposition{
				\llparenthesis
				\kernelComposition{Y_{x,[m]}}{M_x}
				\rrparenthesis_{[x]}
			}{X} \big)(B)
			&= \kernelCompound{X \otimes \kernelCompound{M_x \otimes Y_{x,m}}_x}
			(\topSpace{X}\times\topSpace{M}\times B) \\
			&= \big(\kernelComposition{
				\kernelCompoundConfounded{M_x \otimes Y_{x,m}}_{[x]}
			}{X}\big)(\topSpace{M}\times B)
		\end{align*}
	\end{notation}
	\begin{example}\label{example:kernels_confounding}
		For a causal graph $X \leftarrow Z \rightarrow Y$,
		where $\randomVar{X}$ and $\randomVar{Y}$ are confounded by $\randomVar{Z}$,
		we have (see example \ref{example:kernel_product})
		kernels $Z$, $X_z$ and $Y_{z}$, thus for example:
		\begin{equation*}
			\kernelCompound{X_z \otimes Y_z}_z
			\quad\txt{and}\quad
			\kernelComposition{
				\kernelCompoundConfounded{
					X_z \otimes Y_z
				}_{[z]}
			}{Z}
		\end{equation*}
		Note, that composition by definition is taking the
		joint distribution $Z \otimes X_z \otimes Y_z$
		(an observed confounder) and
		integrates out $Z$ (producing "hidden" confounding)
		$\kernelComposition{
			\kernelCompoundConfounded{
				X_z \otimes Y_z
			}_{[z]}
		}{Z}$,
		while \emph{re-associating} as $Z \otimes \kernelCompound{X_z \otimes Y_{x,z}}_z$
		makes the conditioned (on the observed confounder)
		joint $P(\randomVar{X},\randomVar{Y}|\randomVar{Z})$ given by
		$\kernelCompound{X_z \otimes Y_z}_z$ explicit
		(but retains information of\Slash{}dependence on $\randomVar{Z}$).
		
		We can, from the distribution of $\randomVar{Z}$ and
		the conditioned joint distribution $\kernelCompound{X_z \otimes Y_z}_z$
		reproduce the full $Z \otimes X_z \otimes Y_{x,z}$
		(by associativity),
		but from $\kernelComposition{
			\kernelCompoundConfounded{
				X_z \otimes Y_z
			}_{[z]}
		}{Z}$ we can neither decompose (deconvolute)
		in the sense of undoing
		the composition with $Z$ (thus cannot get back
		$\kernelCompound{X_z \otimes Y_{x,z}}_z$)
		nor can we decompose (disintegrate)
		in the sense of obtaining factors
		to get $X_z$, $Y_{z}$ or $\kernelComposition{Y_{[z]}}{Z}$
		(only $\kernelComposition{X_{[z]}}{Z}$ by marginalizing
		from the right).
	\end{example}
	
	\begin{lemma}[Causality Mandated Transformations]
		\label{lemma:causal_transformations}
		Sparsity enables further operations,
		we write $s$, $t$ etc.\ for tuples of arguments,
		and summarize combined arguments as $s\cup t$ and shared
		arguments as $s\cap t$ in the evident way, similarly
		$X$ and $Y$ are allowed to be products of kernels
		(in the sense of having distinguishable factors that
		might partially overlap $s$ or $t$):
		\begin{enumerate}[label=(\roman*)]
			\item\label{lemma:causal_transformations:commute}
			Products without parent--child relation ($Y$ has no index '$x$') commute
			up to transposition:
			\begin{equation*}
				\txt{If $t \cap x = \emptyset = s \cap y$ then}\halfquad
				\kernelCompound{ X_{s} \otimes Y_{t} }_{s\cup t}
				=
				\kernelCompound{ Y_{t} \otimes X_{s} }_{s\cup t}^t
			\end{equation*}
			\item\label{lemma:causal_transformations:compose_trivial}
			Compositions without parent--child relation are trivial:
			\begin{equation*}
				\txt{If $x \cap t = \emptyset$, then}\halfquad
				(\kernelComposition{Y_{t}}{X_{s}})_{s\cup t}
				= 
				Y_{t}
			\end{equation*}
			\item\label{lemma:causal_transformations:compose_unconfounded}
			Composition
			right-distributes over unconfounded	products (cf.\ (ii)):
			\begin{align*}
				\txt{If $t \cap z = \emptyset$, then}\halfquad
				\kernelComposition{
					\kernelCompound{ X_{s} \otimes Y_{t}
					}_{s\cup (t\setminus x)}
				}{Z_u}
				\halfquad&=\halfquad	
				\kernelCompound{				
					(\kernelComposition{ X_{s} }{Z_u})_{(s\setminus z)\cup u}
					\otimes						
					Y_{t}
				}_{(s\setminus z)\cup u \cup (t\setminus x)}
				\\
				\txt{If $s \cap z = \emptyset$, then}\halfquad
				\kernelComposition{
					\kernelCompound{ X_{s} \otimes Y_{t}
					}_{s\cup (t\setminus x)}
				}{Z_u}
				\halfquad&=\halfquad	
				\kernelCompound{				
					X_s
					\otimes	
					(\kernelComposition{ Y_{t} }{Z_u})_{(t\setminus z)\cup u}
				}_{(t\setminus z\setminus x)\cup u \cup s}
			\end{align*}
			\item\label{lemma:causal_transformations:disint_args}
			Anti-causal disintegrations carry arguments and entries
			of the joint distribution
			excluding visible terms:
			\begin{equation*}
				\txt{For $X_s$, $Y_t$:}
				\halfquad
				(X_s|Y_t)_{y\cup s \cup (t\setminus x)}\txt.
			\end{equation*}
		\end{enumerate}
	\end{lemma}	
	\begin{proof}
		\textbf{Part \ref{lemma:causal_transformations:commute}:}
		Note that the argument $\kernelCompound{X_s\otimes_t}_{s\cup (t\setminus x)}$
		by $t\cap x=\emptyset$ is indeed $s\cup (t\setminus x)=s\cup t$.
		By definition of the $\otimes$-product (Def.\ \ref{def:kernel_product}),
		using that $B=B_X\times B_Y$ form a basis
		of $\sigmaAlgebraBorel{X}\otimes\sigmaAlgebraBorel{Y}$:
		\begin{align*}			
			\forall B_X\times B_Y \in
			\sigmaAlgebraBorel{X}\otimes\sigmaAlgebraBorel{Y}&:&\\
			\kernelCompound{X_s \otimes Y_t}(s\cup t,B_X\times B_Y)
			&=
			\int X(s,\dx) \int Y(t,\dy) 1_{B_X\times B_Y}(x, y)
			&
			\\
			&=
			\int X(s,\dx) 1_{B_X}(x) \int Y(t,\dy) 1_{B_Y}(y)
			&\txt{by $t\cap x=\emptyset$}
			\\
			&=
			X(s,B_X) Y(t,B_Y)
			&\txt{def.\ of $\textstyle\int$}
		\end{align*}
		The problem is symmetric under exchange of $X \leftrightarrow Y$,
		thus the same computation
		shows $\kernelCompound{Y_t \otimes X_s}(s\cup t,B_Y\times B_X)
		=\kernelCompound{Y_t \otimes X_s}^t(s\cup t,B_X\times B_Y)
		=X(s,B_X) Y(t,B_Y)$.
		
		\textbf{Part \ref{lemma:causal_transformations:compose_trivial}:}
		By definition of composition (Def.\ \ref{def:kernel_composition}):
		\begin{align*}			
			\forall B_Y \in\sigmaAlgebraBorel{Y}&:&\\
			Y_t\circ X_s(s\cup t, B_Y)
			&=
			\int X(s,\dx) \int Y(t,\dy) 1_{\topSpace{X}\times B_Y}(x, y)
			&
			\\
			&=
			\int Y(t,\dy) 1_{B_Y}(y) \int X(s,\dx)
			&\txt{by $t\cap x=\emptyset$}
			\\
			&=
			Y(t,B_Y)
			&\txt{$\textstyle\int X=1$; def.\ of $\textstyle\int$}
		\end{align*}
		
		\textbf{Part \ref{lemma:causal_transformations:compose_unconfounded}:}
		By definition of composition (Def.\ \ref{def:kernel_composition}),
		using that $B=B_X\times B_Y$ form a basis
		of $\sigmaAlgebraBorel{X}\otimes\sigmaAlgebraBorel{Y}$
		for the first case ($t \cap z=\emptyset$):
		\begin{align*}			
			&\quad\forall B_X \times B_Y
			\in\sigmaAlgebraBorel{X}\otimes\sigmaAlgebraBorel{Y}:\\
			&\kernelCompound{X_s\otimes Y_t}\circ Z_u(
				(s\setminus z) \cup (t\setminus x), B_X\times B_Y)\\
			&\quad=
			\int Z(u,\dz)\int X(s,\dx) \int Y(t,\dy)
			1_{\topSpace{Z}\times B_X\times B_Y}(z, x, y)
			\\
			&\quad=
			\int Z(u,\dz)\int X(s,\dx)1_{B_X}(x)
			\int Y(t,\dy) 1_{B_Y}(y)
		\end{align*}
		Also by
		definition of composition (Def.\ \ref{def:kernel_composition})
		and definition of integrals:
		\begin{align*}
			\forall B_X \in\sigmaAlgebraBorel{X}&\\
			\int 	
			(X_s\circ Z_u)((s\setminus u)\cup t, \dx)
			1_{B_X}(x)
			&=
			X_s\circ Z_u((s\setminus u)\cup t, B_X)\\
			&=
			\int Z(u,\dz) \int X(s,\dx) 1_{\topSpace{Z}\times B_X}(z,x)
			\\
			&=
			\int Z(u,\dz) \int X(s,\dx) 1_{B_X}(x)
		\end{align*}		
		By definition of integrals, for any measurable mapping
		$f: \topSpace{X}\rightarrow [0,1]$ thus:
		\begin{align*}
			\int 	
			(X_s\circ Z_u)((s\setminus u)\cup t, \dx)
			f(x)
			=
			\int Z(u,\dz) \int X(s,\dx) f(x)
		\end{align*}
		The choice
		$f_t(x) := 1_{B_X}(x) \int Y(t,\dy) 1_{B_Y}(y)$
		is valid for any fixed $t$ and we can fix a $t$,
		because by hypothesis, $t \cap z=\emptyset$, so
		plugging this integral-operator into the previous equation
		we get for all $t$:
		\begin{align*}			
			&\quad\forall B_X \times B_Y
			\in\sigmaAlgebraBorel{X}\otimes\sigmaAlgebraBorel{Y}:\\
			&\kernelCompound{X_s\otimes Y_t}\circ Z_u(
			(s\setminus z) \cup (t\setminus x), B_X\times B_Y)\\
			&\quad=
			\int 	
			(X_s\circ Z_u)((s\setminus u)\cup t, \dx)
			1_{B_X}(x)
			\int Y(t,\dy) 1_{B_Y}(y)\\
			&\quad=
			\bigKernelCompound{(X_s\circ Z_u) \otimes Y_t}\big((s\setminus z)\cup u \cup (t\setminus x), B_X\times B_Y\big)
		\end{align*}
		The second case works similar (but is actually easier,
		because we can move the integral over $Z$ inside the integral over $X$).
		
		\textbf{Part \ref{lemma:causal_transformations:disint_args}:}
		The joint distribution
		$\kernelCompound{X_s\otimes Y_t}^t_{s\cup (t\setminus x)}$
		carries arguments $s\cup (t\setminus x)$ by construction
		(Def.\ \ref{def:kernel_product}).
		Disintegration adds an argument $y$ for $Y_t\circ X_s$
		by definition (Def.\ \ref{lemma:disintegration}),
		so in general $(X_s|Y_t)_{y\cup s \cup (t\setminus x)}$.
		
		\emph{Remark:} We only need (in the proof of Lemma \ref{lemma:properties_c_components}) that $(X_s|Y_t)$ does not
		have \emph{more than} these arguments.
		Considering that these are basically all plausible
		arguments that could be added this is result is not actually
		surprising.
		It is, however, in general not possible to remove any of these
		arguments, which further emphasizes that anti-causal
		(as opposed to causal-direction) disintegrations
		are really properties of the joint distribution and
		in general as complicated as joint distributions
		from the causal perspective.
	\end{proof}

	\subsection{Symmetries}
	\label{apdx:symmetries}
	
	As the attentive reader may have noticed, we do not
	really rely on symmetries described as group-actions,
	rather on the orbit-sets of these group-actions.
	Indeed we could easily reformulate
	our formalism as follows.
	Change Def.\ \ref{def:mechanism} to:
	
	\begin{definition}[Mechanism, Alternative Version]
		\label{def:mechanism_alt}
		A mechanism is a probability kernel
		$f$ from
		$\topSpace{X}^{\Pa}=\prod_{k=1}^\kappa\topSpace{X}^{(k)}$
		(with $\kappa$ arguments) to
		$\topSpace{Y}$ (both standard Borel, Ass.\ \ref{ass:standard_borel})
		together with
		a region of applicability $J\subset I$,
		and a mapping
		$\Pa : J \rightarrow I^\kappa$ of relative parents
		such that
		(\textbf{no sub-group or equivariance-condition}):
		
		Denoting the $k$th parent by $\PaIdx{k}$,
		require $\forall k: \PaIdx{k}(j)\neq j$ and
		$k\neq k'$ $\Rightarrow$ $\forall j\in J$:
		$\PaIdx{k}(j) \neq \PaIdx{k'}(j)$.
		
		\textbf{Additionally:}
		$\forall k \: \PaIdx{k}$ is injective.
		
		Finally,
		$\forall j\in J$,
		$\topSpace{Y} \subset \topSpace{X}_j$ and
		$\forall k$:
		$\topSpace{X}_{\PaIdx{k}(j)}\subset\topSpace{X}^{(k)}$.
	\end{definition}
	\begin{rmk}\label{rmk:symm_alt_more_general}
		A mechanism in the sense of Def.\ \ref{def:mechanism}
		is a mechanism in the sense of Def.\ \ref{def:mechanism_alt}.
	\end{rmk}
	
	The additional condition immediately restores freeness
	of the direct embedding (Lemma \ref{lemma:direct_embedding}),
	so if we delete the equivariance (rigidity", 
	\ref{def:local_graph_embedding_family:rigidity})
	in Def.\ \ref{def:local_graph_embedding_family}:		
	
	\begin{definition}[Families of Embeddings,
		Alternative Version]\label{def:local_graph_embedding_family_alt}
		A family of local graph embeddings
		\begin{equation*}
			(J_0, \{\psi_j\}_{j\in J_0}, \graph, y_0)
		\end{equation*}
		is a collection of local graph embeddings
		$\psi_j: \graph\hookrightarrow I$
		of a single model-aligned local graph $\graph$.
		Further there is a fixed element $y_0 \in \nodes$,
		which will be called the anchor\allowbreak{}\mbox{(-}node)
		 and a range of applicability $J_0\subset I$.
		For $n\in\nodes$ denote:
		\begin{equation*}
			\psi_*(n): J_0 \rightarrow I,
			\halfquad
			j\mapsto \psi_j(n)
			\txt.
		\end{equation*}
		Finally, we will
		require the following conditions to be satisfied:
		\begin{enumerate}[label=(\Roman*)]
			\item\label{def:local_graph_embedding_family:trivial_anchor:alt_symm}
			Trivial on Anchor:
			$\forall j\in J_0$: $\psi_j(y_0) = j$.
			\item\label{def:local_graph_embedding_family:rigidity:alt_symm}
			(\textbf{no rigidity condition})
			\item\label{def:local_graph_embedding_family:freeness:alt_symm}
			Freeness:
			$\forall n\in\nodesInnerProper$:
			$\psi_*(n)$
			is injective.
		\end{enumerate}
	\end{definition}	
	
	This actually leaves us in a place where Thm.\ \ref{thm:extract_from_backdoor_complete}  still applies!
	Indeed the proof of Thm.\ \ref{thm:extract_from_backdoor_complete},
	while painstakingly going through all the properties of the definitions,
	does never reference rigidity.
	So everything else works exactly as before.

	\paragraph{Motivation:}
	Considering Rmk.\ \ref{rmk:symm_alt_more_general}
	this alternative formulation is strictly more general,
	while at the same time producing the same output.
	There are, however very good reasons for an approach by
	group-actions.
	
	\begin{itemize}
		\item
		Group-actions are the standard approach to symmetries
		in physics and other natural sciences.
		Thinking about a model for example as time-translation invariant
		(or $\lambda$-periodic etc.) is much more intuitive than
		thinking about symmetries in terms of their induced orbit-spaces.
		Also, the gain in generality only matters, if an application
		would specify a symmetry that is \emph{not induced} by
		Rmk.\ \ref{rmk:symm_alt_more_general}.
		\item 
		Group-actions allow for an elegant structure.
		For example observations as in Rmk.\ \ref{rmk:families_and_symmetries}
		are intuitive and immediate in the language of group-theory.
		More generally, having a sub-group hierarchy paralleling
		constructions has been of great relevance in many
		topics, from (classical) Galois-theory to classifications
		of covering-maps.
		\item 
		Group-actions allow for efficient encoding.
		Going forward, at some point, we will certainly want to
		\emph{find} models from data
		in some form of structure-discovery (§\ref{apdx:structure_discovery}).
		Clearly searching all possible structures
		is not a viable option. Here group-actions should help
		to efficiently encode both models and assumptions
		about priors over models.
		\item 
		This paper is targeted at an audience in machine-learning
		and applied sciences, and for that purpose is already on
		the rather abstract side.
		Making it yet more abstract without substantial benefit
		seems unwise.
		The pure mathematician reading this is of course invited
		to investigate simpler and more abstract ideas.
	\end{itemize}

	\newpage
	
	\tableofcontents
	
\end{document}